\documentclass[12pt]{report}

\usepackage{textcomp}
\usepackage{ifthen}
\usepackage{uga}
\let\appendix\appendices % To make uga.sty and hyperref.sty both happy.
\usepackage{mathrsfs}
\usepackage[centertags]{mathtools}
\usepackage{amsthm,amssymb}
\usepackage{enumitem}
\usepackage{units}
\usepackage{graphicx,overpic}
\usepackage{varioref}
\usepackage{float}
\usepackage[footnotesize,bf]{caption}
\usepackage{subfig}
\usepackage{xcolor}
\usepackage[final,dvips,bookmarks=true,bookmarksnumbered=true,breaklinks,colorlinks,linkcolor=blue!70!black,anchorcolor=blue!70!black,citecolor=red!70!black,filecolor=blue!70!black,menucolor=green!70!black,pagecolor=green!70!black,urlcolor=blue!70!black,hyperindex=false]{hyperref}
\usepackage[author-year]{amsrefs}
\usepackage[matrix,arrow,curve]{xy}
\usepackage{tikz}
\usetikzlibrary{snakes}
\usetikzlibrary{shapes}
\usepackage{ltablex}
\usepackage{index}
\usepackage{varindex}
\usepackage{hhline}
\definecolor{refkey}{rgb}{0,0,0.8}
\definecolor{labelkey}{rgb}{0.8,0,0.8}
\mathtoolsset{centercolon,mathic}

\labelformat{chapter}{Chapter~#1}
\labelformat{section}{Section~#1}
\labelformat{subsection}{Subsection~#1}
\labelformat{equation}{Equation~(#1)}
\labelformat{table}{Table~#1}
\labelformat{figure}{Figure~#1}
\renewcommand{\reftextlabelrange}[2]{\ref{#1} through~\ref{#2}}

\title{Exploring Continuous Tensegrities}
\author{Edward Bruce ``Ted'' Ashton}
\previousdegrees{B.S.E., Walla Walla College, 1993 \\
                 B.S., Southern Adventist University, 1999}
\thisdegree{Doctor of Philosophy}
\thesistype{Dissertation}
\thesismonth{May}
\thesisyear{2007}
\professor{Dr.\ Jason Cantarella}
\fulltitle{Exploring Continuous Tensegrities}
\indexwords{Bar Equivalence, Dissertation, Framework, Infinitesimal Motion,
Infinitesimal Rigidity, Stress, System of Inequalities, Tensegrity, Theorem of
the Alternative}
\dean{Dr.\ Maureen Grasso}

\memberi{Dr.\ Malcolm R. Adams}
\memberii{Dr.\ Edward A. Azoff}
\memberiii{Dr.\ Joseph H. G. Fu}
\memberiv{Dr.\ Robert Varley}

\newcommand{\B}{\mathsf{B}}
\newcommand{\Bb}{{\overline{\B}}}

\newcommand{\C}{\mathsf{C}}
\newcommand{\Cr}{\mathscr{C}}
\newcommand{\cross}{\times}
\newcommand{\ds}{\, ds}
\newcommand{\dt}{\, dt}
\newcommand{\dmu}{\, d\mu}
\newcommand{\E}{\mathscr{E}}

\newcommand{\ev}[1]{\overrightarrow{#1}}
\newcommand{\Gb}{\overline{G}}

\newcommand{\Ib}{\overline{I}}

\newcommand{\mcol}{\mathord{:}\,}
\newcommand{\muh}{\widehat{\mu}}

\newcommand{\one}{\mathbf{1}}
\newcommand{\qrp}{\quotientspace{\RR}{\pi}}
\newcommand{\RR}{\mathbb{R}}
\newcommand{\St}{\mathsf{S}}
\newcommand{\Str}{\mathscr{S}}
\newcommand{\V}{\mathscr{V}}
\newcommand{\Vf}{\mathsc{V\kern-0.1667em f}}
\newcommand{\vvt}{{\quotientspace{(\V \cross \V)}{\sim}}}
\newcommand{\X}{\mathcal{X}}
\newcommand{\zero}{\mathbf{0}}
\newcommand{\ZZ}{\mathbb{Z}}

\newcommand{\capt}[2][null]{\ifthenelse{\equal{#1}{null}}%
  {\caption[#2]{\setlength{\baselineskip}{0.5\baselineskip}#2}}%
  {\caption[#1]{\setlength{\baselineskip}{0.5\baselineskip}#2}}}
\newcommand{\quotientspace}[2]{
  \ensuremath{#1}\kern-0.3em\lower0.6ex\hbox{\ensuremath{\diagup}}%
    \kern-0.3em\lower1ex\hbox{\ensuremath{#2}}
}

\renewcommand{\ge}{\geqslant}

\renewcommand{\le}{\leqslant}

\newenvironment{dotslist}{\begin{list}{\ldots}{}}{\end{list}}

\theoremstyle{plain}
\newtheorem{theorem}{Theorem}[section]
\newtheorem{lemma}[theorem]{Lemma}
\newtheorem{claim}[theorem]{Claim}
\newtheorem{conjecture}[theorem]{Conjecture}
\newtheorem{prop}[theorem]{Proposition}
\newtheorem{corollary}[theorem]{Corollary}
\newtheorem*{KL}{Lemma \ref{lemma:NewFiveOne}}
\newtheorem*{MTone}{Theorem \ref{theorem:nonnegativeMeasure}}
\newtheorem*{MTtwo}{Theorem \ref{theorem:positiveImplies}}
\newtheorem*{MTthree}{Theorem \ref{theorem:minimallyBarEquivalent}}
\newtheorem*{MTfour}{Theorem \ref{theorem:CountablyCovered}}
\newtheorem*{COSone}{Proposition \ref{prop:COSone}}
\newtheorem*{COStwo}{Proposition \ref{prop:COStwo}}
\newtheorem*{ConjOne}{Conjecture \ref{conj:CountCov}}
\newtheorem*{ConjTwo}{Conjecture \ref{conj:StrictlyPositive}}
\newtheorem*{ConjThree}{Conjecture \ref{conj:minXbarConj}}

\newcommand{\statementNewFiveOne}{%
Let $G(p)$ be a (finite) tensegrity framework with rigidity matrix~$Y$.  Then 
$G(p)$ is bar equivalent if and only if $G(p)$ has a strictly positive
stress\Index{{strictly positive} stress $\Leftrightarrow$ {bar
equivalent}}[1!2!34 2!1!34 4!312].
}
\newcommand{\statementPositiveImplies}{%
If a tensegrity has a strictly positive stress, then it is $\X$-bar equivalent.
\Index{{strictly positive} stress $\Rightarrow$ $\Leftarrow$ {$\X$-bar
equivalent@X-bar equivalent}}[1!2!35 2!1!35 5!412]
}
\newcommand{\statementNonnegativeMeasure}{%
The tensegrity $G(p)$ is partially $\X$-bar equivalent if and only if $G(p)$
has a semipositive stress\Index{semipositive stress $\Leftrightarrow$ partially
$\X$-bar@X-bar equivalent equiv.}[1!2!3457 2!1!3457 456!312 56!4!312].
}
\newcommand{\statementCOSone}{%
There is no locally arclength preserving motion of a circle that increases all
antipodal distances.
}
\newcommand{\statementCOStwo}{%
There is no locally arclength preserving motion of a circle that increases
any antipodal distance without decreasing some other one.
}
\newcommand{\statementMinimallyBarEquivalent}{%
If $G(p)$ is minimally $\X$-bar equivalent, $G(p)$ has a strictly positive
stress.
\Index{minimally {$\X$-bar equivalent@X-bar equivalent} {$\X$-bar equiv.@X-bar
equiv.} $\Rightarrow$ $\Leftarrow$ {strictly positive} stress}[12!467 2!1!467
6!7!513 7!6!513]
}
\newcommand{\statementCountablyCovered}{%
If $G(p)$ is countably covered by subtensegrities that have a strictly
positive stress, then $G(p)$ has a strictly positive stress.%
\Index{{countably covered} $\Rightarrow$ $\Leftarrow$ {strictly positive} 
stress}[1!245 5!4!31 4!5!31]%
}
\newcommand{\statementCountCov}{%
Every $\X$-bar equivalent continuous tensegrity is countably covered
by\Index{countably covered}[ ] minimally $\X$-bar equivalent subtensegrities.
}
\newcommand{\statementStrictlyPositive}{%
Every $\X$-bar equivalent tensegrity has a strictly positive stress.
}
\newcommand{\statementminXbarConj}{%
A tensegrity is minimally $\X$-bar equivalent if and only if it admits only
one semipositive stress, up to scaling.  Furthermore, that stress is strictly
positive.
}

\theoremstyle{definition}
\newtheorem{definition}[theorem]{Definition}

\DeclareMathAlphabet{\mathsc}{OT1}{cmr}{m}{sc}
\DeclareMathOperator{\aut}{Aut}
\DeclareMathOperator{\im}{im}

\DeclareMathOperator{\supp}{supp}

\DeclareMathOperator{\Span}{span}

\DeclareMathOperator{\INT}{int}
\DeclareMathSymbol{:}{\mathop}{operators}{"3A} % Reduce spacing around :

\tikzstyle{every picture}=[>=stealth]
\tikzstyle{strut}=[ultra thick]
\tikzstyle{cable}=[ultra thin,black!82!white,join=round]
\tikzstyle{bar}=[double,thick]
\tikzstyle{vcurve}=[dotted,thick]
\tikzstyle{vertex}=[circle,thin,draw=black,fill=white,scale=0.5]
\tikzstyle{weight}=[circle,color=black,fill=white,scale=0.75]

\makeatletter
\def\toclevel@pseudochapter{0}
\let\ps@plain=\ps@empty
\@addtoreset{figure}{section}
\renewcommand*\l@figure{\@dottedtocline{1}{1.5em}{3.9em}}
\renewcommand{\BibLabel}{%
  \smash{\hyper@anchorstart{cite.\CurrentBib}\hyper@anchorend}%
}
\makeatother

\newcommand{\Index}{\varindex(){\varindextwoScan}{\varindextwo}[]}
\newcommand{\IndexDef}[1]{\Index[|textbf]{#1}}
\newcommand{\IndDefBeg}[1]{\Index[|(textbf]{#1}}
\newcommand{\IndDefEnd}[1]{\Index[|)textbf]{#1}}
\varindexUsePlaceholderAfalse

\makeindex
\indexproofstyle{\tiny\sf}

\begin{document}
\nonfrenchspacing
\pagenumbering{roman}

\begin{abstract}
A discrete tensegrity framework can be thought of as a graph in Euclidean
$n$-space where each edge is of one of three types: an edge with a fixed length
(bar) or an edge with an upper (cable) or lower (strut) bound on its length.
Roth and Whiteley, in their 1981 paper ``Tensegrity Frameworks'', showed that
in certain cases, the struts and cables can be replaced with bars when
analyzing the framework for infinitesimal rigidity.  In that case we call the
tensegrity \emph{bar equivalent}.  In specific, they showed that if there
exists a set of positive weights, called a positive \emph{stress}, on the edges
such that the weighted sum of the edge vectors is zero at every vertex, then
the tensegrity is bar equivalent.  

In this paper we consider an extended version of the tensegrity framework in
which the vertex set is a (possibly infinite) set of points in Euclidean
$n$-space and the edgeset is a compact set of unordered pairs of vertices.
These are called \emph{continuous tensegrities}.  We show that if a continuous
tensegrity has a strictly positive stress, it is bar equivalent and that it has
a semipositive stress if and only if it is partially bar equivalent.  We also
show that if a tensegrity is \emph{minimally bar equivalent} (it is bar
equivalent but removing any open set of edges makes it no longer so), then it
has a strictly positive stress.

In particular, we examine the case where the vertices form a rectifiable curve
and the possible motions of the curve are limited to local isometries of it.
Our methods provide an attractive proof of the following result:
\statementCOStwo
\end{abstract}

\maketitle

\chapter*{Dedication}
\thispagestyle{headings}
\begin{center}
To my wife and best friend, \\ 
{\fontfamily{pzc}\selectfont \Large Heidi}
\end{center}

\pseudochapter{Acknowledgments}
\thispagestyle{headings}
It is rare in our lives that we reach a point of completion, but the end of a
doctorate is certainly such a point.  As I pause to contemplate, I find that
there are many to whom I owe a debt of gratitude.

Certainly the Chief Contributor is the One who hid Tensegrity within the cell
long before we discovered it and who gives us strength even when the stresses
don't appear to be strictly positive.  My thanks is ever to Him.

Among the people who have helped to bring this day about, my family must head
the list.  My wife has given of herself in myriad ways.  She has made a
home for us all and sacrificed much in the pursuit of this degree.  Her worth
is truly far above rubies.  

I'm grateful, also, for the love and support of our children, Jo and Jim, as
well as our ``other children'', Kristen, Jill and Gene.  I have been immensely
grateful for the encouragement of our ten siblings and particularly of our
loving parents, who have supported us emotionally, spiritually and financially. 

Many are the friends and family who have provided aid and comfort.  In
particular, I'm grateful to Aunt Rilla for giving us true vacation times and to
Judy deLay, the Nelsons and the Crosses all for their hospitality.  The Athens
Seventh-day Adventist Church has been a true church home for us.  They welcomed
us with open arms when we arrived and have continued to love and support us
ever since.  When I think of our friends and extended family around the globe,
I am truly amazed at how blessed we are. 

In the mathematical realm, my advisor, Dr.\ Jason Cantarella, leads the pack.
He took me on as an apprentice long before I started my dissertation and poured
many hours, sometimes entire days, into this work.  He has taught me and
advised me and many times provided just the clue I needed.  I am grateful for
his depth of knowledge of the field and his patience in explaining it to me.  I
am doubly honored to have an advisor whom I can also call my friend.

Untold numbers of people have been involved in getting me to where I am now,
far more than I can mention here.  But there are a few who deserve particular
notice.  The members of my committee, Drs.\ Malcolm Adams, Edward Azoff,
Joseph Fu and Robert Varley all gave time to this work.  Dr.\ Adams provided
many valuable editorial notes.  Dr.\ Varley has answered numerous questions,
both mathematical and otherwise.  Dr.\ Azoff and Dr.\ Fu both spent countless
hours answering questions and teaching me such things as Functional Analysis
and Linear Algebra.  Their sage advice and deep knowledge were the well from
which I repeatedly drew in my efforts.  I have had many fine professors here
and I thank them all, but one more who deserves a mention is one from whom I
never took a class.  Dr.\ Ted Shifrin spent his time and energy going over exams
I had designed and giving me advice on teaching and I am a better teacher
because of it.  

However, he was not the first to encourage me in teaching.  My fifth grade
teacher, Mrs.\ Barbara Stanaway, gave me my first taste of math tutoring and I
have loved it ever since.  Mr.\ Morford, one of my high school math teachers,
was the first to put me in front of a class and has always encouraged me in my
pursuit of mathematics.  Years later he again trusted me with his class and has
followed my progress with interest.  

Many others have not only taught me but also cared about me as a student:
Mr.\ Magi, who gave me my first taste of Geometry and Dr.\ Moore, who took me on
from where he left off; Dr.\ Richert, who showed me that Mathematics is a
performance art; Drs.\ Kuhlman and Hefferlin, the world-class physicists who
invited me to be a physics major and kept on investing time in me even when I
turned them down; Dr.\ Haluska, who taught me to write and Dr.\ Sylvia
Nosworthy, who taught me (amongst many things) that my writing needs more
fluff; Dr.\ Masden, who not only demanded great things of his students, but
always went the distance with us; Dr.\ Barnett, who opened to me the world of
simulation and modeling; Dr.\ Rob Frohne who showed me the miracle of feedback
and control; Ralph Stirling, who was not only a fantastic resource for an
engineering student, but also a friend; Dr.\ Cross and his family, who not only
provided good education, but opened their home to us; Dr.\ Wiggins, who started
me into Numerical Analysis; Dr.\ Tim Tiffin, who gave me my first serious taste
of mathematical research; and Dr.\ Tom Thompson, who taught and encouraged me
not only in Mathematics, but also in learning \LaTeX, knowledge I have put to
good use.

My thanks goes to these and the other professors I have had who have chosen to
give their lives in service and who care deeply about their students.  In
particular, I am grateful to Drs.\ Wiggins, Thompson, Richert and Moore for
taking a lonely graduate student under their wings at the National Math
Meetings and providing for Sabbath away from home.  

My thanks also goes to John Pardon, for an illuminating conversation about
nonconvex curves, to Branko Gr\"unbaum for providing a copy of ``Lost
Mathematics'' and to Kenneth Snelson for his gracious reply to this graduate
student's email about his sculptures.  I appreciate also the camaraderie of my
fellow graduate students: Jeremy Francisco, Becca Murphy, Amy Kelly, Valerie
Hower, Alan Thomas, Xander Faber, Bryan Gallant, David and Alicia Beckworth and
my officemates Kenny Little and Chad Mullikin and also of the (then) postdocs
Aaron Abrams and Nancy Wrinkle.  As I look at the piles of library books on my
desk, I thank the people of the University of Georgia Science Library.  Rarely
have I needed a resource that they did not provide.

Finally, I think of Grandmother, who, eighty years ago, was a math major
herself.  She would be proud.  
\tableofcontents
\clearpage
\phantomsection
\listoffigures
\newpage
\pagenumbering{arabic}

\clearpage
\begin{chapter}{Introduction: Tensegrity}
\label{chapter:Intro}
\section{A History}
Stone\Index{stone} is strong.  Good quality stone can withstand amazing
compressive forces\Index[|(]{compression}.  But in
tension\Index[|(]{tension} it is not nearly as strong.  Similarly,
brick\Index{brick} and concrete\Index{concrete} are much stronger
under compression than under tension (\vref*{table:strengths} gives some
representative numbers).  And wood\Index{wood}, while fairly strong under
tension parallel to the grain, can't compare to the durability of stone.  So
for millenia, buildings were designed to be continuous frameworks of members in
compression \cites{Pugh,Girvin,Wood}.
\begin{table}[ht]
\capt[Ultimate stresses for certain materials.]{Ultimate stresses for certain
materials.  Numbers have been converted to consistent units and rounded.
Bending stress given for wood is actually modulus of rupture.  Values for bone
are from \cite{STM} and are given for interest, since tensegrity theory is used
to understand the musculoskeletal system (see, for example, \ocite{Arthritis}).
Values marked $^\text{a}$ are from \ocite{Wood}; those marked $^\text{b}$ are
from \ocite{Girvin}, and those marked $^\text{c}$ come from \ocite{MatWeb}.
}
\label{table:strengths}
\hfil
\begin{tabular}{|l|r|r|r|} \hhline{~|---|}
\multicolumn{1}{c|}{} & \multicolumn{3}{c|}{\small \bfseries Average Ultimate
Stress (\unitfrac{lb}{in$^2$})} \\ \hline
\textbf{Material} & \textbf{Compression} & \textbf{Tension} & \textbf{Bending}
\\ \hline
% Granite & 12,000 & 1,200 & 1,600 \\
Granite\Index{granite} & 14,000--45,000$^\text{c}$ & 1,020--3,630$^\text{c}$ &
1,600$^\text{b}$ \\
Brick & 10,000$^\text{b}$ & 200$^\text{b}$ & 600$^\text{b}$ \\
Steel\Index{steel} & 60,000-72,000$^\text{b}$ & 60,000-72,000$^\text{b}$ &
60,000-72,000$^\text{b}$ \\
Concrete\Index{concrete} & 2,030--10,200$^\text{c}$ & ---$^{\phantom{a}}$ &
---$^{\phantom{a}}$ \\
Sugar Maple\Index1{maple, sugar} (para. to grain) & 8,000$^\text{a}$ &
15,700$^\text{a}$ & 16,000$^\text{a}$ \\
Pin Oak\Index1{oak, pin} (para. to grain) & 7,000$^\text{a}$ &
16,300$^\text{a}$ & 14,000$^\text{a}$ \\ \hline
Human bone\Index{human bone}[2,1 ] (longitudinal) & 28,000$^{\phantom{a}}$ &
19,000$^{\phantom{a}}$ & 30,000$^{\phantom{a}}$ \\ \hline
% Steel (Bridges and Bldgs) (kips/square inch) & 60-72 & 60-72 & 60-72 \\
% Human bone (MPa) (longitudinal) & 195 & 133 & 208.6 \\
% Human bone (MPa) (transverse) & 133 & 51 & 208.6 \\
\end{tabular}
\hfil
\end{table}

The middle of the nineteenth century saw a number of major advances in the
creation of steel\Index{steel}, including the discovery of its microstructure
and the resulting development of the science of metallurgy\Index{metallurgy}
\cite{Timetables}.  This provided a building material with excellent tensile
strength (see, for example, \ocite{CRC}).

In 1965, Kenneth D.\ Snelson\Index[|(]_,{Kenneth D. Snelson}
patented\Index{tensegrity patent}[2,1 ] what he called ``Continuous Tension,
Discontinuous Compression Structures'' \cite{TensegrityPatent} (see
\ref{fig:patent} and \ref{fig:SnelsonArt}).  His new creations had cables under
tension running throughout and then, at intervals along the cables, steel
struts under compression.
\begin{figure}[ht]
\hspace*{\fill}
\begin{overpic}{pat1.epsf}
\end{overpic}
\hspace*{\fill}
\capt[A figure from Kenneth Snelson's tensegrity patent.]{A figure from patent
3,169,611 by Kenneth Snelson \ycite{TensegrityPatent}.} 
\label{fig:patent}
\end{figure}

Buckminster Fuller\Index_,{Buckminster Fuller} recognized in
Snelson's\Index[|)]_,{Kenneth D. Snelson} idea new possibilities in Architecture
that would allow for much lighter, more efficient structures.  He termed the
concept ``Tensional Integrity''\IndexDef{tensional integrity}[ ] or
``Tensegrity''\Index{tensegrity {origin of term}}[1!2].  Since then, Tensegrity
has been used in fields as diverse as Architecture\Index{Architecture},
Cellular Biology\Index{Cellular Biology}[2,1 ], Dairy Science\Index{Dairy
Science}[ ], Dance\Index{Dance}, Ornithology\Index{Ornithology}, Robot
Kinematics\Index{Robot Kinematics}[ ] and the studies of sleep
disorders\Index{sleep disorders}[ ] and the musculoskeletal
system\Index{musculoskeletal tensegrity}[2!1 ].
\cites{Pugh,Slate,Volokh,Dairy,Dance,Macrophage,Lung,RobotKinematics,Sleep,Arthritis}.

\section{Tensegrity in Mathematics}
In 1981, Ben Roth and Walter Whiteley\Index_,{Ben Roth}\Index{Walter
Whiteley}[2,1] wrote ``Tensegrity Frameworks'' \cite{MR610958} in which they
showed how to extend the tools of rigidity analysis on bar frameworks
($n$-dimensional graphs in which the edge lengths are fixed) to certain
tensegrity frameworks.  This was one of the earliest works treating
tensegrities in the mathematical realm and will be a central motivator in what
we do here.

For us, a \emph{tensegrity}\Index{tensegrity} will consist of four sets and a
map.  First, there will be  a set $\V$ of
vertices\Index[|see{$\V$}]{vertex set}[2!1 ]\IndexDef{$\V$@V
(vertices)}[ ].  We will want these to be points in Euclidean $n$-space, but
we'll also want to be able to think about these points moving, so we'll take
$\V$ to be an abstract set accompanied by a 1-1 continuous map $p\mcol \V
\to \RR^n$\IndexDef{$p$@p (map from $\V$ to $\RR^n$)} that places the elements
of $\V$ in space.  Next, there are three sets of edges connecting various pairs
of vertices, with the understanding that at most one edge connects any given
pair.
\label{notation}

These edge sets are the set of \emph{struts}, $\St$\IndexDef{$\St$@S(0)
(struts)}[ ], the set of \emph{cables}, $\C$\IndexDef{$\C$@C(0) (cables)}[ ],
and the set of \emph{bars}, $\B$\IndexDef{$\B$@B (bars)}[ ].  Two vertices
connected by a strut\IndexDef{strut} are never allowed to move any closer
together than they already are, though they may move farther apart.  Two
vertices connected by a cable\IndexDef{cable} cannot move farther apart but may
move nearer.  And two vertices connected by a bar\IndexDef{bar} must stay
exactly as far apart as they currently are.

If the only edges a tensegrity has are bars, we may call the tensegrity a 
\emph{bar framework}\Index{bar framework}[2,1 ].

To date, tensegrities have been finite affairs, with a finite number of finite
elements.  But the time has come to look outside that realm.  We will take $\V$
to be an arbitrary set.  That will open the door to having continuous families
of edges and so we will call our new constructions ``continuous
tensegrities''.\IndexDef{continuous tensegrity}[2!1 ]

\section{Prestressed Concrete}
\label{section:prestressed}
Concrete, like stone, is fairly strong in compression and not nearly so
strong in tension \cite{Girvin}*{p.\ 213}.  Unfortunately, a large slab of
concrete supported only at the edges, such as might form the floor of a parking
garage or the roadway of a bridge, experiences a bending load that puts the
lower portion of the concrete in tension (see \vref{fig:bentConcrete})
\cite{Girvin}*{pp.\ 81--82}.
\begin{figure}[ht]
\hspace*{\fill}
\begin{tikzpicture}
\tikzstyle{every node}=[circle,scale=0.6]
\draw[line width=1ex,red!40!black] 
  (-3.473,-19.696) node (a) {} arc (260:280:20cm) node (b) {};
\draw[->] (-1.726,-19.725) arc (265:269:19.8cm);
\draw[->] ( 1.726,-19.725) arc (275:271:19.8cm);
\draw[<-] (-1.76,-20.123) arc (265:269:20.2cm);
\draw[<-] ( 1.76,-20.123) arc (275:271:20.2cm);
\draw[line width=1ex,black!30] (a) -- +(0,-1);
\draw[line width=1ex,black!30] (b) -- +(0,-1);
\end{tikzpicture}
\hspace*{\fill}
\capt[A concrete slab supported at the ends.]{A concrete beam or slab supported
only at the ends will experience tension below and compression above.}
\label{fig:bentConcrete}
\end{figure}

Prestressed concrete\Index{prestressed concrete}[2!1 ] is a method of dealing
with that problem.  Steel rods or cables are installed under tension in the
concrete.  Their tension puts the concrete in compression (see
\vref{fig:concrete}).  Now, up to a certain load, the concrete remains under
compression, functioning in its ideal realm \cite{wiki:Psc}.

\begin{figure}[ht]
\hspace*{\fill}
\begin{tikzpicture}
\fill[black!40!white] (0,0) -- (7,0) -- (7,2) -- (0,2) -- cycle;
\draw (0,0) -- (0,2) (7,0) -- (7,2);
\foreach \x in {0.5,1,1.5} {
  \draw[thick] (0,\x) -- (7,\x);
  \draw[->] (-1,\x) -- (0,\x);
  \draw[->] (8,\x) -- (7,\x);
}
\foreach \x in {0.25,0.75,1.25,1.75} {
  \draw[->] (0.75,\x) -> (0,\x);
  \draw[->] (6.25,\x) -> (7,\x);
}
\end{tikzpicture}
\hspace*{\fill}
\capt[Forces in a prestressed concrete slab.]{Forces in a prestressed concrete
slab.  The tension in the steel rods pulls inward on the faceplates, while the
compression in the concrete pushes outward, placing them at equilibrium.  As a
tensegrity, the concrete acts like a strut and the rods like a cable, while the
faceplates serve as vertices.  This diagram is of post-tensioned prestressed
concrete, where the tension is added to the rods after the concrete has
hardened.  In pre-tensioning, the concrete is poured around and bonds to rods
that are already in tension.  The forces then act all along the rods instead of
only at the ends.} 
\label{fig:concrete}
\end{figure} 
The prestressed concrete is a tensegrity framework with a
\emph{stress}\Index{stress}, that is, a set of
compressions\Index[|)]{compression} on its struts (the concrete) and
tensions\Index[|)]{tension} on its cables (the rods) that gives a net zero force
at every vertex (faceplate).

A main result in ``Tensegrity Frameworks'' \cite{MR610958} is that if (and
only if) a tensegrity can be stressed in this fashion, then it functions no
longer as if it had cables and struts but rather as if all of its elements
were bars\Index{{strictly positive} stress $\Leftrightarrow$ {bar
equivalent}}[1!2!34 2!1!34 4!312].  Moving this theorem into the realm of
continuous tensegrities is our primary goal.

\section{Mathematical Context}
\label{section:MathContext}
Let's look, for a moment, at a few places where this type of theorem has been
used to good effect.

In recent years, two teams have settled the Carpenter's Rule
Problem\Index{Carpenter's Rule Problem}[ ], which
asks, given a polygonal arc with fixed-length edges in the plane, whether it is
always possible to straighten the arc without it intersecting itself during
the process \cites{CDR,MR2038499}.

Connelly\Index[|(]_,{Robert Connelly}, Demaine\Index_,{Erik D.  Demaine} and
Rote\Index_,{G\"unter@Guenter Rote} \ycite{CDR} not only settle the question in
the affirmative but also answer the related question of whether a closed,
simple polygon with fixed-length edges in the plane can be made convex through
some motion that, again, avoids self intersections.  In fact, they show that
there is a flow of the polygon or arc that accomplishes the purpose and during
which all of the distances between nonadjacent vertices are strictly
increasing.

Their work uses a theorem very much like Roth \& Whiteley's.  They turn the
polygon into a tensegrity and show that if it has no ``weak stress'', then it
has a strictly expansive motion\Index{semipositive stress
$\Leftrightarrow$ partially $\X$-bar@X-bar equivalent equiv.}[1!2!3457
2!1!3457 456!312 56!4!312].  John Pardon\Index_,{John Pardon} \ycite{Pardon}
has made significant progress in extending this to smooth curves, showing that
for a simple, closed, nonconvex curve, there exists a homotopy that takes it to
a convex curve, during which all of the self-distances on the curve are
nondecreasing.

In our second example, Connelly applies Roth \& Whiteley's theorem to provide
an additional proof for the following lemma, which he credits to
Cauchy\Index_,{Augustin Louis Cauchy} \cite{MR652643} (see \ref{fig:ConCau}).
\begin{lemma}
If, in a convex planar or spherical polygon $ABCDEF$, all the sides $AB, BC,
CD, \dotsc, FG$, with the exception of only $AG$, are assumed invariant, one
may increase or decrease simultaneously the angles $B, C, D, E, F$ included
between these same sides; the variable side $AG$ increases in the first case,
and decreases in the second.
\end{lemma}
\begin{proof}
See \cite{MR652643}*{p.\ 30}.
\end{proof}
\begin{figure}[ht]
\hfil
\begin{tikzpicture}
\draw (0,0) node [below left] {$A$} -- 
              node [circle,pos=0.3] (a1) {} node [circle,pos=0.7] (g2) {}
      (3,0) node [below right] {$G$} --
              node [circle,pos=0.3] (g1) {} node [circle,pos=0.7] (f2) {}
      (5,1) node [below right] {$F$} -- 
              node [circle,pos=0.3] (f1) {} node [circle,pos=0.7] (e2) {}
      (6,3) node [right] {$E$} --
              node [circle,pos=0.3] (e1) {} node [circle,pos=0.7] (d2) {}
      (4,4) node [above right] {$D$} -- 
              node [circle,pos=0.3] (d1) {} node [circle,pos=0.7] (c2) {}
      (0,4) node [above left] {$C$} --
              node [circle,pos=0.3] (c1) {} node [circle,pos=0.7] (b2) {}
      (-2,2) node [left] {$B$} -- 
              node [circle,pos=0.3] (b1) {} node [circle,pos=0.7] (a2) {}
      (0,0);
% \draw [<->,yshift=-0.1cm] (0.5,0) -- (2.5,0);
% \draw [<->] (b1) -- (b2);
% \draw [<->] (c1) -- (c2);
% \draw [<->] (d1) -- (d2);
% \draw [<->] (e1) -- (e2);
% \draw [<->] (f1) -- (f2);
\end{tikzpicture}
\hfil
\capt[A lemma of Cauchy reproven by Connelly \ycite{MR652643}.]{Increasing
(resp.  decreasing) angles $B$ through $F$ lengthens (resp. shortens)
side~$AG$.  A lemma of Cauchy reproven by Connelly \ycite{MR652643}.}
\label{fig:ConCau}
\end{figure}

\label{SchursTheorem}
In his discussion, Connelly\Index[|)]_,{Robert Connelly} notes that this is the
polygonal equivalent of Schur's Theorem\Index{Schur Theorem}[2!1], which says
that if we have a $C^2$ convex plane arc and we decrease its curvature at some
place or places and do not increase it anywhere, then the distance between its
ends increases (see, for example, \ocite{doCarmo}*{p.\ 406}).

\label{GrunbaumPolygons}
Our final example comes from Roth \& Whiteley themselves.
Gr\"unbaum\Index_,{Branko Gr\"unbaum@Gruenbaum} and
Shephard\Index_,{Geoffrey C.  Shephard}, in their ``Lectures on Lost
Mathematics'' \cite{LostMath} introduce a class of tensegrities that Roth \&
Whiteley call the ``Gr\"unbaum polygons''\Index{Gr\"unbaum@Gruenbaum polygons}[
].  These are convex, planar (not necessarily regular) polygons made of struts
that have cables joining a select vertex to all nonadjacent vertices and one
further cable joining the two vertices adjacent to the chosen one (see
\vref{fig:GP}).  
\begin{figure}[ht]
\hspace*{\fill}
\begin{tikzpicture}[scale=1.732]
\tikzstyle{every node}=[vertex]
\draw[cable] (0,0) -- (2,2) (2,0) -- (0,2);
\draw[strut] (0,0) node {} -- (2,0) node {} -- (2,2) node {} --
  (0,2) node {} -- cycle;
\end{tikzpicture}
\hspace*{\fill}
\begin{tikzpicture}
\tikzstyle{every node}=[vertex]
\draw[cable]
  (18:2cm) -- (162:2cm)
  (234:2cm) -- (90:2cm) -- (306:2cm);
\foreach \x in {90,162,...,380} {
  \draw[strut] (\x:2cm) node {} -- (\x+72:2cm) node {};
}
\end{tikzpicture}
\hspace*{\fill}
\begin{tikzpicture}
\tikzstyle{every node}=[vertex]
\draw[cable] 
  (1,1.732) -- (-2,0) -- (1,-1.732) 
  (2,0) -- (-2,0) 
  (-1,1.732) -- (-1,-1.732);
\draw[strut] 
  (2,0) node {} -- (1,1.732) node {} -- (-1,1.732) node {} -- 
  (-2,0) node {} -- (-1,-1.732) node {} -- (1,-1.732) node {} -- cycle;
\end{tikzpicture}
\hspace*{\fill}
\capt[The first three regular Gr\"unbaum polygons.]{The first three regular
Gr\"unbaum polygons, introduced in \cite{LostMath}.  Struts are represented by
thick lines, cables by thin ones, vertices by circles.}
\label{fig:GP}
\end{figure}

Roth and Whiteley \ycite{MR610958}*{p.\ 437} show that the Gr\"unbaum polygons
are infinitesimally rigid\Index{infinitesimally rigid}[ ] in~$\RR^2$.  As the
number of edges in the Gr\"unbaum polygons increases, they seem to approach a
smooth convex curve.  So it seems not unreasonable that the limit of such a
family would be \ldots
\begin{dotslist}
\item a convex curve in $\RR^2$ that is not allowed to shrink locally to first
      order (but may stretch) and  \ldots
\item which has a continuous family of upper bounds on the distances
      between one distinguished point on the curves and all the others and
      \ldots
\item some sort of curvature requirement at that point,
\end{dotslist}
and that this limiting curve would be infinitesimally rigid.

Here we have three situations in which Tensegrity Theory has been used to
solve problems in Discrete Geometry, problems that have analogs in
Differential Geometry.  Furthermore Cantarella\Index_,{Jason Cantarella},
Fu\Index{Joseph H. G. Fu}[4,123 ], Kusner\Index_,{Robert B. Kusner},
Sullivan\Index_,{John M. Sullivan} and Wrinkle\Index_,{Nancy Wrinkle}
\ycite{CFKSW} borrow terms and ideas from Tensegrity Theory to answer
questions in Thick Knot Theory, distinctly a Differential Geometry field.  

One can hardly avoid the feeling that these and other problems would benefit
from having tensegrity tools brought from the discrete realm into the
continuous, and we begin that project with this work.

\section{A Road Map}
Here in \ref{chapter:Intro}, we've given a history, both structural and
mathematical, placing this work in context.  In \ref{chapter:Roth} we
build up the background for the main theorems.  \ref{section:RWex} uses
the setting of a specific example tensegrity to introduce the notation and
terms we will need in the rest of the paper.  For example, we here meet
the \emph{variations}\Index{variation} (Definition
\vref{definition:variation}), the possible ways of moving the vertices; the
(infinitesimal) \emph{motions}\Index{motion} (Definition
\vref{definition:motion}), the variations that respect all of the edge
constraints; and the \emph{loads}\Index{load} (Definition
\vref{definition:load}), the effects of variations on the edges of the
tensegrity.  We'll discover that there are \emph{strictly positive}
motions\Index{{strictly positive} motion}, which increase the lengths of all
struts and decrease the lengths of all cables and the \emph{semipositive}
ones\Index{semipositive motion}, which increase the length of at least one
strut or decrease the length of at least one cable.

We also find out about \emph{infinitesimally rigid}\Index{infinitesimally
rigid}[2!1 ] tensegrities, whose only motions are rigid motions of
space\Index{rigid motion of space}[2!1,34 1!234] and about \emph{bar
equivalent}\Index{bar equivalent}[ ] tensegrities, which act as if they were
made entirely of bars, and \emph{partially bar equivalent}\Index{partially bar
equivalent}[ ] tensegrities, which have bar equivalent pieces (Definitions
\vref{definition:IRBE} and \vref{definition:partiallyBE}).  

We will take, as the variations of interest, a subspace, $\X$\Index{$\X$@X
(design variations)}[ ], of the space of variations, called the \emph{design
variations}\Index{design variation}[2!1 ] (Definition \vref{definition:X}).  If
$\X$ doesn't contain all of the variations, we'll use the terms
\emph{$\X$-infinitesimally rigid}\Index{$\X$-infinitesimally@X-infinitesimally
rigid}[2!1 ], \emph{$\X$-bar equivalent}\Index{$\X$-bar@X-bar equivalent}[ ]
and \emph{partially $\X$-bar equivalent}\Index{partially $\X$-bar@X-bar
equivalent}[ ].

Finally, we'll meet the \emph{stresses}\Index{stress} (Definition
\vref{definition:stress}), sets of tensions\Index{tension} on the cables and
compressions\Index{compression} on the struts that, like those in the
prestressed concrete\Index{prestressed concrete}[2!1 ], put the tensegrity in
equilibrium.  These also come in strictly positive\Index{{strictly positive}
stress} (positive values on all the edges) and semipositive\Index{semipositive
stress} (a positive value on at least one edge).

In \ref{section:WorkingRWex} we work the example we have set up in
\ref{section:RWex}.  \ref{section:Mechanics} covers the mathematical background
needed for the proofs in \ref{section:KeyLemma}.  Then in
\ref{section:KeyLemma}, we state and prove Roth \& Whiteley's key lemma, which
we reformulate as follows:
\begin{KL}[p.~\pageref{lemma:NewFiveOne}]
\statementNewFiveOne
\end{KL}

\ref{chapter:Main} holds our main results.  In
\ref{section:ToA} we talk about a class of theorems called ``Theorems of the
Alternative''\Index{Theorem of the Alternative}[1!234].  Our main theorems
resemble these.  In fact, in \ref{appendix:Motzkin}, we will show how to prove
a weaker version of our first main theorem using a Theorem of the Alternative
by Theodore S. Motzkin\Index_,{Theodore S. Motzkin}\Index{Motzkin
Transposition Theorem}[3!12 ].

Next we launch into the task of moving the things we have defined in
\ref{chapter:Roth} into the continuous realm.  \ref{section:E} builds a metric
and a topology on the vertex set $\V$ and the \emph{edgeset}\Index{$\E$@E
(edgeset)}[ ]\Index[|see{$\E$}]{edge set}[2!1]\Index[|see{$\E$}]{edgeset}~$\E$.
\ref{section:CE} takes on $C(\E)$ and $C^*(\E)$ (the set of continuous
functions on $\E$ and its topological dual), giving them norms and topologies
and then exploring their natures in some depth.  In \ref{section:VfVY} we
address $\Vf(\V)$ and~$Y$.  

Then, in \ref{section:FMT} we reach our first main theorem:
\begin{MTone}[p.~\pageref{theorem:nonnegativeMeasure}]
\statementNonnegativeMeasure
\end{MTone}
We explore the theorem with a couple of examples.  Then we move, in 
\ref{section:SMT}, to our second main theorem:
\begin{MTtwo}[p.~\pageref{theorem:positiveImplies}]
\statementPositiveImplies
\end{MTtwo}
But we only have one direction of the desired theorem and we do some looking
at why the other direction is difficult.

\ref{section:TwoExamples} is an excursion into example land.  We take
on a couple of examples of tensegrities in which $\X$ is not the whole
of~$\Vf(\V)$.  

We discover that the example called the ``circle of struts'' is bar equivalent.
\Index{circle of struts}[ ]
In this example, a circle of vertices is filled with antipodal struts (see
\ref{fig:COS}) and we only allow variations which are local isometries on
the curve.
\begin{figure}[ht]
\hspace*{\fill}
\begin{tikzpicture}
\foreach \x in {0,20,...,181} {
  \draw[strut] (\x:2cm) -- (\x+180:2cm);
}
\draw[vcurve] (0,0) circle (2cm);
\end{tikzpicture}
\hspace*{\fill}
\capt{The circle-of-struts example of \ref{section:TwoExamples}.}
\label{fig:COS}
\end{figure}
When we apply the first main theorem to this example, we get the following
proposition:
\begin{COSone}[p.~\pageref{prop:COSone}]
\statementCOSone
\end{COSone}
When we apply the second main theorem, that becomes:
\begin{COStwo}[p.~\pageref{prop:COStwo}]
\statementCOStwo
\end{COStwo}
The circle acts like a water balloon.  Pulling it out somewhere results in it
being pulled in elsewhere.

In \ref{section:MBE} we again attempt to provide the ``other half'' of the
second main theorem.  By restricting our attention to tensegrities that are
$\X$-bar equivalent and that cease to be so upon the removal of any open set
of edges, we get
\begin{MTthree}[p.~\pageref{theorem:minimallyBarEquivalent}]
\statementMinimallyBarEquivalent
\end{MTthree}
After proving the theorem, we talk about what minimally $\X$-bar equivalent
tensegrities are like.  Then, in \ref{section:Covered}, we show that if a
tensegrity is (at least mostly) covered by subtensegrities which have strictly
positive stresses, then the tensegrity itself has a strictly positive stress:
\begin{MTfour}[p.~\pageref{theorem:CountablyCovered}]
\statementCountablyCovered
\end{MTfour}
We end by conjecturing that either
\begin{ConjOne}[p.~\pageref{conj:CountCov}]
\statementCountCov
\end{ConjOne}
or at least
\begin{ConjTwo}[p.~\pageref{conj:StrictlyPositive}]
\statementStrictlyPositive
\end{ConjTwo}

\ref{chapter:Examples} is full of examples of continuous tensegrities.  In
\ref{section:rectangle} we analyze a tensegrity that is not only bar equivalent
but also infinitesimally rigid.  \ref{section:onACircle} through
\ref{section:CrossedMore} are a theme and variations.  First we start by taking
a circle of vertices and setting struts antipodally on it and cables connecting
vertices a fixed skip around.  That works so well that we move up one dimension
and build the equivalent tensegrity on the sphere, thus providing ourselves
with an example of a tensegrity whose vertex set is two-dimensional.
Then we try to relax the circle into other vertex curves, with varying amounts
of success.  In particular, \ref{section:CrossedMore} shows that bar
equivalence can survive affine transformations of the tensegrity.

Our last example (\ref*{section:corner}) is designed to remind us that the
intuition we get about subspaces and closed convex cones from thinking about
$\RR^2$ and $\RR^3$ may be faulty and to make a little clearer why showing that
every bar equivalent tensegrity has a strictly positive stress is so difficult.

Finally, in \ref{chapter:Future}, we wrap things up.  We remind the reader of
the situations we mentioned above where continuous tensegrities might bring
some of the successes of Discrete Geometry into Differential Geometry and we
also point out some other places for further exploration.  We close with
two appendices.  \ref{appendix:rigidMotions} covers rigid motions of space and
the variations which are induced from them, and \ref{appendix:Motzkin} proves
a special-case version of the first main theorem using Motzkin's Transposition
Theorem.
\begin{figure}[T]
\setlength{\unitlength}{0.5in}
\begin{picture}(11,16.5)(-0.5,0)
\put(-0.5,1){\makebox(0,0)[bl]{
  \begin{overpic}{4way1963.epsf}
  \end{overpic}
}}
\put(4.5,17){\makebox(0,0)[tl]{
  \begin{overpic}{e_c_col4x5albrknox.epsf}
  \end{overpic} 
}}
\put(6.5,8.2){\makebox(0,0){
  \begin{overpic}{xtend.epsf}
  \end{overpic}
}}
\put(0,16){\makebox(0,0)[tl]{
  \begin{overpic}{octahelix.epsf}
  \end{overpic}
}}
\put(6,0){\makebox(0,0)[bl]{
  \begin{overpic}{xplanar.epsf}
  \end{overpic}
}}
\put(4.3,5){\makebox(0,0){
  \begin{overpic}{easy_ksonsbk71.epsf}
  \end{overpic}
}}
\put(4.3,2){\makebox(0,0){
  \begin{overpic}{double_sqs_lt_dk_bmboo.epsf}
  \end{overpic}
}}
\end{picture}
\capt[Works of Kenneth Snelson.]{A collection of pictures of the works of
Kenneth D. Snelson\Index_,{Kenneth D. Snelson} from his website
\url{www.kennethsnelson.net}, used by permission.}
\label{fig:SnelsonArt}
\end{figure}
\end{chapter}

\chapter{The Roth-Whiteley Theorem}
\label{chapter:Roth}
\Index[|(]_,{Ben Roth}
\Index[|(]{Walter Whiteley}[2,1]
\section{Setting up a Roth \& Whiteley Example}
\label{section:RWex}
We need to start by understanding Roth \& Whiteley's original theorem.  To that
end, we'll first look at an example, and then work through their proof.
However, we have another task to accomplish as well.  There are numerous
symbols and terms we'll want to have available when we get into the thick of
the mathematics.  As we work through the example, we'll define the terms and
symbols, relating them to the example.  Afterward we'll give a summary list of
the symbols we have defined.

Let's get a little general notation out of the way first.  In what follows,
$\zero$\Index{$\zero$@zero vs. $0$@zero}[ ] will denote the origin of the
various vector spaces we encounter.  The symbol $0$ will be reserved for the
real additive identity.  $\langle \cdot, \cdot \rangle$\IndexDef{{$\langle
\cdot, \cdot \rangle$}@*01 (inner product)}[ ] will denote the inner product
for any vector space.  If the space is Euclidean $n$-space (denoted
$\RR^n$\IndexDef{$\RR^n$@Rn (Euclidean $n$-space)}[ ]), we may use dot product
notation (e.g. $v \cdot w$)\IndexDef{{$v \cdot w$}@*02 (dot product)}[ ] for
the same purpose.

If $X$ and $Y$ are two sets, then $X \setminus Y$ will denote the set of
elements that are in $X$ but not in~$Y$.  That is,
$$X \setminus Y \coloneqq \{x \in X : x \notin Y\}.$$\IndexDef{{$X \setminus
Y$}@*03 (set-theoretic difference)}[ ]

We'll also want some idea of positivity for vectors and vector fields
in~$\RR^n$.
\begin{definition}
A vector $v \in \RR^n$ is said to be \emph{strictly
positive}\IndexDef{{strictly positive} vector}, written $v > \zero$, if every
coordinate of $v$ is strictly positive.  $v$ is said to be
\emph{nonnegative}\IndexDef{nonnegative vector}, written $v \ge \zero$, if
every coordinate of $v$ is nonnegative.  Finally, $v$ is said to be
\emph{semipositive}\IndexDef{semipositive vector}, written $v \gneqq \zero$ if
$v \ge \zero$ and $v \ne \zero$.

A vector field $V$ is termed \emph{strictly positive}\IndexDef{{strictly
positive} {vector field}} (resp.\ \emph{nonnegative}\IndexDef{nonnegative
{vector field}}), also written $V > \zero$ (resp.\ $V \ge \zero$), if every
value of $V$ is strictly positive (resp.\ nonnegative) as a vector.  $V$ is
termed \emph{semipositive}\IndexDef{semipositive {vector field}} (also written
$V \gneqq \zero$) if $V \ge \zero$ and $V$ is not the ``zero vector field'',
that is, $V \ne \zero$.
\end{definition}

Now we can turn our attention to the example.
\ref{fig:RWexample}\subref{subfig:csq} shows a simple
tensegrity.\IndDefBeg{crossed square}[ ]
\begin{figure}[ht]
\hfil
\subfloat[The crossed square.]{
\begin{tikzpicture}
\useasboundingbox (-2,0) rectangle (4,2);
\draw[cable] (0,0) -- (2,0) -- (2,2) -- (0,2) -- cycle;
\draw[strut] 
  (0,0) node [vertex] {1} -- 
  (2,2) node [vertex] {3}
  (0,2) node [vertex] {4} -- 
  (2,0) node [vertex] {2};
\end{tikzpicture}
\label{subfig:csq}
}
\hfil
\subfloat[The crossed square's bar framework.]{
\begin{tikzpicture}
\useasboundingbox (-2,0) rectangle (4,2);
\draw[bar] (0,0) -- (2,0) -- (2,2) -- (0,2) -- cycle;
\draw[bar] 
  (0,0) node [vertex] {1} -- 
  (2,2) node [vertex] {3}
  (0,2) node [vertex] {4} -- 
  (2,0) node [vertex] {2};
\end{tikzpicture}
\label{subfig:csqb}
}
\hfil
\capt[The crossed-square example with its bar framework.]{The crossed-square
example with its bar framework.  Note that the two struts (thick lines in
\subref{subfig:csq}) and the two center bars (double lines crossing in 
\subref{subfig:csqb}) pass through each other without touching.}
\label{fig:RWexample}
\end{figure}
In our illustrations, struts will be shown as thick lines, cables as thin ones
and bars as double lines.  Vertices may be shown as circles, as in
\ref{fig:RWexample}, or as dotted lines, as in \fullref{fig:COS}.  
\Index{{\mbox{\tikz{\draw[bar] (0,0) -- (1cm,1ex);}}}@*10 (bar)}[ ]
\Index{{\mbox{\tikz{\draw[cable] (0,0) -- (1cm,1ex);}}}@*11 (cable)}[ ]
\Index{{\mbox{\tikz{\draw[strut] (0,0) -- (1cm,1ex);}}}@*12 (strut)}[ ]
\Index{{\mbox{\tikz{\draw (0,0) node[vertex] {};}}}@*13 (vertex)}[ ]
\Index{{\mbox{\tikz{\draw[vcurve] (0,0) -- (1cm,1ex);}}}@*14 (vertices)}[ ]

Using the notation introduced on page \pageref{notation}, we see that the
crossed square has
\begin{align*}
\V & = \{1,2,3,4\} &
\St & = \{\{1,3\},\{2,4\}\} \\
\B & = \varnothing &
\C & = \{\{1,2\},\{2,3\},\{3,4\},\{1,4\}\}
\end{align*}
We need to place the elements of $\V$ in space, so we define $p\mcol \V \to
\RR^2$ by 
$$p(v) = \begin{cases}
(0,0), & v = 1 \\
(1,0), & v = 2 \\
(1,1), & v = 3 \\
(0,1), & v = 4
\end{cases}$$
We'll call the tensegrity $G(p) = \{\V; \St, \C, \B; p\}$\Index{$G(p)$@G(p)
(tensegrity)}[ ].  \ref{fig:RWexample}\subref{subfig:csqb} shows the the same
tensegrity but with all the struts and cables replaced by bars.  The bar
framework constructed from tensegrity $G(p)$ in this fashion will be called
$\Gb(p)$\Index{$\Gb(p)$@G-(p) ($G(p)$'s bar framework)}[ ].

Bars are a problem.  They give us equations instead of inequalities.  One might
argue that we prefer equations to inequalities and perhaps that's true.
However, in general, we can't turn inequalities into equalities.  But we
\emph{can} turn an equation into a pair of inequalities by thinking of
a bar as a strut-cable pair\Index{bar {as strut-cable pair}}[1!2].  That way
the cable in the pair keeps the vertices from moving apart and the strut keeps
them from moving together.

To make it clear that we are thinking of bars no longer as bars but as
strut-cable pairs, we'll define the sets of ``functional struts'' and
``functional cables'':\IndexDef{functional struts and cables}[ ]
\begin{align*}
\Str & \coloneqq \St \cup \B & \Cr & \coloneqq \C \cup \B
\end{align*}
and deliberately refer to the elements of $\Str$ as ``struts'' and the elements
of $\Cr$ as ``cables''.\IndexDef{$\Str$@S(1) (functional struts, $\St \cup
\B$)}[ ]\Index{$\St$@S(0) (struts)}[ ]\IndexDef{$\Cr$@C(1) (functional cables,
$\C \cup \B$)}[ ]\Index{$\C$@C(0) (cables)}[ ]
We're ready now to give a name to the set of edges, which will contain no
bars, only (functional) struts and cables:
\begin{definition}
\label{definition:E}
\IndexDef{$\E$@E (edgeset)}[ ]
The \emph{edgeset} $\E$ is defined as a set by
$$\E = \Str \sqcup \Cr$$
(where $\sqcup$\Index{$\sqcup$@*17 (disjoint union)}[ ] indicates disjoint
union as opposed to the $\cup$\Index{$\cup$@*16 (union)}[ ] indicating union
above).  $\E$ is endowed with a topology which we will describe in
\vref{section:E} and required to be compact.
\end{definition}

Here is a quick preview of how we make $\V$ and $\E$ into topological spaces.
We induce a metric on $\V$ from $\RR^n$ via~$p$.  We then take the maximum of
two distances in $\V$ to give us a distance in $\V \cross \V$.  We use that
metric on $\V \cross \V$ to create a metric on the identification space $\vvt$
of unordered pairs of elements of~$\V$.  Finally we create a metric on $\E$
from the metric on $\vvt$.  

Now that we have metrics for $\V$ and $\E$, we put the associated metric
topologies on them.  These turn out to be familiar topologies
(\vref*{subsection:familiar}) and when the two sets are finite, both of the
topologies are discrete\Index{discrete topology}[2!1 ] (Proposition \vref{prop:discrete}).  Of course, in that
case both sets are compact simply because they are finite, but in the general
case we require $\E$ to be compact and discover that $\V$ might as well be, too
(\vref*{subsection:compactV}).

The possible motions of $G(p)$ can be described by putting vectors on the
vertices.  We'll call the space of all continuous vector fields on the vertices
$\Vf(\V)$ and note that in this case it is~$\RR^8$.  Of course, only some of
the elements of $\Vf(\V)$ respect the edge constraints.  We'll reserve the term
``motion'' to apply to those.  If we want to talk generally about the elements
of $\Vf(\V)$, we'll refer to them as ``variations'' or just call them vector
fields.
\begin{definition}
\label{definition:variation}
The \emph{variations}\IndexDef{variation} of $G(p)$ are the continuous vector
fields $V\mcol \V \to \RR^n$ and the space of variations is denoted
$\Vf(\V)$\IndexDef{$\Vf(\V)$@Vf(V) (variations)}[ ].  
\end{definition}
  
Now the variations are of different kinds.  Some of the vector fields in
$\Vf(\V)$ can be induced from rigid motions\Index[|see{Euclidean motion}]{rigid
motion of $\RR^n$}[2!1,34 1!234] of~$\RR^n$.
\label{rigidMotionMention}
They do things like translating $G(p)$ or rotating it, but they do not change
its shape (see \fullref{appendix:rigidMotions} for more on rigid motions and
the vector fields that are induced by them).
\begin{definition}
$T(p)$\IndexDef{$T(p)$@T(p) (Euclidean motions)}[ ] is the set of all elements
of $\Vf(\V)$ that can be induced from rigid motions of~$\RR^n$.  We'll term
these \emph{Euclidean motions}\IndexDef{Euclidean motion}.
\end{definition}

We need to know how a given vector field interacts with the constraints given
by the edges.  If the positions of the vertices vary differentiably with time
and $\frac{\partial}{\partial t} p(v,t) = V(v)$ for a given vector field $V$,
then $V$ yields a vector of real numbers, which are the values of the
derivatives of the lengths of the edges.  

Here we have 6 edges, so that vector lies in $\RR^6$, but since $\E$ has the
discrete topology on it, we can (and will want later to) call it a ``continuous
function on $\E$'' and the space it inhabits $C(\E)$\IndexDef{$C(\E)$@C(E)
(cont. func. on $\E$)}[ ].  We'll call the functions we can generate this way
\emph{loads}\Index{load}.  In a moment we'll want to talk about the nonnegative
orthant of~$C(\E)$.  We'll name that $C(\E)^+$\IndexDef{$C(\E)^+$@C(E)+ and
$C(\E)^-$@C(E)- (orthants of $C(\E)$)}[ ], and the nonpositive orthant we'll
call $C(\E)^-$.  

Formally, loads are computed via the the \emph{rigidity
operator}\IndexDef{rigidity operator}[ ]:
\begin{definition}
\label{definition:Y}
\IndexDef{$Y$@Y (rigidity operator)}[ ]
$Y \mcol  \Vf(\V) \to C(\E)$ is given by
\begin{equation}
(YV)(\{v_1,v_2\}) \coloneqq \begin{cases}
(V(v_1) - V(v_2)) \cdot (p(v_1) - p(v_2)), & \{v_1,v_2\} \in \Str \\
-(V(v_1) - V(v_2)) \cdot (p(v_1) - p(v_2)), & \{v_1,v_2\} \in \Cr.
\end{cases}
\end{equation}
\end{definition}

The signs here were chosen so that $YV > \zero$ when $V$ expands struts and
contracts cables.  It would certainly be possible to do the work that follows
using a $Y$ which is defined the same for all parts of the edgeset.  But we
would do so at the cost of working in an orthant of $C(\E)$ that has no natural
name and of repeatedly saying ``where $YV(e) \ge 0$ on struts and $YV(e) \le 0$
on cables''.

In some cases we will not want to consider all possible vector fields.  For
example, we may have a curve of vertices and be interested only in variations
that are local isometries of the curve.   We will call the variations of 
interest \emph{design variations}.
\begin{definition}
\label{definition:X}
The \emph{design variations}\IndexDef{design variation}[2!1 ] form a subspace
$\X \subset \Vf(\V)$\IndexDef{$\X$@X (design variations)}[ ] with the one
requirement that the Euclidean motions be included in the design variations.
That is, $T(p) \subset \X$\Index{$T(p)$@T(p) $\subset$@subset $\supset$@supset
$\X$@X {(Euclidean motions)} {(design variations)}}[15!24 46!31].
\end{definition}

With $Y$ defined as above, for $V \in \X$, the statement $YV \in C(\E)^+$ means
that $V$ is a variation on $\V$ that respects all of its constraints.  As these
are variations, they are really infinitesimal
motions\Index[|see{motion}]{infinitesimal motion}[ ], but by a slight abuse of
the language, we'll call them simply \emph{motions}.
\begin{definition}
\label{definition:motion}
A \emph{motion}\IndexDef{motion} is a variation $V \in \X$ such that $YV \in
C(\E)^+$.  The set of all motions is denoted $I(p)$\IndexDef{$I(p)$@I(p)
(motions)}[ ].  A motion is termed \emph{strictly positive}\IndexDef{{strictly
positive} motion} or \emph{semipositive}\IndexDef{semipositive motion} as $YV$
is strictly positive or semipositive.
\end{definition}

So a strictly positive motion would increase the lengths of all struts and
decrease the lengths of all cables, while a semipositive one would change the
length of at least one edge.

And now we're prepared to formally define the term ``load''.
\begin{definition}
\label{definition:load}
The \emph{loads}\IndexDef{load} for a given tensegrity are the elements of the
set $Y(\X) \coloneqq \{YV : V \in \X\}$\IndexDef{$Y(\X)$@Y(X) (loads)}[ ].  A
load is \emph{strictly positive}\IndexDef{{strictly positive} load} if $YV >
\zero$ and \emph{semipositive}\IndexDef{semipositive load} if $YV \gneqq
\zero$.
\end{definition}

A quick note on sign.  This definition of ``load'' is related to the idea of
loads in Engineering, but it differs in sign.  For example, an engineer would
consider the edges in \vref{fig:truck} to be bearing (positive) loads.  We 
agree that they are bearing loads, but, for historical reasons, give them a
negative sign.
\begin{figure}[ht]
\hfil
\begin{tikzpicture}
\begin{scope}[yshift=1.732cm,scale=0.2]
  \draw (-1,0.5) circle (0.5cm)
        (1,0.5) circle (0.5cm);
  \draw (-1.6,0.5) arc (180:0:0.6cm) -- (0.4,0.5) arc (180:0:0.6cm) -- 
        (2,0.5) -- (2,1.3) -- (-0.4,1.3) -- (-0.4,2.3) -- (-1.6,2.3) --
        (-1.7,1.4) -- (-2.3,1.3) -- (-2.3,0.5) -- cycle;
\end{scope}
\draw[cable] (-1,0) -- (1,0);
\draw[strut] (-1,0) node[vertex] { } -- (0,1.732) node[vertex] { } --
             (1,0) node[vertex] { };
\draw[->] (0,1.732) -- (0,1);
\draw[->] (-1,0) -- (-1.4,0);
\draw[->] (1,0) -- (1.4,0);
\end{tikzpicture}
\hfil
\capt{Edges bearing negative loads.}
\label{fig:truck}
\end{figure}

Some tensegrities possess motions that, while not being Euclidean motions, 
still don't change the length of any edges (see \vref{fig:twoBars} for an
example of such a tensegrity).  Since such a motion changes no edge lengths,
it respects the constraints not only of $G(p)$ but also of~$\Gb(p)$.  

So, similar to $I(p)$, we'll define $\Ib(p)$\IndexDef{$\Ib(p)$@I-(p) (motions
of $\Gb(p)$)}[ ] to consist of the vector fields in $\X$ that respect the
constraints of $\Gb(p)$.  Because the constraints of $\Gb(p)$ are all
equalities, any element $V \in \Ib(p)$ must give us $YV = \zero$.\footnote{In
passing, we note that since $Y$ is linear in $V$, any positive multiple of an
element of $I(p)$ still maps by $Y$ into $C(\E)^+$ and the same is true of
\emph{convex} combinations of elements of $I(p)$, so $I(p)$ is a convex cone.  

In contrast, all \emph{linear} combinations of elements of $\Ib(p)$ must still
map to $\zero$, so $\Ib(p)$ is a subspace of~$\Vf(\V)$.} 
Clearly, $T(p) \subset \Ib(p) \subset I(p)$\Index{$T(p)$@T(p) $\subset$@subset
$\supset$@supset $\Ib(p)$@I-(p) $I(p)$@I(p) {(Euclidean motions)} {(motions)}
{(motions of $\Gb(p)$)}}[16!2425 57!3431 48!25 48!31], but the reverse
inclusions are not guaranteed, as \ref{fig:simpTens} and \ref{fig:twoBars}
demonstrate.
\begin{figure}[ht]
\hspace*{\fill}
\begin{tikzpicture}
\draw[strut] (0,0) node [vertex] {} -- (2,0) node [vertex] {};
\draw[->] (0,0) -- (-0.5,0);
\draw[->] (2,0) -- (2.5,0);
\draw[bar] (4,0) node [vertex] {} -- (6,0) node [vertex] {};
\end{tikzpicture}
\hspace*{\fill}
\capt[Evidence that $I(p) \ne \Ib(p)$ in general.]{Evidence that $I(p) \ne
\Ib(p)$ in general.  The tensegrity on the left has a strut between the
vertices, which allows them to move away from each other.  The associated bar
framework on the right does not allow that motion.\Index{$I(p)$@I(p)
$\ne$@notequal $\Ib(p)$@I-(p) {(motions)} {(motions of $\Gb(p)$)}}[14!23
35!21]}
\label{fig:simpTens}
\end{figure}
\begin{figure}[ht]
\hspace*{\fill}
\begin{tikzpicture}
\draw[bar] (0,0) node [vertex] { } -- (2,1) node [vertex] { } --
  (4,0) node [vertex] { };
\draw[->] (0,0) -- (-0.75,0.5);
\draw[->] (2,1) -- (2,0);
\draw[->] (4,0) -- (4.75,0.5);
\end{tikzpicture}
\hspace*{\fill}
\capt[A two-bar tensegrity with a motion that is in $\Ib(p)$ but not
in~$T(p)$.]{A two-bar tensegrity with a motion that is in $\Ib(p)$ but not
in~$T(p)$.  Since this tensegrity is made only of bars, $I(p) = \Ib(p)$
naturally.\Index{$\Ib(p)$@I-(p) $\ne$@notequal $T(p)$@T(p) {(Euclidean
motions)} {(motions of $\Gb(p)$)}}[15!23 34!21]}
\label{fig:twoBars}
\end{figure}

Now, for some tensegrities, we \emph{do} have that $I(p) = \Ib(p)$ or $\Ib(p) =
T(p)$ or even $I(p) = T(p)$.  If $I(p) = T(p)$, then the only motions of $G(p)$
are the Euclidean motions (and similarly for $\Ib(p) =
T(p)$).  There are no motions that change the shape of $G(p)$ (or $\Gb(p)$).
On the other hand, if $I(p) = \Ib(p)$, then the only motions of $G(p)$ are also
motions of~$\Gb(p)$.  That is, the tensegrity acts as if it were made only of
bars.
\begin{definition}
\label{definition:IRBE}
If $I(p) = T(p)$ (resp.\ $\Ib(p) = T(p)$), then we say that $G(p)$ (resp.\
$\Gb(p)$) is \emph{infinitesimally rigid with respect to
$\X$}\Index[|see{$\X$-infinitesimally rigid}]{{infinitesimally rigid} {with
respect to $\X$}}[1!2] or \emph{$\X$-infinitesimally
rigid}\Index{$\X$-infinitesimally@X-infinitesimally rigid}[ ].

If $I(p) = \Ib(p)$, we call $G(p)$ \emph{bar equivalent with respect to
$\X$}\Index[|see{$\X$-bar equivalent}]{{bar equivalent} {with respect to
$\X$}}[1!2] or \emph{$\X$-bar equivalent}\Index{$\X$-bar@X-bar equivalent}[ ].

In the case where $\X=\Vf(\V)$ (as it is in Roth \& Whiteley's work), we
can drop the $\X$ and use the terms \emph{infinitesimally
rigid}\Index{infinitesimally rigid}[2!1 ] and \emph{bar equivalent}\Index{bar
equivalent}[ ].
\end{definition}

This leads to a quick proposition.
\begin{prop}
\label{prop:IRBE}
\Index{$\X$-infinitesimally@X-infinitesimally rigid $\Rightarrow$ $\Leftarrow$
$\X$-bar@X-bar equivalent}[12!356 56!412]
Any tensegrity that is infinitesimally rigid with respect to $\X$ is bar
equivalent with respect to~$\X$.
\end{prop}
\begin{proof}
$\X$-infinitesimal rigidity means that $I(p) = T(p)$, but it is always true 
that $T(p) \subset \Ib(p) \subset I(p)$, so we must have $I(p) = \Ib(p)$, that
is, that $G(p)$ is $\X$-bar equivalent.
\end{proof}

There are three seemingly different ways of saying ``bar equivalent''.  They
really all mean the same thing, however.
\begin{prop}
Let $G(p)$ be a tensegrity with motions $I(p)$, bar framework $\Gb(p)$ and
$\Ib(p)$ the motions of~$\Gb(p)$.  Then the
following are equivalent:
\begin{enumerate}
\item $G(p)$ has no semipositive motions.
\item $Y(\X) \cap C(\E)^+ = \{ \zero \}$.
\item $I(p) = \Ib(p)$.
\end{enumerate}
\end{prop}
\begin{proof}
$(1 \Leftrightarrow 2)$. Let $V \in \X$.  Then $V$ is a semipositive motion
if and only if $YV \in C(\E)^+$ and $YV \ne \zero$.

$(2 \Leftrightarrow 3)$. $YV \in Y(\X) \cap C(\E)^+$ and $YV \ne \zero$
if and only if $V \in I(p)$ and $V \notin \Ib(p)$.
\end{proof}

There is notion which is similar to, though weaker than, ``bar equivalence''.
Some tensegrities have semipositive motions, but no strictly positive ones
(that is, $Y(\X) \cap \INT C(\E)^+ = \varnothing$, see Lemma
\vref{lemma:intCEp}).  
\begin{definition}
\label{definition:partiallyBE}
If $G(p)$ has semipositive motions but no strictly positive motions, we call
$G(p)$ \emph{partially $\X$-bar equivalent}\IndexDef{partially $\X$-bar@X-bar
equivalent}[123 23!1].
\end{definition}
In Corollary \vref{corollary:partiallyBE} we show that if $G(p)$ is partially
$\X$-bar equivalent, then some portion of $G(p)$ is $\X$-bar equivalent.

There's one more concept we'll need and that is the concept of
stress\Index{stress}.  Remembering the prestressed concrete\Index_{prestressed
concrete} of \vref{section:prestressed}, we can see that this is a set of
tensions\Index{tension} on the cables and compressions\Index{compression} on
the struts that give a net zero force at each vertex.  That works well for
discrete tensegrities, but we'll need a slightly more general definition in the
future, so we'll define stresses to be certain elements of
$C^*(\E)$\IndexDef{$C^*(\E)$@C*(E) (top.  dual of $C(\E)$)}, the topological
dual of~$C(\E)$.
\begin{definition}
\label{definition:topologicalDual}
If $X$ is a vector space, then its \emph{topological dual}\IndexDef{topological
dual}[2!1 ],\footnote{This dual is ``topological'' in that it uses continuity,
which requires a topology, in its definition.  The \emph{algebraic
dual}\IndexDef{algebraic dual}[2!1 ], by contrast, is the set of all linear
functionals on~$X$.} denoted $X^*$, is the space of all continuous linear
functionals on~$X$. 
\end{definition}
The nonnegative orthant of $C^*(\E)$ will play a role later, so we'll give
it the name~$C^*(\E)^+$\IndexDef{$C^*(\E)^+$@C*(E)+ (nonneg. orth. of
$C^*(\E)$)}[ ].
\begin{definition}
If $S$ is a set of elements in the vector space $X$, then the
\emph{annihilator}\IndexDef{annihilator} of $S$ is the set 
$$S^\perp = \{\mu \in X^* : \mu(s) = 0, \forall s \in S\}.$$
\IndexDef{$S^\perp$@*18 (annihilator of $S$)}[ ]
\end{definition}
\begin{definition}
If $S$ is a set of elements in the vector space $X$, then the
\emph{dual cone}\IndexDef{dual cone}[2,1 ] of $S$ is the set
$$S^* = \{\mu \in X^* : \mu(s) \ge 0, \forall s \in S\}.$$
\IndexDef{$S^*$@*19 (dual cone of $S$)}[ ]
\end{definition}
\begin{definition}
\label{definition:stress}
A \emph{stress}\IndexDef{stress} is an element $\mu \in C^*(\E)$ such that $\mu
\in Y(\X)^\perp$.  A \emph{strictly positive stress}\IndexDef{{strictly
positive} stress} is a stress $\mu$ with $\mu(f) > 0$ for all nonzero $f \in
C(\E)^+$.  A \emph{semipositive stress}\IndexDef{semipositive stress} is a
nonzero stress $\mu$ with $\mu(f) \ge 0$ for all $f \in C(\E)^+$.
\end{definition}
This choice of definition for ``strictly positive'' and ``semipositive'' is
explained more fully by Proposition \vref{prop:Positive}.

Let's take a moment to see that this definition of stress matches our
intuitive understanding in the finite case.
\begin{prop}
\label{prop:stress}
If $G(p)$ is a finite tensegrity with rigidity operator $Y$ and $\X =
\Vf(\V)$, and $\mu \in C^*(\E)$, then $\mu$ is a stress (in the sense of
Definition \ref{definition:stress}) if and only if $Y^\top \mu = \zero$.
\end{prop}
\begin{proof}
Suppose, first, that $\mu$ is a stress.  That is, $\mu(YV) = 0$ for all $V \in
\X = \Vf(\V)$.  Then $(Y^\top \mu) V = 0$ for all~$V$.  But for $Y^\top \mu$ to
have zero dot product with all $V \in \Vf(\V)$, it must be the zero vector.

Conversely, suppose $Y^\top \mu = \zero$.  Then $\mu(YV) = (Y^\top \mu)V =
\zero \cdot V = 0$.
\end{proof}

With our definitions in place, we can say concisely what Roth and Whiteley
proved.  The essence of their result is that having a strictly positive stress
is the same as being bar equivalent (see Lemma
\vref{lemma:NewFiveOne})\Index{{strictly positive} stress $\Leftrightarrow$
{bar equivalent}}[1!2!34 4!312].

\begin{definition} Here is a list of all the entities we have defined, in 
alphabetical order.
\label{def:tensegrity}
\begin{tabularx}{\linewidth}{|l|>{\setlength{\baselineskip}{0.5\baselineskip}}X|} \hline
\endfirsthead
\hline
\multicolumn{2}{|l|}{\small\sl continued from the previous page}\\
\hline
\endhead

\multicolumn{2}{|r|}{\small\sl continued on the next page}\\
\hline
\endfoot
\endlastfoot
$\B$\IndexDef{$\B$@B (bars)}[ ], $\C$\IndexDef{$\C$@C(0) (cables)}[ ], and
$\St$\IndexDef{$\St$@S(0) (struts)}[ ] & 
  Three pairwise-disjoint sets of two-element subsets of~$\V$.  Elements of
  $\B$ (\emph{bars}) require their vertices to stay a fixed distance apart.
  Elements of $\St$ (\emph{struts}) and $\C$ (\emph{cables}) serve as lower and
  upper bounds on their vertices' distances respectively.  \\\hline
$\Cr$ and $\Str$\IndexDef{$\Cr$@C(1) (functional cables, $\C \cup \B$)}[
  ]\IndexDef{$\Str$@S(1) (functional struts, $\St \cup \B$)}[ ] & $\Cr
  \coloneqq \C \cup \B$ and $\Str \coloneqq \St \cup \B$. \\\hline
$C(\E)$\IndexDef{$C(\E)$@C(E) (cont. func. on $\E$)}[ ] & 
  Continuous functions from $\E$ to~$\RR$.  The \emph{loads}, $Y(\X)$, are a
  subspace of~$C(\E)$. \\\hline
$C(\E)^+$\IndexDef{$C(\E)^+$@C(E)+ and $C(\E)^-$@C(E)- (orthants of $C(\E)$)}[
  ] & $\{f \in C(\E) : f(e) \ge 0, \forall e \in \E\}$ \\\hline
$C(\E)^-$ & $\{f \in C(\E) : f(e) \le 0, \forall e \in \E\}$ \\\hline
$C^*(\E)$\IndexDef{$C^*(\E)$@C*(E) (top. dual of $C(\E)$)}[ ] & 
  Continuous linear functionals on $C(\E)$, tensions\Index{tension} and
  compressions\Index{compression} on the edges. \\\hline
$C^*(\E)^+$\IndexDef{$C^*(\E)^+$@C*(E)+ (nonneg. orth. of $C^*(\E)$)}[ ] & 
  $\{\mu \in C^*(\E) : \mu(f) \ge 0, \forall f \in C(\E)^+\}$.\\\hline
$\E$\IndexDef{$\E$@E (edgeset)}[ ] & 
  $\Str \sqcup \Cr$.  The edgeset.  It receives its topology
  from $\RR^n$ via $\V$, $\V \cross \V$ and $\vvt$ (see \vref{section:E}) and
  is required to be compact. \\\hline
$G(p)$\IndexDef{$G(p)$@G(p) (tensegrity)}[ ] & 
  $\{\V;\St,\C,\B;p\}$.   A tensegrity. \\\hline
$\Gb(p)$\IndexDef{$\Gb(p)$@G-(p) ($G(p)$'s bar framework)}[ ] & 
  $\{\V;\Bb;p\}$.  The bar framework created by replacing all of $G(p)$'s edges
  with bars.  That is, elements of $\Bb$\IndexDef{$\Bb$@B- (bars in
  $\Gb(p)$)}[ ] are bars and $\Bb = \St \cup \C \cup \B$. \\\hline
$I(p)$\IndexDef{$I(p)$@I(p) (motions)}[ ] & 
  $\{V \in \Vf(\V) : YV \in C(\E)^+\}$.  Motions of~$G(p)$.  \\\hline
$\Ib(p)$\IndexDef{$\Ib(p)$@I-(p) (motions of $\Gb(p)$)}[ ] & 
  $\{V \in \Vf(\V) : YV = \zero\}$.  Motions of~$\Gb(p)$.  \\\hline
$p$\IndexDef{$p$@p (map from $\V$ to $\RR^n$)}[ ] & 
  A 1-1 map from $\V$ to~$\RR^n$. \\\hline
$T(p)$\IndexDef{$T(p)$@T(p) (Euclidean motions)}[ ] & 
  The Euclidean motions.  These can be induced from rigid motions of $\RR^n$
  and do not change the shape of the tensegrity (see
  \ref{appendix:rigidMotions}). \\ \hline
$\V$\IndexDef{$\V$@V (vertices)}[ ] & 
  The \emph{vertices} of~$G(p)$.  A set with the topology induced from $\RR^n$
  via~$p$.  When $\V$ is a curve rather just a general set, we may call it
  $\gamma$\IndexDef{$\gamma$@gamma (synonym for $\V$)}[ ]. \\\hline
$\Vf(\V)$\IndexDef{$\Vf(\V)$@Vf(V) (variations)}[ ] & 
  Continuous vector fields on $\V$, the \emph{variations}.  Those variations
  which produce nonnegative loads are called \emph{motions}.  \\\hline
$\X$\IndexDef{$\X$@X (design variations)}[ ] & 
  The \emph{design variations}.  A subspace of $\Vf(\V)$ containing~$T(p)$.  \\
  \hline
$Y\mcol \Vf(\V) \to C(\E)$\IndexDef{$Y$@Y (rigidity operator)}[ ] &  
  $\pm(V(v_1)-V(v_2)) \cdot (p(v_1)-p(v_2))$ with the sign being $+$ for struts
  and $-$ for cables.  The rigidity operator.  Used to calculate the load for a
  given variation.  \\\hline
\end{tabularx}
\end{definition}
\label{enddef:tensegrity}
\IndDefEnd{crossed square}[ ]

\section{Working the Crossed-square Example}
\label{section:WorkingRWex}
\vref{fig:RWexampleNum} has our example again, this time with the edges
numbered.  The time has come to see what we can learn about
it.\Index[|(]{crossed square}[ ]
\begin{figure}[ht]
\hspace*{\fill}
\begin{tikzpicture}
\draw[cable] (0,0) -- node [weight] {1} (2,0) -- node [weight] {4}
  (2,2) -- node [weight] {6} (0,2) -- node [weight] {3} (0,0);
\draw[strut] 
  (0,0) node [vertex] {1} -- node [pos=0.25,weight] {2} (2,2) node [vertex] {3}
  (0,2) node [vertex] {4} -- node [pos=0.25,weight] {5} (2,0) node [vertex] {2};
\end{tikzpicture}
\hspace*{\fill}
\capt{The crossed-square example with its edges numbered.}
\label{fig:RWexampleNum}
\end{figure}
We start by building the rigidity operator and choosing a stress:
$$Y = \begin{bmatrix}
 1 &  0 & -1 &  0 &  0 &  0 &  0 &  0 \\
-1 & -1 &  0 &  0 &  1 &  1 &  0 &  0 \\
 0 &  1 &  0 &  0 &  0 &  0 &  0 & -1 \\
 0 &  0 &  0 &  1 &  0 & -1 &  0 &  0 \\
 0 &  0 &  1 & -1 &  0 &  0 & -1 &  1 \\
 0 &  0 &  0 &  0 & -1 &  0 &  1 &  0 
\end{bmatrix}
\quad \text{ and } \quad
\mu = \begin{bmatrix} 1 \\ 1 \\ 1 \\ 1 \\ 1 \\ 1 \end{bmatrix}
\label{eq:Y}$$
Then $Y^\top \mu = \zero$, so (by Proposition \vref{prop:stress}) $\mu$ is a
strictly positive stress.  By Lemma \vref{lemma:NewFiveOne}, then, $G(p)$ is
bar equivalent.  It will be infinitesimally rigid or not as $\Gb(p)$ is.

Even though our main concern here is with bar equivalence, this example is
simple enough that we can calculate whether $\Gb(p)$ is infinitesimally rigid
or not.  Since the primary value of bar equivalence is tied to knowing whether
the bar framework is infinitesimally rigid, it seems worth the time to
calculate it for this example.  In \vref{section:CrossedMore} we dig even 
deeper, showing that for this example, the image of $Y$ is a hyperplane
in~$C(\E)$.

Suppose that $V$ is a motion of the bar framework $\Gb(p)$, that is, $V \in
\Ib(p)$.  Then $YV = \zero$.  But we have
\begin{align*}
YV & = \begin{bmatrix}
 1 &  0 & -1 &  0 &  0 &  0 &  0 &  0 \\
-1 & -1 &  0 &  0 &  1 &  1 &  0 &  0 \\
 0 &  1 &  0 &  0 &  0 &  0 &  0 & -1 \\
 0 &  0 &  0 &  1 &  0 & -1 &  0 &  0 \\
 0 &  0 &  1 & -1 &  0 &  0 & -1 &  1 \\
 0 &  0 &  0 &  0 & -1 &  0 &  1 &  0 
\end{bmatrix}
\begin{bmatrix}
v_{x,1} \\ v_{y,1} \\ v_{x,2} \\ v_{y,2} \\ 
v_{x,3} \\ v_{y,3} \\ v_{x,4} \\ v_{y,4} 
\end{bmatrix} =
\begin{bmatrix}
v_{x,1} - v_{x,2} \\ - v_{x,1} - v_{y,1} + v_{x,3} + v_{y,3} \\
v_{y,1} - v_{y,4} \\ v_{y,2} - v_{y,3} \\
v_{x,2} - v_{y,2} - v_{x,4} + v_{y,4} \\ -v_{x,3} + v_{y,4}
\end{bmatrix}
\end{align*}
A little calculation gives us that any motion of the bar framework must look
like
$$V = \begin{bmatrix} 
v_{x,1} & v_{y,1} & v_{x,1} & v_{y,2} &
(v_{x,1} + v_{y,1} - v_{y,2}) & v_{y,2} & 
(v_{x,1} + v_{y,1} - v_{y,2}) & v_{y,1} 
\end{bmatrix}^\top$$

We'd like to know whether this variation is a Euclidean motion or whether it
changes the shape of the bar framework.  We can remove any translation from $V$
by setting the sum of the horizontal components to zero and the sum of the
vertical components to zero.  That is we require that
$$
v_{x,1} + v_{x,2} + v_{x,3} + v_{x,4} = 0 \quad \text{ and } \quad
v_{y,1} + v_{y,2} + v_{y,3} + v_{y,4} = 0.
$$

If we do that, we get
$$V = \begin{bmatrix}
v_{x,1} & -v_{x,1} & v_{x,1} & v_{x,1} & -v_{x,1} & v_{x,1} & -v_{x,1} &
-v_{x,1}
\end{bmatrix}^\top,$$
which is simply a rotation around the center of the figure (see
\vref{fig:rotation}).
\begin{figure}[ht]
\hfil
\begin{tikzpicture}
\draw[bar] (0,0) -- (2,0) -- (2,2) -- (0,2) -- (0,0);
\draw[bar] 
  (0,0) node [vertex] {1} -- (2,2) node [vertex] {3}
  (0,2) node [vertex] {4} -- (2,0) node [vertex] {2};
\draw[->] (0,0) -- +( 1,-1) node [scale=0.8,below] {$( v_{x,1},-v_{x,1})$};
\draw[->] (2,0) -- +( 1, 1) node [scale=0.8,above] {$( v_{x,1}, v_{x,1})$};
\draw[->] (2,2) -- +(-1, 1) node [scale=0.8,above] {$(-v_{x,1}, v_{x,1})$};
\draw[->] (0,2) -- +(-1,-1) node [scale=0.8,below] {$(-v_{x,1},-v_{x,1})$};
\end{tikzpicture}
\hfil
\capt[Every motion of the crossed-square bar framework is Euclidean.]{Every
motion of the crossed-square bar framework is a Euclidean motion.}
\label{fig:rotation}
\end{figure}
So every element of $\Ib(p)$ is an element of~$T(p)$.
That makes $\Gb(p)$ infinitesimally
rigid and by Theorem \vref{theorem:RW}, $G(p)$ is infinitesimally rigid.
\Index[|)]{crossed square}[ ]

\section{The Mechanics}
\label{section:Mechanics}
In their proof, Roth and Whiteley make use of two theorems of
Rockafellar\Index[|(]{R. Tyrrell Rockafellar}[3,12 123], but he uses some terms
we haven't seen before, so we'll mention them first.  

Rockafellar uses the term ``relative interior''.  The idea is useful in a
vector space when dealing with objects that would fit within proper subspaces
of the vector space.  For example, when looking at a line segment in $\RR^3$,
it is clear that there are two endpoints and many ``non-end'' points.  However,
a line segment in $\RR^3$ has no 3-dimensional interior.  That is where
``relative interior'' comes in handy.  Those ``non-end'' points form the
relative interior of the line segment (in the case of a single point, the point
itself is the relative interior).

Let's define that more rigorously.  Our first two definitions differ in the 
bounds on~$\lambda$.
\begin{definition}[\ocite{Rockafellar}*{p.\ 10}]
A subset $C$ of $\RR^n$ is called a \emph{convex set}\Index{convex set} if
$(1 - \lambda)x + \lambda y \in C$ whenever $x, y \in C$ and $\lambda \in
(0,1)$.
\end{definition}
\begin{definition}[\ocite{Rockafellar}*{p.\ 3}]
A subset $M$ of $\RR^n$ is called an \emph{affine set}\Index{affine set} if
$(1 - \lambda)x + \lambda y \in M$ whenever $x, y \in M$ and $\lambda \in \RR$.
\end{definition}
\begin{definition}[\ocite{Rockafellar}*{p.\ 6}]
The \emph{affine hull}\Index{affine hull} of a set $S \subset \RR^n$ is the
intersection of the collection of affine sets $M$ such that $M \supset S$.
\end{definition}
\begin{definition}[\ocite{Rockafellar}*{p.\ 44}]
The \emph{relative interior}\Index_{relative interior} of a convex set $C
\subset \RR^n$ is the set of points $x$ in the affine hull of $C$ for which
there exists an $\varepsilon > 0$ such that $y \in C$ whenever $y$ is in the
affine hull of $C$ and $d(x,y) \le \varepsilon$.
\end{definition}

For the sake of understanding, let's apply these definitions to the line 
segment $S$ running from $(0,0,0)$ to $(0,0,1)$ in~$\RR^3$. 
\begin{prop}
The relative interior of the line segment $\{(0,0,z) : z \in [0,1]\}$ is
the open set $\{(0,0,z) : z \in (0,1)\}$.
\end{prop}
\begin{proof}
The affine sets that contain $S$ are the $z$-axis, any plane containing the
$z$-axis and~$\RR^3$.  The affine hull of $S$, then, is the intersection of
those set, or simply the $z$-axis.  Then, for any point $(0,0,z)$, one of three
things is true:
\begin{enumerate}
\item $z > 1$ or $z < 0$.  In this case $(0,0,z)$ itself is outside $S$, so no
epsilon neighborhood of it in the affine hull could be completely contained
within~$S$.
\item $z = 0$ or $z = 1$.  Here $(0,0,z)$ is inside $S$, but for no positive
$\varepsilon$ is the neighborhood completely inside $S$, so these points are
also not in the relative interior.
\item $z > 0$ and $z < 1$.  In this case, setting $\varepsilon =
\min\{1-z, z-0\}$ gives us a neighborhood in the affine hull that is also
in~$S$.  These are the points of the relative interior. \qedhere
\end{enumerate}
\end{proof}
While we're on the topic, let's support that offhand comment about the relative
interior of a point.
\begin{prop}
The relative interior of a point is the point itself.
\end{prop}
\begin{proof}
Let $x \in \RR^n$.  Now the smallest (in terms of set inclusion) affine set
containing $x$ (and thus the affine hull of $x$) is $x$ itself.  But then, if
$y$ is in the affine hull, $y$ must be $x$, so $d(x,y) = 0 \le \varepsilon$ for
any $\varepsilon > 0$.
\end{proof}

Let's get back to definitions.
\begin{definition}[\ocite{Rockafellar}*{p.\ 170}]
A \textit{polyhedral set}\Index_{polyhedral set} set is one that is the
intersection of a finite number of closed half-spaces.
\end{definition}
\begin{definition}[\ocite{Rockafellar}*{p.\ 162}]
A \emph{face}\Index1{face (of a convex set)} of a convex set $C$ is a convex
subset $C'$ of $C$ such that every (closed) line segment in $C$ with a relative
interior point in $C'$ has both endpoints in~$C'$.
\end{definition}
\begin{prop}
For a cube in $\RR^3$, the faces are: the 6 squares that form the boundary, 
the 12 edges of the squares and the 8 corners as well as the entire cube and
the empty set.
\end{prop}
\begin{proof}
Clearly the empty set is trivially a face, as it contains no line segments.  
Similarly, the entire cube is a face as any closed line segment in the cube
has both ends in the cube.

The 8 corners are faces in much the same fashion.  They can hold endpoints of
line segments, but the only way they can hold relative interior points of line
segments is if those line segments are the zero-length ones which are the
corners themselves, and in that case both endpoints are certainly in.
No other single point can be a face, though, as any other point can be a 
relative interior point of some line segment that passes through it.

The 12 edges are faces, as any line segment that has a relative interior
point in an edge must lie entirely in that edge.  No other segment in the
cube can be a face since it will be crossed by some other segment that it does
not contain.

Finally, any segment with a relative interior point in one of the boundary
squares of the cube must lie entirely in that square, but any other plane
passing through the cube can be intersected transversely with some other
segment in the cube, so those squares are the only 2-dimensional faces.
\end{proof}

In considering unbounded convex sets, Rockafellar introduces a type of ``point
at infinity''.  An unbounded convex set (not necessarily conic), he notes, must
contain some entire half-line or it wouldn't be an unbounded set.  However,
that half-line might point in any direction.  So the directions in which the
half-lines lie can be thought of as ``horizon points'' or points at infinity.
Formally,
\begin{definition}[\ocite{Rockafellar}*{p.\ 60}]
A \emph{direction}\Index{direction} of $\RR^n$ is an equivalence class of the
collection of all closed half-lines of $\RR^n$ under the equivalence relation
``half-line $L_1$ is a translate of half-line $L_2.$''
\end{definition}
\begin{definition}[\ocite{Rockafellar}*{p.\ 12}]
The \emph{convex hull}\Index{convex hull} of a set $S \subset \RR^n$ is the
intersection of all convex sets containing~$S$.
\end{definition}
\begin{definition}[\ocite{Rockafellar}*{p.\ 170}]
A \emph{finitely generated convex set}\Index{{finitely generated} convex
set}[123 3!12 2!3!1]
is the convex hull of a finite set of points and directions.
\end{definition}

As an example, the convex hull in $\RR^2$ of the two points $(0,0)$ and $(1,0)$
and the direction of the vector $(1,1)$ is the area shown in
\vref{fig:convexHull}.
\begin{figure}[ht]
\hspace*{\fill}
\begin{tikzpicture}
\foreach \x in {0,90,180,270} {
  \draw[->] (0,0) -- (\x:2cm);
}
\fill[blue!20!white] (1.5,1.5) -- (0,0) -- (1,0) -- (2,1) 
   [snake=snake,segment amplitude=0.5mm,segment length=0.1414cm] -- (1.5,1.5);
\draw[ultra thin] (1.5,1.5) -- (0,0) -- (1,0) -- (2,1);
\end{tikzpicture}
\hspace*{\fill}
\capt[Convex set generated by $(0,0)$, $(1,0)$ and the direction $(1,1)$.]{The
convex set generated by the points $(0,0)$ and $(1,0)$ and the direction of the
vector $(1,1)$.}
\label{fig:convexHull}
\end{figure}

Now we are ready for the theorems.  The first is a separation theorem and the
second gives an equivalence between attributes of a convex set.
\begin{theorem}
\label{theorem:Separation}
Let $X$ be a closed convex cone in $\RR^n$ and $x_0 \in \RR^n \setminus X$.
Then, since the relative interiors of $X$ and $\{x_0\}$ are disjoint, there
exists a hyperplane separating $X$ and $\{x_0\}$
properly.
\end{theorem}
% Statement of the theorem from the book at Google:
% Let C1 and C2 be non-empty convex sets in R^n.  IN order that there exist a
% hyperplane separating C1 and C2 properly, it is necessary and sufficient that
% ri C1 and ri C2 have no point in common.
%
% It might also be valuable to consider the theorem from Covex Polytopes by 
% Branko Grunbaum, p. 11, Theorem 1.
\begin{proof}
See \ocite{Rockafellar}*{Theorem 11.3, p. 97}.
\end{proof}
\begin{theorem}
\label{theorem:Rockafellar}
The following properties of a convex set $C$ are equivalent:
\begin{enumerate}[label=(\alph*)]
\item $C$ is polyhedral;
\item $C$ is closed and has only finitely many faces;
\item $C$ is finitely generated.
\end{enumerate}
\end{theorem}
\begin{proof}
See \ocite{Rockafellar}*{Theorem 19.1, p. 171}.
\end{proof}\Index[|)]{R. Tyrrell Rockafellar}[3,12 123]

We also have two more objects that will be of value to us.  Suppose we have
some set $X$ of vectors in~$\RR^n$.  It seems natural to consider the set of
vectors that are normal to all of our vectors.  This is the finite-dimensional
analogue of the annihilator\Index{annihilator}, so we use the same notation:
$$X^\perp = \{ v \in \RR^n : v \cdot x = 0 \text{ for all } x \in X\}$$
\Index{$S^\perp$@*18 (annihilator of $S$)}[ ]
Similarly, the set of vectors that ``point in the same general
direction'' as ours is the analog of the dual cone\Index{dual cone}[2,1 ]:
$$X^* = \{ v \in \RR^n : v \cdot x \ge 0 \text{ for all } x \in X\}.$$
\Index{$S^*$@*19 (dual cone of $S$)}[ ]

\begin{lemma}
\label{lemma:above}
For any set $X \subset \RR^n$, $X^*$ is a closed convex cone.
\end{lemma}
\begin{proof}
If $x_1,x_2 \in X^*$, then, for any $x \in X$,
$$((1-\lambda) x_1 + \lambda x_2) \cdot x = (1-\lambda) x_1 \cdot x + \lambda
x_2 \cdot x \ge 0$$
whenever $\lambda \in [0,1]$.  Also, $\alpha x_1 \cdot x \ge 0$ whenever
$\alpha > 0$, so $X^*$ is a convex cone.  If $y$ lies outside $X^*$, then by
the definition of $X^*$ there must be some $x_- \in X$ such that $y \cdot x_- <
0$.  Let $z \in \RR^n$ with $\|y - z\|< \frac{- y \cdot x_-}{\|x_-\|}$.  Then
we have
\begin{align*}
z \cdot x_- & = (z - y + y) \cdot x_- 
  = (z - y) \cdot x_- + y \cdot x_- \\
& \le (z - y) \cdot (z - y) \frac{\|x_-\|}{\|z - y\|} + y \cdot x_- 
  = \|y - z\| \|x_-\| + y \cdot x_- \\
& < -y \cdot x_- + y \cdot x_- = 0.
\end{align*}
That is, $z$ is also outside~$X^*$.  So $X^*$ is closed.
\end{proof}

Here's a corollary to Theorem \ref{theorem:Separation} that will prove useful.
\begin{corollary}
\label{corollary:Separation}
If $X$ is a closed convex cone in $\RR^n$ and $x_0 \in \RR^n \setminus X$, then
there exists $\mu \in X^*$ with $\mu \cdot x_0 < 0$.
\end{corollary}
\begin{proof}
By Theorem \ref{theorem:Separation} there exists a hyperplane, $h$, with 
$x_0$ lying strictly on one side of it and $X$ in the half-space on the other
side.  The normal to $h$ pointing into the half-space with $X$ has nonnegative
dot product with all the elements of $X$ and hence is in $X^*$, but it has
negative dot product
with~$x_0$.
\end{proof}

Here are a few more attributes of $X^*$:
\begin{prop}~
\label{prop:Attributes}
\begin{enumerate}[label=(\roman*)]
\item for $Z \subset \RR^n$ and $X$ a closed convex cone in $\RR^n$, $Z \subset
X$ if and only if $X^* \subset Z^*$.~
\label{item:Containment}
\item for $Z \subset \RR^n$, $Z^{**} = (Z^*)^*$ is the smallest closed convex
cone in $\RR^n$ that contains~$Z$.
\label{item:smallestCone}
\item if $Z = \{z_1, \dotsc, z_k\}$ is a finite subset of $\RR^n$, then the
convex cone $$\left\{\sum_{i=1}^k \lambda_i z_i : \lambda_i \ge 0 \text{ for }
1 \le i \le k\right\}$$ 
generated by $Z$ is closed in $\RR^n$ (and therefore equals $Z^{**}$).
\label{item:Rock3}
\end{enumerate}
\end{prop}
\begin{proof}~
\begin{enumerate}[label=(\roman*)]
\item Suppose that $X$ is some closed convex cone in $\RR^n$ and we have some
$Z \subset \RR^n$.  If $Z \subset X$, then certainly any vector that has a
positive dot product with everything in $X$ must have a positive dot product
with everything in $Z$, hence $X^* \subset Z^*$.  What about the
other way around?  If there is some $z \in Z, z \notin X$, then by Corollary
\vref{corollary:Separation}, there is some $\mu \in X^*$ such that $\mu
\cdot z < 0$.  But that means that $\mu \notin Z^*$.

\item $Z^{**}$ is a closed convex cone by Lemma \ref{lemma:above}.
Furthermore, if $z \in Z$, then $z$ has a positive dot product with everything
in $Z^*$, so $Z \subset Z^{**}$.  Suppose that some closed convex cone $C
\supset Z$ but that there were some $\hat{z} \in Z^{**}$ and $\hat{z} \notin
C$.  Then, by Corollary \ref{corollary:Separation}, there must be an element
$\mu$ of $C^*$ that separates $\hat{z}$ and~$C$.  Now by Item
\ref{item:Containment}, $C^* \subset Z^*$, so $\mu \in Z^*$ and yet $\mu \cdot
\hat{z} < 0$.  That contradicts the definition of~$Z^{**}$.  Hence $Z^{**}$
must be the smallest closed convex cone containing~$Z$.

\item Finally, if $Z$ is a finite subset of $\RR^n$, then the convex cone
$$Z_c = \left\{ \sum_{i=1}^k \lambda_i z_i : \lambda_i \ge 0 \text{ for } 1 \le
i \le k \right\}$$ 
is finitely generated (being generated by the origin and the directions of
the $z_i$) and hence, by Theorem \vref{theorem:Rockafellar}, it is closed.  Now
any positive linear combination of the $z_i$ must have positive dot product
with any vector in $Z^*$, so $Z_c \subset Z^{**}$.  Since it is a closed convex
cone it must, by Item \ref{item:smallestCone}, equal~$Z^{**}$. \qedhere
\end{enumerate}
\end{proof}

One lemma remains that we will want in the next section.
\begin{lemma}
\label{lemma:subspaceClaim}
Let $X = \{x_1, \dotsc, x_k\} \in \RR^n$.  Then $X^* = X^\perp$ implies that 
$X^{**}$ is a subspace.
\end{lemma}
\begin{proof}
Let $\alpha \in \RR$ and $v_1, v_2 \in X^{**}$.  We want to show
that $v_1 + v_2 \in X^{**}$ and that $\alpha v_1 \in X^{**}$.
\begin{align*}
v_1, v_2 \in X^{**} & \Rightarrow
v_1 \cdot x^{*} \ge 0 \text{ and } v_2 \cdot x^{*} \ge 0, & \forall
x^{*} \in X^{*} \\
& \Rightarrow (v_1 + v_2) \cdot x^* = v_1 \cdot x^* + v_2 \cdot
x^* \ge 0, & \forall x^* \in X^* \\
& \Rightarrow (v_1 + v_2) \in X^{**}
\end{align*}
Now if $\alpha \ge 0$, then we easily have that $(\alpha v_1) \cdot x^* =
\alpha v_1 \cdot x^* \ge 0$ for all $x^* \in X^*$.  But if $\alpha < 0$, we
need to be a little more clever.  In that case, we have
$\alpha v_1 \cdot x^* = |\alpha| v_1 \cdot -x^*$, so we are fine if $x^* \in
X^* \Rightarrow -x^* \in X^*$.  But
\begin{align*}
x^* \in X^* & \Rightarrow x^* \cdot x = 0, \forall x \in X,
\text{ because $X^* = X^\perp$ by hypothesis} \\
& \Rightarrow -x^* \cdot x = 0, \forall x \in X \\
& \Rightarrow -x^* \in X^*
\end{align*}
and we have that $X^{**}$ is a subspace.
\end{proof}
\section{Key Lemma and Theorem}
\label{section:KeyLemma}
We now reach the key lemma \cite{MR610958}*{Lemma 5.1, p. 426}.  This is the
engine that Roth \& Whiteley use to establish bar equivalence\Index{bar
equivalent}[ ] for a tensegrity so that they can analyze the infinitesimal
rigidity\Index{infinitesimal rigidity}[2,1 ] of the bar framework instead of
the tensegrity.  It is the predecessor to our main theorem.  The proof we give
here is essentially that of Roth \& Whiteley.

\begin{lemma}
\label{lemma:LemmaFiveOne}
\Index[|(]{key lemma (Roth \& Whiteley)}[ ]
Suppose $X = \{x_1, \dotsc, x_k\} \subset \RR^n$.  Then $X^* = X^\perp$ if
and only if there exist positive scalars $\lambda_1, \dotsc, \lambda_k$ with 
$\sum_{i=1}^k \lambda_i x_i = \zero$.
\end{lemma}
\begin{proof}  
By Item \ref{item:Rock3} of Proposition \ref{prop:Attributes}, $X^{**} =
\left\{\sum_{i = 1}^k \lambda_i x_i : \lambda_i \ge 0 \text{ for } 1 \le i \le
k \right\}$.

Suppose that $X^* = X^\perp$.  Then, by Lemma \ref{lemma:subspaceClaim},
$X^{**}$ is a subspace.  Since every $x_j \in X$ is also in $X^{**}$ we can
write
$$-x_j = \sum_{i=1}^k \lambda_i x_i$$ 
for some set of $\lambda_i \ge 0$.  But then we have
$$\lambda_1 x_1 + \dotsb + (1 + \lambda_j) x_j + \dotsb + \lambda_k x_k =
\zero.$$ 
Adding up $k$ such expressions gives a linear dependency of $x_1, \dotsc, x_k$
with all positive coefficients.

Conversely, if there is some set $\lambda_i > 0$ for which $\sum_{i = 1}^k
\lambda_i x_i = \zero$ and if $x^* \in X^*$, then 
$$0 = x^* \cdot \zero = x^* \cdot \sum_{i=1}^k \lambda_i x_i = \sum_{i=1}^k
\lambda_i (x^* \cdot x_i).$$
Since $x^* \cdot x_i \ge 0$ and $\lambda_i > 0$ for all $i$, we must have 
$x^* \cdot x_i = 0$ for all $i$, so $x^* \in X^\perp$.
\end{proof}

Roth \& Whiteley's lemma is very much in the world of finite-dimensional 
Euclidean spaces.  We need to reformulate it in preparation for the move
to more general spaces.  
\begin{lemma}
\label{lemma:MidFiveOne}
Let $Y$ be an $m \cross n$ matrix with elements in~$\RR$.  Then exactly one of
the following is true:
\begin{enumerate}[label=(\Roman*)]
\item There exists $V \in \RR^n$ such that $YV \gneqq \zero$.
\label{item:LemmaFirstCase}
\item There exists $\mu \in \RR^m$ with $\mu > \zero$ such that $Y^\top \mu =
\label{item:LemmaSecondCase}
\zero$.
\end{enumerate}
\end{lemma}
\begin{proof}
Let $X = \{ x_i \}$ be the set of row vectors of~$Y$.  Then those $V$ for which
$YV = \zero$ make up $X^\perp$ and those $V$ for which $YV$ has no negative
coordinates make up~$X^*$. 

Case \ref{item:LemmaFirstCase} is false if and only if $X^* = X^\perp$.  
By Lemma \ref{lemma:LemmaFiveOne}, $X^* = X^\perp$ if and only if 
there exist positive scalars $\lambda_1, \dotsc, \lambda_m$ such that
$\sum_{i=1}^m \lambda_i x_i = \zero$.  And by defining $\mu = \left[
\begin{smallmatrix} \lambda_1 \\ \vdots \\ \lambda_m \end{smallmatrix} \right]$
we have that the $\lambda_i$'s exist if and only if Case
\ref{item:LemmaSecondCase} is true.
\end{proof}

Lemma \ref{lemma:MidFiveOne} resembles the ``Theorem of the
Alternative''\Index{Theorem of the Alternative}[1!234] we talk about in 
\ref{section:ToA}, but it doesn't provide a lot of geometric insight.  
Let's recast it yet one more time.

\begin{lemma}
\label{lemma:NewFiveOne}
\statementNewFiveOne
\end{lemma}
\begin{proof}
Let $G(p)$ be bar equivalent.  Then there is no $V$ for which $YV \gneqq \zero$,
so we are not in Case \ref{item:LemmaFirstCase} of Lemma
\ref{lemma:MidFiveOne}.  Hence we must be in Case \ref{item:LemmaSecondCase},
i.e., $G(p)$ has a strictly positive stress.

Conversely, let $G(p)$ have a strictly positive stress.  We are again in
Case \ref{item:LemmaSecondCase}, so there can be no $V$ like that described in
Case \ref{item:LemmaFirstCase}, and $G(p)$ is bar equivalent.
\end{proof}
\Index[|)]{key lemma (Roth \& Whiteley)}[ ]

Finally, as an application of Lemma \vref{lemma:NewFiveOne}, we prove Roth and
Whiteley's theorem \ycite{MR610958}*{Theorem 5.2, p. 427}.
\begin{theorem}
\label{theorem:RW}
Suppose $G(p)$ is a (finite) tensegrity framework in $\RR^n$ and $\Gb(p)$ the
associated bar framework.  Then $G(p)$ is infinitesimally rigid if and only if
$\Gb(p)$ is infinitesimally rigid and there exists a strictly positive stress
of~$G(p)$.
\end{theorem}
\begin{proof} 

$(\Rightarrow)$.  Suppose that $G(p)$ is infinitesimally rigid.  Then $I(p) =
\Ib(p) = T(p)$ by Definition \vref{definition:IRBE}.  So $\Gb(p)$ is
infinitesimally rigid (again by Definition \ref{definition:IRBE}) and $G(p)$
is bar equivalent by Proposition \ref{prop:IRBE}.  By Lemma
\ref{lemma:NewFiveOne}, $G(p)$ has a strictly positive stress.

$(\Leftarrow)$.  Conversely, suppose that $\Gb(p)$ is infinitesimally rigid and
that $G(p)$ has a strictly positive stress.  By Lemma \ref{lemma:NewFiveOne},
$G(p)$ is bar equivalent.  Infinitesimal rigidity of $\Gb(p)$ gives us $\Ib(p)
= T(p)$ and bar equivalence of $G(p)$ gives us $I(p) = \Ib(p)$, so $I(p) =
T(p)$ and thus $G(p)$ is infinitesimally rigid.
\end{proof}
\Index[|)]_,{Ben Roth}
\Index[|)]{Walter Whiteley}[2,1]

\begin{chapter}{Main Theorems}
\label{chapter:Main}
\section{Under the Hood: Theorems of the Alternative}
\label{section:ToA}
Having seen what Roth \& Whiteley accomplished, 
\IndDefBeg{Theorem of the Alternative}[1!234]
we'd like to get a better grasp of what makes their theorem work, so that we
can extend it to continuous tensegrities.  Let's take another look at Lemma
\vref{lemma:MidFiveOne}.  The statement that $Y^\top \mu = \zero$ implies that
$(Y^\top \mu) V = 0$ for all $V \in \Vf(\V)$.  That's the same as saying that
$\mu^\top YV = 0$ for all $V$ or that $\mu \in (\im Y)^\perp$.

With that understanding, we can see Lemma \ref{lemma:MidFiveOne} as saying,
``either $\im Y$ intersects the nonnegative orthant somewhere other than the
origin, or else $(\im Y)^\perp$ intersects the interior of the nonnegative
orthant, but not both'' (see \vref{fig:Stiemke}).
\begin{figure}[ht]
\hspace*{\fill}
\begin{tikzpicture}
\path [fill=blue!20!white] (0,0) -- (0,3) -- (3,3) -- (3,0) -- cycle;

\draw [->] (0,0) -- (3,0) node [right] {$x$};
\draw [->] (0,0) -- (0,3) node [above] {$y$};
\draw [->] (0,0) -- (-3,0);
\draw [->] (0,0) -- (0,-3);

\draw [thick] (-3,1) -- node [pos=0.1,anchor=south west] {$\im Y$} (3,-1);
\draw [thick] (-1,-3) -- node [pos=0.9,anchor=north west] {$(\im Y)^\perp$}
              (1,3);
\end{tikzpicture}
\hspace*{\fill}
\capt{A Theorem of the Alternative in~$\RR^2$.}
\label{fig:Stiemke}
\end{figure}

It turns out that this is one of many ``Theorems of the Alternative'' (also
sometimes called ``Theorems on the Alternative''), in which exactly one of two
alternative systems of inequalities or equalities is shown to have a solution.
In his \textit{Nonlinear Programming}, Olvi Mangasarian\Index_,{Olvi
Mangasarian} has a wonderful chapter where he catalogs no fewer than 11 of
these theorems \cite{MR0252038}*{ch. 2}.  Mangasarian indicates that the one we
have in hand is the one from Stiemke's\Index_,{Erich Stiemke} 1915 paper
\cite{Stiemke}.\footnote{The papers of Stiemke, Gordan and Motzkin (see
\ref{appendix:Motzkin}) are all in German, so, due to time constraints, I have
not verified Mangasarian's analysis of them.} It is very similar to a theorem
he credits to Paul Gordan\Index_,{Paul Gordan} \ycite{Gordan73}.  Here are
Mangasarian's statements of the two, using our notation.
\begin{theorem}[Stiemke 1915]
\IndexDef{Stiemke Theorem}[2!1]
\label{theorem:Stiemke}
For a given matrix $Y$, either 
\begin{enumerate}[label=(\Roman*)]
\item $YV \gneqq \zero$ has a solution $V$ \\ or
\item $Y^\top \mu = \zero$, $\mu > \zero$ has a solution $\mu$ 
\end{enumerate}
but never both.
\end{theorem}
\begin{theorem}[Gordan 1873]
\IndexDef{Gordan Theorem}[2!1]
\label{theorem:Gordan}
For a given matrix $Y$, either 
\begin{enumerate}[label=(\Roman*)]
\item $YV > \zero$ has a solution $V$ \\ or
\item $Y^\top \mu = \zero$, $\mu \gneqq \zero$ has a solution $\mu$ 
\end{enumerate}
but never both.
\end{theorem}
The difference between those two is in which case is true when both $\im Y$ and
$(\im Y)^\perp$ lie on the boundar of the nonnegative orthant.  That difference
will appear again later in our main theorems.  Before we get there, though, we
need to go back through the various entities listed in Definition
\ref{def:tensegrity} \vpagerefrange{def:tensegrity}{enddef:tensegrity} and see
how they change during the move to the continuous world.
\IndDefEnd{Theorem of the Alternative}[1!234]

\section{Turning Sets into Spaces: \texorpdfstring{$\V$}{V} and \texorpdfstring{$\E$}{E}}
\label{section:E}
\subsection{New Topologies via Metrics}
We will need $\E$ to be a compact metric space, so we set out to provide a
metric on the set.  Now $\E$ is, by definition, a disjoint union of two sets.
Each of those sets is a subset of the quotient space $\vvt$ formed from the 
direct product\Index{direct product}[2!1 ] $\V \cross \V$ by the identification
$(v_1,v_2) \sim (v_2,v_1)$.  We'll work our way through these spaces in
creating the desired metric.\Index[|(]{$\E$@E (edgeset)}[ ]

\begin{definition}
\IndexDef{$\V$@V (vertices)}[ ]
\IndexDef{$d_\V$@dV (metric on $\V$)}[ ]\IndexDef{metric on $\V$}[1!23]
The space $\V$ is the set of vertices, the metric induced by its embedding into
$\RR^n$:
$$d_\V(v_1,v_2) \coloneqq \|p(v_1) - p(v_2)\|$$
and the topology defined by $d_\V$.\footnote{That is, the topology whose basic
open sets are the open balls $d_\V(v_1,v_2) < \varepsilon$ for all $v_1,v_2 \in
V$, $\varepsilon > 0$.  See, for example \ocite{Munkres}*{p.\ 119}.}
\end{definition}

Now, a metric has three attributes.  It must be positive definite and
symmetric, and it must satisfy the triangle inequality.  Of course $d_\V$ has
those attributes naturally, but we'll need to check them for the metrics yet to
come.

We give $\V \cross \V$ the metric
$$d_{\V \cross \V}((v_1,v_2),(v_3,v_4)) \coloneqq \max \{d_\V(v_1,v_3),
d_\V(v_2,v_4)\}.$$
This has the decided advantage that the basic open sets are the basic open
sets of the product topology.  For example, if we take the basic open set
of size $\varepsilon$ around $(v_1,v_2)$, then the points inside it are the
edges $(v_3,v_4)$ where $d_\V(v_1,v_3) < \varepsilon$ and $d_\V(v_2,v_4) < 
\varepsilon$\IndexDef{$d_{\V \cross \V}$@dVcV (metric on $\V \cross \V$)}[
]\IndexDef{metric on {$\V \cross \V$}}[1!23].

\begin{prop}
$d_{\V \cross \V}$ is a metric.
\end{prop}
\begin{proof}
Now $d_{\V \cross \V}((v_1,v_2),(v_3,v_4))$ is clearly nonnegative and is zero
only when $v_1 = v_3$ and $v_2 = v_4$.  Its symmetry comes from the nature of 
$\max$.  So we only need to check the triangle inequality.

If $(v_1,v_2), (v_3,v_4), (v_5,v_6) \in \V \cross \V$, then 
\begin{align*}
d_{\V \cross \V}((v_1,v_2),(v_5,v_6)) 
& = \max \{ d_\V(v_1,v_5), d_\V(v_2,v_6) \}.
\end{align*}
Now, if $d_\V(v_1,v_5) \ge d_\V(v_2,v_6)$, we have
\begin{align*}
d_\V(v_2,v_6) & \le d_\V(v_1,v_5) \le d_\V(v_1,v_3) + d_\V(v_3,v_5) \\
& \le \max \{ d_\V(v_1,v_3), d_\V(v_2,v_4) \} + \max \{ d_\V(v_3,v_5),
d_\V(v_4,v_6) \}.
\end{align*}
Otherwise, we have
\begin{align*}
d_\V(v_1,v_5) & \le d_\V(v_2,v_6) \le d_\V(v_2,v_4) + d_\V(v_4,v_6) \\
& \le \max \{ d_\V(v_1,v_3), d_\V(v_2,v_4) \} + \max \{ d_\V(v_3,v_5),
d_\V(v_4,v_6) \}.
\end{align*}
So, either way, we have 
$$d_{\V \cross \V}((v_1,v_2),(v_5,v_6)) \le d_{\V \cross
\V}((v_1,v_2),(v_3,v_4)) + d_{\V \cross \V}((v_3,v_4),(v_5,v_6)).$$
Thus $d_{\V \cross \V}$ is a metric.
\end{proof}

Now we come to $\vvt$.
\begin{prop}
The function
$$
d_\vvt(\{v_1,v_2\},\{v_3,v_4\}) \coloneqq \min_{i \ne j \in \{3,4\}}
\left\{ d_{\V \cross \V}((v_1,v_2),(v_i,v_j)) \right\}
$$
is a metric for $\vvt$.%
\IndexDef{$d_\vvt$@dVVt (metric on $\vvt$)}[ ]%
\IndexDef{metric on {$\vvt$}}[1!23]
\end{prop}
\begin{proof}
Let $\{v_1,v_2\}$, $\{v_3,v_4\}$, and $\{v_5,v_6\}$ be in $\vvt$.
\begin{enumerate}[label=(\alph*)]
\item $d_\vvt(\{v_1,v_2\},\{v_3,v_4\})$ is the smaller of two nonnegative
numbers and hence nonnegative.  Furthermore, as $d_{\V \cross \V}$ is a metric,
$d_\vvt$ will only return $0$ when either $(v_1,v_2) = (v_3,v_4)$ or $(v_1,v_2)
= (v_4,v_3)$.  That is, only when the two are in the same equivalence class in
$\vvt$.
\item We need symmetry of the metric: \begin{align*}
d_\vvt(\{v_1,v_2\},\{v_3,v_4\}) & = \min_{i \ne j \in \{3,4\}}
  \left\{ d_{\V \cross \V}((v_1,v_2),(v_i,v_j)) \right\} \\
& = \min_{i \ne j \in \{3,4\}}
  \left\{ d_{\V \cross \V}((v_i,v_j),(v_1,v_2)) \right\} \\
& = \min_{i \ne j \in \{1,2\}}
  \left\{ d_{\V \cross \V}((v_3,v_4),(v_i,v_j)) \right\} \\
& = d_\vvt(\{v_3,v_4\},\{v_1,v_2\})
\end{align*}
\item Finally, we need to check the triangle inequality:
\begin{align*}
d_\vvt&(\{v_1,v_2\},\{v_5,v_6\}) = \min_{ i \ne j \in \{5,6\} }
\left\{ d_{\V \cross \V}((v_1,v_2),(v_i,v_j)) \right\} \\
& \le \min_{\substack{
s \ne t \in \{3,4\} \\
i \ne j \in \{5,6\}
}} \left\{ 
d_{\V \cross \V}((v_1,v_2),(v_s,v_t)) + d_{\V \cross \V}((v_s,v_t),(v_i,v_j)) 
\right\} \\
& = \min_{s \ne t \in \{3,4\}} \left\{ d_{\V \cross
\V}((v_1,v_2),(v_s,v_t)) \right\} + 
\min_{i \ne j \in \{5,6\}} \left\{ d_{\V \cross \V}((v_3,v_4),(v_i,v_j))
\right\} \\
& = d_\vvt(\{v_1,v_2\},\{v_3,v_4\}) + d_\vvt(\{v_3,v_4\},\{v_5,v_6\}).
\end{align*}
\end{enumerate}
So $d_\vvt$ is a metric for $\vvt$ and also for $\Str$ and $\Cr$ as subsets of
$\vvt$.  
\end{proof}

We've nearly achieved our goal, but we have a decision to make.  How far apart
should a strut be from a cable?  Both $\Str$ and $\Cr$ are
subsets of $\vvt$, so we could simply define the distance from strut $s$ to 
cable $c$ to be $d_\vvt(s,c)$.  

At first glance this seems a wise idea.  After all, suppose we have a
tensegrity in which the struts or the cables alone do not form a compact
set, but the two together do (such a construction may be seen in
\vref{fig:openstruts}).\footnote{Since $\E$ is Hausdorff\Index{$\E$@E (edgeset)
Hausdorff}[12!3], we know that if $\Str$ and $\Cr$ were compact, they would be
closed (see, for example, \ocite{Munkres}*{p.\ 165}).  But the sequence of
struts at angles $\pi - \frac{1}{2^n}$ and the sequence of cables at angles
$-\frac{1}{2^n}$ show that we would need a strut at $\pi$ and a cable at $0$
for closure.  

On the other hand, $\V$ is compact as a $p(\V)$ is a closed, bounded subset of
$\RR^2$, so $\V \cross \V$ is compact.  Now $\Str \cup \Cr$ is closed, so if
$\pi\colon \V \cross \V \to \vvt$ is the quotient map, then $\pi^{-1}(\Str \cup
\Cr)$ is a closed subset of $\V \cross \V$ and hence compact.  But then $\Str
\cup \Cr$ is the continuous image of a compact set and thus also compact
\cite{Munkres}*{p.\ 165--167}.
\label{footnote:compact}
}
\begin{figure}[ht]
\hfil
\begin{tikzpicture}
\draw[vcurve] (0,0) circle (1cm) circle (3cm);
\foreach \x in {0,10,...,179} {
  \draw[strut] (\x:1cm) -- (\x:3cm);
  \draw[cable] (\x:-1cm) -- (\x:-3cm);
}
\draw  (3.2,0) node [circle,draw=black,fill=black,scale=0.5] (a) {} 
       arc (0:180:3.2cm) 
               node [circle,draw=black,fill=white,scale=0.5] (b) {}
      (-3.5,0) node [circle,draw=black,fill=black,scale=0.5] (c) {} 
       arc (180:360:3.5cm)
               node [circle,draw=black,fill=white,scale=0.5] (d) {};
\draw (0,3.2) node [above] {Struts}
      (0,-3.4) node [below] {Cables};
\end{tikzpicture}
\hfil
\capt{A tensegrity with noncompact $\Str$ and $\Cr$ but compact $\Str \cup
\Cr$.}
\label{fig:openstruts}
\end{figure}
If we make the distance between struts and cables simply the $d_\vvt$ distance,
then this edgeset will be compact and our results here will apply to it.

However, if we do that and $\B$ is nonempty, we run into trouble.  If $b_s$ is
a bar in $\Str$ and $b_c$ is the same bar in $\Cr$, then
$d_\vvt(b_s,b_c)$ would be zero.  Putting the metric topology\Index{metric
topology}[2!1 ] on $\E$ in this situation would mean that there would be no
open sets which would separate $b_s$ and~$b_c$.  $\E$ would no longer be
Hausdorff.\Index{$\E$@E (edgeset) Hausdorff}[12!3]  As we certainly need $\E$
to be Hausdorff, we'll look for other options (and we'll circle around later
and talk for a bit about why excluding this example doesn't trouble us
greatly).

Since having a strut and cable too close together was a problem, perhaps we
should  define them to be far apart.  To do this, we could find out how far
apart any two struts or cables are and then make all struts farther away from
all cables than that.  This is not difficult, but it seems clearer to use what
\ocite{MR1912709} call an ``extended metric'' (they credit this to
\ocite{MR0263062}, who pass the credit on to \ocite{MR0124554}).

An extended metric\Index{extended metric}[2!1 ] differs from a normal metric
only in that it is allowed to take the value $\infty$ (to readers interested in
the properties of metrics are recommended \ocite{wiki:Metric} and
\ocite{Munkres}*{p.\ 119}).

\begin{prop}
The set $\Str \sqcup \Cr$ is a metric space with (extended)
metric $d_\E$ defined as:
$$
d_\E(e_1,e_2) = \begin{cases}
d_\vvt(e_1,e_2), & \text{if } e_1,e_2 \in \Str \text { or }
                              e_1,e_2 \in \Cr \text{, and} \\
\infty & \text{otherwise.}
\end{cases}
$$
\IndexDef{$d_\E$@dE (metric on $\E$)}[ ]\IndexDef{metric on $\E$}[1!23]
\end{prop}
\begin{proof}
Let $e_1$, $e_2$ and $e_3$ be in $\Str \sqcup \Cr$.  Then, if $e_1$ and $e_2$
are both struts or both cables, we get that $d_\E(e_1,e_2) \ge 0$ with equality
only if $e_1 = e_2$ and that $d_\E(e_1,e_2) = d_\E(e_2,e_1)$, both from the
nature of $d_\vvt$.  If they are different, then $d_\E(e_1,e_2) = d_\E(e_2,e_1)
= \infty > 0$.  So we are left with the triangle inequality to check.  

Again, if all three edges are of the same type, 
\begin{equation}
d_\E(e_1,e_2) + d_\E(e_2,e_3) \ge d_\E(e_1,e_3)
\label{eq:triangle}
\end{equation}
due to  that being true for $d_\vvt$.  Otherwise, at least two differ in type.
If $e_1$ differs from $e_3$, then the right side of \ref{eq:triangle} is
$\infty$, but one of the terms on the left must be as well ($e_1$ and $e_3$ are
different, so $e_2$ must differ from one or the other).  On the other hand, if
$e_1$ and $e_3$ are of the same type, the right side is finite and both of
the terms on the left are infinite.  We see that the triangle inequality holds
for~$d_\E$.

So $d_\E$ is a metric on $\Str \sqcup \Cr$.  
\end{proof}

\begin{definition}
The \emph{edgeset} $\E$\IndexDef{$\E$@E (edgeset)}[ ] is the space formed by
endowing the set $\Str \sqcup \Cr$ whith the metric $d_\E$ and the associated
metric topology\Index{metric topology}[2!1 ].
\end{definition}

\subsection{Familiar Topologies}
\label{subsection:familiar}
Now that we have topologies on $\V$ and $\E$, let's see what we can learn about
them.

\begin{prop}
The topology on $\V$, $\mathcal{T}_\V$\IndexDef{$\mathcal{T}_\V$@TV (topology
on $\V$)}[ ] is same as the topology,
$\mathcal{T}_{p(\V)}$\IndexDef{$\mathcal{T}_{p(\V)}$@TpV (alt. topology on
$\V$)}[ ], which $\V$ would have if we put the subspace topology on $p(\V)$
and then mapped all of the open sets back to $\V$ via~$p$.
\end{prop}
\begin{proof}
We need to show that the basic open sets in each topology are open in the other
topology.

Basic open sets in $\mathcal{T}_\V$ are induced by the intersection with
$p(\V)$ of open balls in $\RR^n$ \emph{centered at points of $p(\V)$}.  Basic
open sets in $\mathcal{T}_{p(\V)}$ are induced by the intersection with $p(\V)$
of open balls in $\RR^n$ \emph{centered anywhere}, so $\mathcal{T}_\V \subset
\mathcal{T}_{p(\V)}$.

To show the converse, let $B$ be a basic open set in $\mathcal{T}_{p(\V)}$
containing the vertex $v$, and let $r \in (0,\infty)$ and $c \in \RR^n$ be the
radius and center, respectively, of an $n$-ball which, when intersected with
$p(\V)$ gives~$p(B)$.  Since $v \in B$, we know that $\|c - p(v)\| < r$, so
the intersection with $p(\V)$ of the open $n$-ball of radius $r - \|c -
p(v)\|$ centered at $p(v)$ induces an open set in $\mathcal{T}_\V$ (by
definition) which contains $v$ (by design) and is contained in $B$ because if
$x \in \RR^n$ with $\|p(v) - x\| < r - \|c - p(v)\|$ we have
\begin{align*}
\|c - x\| = \|c - p(v) + p(v) - x\| \le \|c - p(v)\| + \|p(v) - x\| < r
\end{align*}
(see \vref{fig:openSetV}).
\begin{figure}[ht]
\hfil
\begin{tikzpicture}
\draw[dashed,fill=black!40] (2,2) node {$\bullet$} node [above] {$p(v)$} 
  circle (0.75cm);
\draw[vcurve] plot[smooth] coordinates{(0,0) (2,2) (4,1) (6,3)};
\draw[dashed] (2,3) node {$\bullet$} node [above left] {$c$} 
  circle (1.75cm);
\draw[yshift=0.1cm,snake=brace,segment amplitude=0.2cm] 
  (2,3) -- node [above,yshift=0.1cm] {$r$} (3.75,3);
\end{tikzpicture}
\hfil
\capt{Open sets in $\mathcal{T}_{p(\V)}$ are open in $\mathcal{T}_\V$.}
\label{fig:openSetV}
\end{figure}
\end{proof}

So the topology in $\V$ is familiar.  What about the one on $\E$?  Let
$\mathcal{T}_\Str$\IndexDef{$\mathcal{T}_\Str$@TStr (topology on $\Str$)}[
] and $\mathcal{T}_\Cr$\IndexDef{$\mathcal{T}_\Cr$@TCr (topology on $\Cr$)}[ ]
be the topologies on $\Str$ and $\Cr$ which come from seeing each of those
spaces as a subspace of a quotient space of $\V \cross \V$.  Then
let\IndexDef{$\mathcal{T}_{sqp}$@Tsqp (alt. topology on $\E$)}[ ]
$$\mathcal{T}_{sqp} = \left\{U_\Str \sqcup U_\Cr \colon U_\Str \in
\mathcal{T}_\Str, U_\Cr \in \mathcal{T}_\Cr \right\}.$$
\begin{prop}
$\mathcal{T}_{sqp}$ is a topology.
\end{prop}
\begin{proof}
We need to show that $\mathcal{T}_{sqp}$ contains the empty set, the set $\E$,
and that it is closed under arbitrary unions and finite intersections.

We note that the empty set is in $\mathcal{T}_{sqp}$ as it is the union of 
the empty set in $\mathcal{T}_\Str$ and the empty set in $\mathcal{T}_\Cr$.  
Similarly, the set $\E$ is in $\mathcal{T}_{sqp}$ as it is the union of the
sets $\Str \in \mathcal{T}_\Str$ and $\Cr \in \mathcal{T}_\Cr$.
If $\{U_{\Str\alpha} \sqcup U_{\Cr\alpha}\}$ is a collection of elements of
$\mathcal{T}_{sqp}$, then 
$$\bigcup \{U_{\Str\alpha} \sqcup U_{\Cr\alpha}\} = \bigcup U_{\Str\alpha}
\sqcup \bigcup U_{\Cr\alpha}$$
is also and $\mathcal{T}_{sqp}$.  Likewise, if $\{U_{\Str 1} \sqcup U_{\Cr 1},
\dotsc, U_{\Str n} \sqcup U_{\Cr n} \}$ is a finite collection of elements of 
$\mathcal{T}_{sqp}$, then
$$\bigcap \{U_{\Str i} \sqcup U_{\Cr i}\} = \bigcap U_{\Str i} \sqcup \bigcap
U_{\Cr i}$$
is in $\mathcal{T}_{sqp}$ as well.  So $\mathcal{T}_{sqp}$ is a topology
for~$\E$.
\end{proof}
\begin{prop}
$\mathcal{T}_\E$\IndexDef{$\mathcal{T}_\E$@TE (topology on $\E$)}[ ] equals
$\mathcal{T}_{sqp}$\IndexDef{$\mathcal{T}_{sqp}$@Tsqp (alt. topology on $\E$)}[
].
\end{prop}
\begin{proof}
We note that the basic open sets of $\mathcal{T}_{sqp}$ consist of all edges
in $\E$ whose ends fall within open sets in $\V$, $U_1$ and $U_2$.  We'll call
such a set $\nicefrac{U_1 \cross U_2}{\sim}$.

Let $U_1$ be an open set in $\mathcal{T}_\E$ containing edge $e_1 =
\{v_{11},v_{12}\}$.
Without loss of generality, we can take $e_1 \in \Str$.  Then there is a
basic open set, $B_1$, centered at $e_1$ and of radius $0 < r_1 < \infty$ that
is contained in~$U_1$.  Since $r_1$ is finite, we know that the elements of
$B_1$ are also in $\Str$.

Let $B_{11}$ and $B_{12}$ be open balls in $\V$ centered at $v_{11}$ and
$v_{12}$, respectively, both of radius $r_1$.  Then, if $e_2 =
\{v_{21},w_{22}\} \in \nicefrac{B_{11} \cross B_{12}}{\sim}$, we have 
\begin{align*}
d_\E(e_1,e_2) & = d_\vvt(e_1,e_2) \\
& = \min_{i \ne j \in \{1,2\}} \left\{ 
    d_{\V \cross \V}((v_{11},v_{12}),(v_{2i},v_{2j})) \right\} \\
& \le d_{\V \cross \V}((v_{11},v_{12}),(v_{21},v_{22})) \\
& = \max \{ d_\V(v_{11},v_{21}), d_\V(v_{12},v_{22}) \} < r_1.
\end{align*}
so $e_2 \in B_1$.  That is, the open set $\nicefrac{B_{11} \cross
B_{12}}{\sim} \cap \E$ is contained in~$B_1$.  So $B_1$ is open
in~$\mathcal{T}_{sqp}$.   

Conversely, let $U_{sqp}$ be an open set in $\mathcal{T}_{sqp}$ and let $e_1 =
\{v_{11},v_{12}\} \in U_{sqp}$.  Once again, we take $e_1$ to be in $\Str$.
Then there is some basic open set $\nicefrac{U_{11} \cross U_{12}}{\sim}$,
centered at $e_1$, that lies inside~$U_{sqp}$.  Let, $B_{11}$ (centered at
$v_{11}$) and $B_{12}$ (centered at $v_{12}$) be open balls in $\V$ of equal
radius,~$r$, contained within $U_{11}$ and $U_{12}$ respectively.

Let $e_2 = \{v_{21},v_{22}\} \in \E$ with $d_\E(e_1,e_2) < r$ (which of course
means that $e_2$ is also in $\Str$).  Then 
\begin{align*}
r & > d_\E(e_1,e_2) \\
& = d_\vvt(e_1,e_2) \\
& = \min_{i\ne j \in \{1,2\}} 
    d_{\V \cross \V}((v_{11},v_{12}),(v_{2i},v_{2j})) \\
& = \min_{i\ne j \in \{1,2\}}
    \max \left\{ d_\V(v_{11},v_{2i}), d_\V(v_{12},v_{2j}) \right\}.
\end{align*}
That means that either $v_{21} \in B_{11}$ and $v_{22} \in B_{12}$ or the other
way around (or both).  Either way, $e_2 \in \nicefrac{B_{11} \cross
B_{12}}{\sim}$, so there is an open $d_\E$ ball completely contained within
$U_{sqp}$ around every point in~$U_{sqp}$.  Thus $U_{sqp}$ is open
in~$\mathcal{T}_\E$.
\end{proof}

That's good news because it wasn't at all obvious that we could simply take two
open sets on $\V$ and know that the edges whose ends fall into those open sets
themselves form open sets on $\Str$ and $\Cr$.

\subsection{Compactness of $\V$}
\label{subsection:compactV}
Now that we have topologies on $\V$ and $\E$, let's talk for a moment
about~$\V$\Index{$\V$@V (vertices) compactness}[12!3].  Consider the tensegrity
shown in \ref{fig:CSplus}.
\begin{figure}[ht]
\hfil
\begin{tikzpicture}
\draw[cable] (0,0) -- (2,0) -- (2,2) -- (0,2) -- (0,0) -- (2,0);
\draw[strut] (0,0) node [vertex] {} -- (2,2) node [vertex] {};
\draw[strut] (0,2) node [vertex] {} -- (2,0) node [vertex] {};
\draw (-1,1) node [vertex] {} (1,3) node [vertex] {}
      (3,1) node [vertex] {} (1,-1) node [vertex] {};
\end{tikzpicture}
\hfil
\capt{A tensegrity with ``free vertices''.}
\label{fig:CSplus}
\end{figure}
This tensegrity is our (bar equivalent) crossed-square\Index{crossed square}[ ]
plus four extra vertices.  What effect do those extra vertices have on the bar
equivalence of this tensegrity?  None whatsoever.
\begin{claim}
Adding or removing vertices which have no edges connected to them does not 
change the $\X$-bar equivalence of a tensegrity.
\end{claim}
\begin{proof}
$\X$-bar equivalence and partial $\X$-bar equivalence are conditions
about~$Y(\X)$.  Adding or removing vertices, so long as it does not change the
edgeset, cannot change~$Y(\X)$.
\end{proof}
So, when it is convenient, we will consider $\V$ to have only vertices to 
which edges are attached.  That will give us an additional attribute for~$\V$.
\begin{claim}
If $\V$ consists only of vertices to which edges in $\E$ are attached, then
$\V$ is compact.
\end{claim}
\begin{proof}
Suppose every vertex in $\V$ has an edge attached to it.  Let $\{U_\alpha\}$ be
an open cover of~$\V$.  Then, by the work we've done above, the sets
$\nicefrac{U_{\alpha_1} \cross U_{\alpha_2}}{\sim} \cap \E$ (where
$U_{\alpha_1}$ and $U_{\alpha_2}$ are in $\{U_\alpha\}$) are open sets on
$\Str$ and $\Cr$.  Since every edge has to connect to vertices at both ends,
these open sets cover~$\E$.  Since $\E$ is compact, a finite number of them,
$\nicefrac{U_{\alpha_{11}} \cross U_{\alpha_{21}}}{\sim}, \dotsc,
\nicefrac{U_{\alpha_{1n}} \cross U_{\alpha_{2n}}}{\sim}$ will suffice to
cover~$\E$. 

We want to show that the open sets $U_{\alpha_{11}}, \dotsc, U_{\alpha_{1n}},
U_{\alpha_{21}}, \dotsc, U_{\alpha_{2n}}$ cover~$\V$.  Suppose $v_1 \in \V$.
Then, since all vertices of $\V$ have edges connected, there must be some $e =
\{v_1,v_2\}$ in~$\E$.  But then there must be some $\nicefrac{U_{\alpha_{1i}}
\cross U_{\alpha_{2i}}}{\sim}$ which contains~$e$.  So either $v_1 \in
U_{\alpha_{1i}}$ or $v_1 \in U_{\alpha_{2i}}$.  

So the $U_{\alpha_{11}}, \dotsc, U_{\alpha_{1n}}, U_{\alpha_{21}}, \dotsc,
U_{\alpha_{2n}}$ cover~$\V$.  Hence $\V$ is compact.
\end{proof}

Thus, while we never require that $\V$ be compact, we discover that if $\E$ is
compact, then $\V$ is effectively compact as well.

\subsection{The New Topologies made Visible}
Let's take a moment to look at a specific example.  Suppose we have the
tensegrity shown in \vref{fig:stadium}\subref{subfig:stadB}, whose the vertex
curve is shown in \ref{fig:stadium}\subref{subfig:stadA}.  The vertex curve
consists of two connected components: a central segment running from $(0,0)$ to
$(0,2)$ and a stadium curve running at distance $2$ around that segment.
Struts are connected orthogonally to the curve.  Cables are connected in the
piecewise linear fashion shown in \vref{fig:ETop} (essentially, endpoints on
the curves move with one fixed speed and endpoints on the straight segments
with another).
\begin{figure}[ht]
\hfil
\subfloat[The vertex curve, with basic open sets $a$ and $b$ of $\RR^n$ shown
(dashed lines), delineating basic open sets of $\V$ (shaded regions).  The 
edge (strut/cable) on which the open sets in \ref{fig:ETop} will be
centered is shown as a ``double helix''.]{
  \begin{tikzpicture}[rotate=90]
    \begin{scope}
      \draw (0,0) node [below,scale=0.5] {$(0,0)$} 
            (2,0) node [above,scale=0.5] {$(2,0)$};
      \useasboundingbox (-2.5,-5.5) rectangle (4,3.5);
      \draw[dashed,blue] (0,0.5) node[left] {$b$} arc (90:270:2.5cm) 
        arc (-90:90:2.5cm);
      \clip (0,0.5) arc (90:270:2.5cm) arc (-90:90:2.5cm);
      \draw[line width=0.75ex,blue!30!white,cap=rounded] (0,0) -- (2,0) 
        (2,-2) arc (-90:90:2cm) -- (0,2) arc (90:270:2cm) -- cycle;
      \draw[dashed,red!70!black] (0,-3) node [right] {$a$} arc (-90:90:1cm)
        arc (90:270:1cm);
      \clip (0,-1) arc (90:270:1cm) arc (-90:90:1cm);
      \draw[line width=1.5ex,red!70!black!30!white,cap=rounded] (0,0) -- (2,0) 
        (2,-2) arc (-90:90:2cm) -- (0,2) arc (90:270:2cm) -- cycle;
    \end{scope}
    \draw[vcurve] (0,0) -- (2,0) 
      (2,-2) arc (-90:90:2cm) -- (0,2) arc (90:270:2cm) -- cycle;
    \draw[->,yshift=0.2cm,thick] (0,2) node {$\bullet$} -- (0.75,2);
    \draw[snake=coil,segment aspect=0,segment amplitude=2pt,green!40!black] (0,0) -- (0,-2);
    \draw[snake=coil,segment aspect=0,segment amplitude=2pt,violet] (0,-2) -- (0,0);
  \end{tikzpicture}
  \label{subfig:stadA}
}
\hfil
\subfloat[The tensegrity with struts and cables.  Edges in the open sets shown
          in \ref{fig:ETop} are shaded.]{
  \begin{tikzpicture}[rotate=90]
    \useasboundingbox (-2,-2.5) rectangle (4,2.5);
    % Edges in the open set
    \filldraw[line width=1ex,green!70!black!30,cap=round,join=round]
      (0,0) -- (2,0) -- +(277:2cm) arc (277:270:2cm) -- 
      (0,-2) arc (270:202:2cm) -- cycle;
    \draw[line width=1ex,violet!30,cap=round]
      (22:-2cm)        -- (-90:-2cm)
      (26:-2cm)        -- (-86:-2cm)
      (30:-2cm)        -- (-82:-2cm)
      (34:-2cm)        -- (-78:-2cm)
      (38:-2cm)        -- (-74:-2cm)
      (42:-2cm)        -- (-70:-2cm)       
      (46:-2cm)        -- (-66:-2cm)
      (50:-2cm)        -- (-62:-2cm)
      (54:-2cm)        -- (-58:-2cm)
      (58:-2cm)        -- (-54:-2cm)
      (62:-2cm)        -- (-50:-2cm)
      (66:-2cm)        -- (-46:-2cm)
      (70:-2cm)        -- (-42:-2cm)
      (74:-2cm)        -- (-38:-2cm)       
      (78:-2cm)        -- (-34:-2cm)
      (82:-2cm)        -- (-30:-2cm)
      (86:-2cm)        -- (-26:-2cm)
      (0,-2)           -- (-22:-2cm)
      (0.075,-2)       -- (-18:-2cm)
      (0.15,-2)        -- (-14:-2cm)
      (0.225,-2)       -- (-10:-2cm)
      (0.3,-2)         -- (-6:-2cm)         
      (0.375,-2)       -- (-2:-2cm)
      (0.45,-2)        -- ( 2:-2cm)
      (0.525,-2)       -- ( 6:-2cm)
      (0.6,-2)         -- (10:-2cm)
      (0.675,-2)       -- (14:-2cm)
      (0.75,-2)        -- (18:-2cm)
      (0.825,-2)       -- (22:-2cm)
      (0.9,-2)         -- (26:-2cm)         
      (0.975,-2)       -- (30:-2cm)
      (1.05,-2)        -- (34:-2cm)
      (1.125,-2)       -- (38:-2cm)
      (1.2,-2)         -- (42:-2cm)
      (1.275,-2)       -- (46:-2cm)
      (1.35,-2)        -- (50:-2cm)
      (1.425,-2)       -- (54:-2cm)
      (1.5,-2)         -- (58:-2cm)         
      (1.575,-2)       -- (62:-2cm)
      (1.65,-2)        -- (66:-2cm)
      (1.725,-2)       -- (70:-2cm)
      (1.8,-2)         -- (74:-2cm)
      (1.875,-2)       -- (78:-2cm)
      (1.95,-2)        -- (82:-2cm)
      (2.0,-2)         -- (86:-2cm) % hmm
      (2,0) +(-86:2cm) -- (0,-2);
    \draw[vcurve] (0,0) -- (2,0) 
      (2,-2) arc (-90:90:2cm) -- (0,2) arc (90:270:2cm) -- cycle;
    % Struts along sides
    \foreach \x/\y in {0/2.0,0.3/1.7,0.6/1.4,0.9/1.1,1.2/0.8,1.5/0.5,1.8/0.2} {
      \draw[strut] (\x,0) -- (\x,-2);
      \draw[strut] (\y,0) -- (\y,2);
    }
    % Struts around ends
    \foreach \x in {-86,-70,-54,-38,-22,-6,10,26,42,58,74} {
      \draw[strut] (2,0) -- +(\x:2cm);
      \draw[strut] (0,0) -- +(\x:-2cm);
    }
    \draw[cable] 
      (2,0) +(-86:2cm) node (a) {} -- (0,-2)    
      (2,0) +(-70:2cm) node (b) {} -- (0.3,-2)
      (2,0) +(-54:2cm) node (c) {} -- (0.6,-2)  
      (2,0) +(-38:2cm) node (d) {} -- (0.9,-2)
      (2,0) +(-22:2cm) node (e) {} -- (1.2,-2)  
      (2,0) +(-6:2cm)  node (f) {} -- (1.5,-2)
      (2,0) +(10:2cm)  node (g) {} -- (1.8,-2)   
      (2,0) +(26:2cm)  node (h) {} -- (a)
      (2,0) +(42:2cm)  node (i) {} -- (b)
      (2,0) +(58:2cm)  node (j) {} -- (c)
      (2,0) +(74:2cm)  node (k) {} -- (d)
      (2,2)            node (l) {} -- (e)
      (1.7,2)          node (m) {} -- (f)
      (1.4,2)          node (n) {} -- (g)
      (1.1,2)          node (o) {} -- (h)
      (0.8,2)          node (p) {} -- (i)
      (0.5,2)          node (q) {} -- (j)
      (0.2,2)          node (r) {} -- (k)
      (-86:-2cm)       node (s) {} -- (l)
      (-70:-2cm)       node (t) {} -- (m)
      (-54:-2cm)       node (u) {} -- (n)
      (-38:-2cm)       node (v) {} -- (o)
      (-22:-2cm)       node (w) {} -- (p)
      (-6:-2cm)        node (x) {} -- (q)
      (10:-2cm)        node (y) {} -- (r)
      (26:-2cm)        node (z) {} -- (s)
      (42:-2cm)       node (aa) {} -- (t)
      (58:-2cm)       node (ab) {} -- (u)
      (74:-2cm)       node (ac) {} -- (v)
      (0,-2)                       -- (w)
      (0.3,-2)                     -- (x)
      (0.6,-2)                     -- (y)
      (0.9,-2)                     -- (z)         
      (1.2,-2)                     -- (aa)
      (1.5,-2)                     -- (ab)         
      (1.8,-2)                     -- (ac);
  \end{tikzpicture}
  \label{subfig:stadB}
}
\hfil
\caption{The stadium-curve tensegrity.}
\label{fig:stadium}
\end{figure}

In \ref{fig:ETop}, we show the sets $\Str$ and $\Cr$ for the 
stadium curve tensegrity.  These two sets are subsets of the space $\vvt$.  If
we plot $\V \cross \V$ as a square, the space $\vvt$ can be viewed as either
the triangle above the diagonal or the one below.  We've used one of each 
for the sake of space.  

The parameter for $\V$, as it runs along the edge of the square, first sweeps
out the component whose vertices lie along the center segment (from bottom to
top) and then covers the outer curve, starting at the point marked in
\ref{fig:stadium} and running in the direction of the arrow around the outside.

The open sets $a$ and $b$ of $\V$ are shown marked along the edges of the
square.  The struts are plotted in the upper triangle, while the cables are
shown in the lower one.  We've marked two basic open sets (in both cases, a
ball of radius $\sqrt{5}$ around the edge shown in
\ref{fig:stadium}\subref{subfig:stadA}).

\begin{figure}[ht]
\hspace*{\fill}
\begin{tikzpicture}[scale=0.5]
\begin{scope}[xshift=-7cm]
  \clip (-1,0) rectangle (0,18.566);
  \draw[line width=0.75ex,blue!30!white,cap=rounded] 
    (0,14.984) -- (0,9.779) (0,0) -- (0,1.5);
  \draw[line width=1.5ex,red!70!black!30!white,cap=rounded] 
    (0,13.294) -- (0,11.283);
\end{scope}
\begin{scope}[xshift=-7cm]
  \clip (0,18.566) rectangle (18.566,19.566);
  \draw[line width=0.75ex,blue!30!white,cap=rounded] 
    (14.984,18.566) -- (9.779,18.566) (0,18.566) -- (1.5,18.566);
  \draw[line width=1.5ex,red!70!black!30!white,cap=rounded] 
    (13.294,18.566) -- (11.283,18.566);
\end{scope}
\begin{scope}[xshift=-7cm]
  \draw[dashed,draw=green!40!black,fill=green!70!black!30]
    (0,0) -- (0,2) -- (1,2) -- (1,1) -- cycle
    (0,2) rectangle (1,3)
    (0,10.047) -- (0,18.566) -- (1,18.566) -- (1,14.656) -- (2,14.656) --
    (2,10.047) -- cycle
    (2,10.047) rectangle (3,14.656)
    (10.047,11.283) -- (10.047,18.566) -- (14.656,18.566) -- (14.656,14.656) --
    (11.283,11.283) -- cycle;
  \draw[dotted] (0,0) -- (18.566,18.566) -- (0,18.566) -- cycle;
  \draw[dotted] (0,2) -- (2,2) -- (2,18.566);
  \draw[yshift=0.1cm,snake=brace,segment amplitude=0.2cm] 
    (0,18.566) -- node[above,yshift=0.3cm,scale=0.7] {Middle Curve} (2,18.566);
  \draw[yshift=0.1cm,snake=brace,segment amplitude=0.2cm]
    (2,18.566) -- node[above,yshift=0.3cm,scale=0.7] {Outer Curve} 
    (18.566,18.566);
  \draw[xshift=-0.1cm,snake=brace,segment amplitude=0.2cm]
    (0,0) -- node[above,rotate=90,yshift=0.3cm,scale=0.7] {Middle Curve} (0,2);
  \draw[xshift=-0.1cm,snake=brace,segment amplitude=0.2cm]
    (0,2) -- node[above,rotate=90,yshift=0.3cm,scale=0.7] {Outer Curve} 
    (0,18.566);
  \draw[strut] (0,18.566) -- (0,12.283) -- (2,10.283) -- (2,4) -- 
    node [sloped,above,scale=0.8,pos=0.4] {$\Str$} 
    (0,2);
\end{scope}
\begin{scope}
  \clip (18.566,0) rectangle (19.566,18.566);
  \draw[line width=0.75ex,blue!30!white,cap=rounded] 
    (18.566,14.984) -- node [pos=0.2,right,blue] {$b$}  
                       node [pos=0.8,right,blue] {$b$} (18.566,9.779) 
    (18.566,0) -- node [right,blue] {$b$} (18.566,1.5);
  \draw[line width=1.5ex,red!70!black!30!white,cap=rounded] 
    (18.566,13.294) -- node [right,red!70!black] {$a$} (18.566,11.283);
\end{scope}
\begin{scope}
  \clip (0,-1) rectangle (18.566,0);
  \draw[line width=0.75ex,blue!30!white,cap=rounded] 
    (14.984,0) -- node [pos=0.2,below,blue] {$b$}  
                  node [pos=0.8,below,blue] {$b$} (9.779,0) 
    (0,0) -- node [below,blue] {$b$} (1.5,0);
  \draw[line width=1.5ex,red!70!black!30!white,cap=rounded] 
    (13.294,0) -- node [below,red!70!black] {$a$} (11.283,0);
\end{scope}
\begin{scope}
  \draw[dashed,draw=violet,fill=violet!30]
    (0,0) -- (2,0) -- (2,1) -- (1,1) -- cycle
    (2,0) rectangle (3,1)
    (10.047,0) -- (18.566,0) -- (18.566,1) -- (14.656,1) -- (14.656,2) --
    (10.047,2) -- cycle
    (10.047,2) rectangle (14.656,3)
    (11.283,10.047) -- (18.566,10.047) -- (18.566,14.656) -- (14.656,14.656) --
    (11.283,11.283) -- cycle;
  \draw[dotted] (0,0) -- (18.566,18.566) -- (18.566,0) -- cycle;
  \draw[dotted] (2,0) -- (2,2) -- (18.566,2);
  \draw[cable] (4.187,2) -- (7.910,4) -- 
               node [sloped,above,scale=0.8] {$\Cr$} 
               (10.283,6.421) -- 
               (12.283,10.144) -- (12.422,10.283) -- 
               (16.145,12.283) --
               (18.566,14.704);
  \draw[cable] (14.704,2) -- node [sloped,above,scale=0.8] {$\Cr$} 
               (18.427,4) -- (18.566,4.139);
\end{scope}
\end{tikzpicture}
\hspace*{\fill}
\capt[The edgeset of the stadium-curve tensegrity]{Here we show the sets $\Str$
and $\Cr$ for the stadium-curve tensegrity of \vref{fig:stadium} with one basic
open set for each connected component of $\E$ (shaded, with dashed outlines).
See the text for a more complete description of this figure.}
\label{fig:ETop}
\end{figure}

\subsection{About non-compact edgesets}
Let's talk again for a moment about the example from \vref{fig:openstruts}.
\Index[|(]{non-compact $\E$@E (edgeset) edgeset}[14 23!1]
Why doesn't it trouble us to exclude this case?  Frankly, because it is
indistinguishable from the case shown in \vref{fig:compactBoth}.
\begin{figure}[ht]
\hfil
\begin{tikzpicture}
\draw[vcurve] (0,0) circle (1cm) circle (3cm);
\foreach \x in {10,20,...,179} {
  \draw[strut] (\x:1cm) -- (\x:3cm);
  \draw[cable] (\x:-1cm) -- (\x:-3cm);
}
\draw[bar] (0:1cm) -- (0:3cm) (0:-1cm) -- (0:-3cm);
\end{tikzpicture}
\hfil
\capt[The tensegrity of \ref{fig:openstruts} ``compactified''.]{The tensegrity
of \ref{fig:openstruts} ``compactified''.  Bars have been added at $0$
and~$\pi$.  Now $\Str$ and $\Cr$ are both compact.}
\label{fig:compactBoth}
\end{figure}
Here, the strut at $0$ and the cable at $\pi$ have been replaced with bars.  
Now, both $\Str$ and $\Cr$ are compact\footnote{As we noted in 
the footnote on page \pageref{footnote:compact}, $\V \cross \V$ is compact and
since $\Str$ and $\Cr$ are now closed, their preimages in 
$\V \cross \V$ are compact and hence they are themselves compact.} and all
of our results hold.  

This tensegrity, however, works exactly as the one in \ref{fig:openstruts} did,
because the elements of $\Vf(\V)$ are continuous.  Any vector field which
expanded the strut at $0$, $s_0$, would have had to expand all edges in an open
set around~$s_0$.  But any such open set contains cables, which cannot expand.  

Thus $s_0$ cannot change length and it might as well be a bar.  By the same
type of argument, the cable at $\pi$ is effectively a bar as well.
\Index[|)]{non-compact $\E$@E (edgeset) edgeset}[14 23!1]

\begin{prop}
For any $V \in \Vf(\V)$, $YV \in C(\E)^+ \setminus \{\zero\}$ if and only if
$Y_c V \in C(\E_c)^+ \setminus \{\zero\}$, where $\E_c$ is the disjoint union
of the closures of $\Str$ and $\Cr$\Index{$\E_c$@Ec ($\overline{\Str} \sqcup
\overline{\Cr}$)}[ ] (that is, $\E_c = \overline{\Str} \sqcup \overline{\Cr}$),
and $Y_c$ is the rigidity operator which results from replacing $\E$
with~$\E_c$.
\end{prop}
\begin{proof}
If $Y_c V \in C(\E_c)^+$, then $YV$, which is the restriction of $Y_c V$ to
$\E$, must be in~$C(\E)^+$.  If $Y_c V \ne \zero$, then there must be some edge
$e_c \in \E_c$ for which $Y_c V(e_c) > 0$.  But because $Y_c V$ is continuous,
there is some open set on which $Y_c V > 0$, and that open set must intersect
$\E$ (as $e_c$ is in the closure of either $\Str$ or $\Cr$).  So $YV$ must 
also be nonzero.  That's one direction.  

Conversely, assume that $YV \in C(\E)^+$.  Then its image must be contained in
the interval $[0,\infty)$.  We want to show that the $\im Y_c V$ is also
contained in that interval.  Suppose, to the contrary, there is some edge $e_n$
such that $Y_c V(e_n) < 0$.  Now $e_n$ is either in $\overline{\Str}$ or in
$\overline{\Cr}$.  Without loss of generality, take it to be in
$\overline{\Str}$.

Since $Y_c V$ is continuous, there must be an open set $U_-$ containing $e_n$
on which $Y_c V$ is strictly negative.  But since $e_n \in \overline{\Str}$, we
must have $U_- \cap \Str \ne \varnothing$, and $Y_c \ge 0$ on $\Str$, so it
cannot be strictly negative on $U_-$.  That contradiction shows that $\im Y_c V
\in [0,\infty)$.

Also, if $YV \ne \zero$, then on some $e \in \E$, $YV(e) \ne 0$ and hence
$Y_cV(e) \ne \zero$, so $Y_c V \ne \zero$.
\end{proof}

\subsection{Discrete Topologies}
Just to finish things off well, let's end this section by showing that in the
finite case, the topologies on $\V$ and $\E$ are discrete.
\Index[|(]{discrete topology}[2!1 ]
\begin{prop}
\label{prop:discrete}
If $\V$ and $\E$ are finite sets, the topologies on them are discrete.
\end{prop}
\begin{proof}
If $\V$ is a finite set, then $p(\V)$ is a finite set of points in~$\RR^n$.  
Since $p$ is 1-1, if $v_1 \ne v_2$, then $d_\V(v_1,v_2) = \|p(v_1)
- p(v_2)\| > 0$, so open balls in $\V$ centered on $v_1$ and $v_2$ and having
radius less than half of $d_\V(v_1,v_2)$ separate $v_1$ and~$v_2$.  For any
given vertex, $v$, the (finite) intersection of all such balls is an open set
which contains $v$ and no other element of~$\V$.  Thus we have the discrete 
topology on~$\V$.

Similarly, if $e_1, e_2 \in \E$, then either they are both struts (or both
cables) or there is one of each.  Since the sets $\St$, $\C$ and $\B$ are
pairwise disjoint, no two struts connect the same pair of vertices and no two 
cables connect the same pair of vertices.  So if $e_1$ and $e_2$ are of the
same type, they are some positive distance apart.  But if they are of different
type, that is even more true.  Using the same process as we did for $\V$, we
see that this topology is also discrete.
\end{proof}
\Index[|)]{$\E$@E (edgeset)}[ ]
\Index[|)]{discrete topology}[2!1 ]

\section{\texorpdfstring{$C(\E)$}{C(E)} and \texorpdfstring{$C^*(\E)$}{C*(E)}}
\label{section:CE}
Since $C(\E)$ and $C^*(\E)$ are no longer finite vector spaces,
the Euclidean norm\Index{Euclidean norm} will no longer serve.%
\IndDefBeg{$C(\E)$@C(E) (cont. func. on $\E$)}[ ]%
\IndDefBeg{$C^*(\E)$@C*(E) (top. dual of $C(\E)$)}[ ]
We'll put the $\sup$
norm\Index{$\sup$@sup norm}[2!1 ] on $C(\E)$\IndexDef{norm on
$C(\E)$@C(E)}[1!23] and the operator norm\Index{operator norm}[2!1 ]
on~$C^*(\E)$\IndexDef{norm on $C^*(\E)$@C*(E)}[1!23].  That is, for $f \in
C(\E)$ and $\mu \in C^*(\E)$, we have
\begin{align*}
\|f\| & = \sup_{e \in \E} |f(e)|\text{ and } \\
\|\mu\| & = \sup_{\|f\| = 1} |\mu f|.
\end{align*}
We can induce a metric on each\IndexDef{metric on
$C(\E)$@C(E)}[1!23]\IndexDef{metric on $C^*(\E)$@C*(E)}[1!23] by considering
$\|x-y\|$ to be the distance between $x$ and $y$ and establish on each the
associated metric topology\Index{metric topology}[2!1 ].  

One quick definition:
\begin{definition}[see, for example, \ocite{MR1681462}*{p.\ 132}]
If $X$ is a topological space, the \emph{support}\IndexDef{support of a
function}[1!234] of a function $f\mcol X \to \RR$, denoted $\supp f$, is the
smallest closed set $S \subset X$ such that $f$ is zero everywhere on the
complement $X \setminus S$. In other words, $S$ is the closure of the set of
all points $x \in X$ for which $f(x) \ne 0$.
\end{definition}
And now we're prepared to talk about the interiors of $C(\E)^+$\Index{interior
of $C(\E)^+$@C(E)+ and $C(\E)^-$@C(E)- {(orthants of $C(\E)$)}}[1!2345
3456!1] and $C(\E)^-$.
\begin{lemma}
\label{lemma:intCEp}
The interiors of the nonnegative and nonpositive orthants of $C(\E)$, that is
$\INT C(\E)^+$ and $\INT C(\E)^-$, consist of the strictly positive and
strictly negative continuous functions from $\E$ to $\RR$, respectively.
\end{lemma}
\begin{proof}
We will prove here that $\INT C(\E)^+$ consists of the strictly positive
elements of~$C(\E)$.  The proof for $\INT C(\E)^-$ is identical except for
sign.  We need to establish our result using the Euclidean norm\Index{Euclidean
norm} for the finite case and the $\sup$ norm\Index{$\sup$@sup norm}[2!1
] for the general case.  

Let $f$ be a strictly positive, continuous function on~$\E$.  Then, since $\E$
is compact, $f$ achieves a (strictly positive) minimum on~$\E$.  Now let $g \in
C(\E)$ with $\|f - g\| = d < \min f$.  Of course $g$ also achieves a minimum on
$\E$ at some edge~$e$.  Then, under either norm,
\begin{align*}
\min g = g(e) = f(e) - (f(e) - g(e)) \ge f(e) - \|f - g\| \ge \min f - d > 0,
\end{align*}
so $g$ is strictly positive.  Thus $f$ is in $\INT C(\E)^+$.

Conversely, let $f \in \INT C(\E)^+$.  Then there exists some $\delta > 0$ such
that $g \in C(\E)^+$ for all $\|f - g\| < \delta$.  Let $e \in \E$ be where
$f$ achieves its minimum.  Suppose $f(e) < \delta$.  

In the finite case, let $g = f$ everywhere except at $e$ and set $g(e) =
\frac{f(e) - \delta}{2}$.  In the infinite case, select an open set $U$
containing $e$ and let $g_e \in C(\E)^+$ with $\sup g_e = g_e(e) = \frac{f(e) -
\delta}{2}$ and $\supp g_e \subset U$, and set $g = f - g_e$.  Then in both
cases, we have $\|f - g\| = \frac{f(e) + \delta}{2} \in (0,\delta)$, but $g(e)
= \frac{f(e) - \delta}{2} < 0$.  That contradiction means that we must have
$f(e) = \min f \ge \delta > 0$, and $f$ is strictly positive.
\end{proof}

\IndDefEnd{$C(\E)$@C(E) (cont. func. on $\E$)}[ ]
Next we turn our attention from $C(\E)$ to its topological dual, $C^*(\E)$,
which, we remember, consists of the continuous linear functionals on $C(\E)$
(Definition \vref{definition:topologicalDual}).  We'll want to relate the
elements of $C^*(\E)$ not only to $C(\E)$ but also to $\E$ itself.

First, we should note that for linear functionals, ``continuous'' (Definition
\vref{definition:continuous}) and ``bounded'' (Definition
\vref{definition:bounded}) are connected:
\begin{theorem}
\label{theorem:ContBounded}
For a linear transformation $\Lambda$ of a normed linear space $X$ into a
normed linear space $Z$, each of the following three conditions implies the
other two:
\begin{enumerate}[label=(\alph*)]
\item $\Lambda$ is bounded.
\item $\Lambda$ is continuous.
\item $\Lambda$ is continuous at one point of~$X$.
\end{enumerate}
\end{theorem}
\begin{proof}
See \ocite{MR924157}*{p.\ 96}.
\end{proof}

Next, a quick definition.
\begin{definition}[\ocite{MR924157}*{p.\ 47}]
A measure is called \emph{regular}\IndexDef{regular measure}[2!1 ] if, for any
measurable set $E$,
$$\mu(E) = \inf \{\mu(U) : E \subset U, U \text{ open}\}$$
and
$$\mu(E) = \sup \{\mu(K) : K \subset E, K \text{ compact}\}.$$
\end{definition}
And now we're ready for a theorem that will allow us to understand the elements
of $C^*(\E)$ not only as continuous linear functionals on $C(\E)$, but also as
measures on~$\E$.
\begin{theorem}[Riesz Representation Theorem]
\label{theorem:Riesz}
If $X$ is a compact Hausdorff space, then every bounded linear functional
$\Phi$ on $C(X)$ is represented by a unique regular Borel measure $\mu$, in
the sense that\IndexDef{Riesz Representation Theorem}[3!12 ]
$$\Phi f = \int_X f \dmu$$
for every $f \in C(X)$.
\end{theorem}
\begin{proof}
See \ocite{MR924157}*{p.\ 130}.  Rudin actually gives the theorem for locally
compact spaces and requires the functions to be in $C_0(X)$ (the continuous 
functions which vanish at infinity), but for a compact set, $C_0(X) = C(X)$.
Also, he returns a regular complex measure $\nu$, but our $\mu = |\nu|$, the
total variation of $\nu$, which is regular since $\nu$ is regular
\cite{MR924157}*{pp.\ 70, 116--117, 130}.
\end{proof}

Since we are dealing in measures, we'll want to talk about ``strictly positive
measures'' and ``semipositive measures'', and since the Riesz Representation
theorem gives it to us, we'll expect regularity.  We'll take the following
definitions.
\begin{definition}
\label{definition:SPMandSM}
A \emph{strictly positive measure}\IndexDef{{strictly positive} measure} is a
regular measure $\mu$ on a set $X$ such that for every nonempty, open $U
\subset X$, $\mu(U) > 0$.  A \emph{semipositive measure}\IndexDef{semipositive
measure} is a nonzero regular measure $\mu$ on a set $X$ such that for every
nonempty, open $U \subset X$, $\mu(U) \ge 0$.

Following general practice, we will use the term \emph{positive
measure}\Index{positive measure}[2!1 ] to refer to what might more accurately
be called ``nonnegative measure''\Index{nonnegative measure}.  The
regular positive measures differ from the semipositive measures only in that
the zero measure is considered a regular positive measure.
\end{definition}

It seems now that we have competing definitions for ``strictly positive'' and
``semipositive'' for those measures that are stresses of our tensegrity (in
Definition \vref{definition:stress} they were defined in terms of $\mu(f)$
rather than $\mu(U)$).  In Proposition \vref{prop:Positive}, we'll show that
the definitions are equivalent, but we have some work to do before we get
there.

First, we need to define the support of a measure.
\begin{definition}[see, for example, \ocite{Morrison}*{p.\ 152}]
\Index{support {of a measure}}[1!2]
If $\mu$ is a measure on a set $X$ (with its associated $\sigma$-algebra), then
the \emph{support} of $\mu$, denoted $\supp \mu$ is the set of $x \in X$ such
that $|\mu(U)| > 0$ for all open neighborhoods $U$ of~$x$.
\end{definition}
There are a few lemmas about supports we will want.  The proof for the first
was found on Wikipedia \ycite{wiki:SupM}.
\begin{lemma}
\label{lemma:closedsupport}
The support of a measure $\mu$ on a topological set $X$ is a closed set.
\Index{support {of a measure} {is closed}}[1!2!3]
\end{lemma}
\begin{proof}
Let $\mu$ be a measure on a topological set $X$ and let $\{x_i \in \supp \mu\}$
be a sequence which approaches some limit $x \in X$.  If $U \subset X$ is any
open set which contains $x$, then, since $x$ is the limit of the $x_i$, $U$
must intersect~$\{x_i\}$.  But that means $U$ is an open set containing some
element of $\supp \mu$ and hence $|\mu(U)| > 0$.  So $x \in \supp \mu$.
\end{proof}
\begin{lemma}
\label{lemma:mustIntersect}
If $\mu$ is a regular positive measure on a topological set $X$ and $A$ is an
open subset of $X$ such that $\supp \mu \cap A = \varnothing$, then $\mu(A) =
0$.
\end{lemma}
\begin{proof}
We note first that if $x \in A$, then $x \notin \supp \mu$, so there must be
some open neighborhood $U_x$ of $x$ with measure~$0$.  If $K$ is any compact
subset of $A$, then $K$ can be covered by a finite number of these~$U_x$.  
So $0 \le \mu(K) \le \sum_{x \in K} \mu(U_x) = 0$.  But since $\mu$ is
regular, $\mu(A) = \sup\{\mu(K) : K \subset A, K \text{ compact}\}$, so 
$\mu(A) = 0$.
\end{proof}

Now we can show how positivity relates to the support of a measure.
\begin{lemma}
\label{lemma:spisE}
$\mu$ is a strictly positive measure on $\E$ if and only if $\mu$ is a
regular positive measure with $\supp \mu = \E$.
\end{lemma}
\begin{proof}
Let $\mu$ be strictly positive and let $e \in \E$.  Then, for every
neighborhood $U$ of $e$, $\mu(U) > 0$ (after all, $\mu$ is strictly positive),
so $e \in \supp \mu$.  

Conversely, let $\mu$ be a regular positive measure on $\E$ with $\supp \mu =
\E$ and let $U \subset \E$ be nonempty and open.  Then $U \cap \supp \mu = U
\ne \varnothing$, so $\mu(U) > 0$.
\end{proof}
\begin{lemma}
\label{lemma:spins}
$\mu$ is a semipositive measure on $\E$ if and only $\mu$ is a regular positive
measure with $\supp \mu \ne \varnothing$.
\end{lemma}
\begin{proof}
Let $\mu$ be semipositive.  Suppose, to the contrary, that $\supp \mu =
\varnothing$.  Then, $\E$ is an open set disjoint from $\supp \mu$, so by Lemma
\ref{lemma:mustIntersect}, $\mu(\E) = 0$ and since $\mu$ is semipositive, we
must have $\mu(U) = 0$ for every $U \subset \E$.  This contradiction shows that
$\supp \mu \ne \varnothing$.

Conversely, let $\mu$ be a regular positive measure on $\E$ with $\supp \mu \ne
\varnothing$ and let $e \in \supp \mu$.  Then, since $\mu$ is positive, $\mu(U)
\ge 0$ for all $U \subset \E$.  Furthermore, $\E$ is an open set which contains
$e$ and hence $\mu(\E) > 0$.
\end{proof}

We'll need to know about the Jordan decomposition of a measure, which allows a
signed measure to be ``split'' into positive and negative parts and the Urysohn
Lemma, which allows us to create certain continuous maps.  Here are the
details:
\begin{definition}[\ocite{MR1681462}*{p.\ 87}]
\label{definition:MutuallySingular}
Two signed measures $\mu$ and $\nu$ on a measurable space $(X,\mathcal{M})$
are \emph{mutually singular} if there exist $E, F \in \mathcal{M}$ such that
$E \cap F = \varnothing$, $E \cup F = X$, $E$ is null for $\mu$ and $F$ is
null for~$\nu$.  We denote this $\mu \perp \nu$.
\end{definition}
\begin{theorem}[Jordan Decomposition Theorem]
If $\nu$ is a signed measure, there exist unique positive measures $\nu^+$ and
$\nu^-$ such that $\nu = \nu^+ - \nu^-$ and $\nu^+ \perp \nu^-$.
\end{theorem}
\begin{proof}
See \ocite{MR1681462}*{p.\ 87}.
\end{proof}
\begin{definition}[\ocite{Munkres}*{p.195}]
A space $X$ is said to be \emph{normal} if, for each pair $A$, $B$ of 
disjoint closed sets of $X$, there exist disjoint open sets containing $A$ and
$B$, respectively.
\end{definition}
\begin{lemma}
\label{lemma:normal}
Every metric space is normal.
\end{lemma}
\begin{proof}
For this proof, we'll take some guidance from a hint in \ocite{MR1984838}*{p.\
52}.  We let $X$ be a metric space with metric $d$ and let $A$ and $B$ be
disjoint closed subsets of $X$.  For each $a \in A$, we know that since $B$ is
closed and $A$ and $B$ are disjoint, $a$ is not in the closure of $B$.  So
there must be some strictly positive infimal distance $d_a$ between $a$ and all
of $B$.  Let $U_a$ be the open ball of radius $\nicefrac{d_a}{2}$ centered at
$a$.  Clearly, $U_a \cap B = \varnothing$.  Furthermore, the collection of all
such $U_a$ covers $A$.  Similarly, we can cover $B$ with $U_b$ where the radius
of a given $U_b$ is half the minimum distance between $b$ and all of $A$.  

Now let $U_A = \bigcup U_a$ and $U_B = \bigcup U_b$.  We'd like to show that
$U_A \cap U_B = \varnothing$.  Suppose, to the contrary, that some element 
$x \in X$ lies in both $U_A$ and $U_B$.  Since $x$ is in a union of sets, it
must be in one of the sets.  That is, there must be some $a \in A$ and some
$b \in B$ such that $x \in U_a \cap U_b$.  Let $d_a$ be the infimal distance
between $a$ and $B$ and $d_b$ the infimal distance between $b$ and $A$.  Then
the distance between the two points $d(a,b)$ has to be at least as large as
$d_a$ and $d_b$.  But by the triangle inequality
$$\max \{ d_a,d_b \} \le d(a,b) \le d(a,x) + d(x,b) < \frac{d_a}{2} +
\frac{d_b}{2} \le 2 \frac{\max \{ d_a,d_b \}}{2} = \max \{ d_a,d_b \}.$$
That contradiction shows us that there can be no such $x$ and hence $U_A$ and
$U_B$ are disjoint.  So every metric space is normal.
\end{proof}
\begin{theorem}[Urysohn Lemma]
\label{theorem:Urysohn}
\Index{Urysohn Lemma}[ ]
Let $X$ be a normal space; let $A$ and $B$ be disjoint closed subsets of~$X$.
Let $[a,b]$ be a closed interval in the real line.  Then there exists a
continuous map $f\mcol X \to [a,b]$ such that $f(x) = a$ for every $x$ in
$A$, and $f(x) = b$ for every $x$ in~$B$.
\end{theorem}
\begin{proof}
See, for example, \ocite{Munkres}*{p.\ 207}.
\end{proof}

\begin{prop}
\label{prop:Positive}
Let $\mu \in C^*(\E)$ with $\mu \ne \zero$.  Then $\mu(f) > 0$ (resp.\ $\mu(f)
\ge 0$) for all nonzero $f \in C(\E)^+$ if and only if $\mu(U) > 0$ (resp.\
$\mu(U) \ge 0$) for all nonempty, open $U \subset \E$.
\end{prop}
\begin{proof}
We'll start by showing that if $\mu(f) \ge 0$ for all $f \in C(\E)^+$, then
$\mu(U) \ge 0$ for all open $U \subset E$.

Let $\mu \in C^*(\E)^+$ such that $\mu \ne \zero$ and $\mu(f) \ge 0$ for all $f
\in C(\E)^+$.  Let $\mu = \mu^+ - \mu^-$ be the Jordan decomposition of~$\mu$.
By Lemma \ref{lemma:closedsupport}, $\supp \mu^+$ is closed, so $U_s = \E
\setminus \supp \mu^+$ is open.  If $U_s = \varnothing$, then by Lemma
\ref{lemma:spisE}, $\mu^+$ is a strictly positive measure.  Since $\mu^+ \perp
\mu^-$ in the Jordan decompositions, then the set $F$ from Definition
\ref{definition:MutuallySingular} must be the empty set.  So $\mu^-$ is the
zero measure and $\mu$ is strictly positive as well and we are done.

Suppose, then, that $U_s \ne \varnothing$ and, further, that there exists some
$U_- \subset \E$ such that $\mu(U_-) < 0$.  If $U_-$ were completely contained
in $\supp \mu^+$, we'd have $\mu^+(U_-) > 0$ and $\mu^-(U_-) = 0$ and thus
$\mu(U_-) > 0$.  Since that isn't the case, we know that $U_-$ reaches (or lies
completely) outside of $\supp \mu^+$.  That is, $U_i = U_- \cap U_s \ne
\varnothing$ and $\mu(U_i) = -\mu^-(U_i) < 0$.  

Since $-\mu^-(U_i) < 0$, $U_i$ must intersect the support of $\mu^-$
nontrivially (by Lemma \vref{lemma:mustIntersect}) , so define $E_i = \supp
\mu^- \cap \overline{U_i}$, which means that $\mu(E_i) < 0$.  Now $E_i$ and
$\supp \mu^+$ are disjoint closed sets in the metric (hence normal, by
Lemma \ref{lemma:normal}) space $\E$, so by the Urysohn Lemma, there is a
continuous function $f$ which takes the value $1$ on $E_i$ and the value $0$ on
$\supp \mu^+$ and whose values lie in the interval $[0,1]$ elsewhere.  Clearly
$f \in C(\E)^+$.

So now we have $\mu(f) = \int_\E f \dmu = \int_{\supp f} f \dmu$.  But by
design, $\supp f \cap \supp \mu^+ = \varnothing$, so we have $\mu(f) =
-\int_{\supp f} f \, d\mu^- \le -\int_{E_i} f \, d\mu^- = \mu^-(E_i) < 0$.

That contradicts our hypothesis that $\mu(f) \ge 0$ for all $f \in C(\E)^+$,
so there can be no such~$U_-$.  Hence $\mu(U) \ge 0$ for all $U \subset \E$.

Next we show that if $\mu(f) > 0$ for all nonzero $f \in C(\E)^+$, then $\mu(U)
> 0$ for all nonempty, open $U \subset \E$.

Let $\mu \in C^*(\E)^+$ such that $\mu(f) > 0$ for all nonzero $f \in C(\E)^+$.
Now we know (by the work we just finished) that $\mu(U) \ge 0$ for all
nonempty, open $U \subset \E$, but suppose that for one such set, $\mu(U) = 0$.
Select any nonnegative, continuous function $f_U$ on $\E$ with $|f_U| = 1$ and
$\supp f_U \subset U$.  

Then, defining $\one$ to be the constant function $\one(e) \equiv
1$\IndexDef{$\one$@one (constant function 1)}[ ] on $\E$, we have
$$0 < \mu(f_U) = \int_\E f_U \dmu = \int_{\supp f_U} f_U \dmu \le \int_{\supp
f_U} \one \dmu \le \mu(U) = 0.$$
That contradiction show that there can be no such~$U$.

Finally, we work the other direction.  The material given here is for the 
strict inequalities, but it works \textit{a fortiori} if they are not strict.

Let $\mu \in C^*(\E)^+$ such that $\mu(U) > 0$ for all nonempty, open $U
\subset \E$ and let $f$ be a nonzero element of~$C(\E)^+$.  Since $f$ is
nonzero, there must be some $e_+ \in \E$ such that $f(e_+) > 0$.  As $f$ is
continuous $U_+ = f^{-1}( (\nicefrac{f(e_+)}{2}, \infty) )$ is nonempty and
open.  Thus 
$$\mu(f) = \int_{\supp f} f \dmu \ge \int_{U_+} f \dmu \ge
\int_{U_+} \frac{f(e_+)}{2} \dmu = \frac{f(e_+)}{2} \mu(U_+) > 0.$$
\end{proof}
\IndDefEnd{$C^*(\E)$@C*(E) (top. dual of $C(\E)$)}[ ]

\section{\texorpdfstring{$\Vf(\V)$}{Vf(V)} and \texorpdfstring{$Y$}{Y}}
\label{section:VfVY}
$\Vf(\V)$ consists of the continuous vector fields from $\V$ to~$\RR^n$.  
\IndexDef{$\Vf(\V)$@Vf(V) (variations)}[ ]
\IndexDef{$Y$@Y (rigidity operator)}[ ]
\IndexDef{rigidity operator}[ ]
We'll want it to be a topological space, so we'll give it the $\sup$
norm\Index{$\sup$@sup norm}[2!1 ].  That is, for $V \in \Vf(\V)$, $$\|V\| =
\sup_{v \in \V} \|V(v)\|.$$ Then we'll use the norm to give us a metric and
endow $\Vf(\V)$ with the associated metric topology.

When we first met $Y$, back on page \pageref{definition:Y}, we defined it as
a map from $\Vf(\V)$ to $C(\E)$, but much of what we have done since then has
thought of $Y$ simply as a matrix.  That was fine for the finite world, but 
we need to get back to thinking of $Y$ as a map.  As a map, $Y$ has some nice
attributes.  
\begin{definition}[\ocite{MR924157}*{p.\ 96}]
\label{definition:bounded}
Consider a linear transformation $\Lambda$ from a normed linear space $X$ into
a normed linear space $Z$, and define its \emph{norm} by $$\|\Lambda\| = \sup
\{ \|\Lambda x \| : x \in X, \|x \| \le 1 \}.$$
If $\|\Lambda\| < \infty$, then $\Lambda$ is called a \emph{bounded linear
transformation}.
\end{definition}

In a moment we'll see that our map $\|Y\|$ is bounded.  First, though, we have
a lemma about edge lengths.
\begin{lemma}
\label{lemma:LenIsCont}
The edge length $\|p(v_1) - p(v_2)\|$ is continuous in $e \in \E$.
\end{lemma}
\begin{proof}
Let $\varepsilon > 0$.  Select $\delta = \nicefrac{\varepsilon}{2}$.  Then,
whenever $e_1 = \{v_{11},v_{12}\}, e_2 = \{v_{21},v_{22}\} \in \E$ with
$d_\E(e_1,e_2) < \delta$, we have either 
\begin{equation}
\|p(v_{11}) - p(v_{21})\| < \delta \text{ and } \|p(v_{12}) - p(v_{22})\| <
\delta
\label{eq:boundsA}
\end{equation}
or
\begin{equation}
\|p(v_{11}) - p(v_{22})\| < \delta \text{ and } \|p(v_{12}) - p(v_{21})\| <
\delta.
\label{eq:boundsB}
\end{equation}
But then
\begin{equation}
\begin{split}
\big| \| p(v_{11}) - p(v_{12}) \| - \| p(v_{21}) - p(v_{22}) \| \big|
& \le \big| \| p(v_{11}) - p(v_{12}) - p(v_{21}) + p(v_{22}) \| \big|\\
& \le \| p(v_{11}) - p(v_{21}) \| + \|p(v_{22}) - p(v_{12})\| 
\end{split}
\label{eq:boundsC}
\end{equation}
and also
\begin{equation}
\begin{split}
\big| \| p(v_{11}) - p(v_{12}) \| - \| p(v_{21}) - p(v_{22}) \| \big|
& = \big| \| p(v_{11}) - p(v_{12}) \| - \| p(v_{22}) - p(v_{21}) \| \big|\\
& \le \big| \| p(v_{11}) - p(v_{12}) - p(v_{22}) + p(v_{21}) \| \big|\\
& \le \| p(v_{11}) - p(v_{22}) \| + \|p(v_{21}) - p(v_{12})\|.
\end{split}
\label{eq:boundsD}
\end{equation}
Since both \ref{eq:boundsC} and \ref{eq:boundsD} are true, and either 
\ref{eq:boundsA} or \ref{eq:boundsB} is true as well, we have
$$\big| \| p(v_{11}) - p(v_{12}) \| - \| p(v_{21}) - p(v_{22}) \| \big| \le
\delta + \delta = \varepsilon$$
and we are done.
\end{proof}
\begin{prop}
\label{prop:Ybounded}
For any rigidity operator $Y$, $\|Y\| < \infty$.
\end{prop}
\begin{proof}
Let $V \in \Vf(\V)$ with $\|V\| \le 1$.  Then $\sup_{v \in \V} \|V(v)\| \le 1$.
So for any given edge, $e = \{v_1,v_2\}$, we have $\|V(v_1) - V(v_2)\| \le 2$.
By Lemma \ref{lemma:LenIsCont}, edge length is a continuous function on
$\E$, and $\E$ is compact, so there is some edge $e_L$ of longest length $L <
\infty$.  Thus, whenever $\|V\|\le 1$, we have
\begin{align*}
\|YV\| & = \sup_{\{v_1,v_2\} \in \E} |(V(v_1) - V(v_2)) \cdot 
  (p(v_1) - p(v_2))| \le 2 L.
\end{align*}
Since $\|Y\|$ is defined as
\begin{align*}
\|Y\| & = \sup \{ \|YV\| : V \in \Vf(\V), \|V\| \le 1\},
\end{align*}
we have $\|Y\| \le 2 L < \infty$.
\end{proof}

We'd like to show that $Y$ is a continuous map from $\Vf(\V)$ to~$C(\E)$.
To do that we'll need a little machinery.  We need a theorem about continuity
in quotient spaces and a lemma about a continuous function on $\V \cross \V$.
\begin{theorem} 
\label{theorem:contQuot}
Let $\pi\mcol X \to W$ be a quotient map.  Let $Z$ be a space and let $g\mcol X
\to Z$ be a map that is constant on each set $\pi^{-1}({w})$, for $w \in W$.
Then $g$ induces a map $f\mcol W \to Z$ such that $f \circ \pi = g$.  The
induced map $f$ is continuous if and only if $g$ is continuous.
$$
\xymatrix{
X \ar[d]_\pi \ar[dr]^g \\
W \ar@{.>}[r]_f & Z
}
$$
\end{theorem}
\begin{proof}
See \ocite{Munkres}*{p.\ 142}.
\end{proof}
\begin{lemma}
\label{lemma:diffCont}
If $f(v)$ is a continuous function from $\V$ to $\RR^n$, then $f(v_1) - f(v_2)$
is a continuous function from $\V \cross \V$ to~$\RR^n$.
\end{lemma}
\begin{proof}
Let $\varepsilon > 0$.  Since $f(v)$ is continuous, there exists some $\delta >
0$ such that whenever $v_1$ and $v_2$ are in $\V$ with $d_\V(v_1,v_2) <
\delta$, we have $\|f(v_1) - f(v_2)\| < \nicefrac{\varepsilon}{2}$.  Now,
whenever $(v_1,v_2)$ and $(v_3,v_4)$ are in $\V \cross \V$ with $d_{\V \cross
\V}((v_1,v_2),(v_3,v_4)) < \delta$, we have $d_\V(v_1,v_3) < \delta$ and
$d_\V(v_2,v_4) < \delta$.  But then
\begin{align*}
\|(f(v_1) - f(v_2)) - (f(v_3) - f(v_4))\|
& = \|f(v_1) - f(v_2) - f(v_3) + f(v_4)\| \\
& = \|(f(v_1) - f(v_3)) + (f(v_4) - f(v_2))\| \\
& \le \|f(v_1) - f(v_3)\| + \|f(v_4) - f(v_2)\| \\
& < \frac{\varepsilon}{2} + \frac{\varepsilon}{2} = \varepsilon.
\end{align*}
\end{proof}

Now we can state and prove our proposition.
\begin{prop}
\label{prop:YVcont}
The rigidity operator $Y$ is a continuous linear map from $\Vf(\V)$ to~$C(\E)$.
\end{prop}
\begin{proof}
We have two things to accomplish.  First, we need to show that continuous 
vector fields get mapped to continuous functions under~$Y$.  Second, we need
to show that $Y$ is a continuous linear map.

Remembering that $YV$ is the pointwise dot product of the functions $V(v_1) -
V(v_2)$ and $p(v_1) - p(v_2)$, we note that 
$$(V(v_1) - V(v_2)) \cdot (p(v_1) - p(v_2)) = 
  (V(v_2) - V(v_1)) \cdot (p(v_2) - p(v_1)).$$
So if $\pi$ is the quotient map from $\V \cross \V \to \vvt$, we see that 
$(V(v_1) - V(v_2)) \cdot (p(v_1) - p(v_2))$ is constant on
$\pi^{-1}(\{v_1,v_2\})$.

Now $V$ is continuous, since it is an element of $\Vf(\V)$, and $p$ is
continuous because the topology on $\V$ is defined by it.  So by Lemma
\ref{lemma:diffCont}, $V(v_1) - V(v_2)$ and $p(v_1) - p(v_2)$ are both
continuous on $\V \cross \V$.  Hence $(V(v_1) - V(v_2)) \cdot (p(v_1) -
p(v_2))$ is continuous on $\V \cross \V$.  But by Theorem
\ref{theorem:contQuot}, then, $YV = (V(v_1) - V(v_2)) \cdot (p(v_1) - p(v_2))$
is continuous on $\vvt$ and thus on~$\E$.

So we've shown that $Y$ maps from $\Vf(\V)$ to~$C(\E)$.  Let's show that it
does so in a continuous fashion.  We could do that by invoking Theorem \vref{theorem:ContBounded}, but doing the work ourselves is not difficult.

Let $U$ be an open set of functions in~$C(\E)$.  If $\im Y \cap U =
\varnothing$, then $Y^{-1} U = \varnothing$, which is open.  Otherwise,
let $V_1 \in Y^{-1} U$.  Now there is a basic open set $B$ around
$YV_1$ contained in $U$ with radius $r > 0$.  Let $V_2 \in
\Vf(\V)$ with $\|V_1 - V_2\| < \nicefrac{r}{\|Y\|}$ (which is well defined, 
thanks to Proposition \vref{prop:Ybounded}).  Then 
\begin{align*}
\|YV_1 - YV_2\| & = \|Y(V_1 - V_2)\| \\
& \le \|Y\|\|V_1 - V_2\| \text{, by the definition of $\|Y\|$} \\
& < \|Y\| \frac{r}{\|Y\|} = r
\end{align*}
So there is an open ball of radius $\nicefrac{r}{\|Y\|}$ around $V_1$ contained
in $Y^{-1} U$.  That is, $Y^{-1} U$ is open.  Thus $Y$ is a continuous map.
\end{proof}

\section{First Main Theorem}
\label{section:FMT}
\subsection{Statement and Proof}
We're almost ready to state and prove our first main theorem.  Before we do,
though, we'll need another of the fundamental results from Functional Analysis,
the Hahn-Banach Theorem.  Morrison, in his well-written book on Functional
Analysis \ycite{Morrison}, works through various versions of the Hahn-Banach
Theorem.  This is the extension form, which allows us to extend a linear
functional on a subspace of
$C(\E)$ to all of~$C(\E)$.
\begin{theorem}[Hahn-Banach Theorem, Extension Form]
\label{theorem:HahnBanach}
\IndexDef{Hahn-Banach Theorem}[2!1 ]
Suppose $X$ is a real linear space and $S$ is a linear subspace of $X$ and
$p\mcol X \to \RR$ is 
\begin{itemize}
\item[] \textbf{Subadditive} (that is, $p(x_1 + x_2) \le p(x_1) + p(x_2)$
for all $x_1, x_2 \in X$) 
\item[] \textbf{Nonnegatively subhomogeneous} (i.e., $p(\lambda x) \le \lambda p(x)$ for all $\lambda \ge 0$ and $x \in X$).
\end{itemize}
Let $f\mcol S \to \RR$ be a linear functional that satisfies $f(s) \le p(s)$
for all $s \in S$.  Then $f$ may be extended in a linear fashion to a
functional $F$ defined on all of $X$ in such a manner that the extension
satisfies $F(x) \le p(x)$ for all $x \in X$.
\end{theorem}
\begin{proof}
See \ocite{Morrison}*{p.\ 65}.
\end{proof}

\begin{theorem}[First Main Theorem]
\label{theorem:nonnegativeMeasure}
\statementNonnegativeMeasure
\end{theorem}
\begin{proof}
Let us begin by assuming that $G(p)$ has a semipositive stress,~$\mu$.  We wish
to show that $G(p)$ has no strictly positive motion.  Let $V \in \X$ be such
that $YV \in C(\E)^+$.  Since $\mu$ is a stress, we know that 
$$0 = \mu(YV) = \int_\E YV \, d\mu = \int_{\supp \mu} YV \, d\mu.$$
So we must have $YV = 0$ almost everywhere on $\supp \mu$, a set which, by 
Lemma \vref{lemma:spins}, is nonempty due to $\mu$ being semipositive.

Conversely, let us assume that $G(p)$ is partially $\X$-bar equivalent.  Then
$Y(\X) \cap \INT C(\E)^+ = \varnothing$ (by definition).  So $Y(\X)$ contains no
strictly positive motions, that is, for any motion $V \in \X$, there must be
some $e_V \in \E$ such that $YV(e_V) = 0$.  We also note that since $\X$ is a
subspace, $Y(\X)$ must be a subspace of $C(\E)$, so $Y(\X) \cap \INT C(\E)^- =
-(-Y(\X) \cap \INT C(\E)^+) = -(Y(\X) \cap \INT C(\E)^+) = \varnothing$.

Let $\mathfrak{s}\mcol C(\E) \to \RR$ by $\mathfrak{s}(f) = \sup f$ for all $f
\in C(\E)$.  We note that if $f_1$ and $f_2$ are in $C(\E)$ and $a \in
(0,\infty)$, then
$$\mathfrak{s}(f_1 + f_2) = \sup_{e \in \E} \left\{ \left(f_1 + f_2\right)(e)
\right\} \le \sup_{e \in \E} \{ f_1(e) \} + \sup_{e \in \E} \{ f_2(e) \} =
\mathfrak{s}(f_1) + \mathfrak{s}(f_2),$$
and 
$$\mathfrak{s}(a f_1) = \sup_{e \in \E} \{ a f_1(e) \} = a \sup_{e \in \E} \{
  f_1(e) \} = a \mathfrak{s}(f_1),$$ 
so $\mathfrak{s}$ is subadditive and nonnegatively subhomogeneous.

\Index[|(]{$\one$@one (constant function 1)}[ ]
Now, by Lemma \vref{lemma:intCEp}, $\one \in \INT C(\E)^+$, so, by hypothesis,
$\one \notin Y(\X)$.  We can define $\muh\mcol \Span \{ Y(\X), \one \} \to \RR$
by $\muh(YV + \alpha \one) = \alpha$.  What can we say about $\muh$?  First, it
is clearly linear.  Secondly, it is certainly zero on all of~$Y(\X)$.  How does
it compare to $\mathfrak{s}$?  Well, we know that any element $YV$ of $Y(\X)$
must be zero somewhere (so as to avoid being in either $\INT C(\E)^+$ or $\INT
C(\E)^-$), so at that point $YV + \alpha \one$ equals~$\alpha$.  That gives us
that $\mathfrak{s}(YV + \alpha \one) \ge \alpha = \muh(YV + \alpha \one)$.
Hence $\muh$ is dominated by~$\mathfrak{s}$.  

By the Hahn-Banach\Index{Hahn-Banach Theorem}[2!1 ] theorem (Theorem
\vref{theorem:HahnBanach}), then,  we can extend $\muh$ to some linear
functional $\mu$ on all of $C(\E)$ that is also dominated by~$\mathfrak{s}$.

Let's calculate 
$$\|\mu\| = \sup_{\|f\| = 1} |\mu(f)|$$
Remembering that $\|f\| = \sup_{e \in \E} |f(e)|$, we see that if $\|f\|
= 1,$ we must have $\mathfrak{s}(f) \le 1$ and $\mathfrak{s}(-f) \le 1$.  Since
$\mathfrak{s}$ dominates $\mu$, then, we must have $\mu(f) \le \mathfrak{s}(f)
\le 1$ and $\mu(-f) \le \mathfrak{s}(-f) \le 1$, giving us (by the linearity of
$\mu$) $\mu(f) = -\mu(-f) \ge -1$.  So $\|\mu\| \le 1$.  Since $\mu$ is linear
and bounded, by the Riesz Representation Theorem\Index{Riesz Representation
Theorem}[3!12] (Theorem \vref{theorem:Riesz}), $\mu$ is a regular measure
on~$\E$.  We need to show that $\mu(f) \ge 0$ for all $f \in C(\E)^+$ and
that $\mu \ne \zero$.

Suppose that $f \in C(\E)^+$.  Then $0 \ge \mathfrak{s}(-f) \ge \mu(-f)$.
But $\mu(-f) \le 0$ implies that $\mu(f) \ge 0$ by the linearity of~$\mu$.
Since $\mu(\one) = 1$, $\mu$ is nonzero, and since $\mu(f) \ge 0$ for all 
$f \in C(\E)^+$, $\mu$ is (at least) a semipositive measure.  Finally, as
$\mu(YV) = 0$ for all $V \in Y(\X)$, $\mu$ is a semipositive stress and our
proof is complete.
\Index[|)]{$\one$@one (constant function 1)}[ ]
\end{proof}

\subsection{An Example}
\label{subsection:FMTExmp}
Perhaps we could benefit from an example.  Suppose we have a tensegrity like
the one in \vref{fig:NoPosStress}.
\renewcommand{\V}{\mathscr{V}}
\begin{figure}[ht]
\hspace*{\fill}
\begin{tikzpicture}[scale=2]
\draw[strut] 
  (0,0) -- node [weight,pos=0.75] {$\nicefrac{\sqrt{2}}{2}$} (1,1) 
  (0,1) -- node [weight,pos=0.25] {$\nicefrac{\sqrt{2}}{2}$} (1,0) 
  (2,0.5) -- node [weight,pos=0.75] {$\nicefrac{\sqrt{2}}{2}$} (3,1.5) 
  (2,1.5) -- node [weight,pos=0.25] {$\nicefrac{\sqrt{2}}{2}$} (3,0.5);
\draw[cable] 
  (0,0) node [vertex] {} -- node [weight] {$1$}
  (0,1) node [vertex] {} -- node [weight] {$1$}
  (1,1) node [vertex] {} -- node [weight] {$0$}
  (2,1.5) node [vertex] {} -- node [weight] {$1$}
  (3,1.5) node [vertex] {} -- node [weight] {$1$}
  (3,0.5) node [vertex] {} -- node [weight] {$1$}
  (2,0.5) node [vertex] {} -- node [weight] {$0$}
  (1,0) node [vertex] {} -- node [weight] {$1$}
  (0,0) node [vertex] {} 
  (1,0) -- node [weight] {$1$} (1,1) 
  (2,1.5) -- node [weight] {$1$} (2,0.5);
\end{tikzpicture}
\hspace*{\fill}
\capt[A tensegrity with a semipositive stress but no strictly positive one.]{A
tensegrity with a semipositive stress (shown) but no strictly positive one.}
\label{fig:NoPosStress}
\end{figure}
In the figure, we have a semipositive stress for this tensegrity, so there is no
motion that strictly increases the lengths of all struts and decreases the
lengths of all cables.  

On the other hand, if we were to put a positive weight on the two edges which
are currently marked $0$, there would be no way to get a zero vector sum at the
vertices they touch.  As we saw in Proposition \vref{prop:stress}, that says
there is no strictly positive stress.  So the tensegrity is partially $\X$-bar
equivalent (each of the squares is $\X$-bar equivalent), but it is not $\X$-bar
equivalent as a whole.

Now, Theorem \ref{theorem:nonnegativeMeasure} doesn't seem particularly helpful
if there are bars in the tensegrity.  After all, one could create a
semipositive stress by putting equal positive weights on the strut and cable
portions of one of the bars and zeros everywhere else.  That is true, but there
are a couple of points worth noting.  First, there are many continuous
tensegrities which do not have bars.  Connelly et al. \ycite{CDR}, for example,
used just such a theorem on a tensegrity which had no bars.

Secondly, in some cases we can narrow $\X$ down to variations which would not
change the bars, and then remove the bars and analyze the tensegrity without 
them.  This has to be done carefully, as removing the bars may destroy
the compactness of $\E$ (it would, for example, for the tensegrity in
\vref{fig:compactBoth}), but it is useful in other cases.
\begin{figure}[ht]
\hspace*{\fill}
\begin{tikzpicture}
\draw[bar] (2,0) -- (0,2);
\draw[cable] (2,0) node [vertex] {} -- (0,0) node [vertex] {} -- 
  (0,2) node [vertex] {};
\draw[->] (0,0) -- (0.5,0.5);
\end{tikzpicture}
\hspace*{\fill}
\capt[A tensegrity with no semipositive stress.]{A tensegrity that has no
semipositive stress (and hence has a strictly positive motion) if $\X$ only
contains the variations that do not change the length of the bar.}
\label{fig:hasAmotion}
\end{figure}

As an example, the tensegrity in \ref{fig:hasAmotion} has no motion that
strictly shortens all cables and lengthens all struts, because that would
require both lengthening and shortening the bar.  But if $\X$ excludes those
variations that would change the length of the bar, the tensegrity has no
semipositive stress.  By Theorem \ref{theorem:nonnegativeMeasure}, then, there
is a strictly positive motion of the other edges.

\section{Second Main Theorem}
\label{section:SMT}
Our original goal was to move Roth \& Whiteley's lemma into the continuous
world.  We haven't quite accomplished it.  The difference between what we
have and what we were after is the same as the difference between Gordan's 
theorem and Stiemke's (shown \vpagerefrange{theorem:Stiemke}{theorem:Gordan}).
Let's take another pass at it.
\subsection{One Direction}
One direction of the theorem is easy to do.
\begin{theorem}[Second Main Theorem]
\label{theorem:positiveImplies}
\statementPositiveImplies
\end{theorem}
\begin{proof}
Let $G(p)$ be a tensegrity.  $G(p)$ having a strictly positive stress means
that there exists some strictly positive measure $\mu \in Y(\X)^\perp$.

Let $V$ be any element of $\X$ such that $YV \in C(\E)^+$.  Since $\mu \in
Y(\X)^\perp$, we must have $\mu(YV) = 0$.  On the other hand, $\mu$ is strictly
positive and $YV \in C(\E)^+$, so $\mu(YV) = 0$ implies $YV = \zero$.
\end{proof}

\subsection{Partially $\X$-bar equivalent}
Let's pause for a moment to defend the choice of ``partially $\X$-bar
equivalent'' as a term.  
\begin{corollary}
\label{corollary:partiallyBE}
If $G(p)$ is partially $\X$-bar equivalent, then some subtensegrity of $G(p)$
is $\X$-bar equivalent.
\end{corollary}
\begin{proof}
Let $G(p)$ be partially $\X$-bar equivalent.  By Theorem
\ref{theorem:nonnegativeMeasure}, $G(p)$ has a semipositive stress~$\mu$.  
Then $\supp \mu$ is nonempty, by Lemma \ref{lemma:spins}, and closed, by
Lemma \ref{lemma:closedsupport}, (and hence compact as the closed subspace
of a compact space, see, for example \ocite{Munkres}*{p.\ 165}) and has
positive measure, since $\mu \ne \zero$. 

So let $S(p)$ be the subtensegrity we get from $G(p)$ by keeping the vertex set
$\V$ intact but by restricting the edgeset to $\supp \mu$.  Lemma
\ref{lemma:spisE} gives us that  $\mu$ is a strictly positive measure
for~$S(p)$, so we need only show that it is a stress.  But since $\mu$ is a
stress for $G(p)$, we already have 
$$0 = \mu(YV) = \int_{\supp \mu} YV,$$
so $\mu$ is a stress for $S(p)$ as well.  Hence, by Theorem
\ref{theorem:positiveImplies}, $S(p)$ is $\X$-bar equivalent.
\end{proof}

We note, in passing, that \vref{fig:openstruts} does not give us a
counterexample to this theorem, as $\E$ in that case is not compact.

\subsection{The Other Direction}
\label{subsection:TheOtherDirection}
We have that having a strictly positive stress implies $\X$-bar equivalence.
Showing that the converse is true for continuous tensegrities is surprisingly
difficult.  Here's a little idea of why.  \vref{fig:notMinimal} shows a
tensegrity which resembles the crossed square\Index[|(]{crossed square}[ ] of
\ref{section:WorkingRWex}.
\begin{figure}[ht]
\hspace*{\fill}
\begin{tikzpicture}
\draw[cable] (0,0) -- (2,0) -- (2,2) -- (0,2) -- (0,0) -- (2,0);
\draw[strut] (0,0) node [vertex] {1} -- (2,2) node [vertex] {3};
\draw[bar] (0,2) node [vertex] {4} -- (2,0) node [vertex] {2};
\end{tikzpicture}
\hspace*{\fill}
\capt{A non-minimal crossed square.}
\label{fig:notMinimal}
\end{figure}
This one differs, however, in that it has a bar in place of one of the struts.
We already know that the crossed square is bar equivalent, what effect does
``extra cable'' of the bar have?  In \ref{section:WorkingRWex}, we found that
giving a weight of $1$ to every edge provided a strictly positive stress.  
By inspecting $Y^\top$ (on page \pageref{eq:Y}), we can discover that elements
of its kernel must have the same weight on every edge, so our $\mu$ is, up to
scaling, the only strictly positive stress for the crossed square.

By changing that strut to a bar, though, we changed~$Y$.  It now has another
row corresponding the new ``cable''.
$$\hat{Y} = \begin{bmatrix}
 1 &  0 & -1 &  0 &  0 &  0 &  0 &  0 \\
-1 & -1 &  0 &  0 &  1 &  1 &  0 &  0 \\
 0 &  1 &  0 &  0 &  0 &  0 &  0 & -1 \\
 0 &  0 &  0 &  1 &  0 & -1 &  0 &  0 \\
 0 &  0 &  1 & -1 &  0 &  0 & -1 &  1 \\
 0 &  0 &  0 &  0 & -1 &  0 &  1 &  0 \\
 0 &  0 & -1 &  1 &  0 &  0 &  1 & -1 \\
\end{bmatrix}$$
A little arithmetic and we find that the kernel of $\hat{Y}^\top$ is now a
family of stresses which look like $\begin{bmatrix} a & a & a & a & b & a &
(b-a) \end{bmatrix}^\top$.  There is now an interval's worth of strictly
positive stresses.  For any given positive value of $b$, $a$ can take on any
value in the interval $(0,b)$.  
\Index[|)]{crossed square}[ ]

In \vref{subsection:minimal}, we'll return to this example and see what that
says about the tensegrity, but at the moment we're more interested in what it
says about~$\mu$.

While the values $a \in (0,b)$ give us strictly positive stresses, the values
$a=0$ and $a=b$ give us semipositive stresses, which lie on the boundary of
$C^*(\E)^+$.  We could turn those semipositive stresses into strictly positive
ones by ``rotating'' the stress into the interior of
$C^*(\E)$\Index[|(]{interior {of $C^*(\E)^+$}}[1!2].  In the continuous case,
Theorem \vref{theorem:nonnegativeMeasure} gives us semipositive stresses, so
why not do something similar and ``rotate'' those stresses into the interior of
$C^*(\E)^+$?  The answer is fairly fundamental: In the infinite case,
$C^*(\E)^+$ has no interior.  

Suppose $X$ is any compact Hausdorff space and let $C(X)$, $C^*(X)$ and
$C^*(X)^+$ be defined for $X$ just as $C(\E)$, $C^*(\E)$ and $C^*(\E)^+$ are
defined for $\E$ for continuous tensegrities.  Then,
\begin{theorem}
\label{theorem:noInterior}
$C^*(X)^+$ has interior if and only if $X$ is finite.
\end{theorem}
\begin{proof}
Suppose $X$ is finite.  The elements of $C(X)$ and $C^*(X)$ are scalar
fields on~$X$.  Let $\mu \in C^*(X)^+$ with $\mu(x) > 0$ for all $x \in X$
and let $x_0$ denote some element of $X$ on which $\mu$ takes its smallest
value.  Now let $\eta \in C^*(X)$ with $\|\eta - \mu\| < \mu(x_0)$, that is,
$$\sup_{\|f\|=1} \left|(\eta - \mu)(f)\right| < \mu(x_0).$$  
Then certainly $\left|(\eta - \mu)(f)\right| < \mu(x_0)$ for all $f$ of norm 1,
including those functions which map one element of $X$ to $1$ and all of the
others to~$0$.  So for every $x \in X$, we have \begin{align*}
|(\eta - \mu)(x)| < \mu(x_0) 
  & \Rightarrow -\mu(x_0) < \eta(x) - \mu(x) < \mu(x_0) \\
  & \Rightarrow \mu(x) - \mu(x_0) < \eta(x) < \mu(x) + \mu(x_0)
\end{align*}
and since $\mu(x) \ge \mu(x_0)$ for all $x$, we have $\eta(x) > 0$ for all
$x \in X$.  Thus in the finite case, $C^*(X)^+$ has nonempty interior --- the
strictly positive elements of $C^*(X)$ form its interior.

Let's consider the general case.  Suppose $\mu \in \INT C^*(X)^+$.  At the very
least that means that $\mu(X) \ge 0$, so by scaling we can assume that $\mu(X)
= 1$.  Since $C^*(X)$ has the metric topology,\Index{metric topology}[2!1 ]
$\mu$ being in the interior of $C^*(X)^+$ also means that there must be some
$\delta > 0$ such that whenever $\eta \in C^*(X)$ with $\|\mu - \eta\| <
\delta$, we have $\eta \in C^*(X)^+$.

Now either there exists some $\varepsilon > 0$ such that $\mu(U) \ge
\varepsilon$ for every nonempty open $U \subset X$, or not.

Let us assume first that such an $\varepsilon$ exists.  Then any collection of
nonempty, pairwise-disjoint open sets on $X$ must have at most $\left\lfloor
\nicefrac{1}{\varepsilon} \right\rfloor$ members (where $\lfloor x \rfloor$
designates the greatest integer less than or equal to $x$).\Index{{$\lfloor
x \rfloor$}@*20 (greatest integer $\le x$)}[ ]

Let $m$ be the maximum number of sets in any such collection and let $U_1,
\dotsc, U_m$ be a collection with that maximal number of members.  Now the
union of the closures of the $U_i$, 
$$\bigcup_{i=1}^m \overline{U_i}$$
must contain all of~$X$.  Otherwise $X \setminus \bigcup \overline{U_i}$ would
be another open set disjoint from each of the $U_i$, contradicting the
maximality of~$m$.

There are only finitely many $U_i$, so if there are infinitely many elements in
$X$, then by the Pigeon-hole Principle (see, for example, \ocite{MR0092794}*{p.
78}), at least one of the $\overline{U_i}$ must contain infinitely many of
them.  Without loss of generality, we'll take that set to be $\overline{U_m}$.

Let $x_1$ and $x_2$ be two distinct elements in $\overline{U_m}$.  Since $X$ is
Hausdorff, there exist disjoint open sets $S_1$ and $S_2$ which separate them.
Since $x_1$ and $x_2$ are in the closure of $U_m$, it must be true that $S_1
\cap U_m \ne \varnothing$ and $S_2 \cap U_m \ne \varnothing$.  Hence, by
removing $U_m$ from the collection and replacing it with $S_1 \cap U_m$ and
$S_2 \cap U_m$, we have a new collection of nonempty, pairwise-disjoint open
sets with $m+1$ elements, again contradicting the maximality of~$m$.  

Thus, if there exists such an $\varepsilon$, there are only a finite number of
elements in~$X$.

On the other hand, suppose that there is no such~$\varepsilon$.  Then we can
select some $U_0$ such that $\mu(U_0) < \nicefrac{\delta}{4}$ and let
$\mu_\delta$ be a semipositive measure such that $\mu_\delta(U_0) =
\nicefrac{\delta}{2}$, but $\mu_\delta(U) = 0$ whenever $U \subset X$ with $U
\cap U_0 = \varnothing$ (for example, since $U_0$ is nonempty, 
$\mu(U_0) > 0$, so select any $x \in U_0$ and put an atom of size
$\nicefrac{\delta}{2}$ on it and zeros everywhere else).  We note that since
$(X \setminus U_0) \cap U_0 = \varnothing$, $\mu_\delta(X) = \mu_\delta(U_0) =
\nicefrac{\delta}{2}$.  And, since $\mu_\delta$ is semipositive,
$$\|\mu_\delta\| = \sup_{\|f\| = 1} \mu_\delta(f) = \mu_\delta(\one) =
\mu_\delta(U_0) = \nicefrac{\delta}{2}.$$

Now consider $\eta \coloneqq \mu - \mu_\delta$.  By construction,
$$\|\mu - \eta\| = \|\mu - (\mu - \mu_\delta)\| = \|\mu_\delta\| =
\nicefrac{\delta}{2} < \delta,$$
so $\eta$ lies with in a ball of radius $\delta$ about~$\mu$.  But 
$$\eta(U_0) = \mu(U_0) - \mu_\delta(U_0) < \nicefrac{\delta}{4} -
\nicefrac{\delta}{2} = -\nicefrac{\delta}{4} < 0$$
so $\eta \notin C^*(X)^+$.

Thus $C^*(X)^+$ has interior if and only if $X$ is finite.
\end{proof}
\Index[|)]{interior {of $C^*(\E)^+$}}[1!2]
As that route to a strictly positive stress has proved unproductive, let's
briefly set aside our quest to show that $\X$-bar equivalent tensegrities
have strictly positive stresses and take a look at the results of what we have
accomplished so far.

\section{Two Examples}
\label{section:TwoExamples}
\subsection{The setup}
We can explore the new theorems with a couple of examples.  In the first, shown
in \vref{fig:CircleOfStruts}, the vertices form a unit circle\IndDefBeg{circle
of struts}[ ] and struts are placed antipodally at all points.  Our design
variations will be only those variations which do not even locally stretch or
shrink the curve of vertices.
\begin{figure}[ht]
\hspace*{\fill}
\begin{tikzpicture}
\foreach \x in {0,20,...,181} {
  \draw[strut] (\x:2cm) -- (\x+180:2cm);
}
\draw[vcurve] (0,0) circle (2cm);
\end{tikzpicture}
\hspace*{\fill}
\capt{The circle-of-struts example.}
\label{fig:CircleOfStruts}
\end{figure}
Since the circle cannot stretch or shrink, this example certainly appears to be
$\X$-infinitesimally rigid and so we would expect it to be $\X$-bar equivalent.
\IndDefEnd{circle of struts}[ ]

We get the second example by removing a little more than half the struts of the
first example (see \vref{fig:AlmostHalfACircleOfStruts}).  We'll leave those
struts that touch the circle from $\varepsilon$ to $\nicefrac{\pi}{2} -
\varepsilon$ for some $\varepsilon \in (0,\nicefrac{\pi}{4})$.  That should
allow enough freedom to establish a strictly positive motion, even with our
requirement about maintaining length locally on the vertex set.
\begin{figure}[ht]
\hspace*{\fill}
\begin{tikzpicture}
\foreach \x in {5,15,...,86} {
  \draw[strut] (\x+180:2cm) -- (\x:2cm);
}
\draw[vcurve] (0,0) circle (2cm);
\end{tikzpicture}
\hspace*{\fill}
\capt{The almost-half-a-circle-of-struts example.}
\label{fig:AlmostHalfACircleOfStruts}
\end{figure}

In both cases, we can parameterize $\V$ and $\E$ by angle (since angle and
arclength are the same on a unit circle).  The values for $\V$ will lie in 
$\quotientspace{\RR}{2\pi}$ and for $\E$ in~$\qrp$.

\subsection{A New $\X$}
\label{subsection:VandVf}
Much of our preparation for this example will center around the set of design
variations,~$\X$.
\renewcommand{\V}{\gamma}
Specifically, we need to identify those variations which are local isometries
of the vertex curve (to remind us that it is a curve, we'll call it $\V$
instead of~$\mathscr{V}$).

For a little while, we'll even forget that our examples have circles
for vertex curves, and see how general we can let $\V$ and $\X$ be.

We need to be able to define length.  To do that, we approximate $\V$ with
inscribed polygons or polygonal lines and take its length to be the supremum of
the lengths of all such approximations.  Here's a more rigorous definition:
\begin{definition}[\ocite{Strichartz95}*{p.\ 611}]
The \emph{length}\IndexDef{length (of a curve)}[ ] of a continuous curve
$\V\mcol [a,b] \to \RR^n$ is the $\sup$ of $\sum_{j=1}^{N} \|\V(t_j) -
\V(t_{j-1})\|$ taken over all partitions $a = t_0 < t_1 < \dotsb < t_N = b$ of
the interval $[a,b]$, where $\|\V(t_j) - \V(t_{j-1})\|$ is the distance between
the points $\V(t_j)$ and $\V(t_{j-1})$ in~$\RR^n$.  If the length is finite, we
say the curve is \emph{rectifiable}\IndexDef{rectifiable}.  
\end{definition}

We'll require that $p(\V)$ be a simple, rectifiable curve with a finite number
of connected components and total length $\ell < \infty$.  

One advantage of a rectifiable curve is that we can give it a nice
parameterization.  Using our definition of length, we can define:
\begin{definition}
An \emph{arclength parameterization}\IndexDef{arclength parameterization}[2,1 ]
of a curve $\V\mcol [a,b] \to \RR^n$ of length $\ell$ is a continuous map
$t(s)\mcol [0,\ell] \to [a,b]$ such that at any point $s \in [0,\ell]$, the
length of $\V$ from $\V(t(0)) = \V(a)$ to $\V(t(s))$ is equal to~$s$.
\end{definition}
\begin{prop}
Every rectifiable curve admits an arclength parameterization.
\end{prop}
\begin{proof}
See, for example, \ocite{Strichartz95}*{p.\ 616} or \ocite{MR1835418}*{p.\ 46}.
\end{proof}

Now suppose we have some arclength parameterized $\V\mcol [0,\ell] \to \RR^n$. 
What more can we say about $\V$?  It turns out we can say some strong things
about continuity.  Here are four different types of continuity, each stronger
than the previous:
\begin{definition}[see, for example, \ocite{Royden}*{p.\ 44}]
\label{definition:continuous}
A function $f(x)$ is \emph{continuous}\IndexDef{continuous} if, for every
$\varepsilon > 0$ and every $x$, there exists a $\delta > 0$ such that whenever
$|x - y| < \delta$, $|f(x) - f(y)| < \varepsilon$.
\end{definition}
\begin{definition}[see, for example, \ocite{Royden}*{p.\ 135}]
A function $f(x)$ is \emph{uniformly continuous}\IndexDef{uniformly
continuous}[2!1 ] if, for every $\varepsilon > 0$ there exists a $\delta > 0$
such that for every $x$, whenever $|x-y| < \delta$, $|f(x) - f(y)| <
\varepsilon$.
\end{definition}
\begin{definition}[see, for example, \ocite{MR1681462}*{p.\ 105}]
\label{definition:AbsCont}
A function $f(x)$ is \emph{absolutely continuous}\IndexDef{absolutely
continuous}[2!1 ] if, for every $\varepsilon > 0$ there exists a $\delta > 0$
such that for every finite collection of pairwise-disjoint intervals
$(x_i,y_i)$ with $\sum |x_i - y_i| < \delta$ we have $\sum |f(x_i) - f(y_i)| <
\varepsilon$.
\end{definition}
\begin{definition}[see, for example, \ocite{MR1681462}*{p.\ 108} or
\ocite{MR1070713}*{p.\ 25}]
A function $f$ is \emph{Lipschitz continuous} (or just
\emph{Lipschitz})\IndexDef{Lipschitz continuous}[2!1 ] if there exists some
constant $k \ge 0$ (called the \emph{Lipschitz constant}\IndexDef{Lipschitz
constant}[ ]) such that for all $x$ and $y$, we have $|f(x) - f(y)| \le k |x -
y|$.
\end{definition}

We claimed that each of the types of continuity was stronger than the previous.
Here's proof of that relationship for the last two:
\begin{prop}
\label{prop:lac}
If a map $f$ is Lipschitz, it is absolutely continuous.
\end{prop}
\begin{proof}
Let $\varepsilon > 0$ be given.  Select $\delta = \nicefrac{\varepsilon}{k}$
where $k$ is the Lipschitz constant of~$f$.  Then, for every 
finite collection of intervals, $(x_i,y_i)$, with $\sum |x_i - y_i| < \delta$,
we have that 
\begin{align*}
\sum |f(x_i) - f(y_i)| & \le \sum k|x_i - y_i| \text{ since $f$ is Lipschitz}\\
& = k \sum|x_i - y_i| < k \delta = \varepsilon.
\end{align*}
\end{proof}
Now the converse of Proposition \ref{prop:lac} is not generally true.  Here's
a case in which it is not.
\begin{claim}
The function $\sqrt{x}$ is absolutely continuous but not Lipschitz continuous
on $[0,1]$.
\end{claim}
\begin{proof}
Consider the function~$\sqrt{x}$.  Let's show that $\sqrt{x}$ is absolutely
continuous.  

For any given $\varepsilon > 0$, we can select $\delta = \varepsilon^2$.  Now
$\sqrt{x}$ has monotonically decreasing derivative on $(0,\infty)$.  So the
worst-case collection of intervals of total length less than $\delta$ will be
the single interval $(0,\rho)$ for any $\rho < \delta$.  But then
$|\sqrt{0} - \sqrt{\rho}| < \sqrt{\delta} = \sqrt{\varepsilon^2} =
\varepsilon$.

On the other hand, $\sqrt{x}$ is not Lipschitz continuous, since the 
expression $\frac{\left|\sqrt{0} - \sqrt{x}\right|}{\left|0 - x\right|} =
\frac{\sqrt{x}}{x} = \frac{1}{\sqrt{x}}$ can be made arbitrarily large by getting $x$ near enough to~$0$.
\end{proof}

John Sullivan, in his survey paper, ``Curves of Finite Total Curvature'', gives
the following proposition.
\begin{prop}
A curve is rectifiable if and only if it admits a Lipschitz parameterization.
\end{prop}
\begin{proof}
See \ocite{math.GT0606007}*{p.\ 3}.
\end{proof}
What he is noting, in the ``only if'' portion of his proof, is that an
arclength parameterized curve, while not necessarily differentiable, is 
Lipschitz (with Lipschitz constant 1).  So we have both Lipschitz continuity
and absolute continuity of $\V$ simply due to its being rectifiable.  What more
can we get?
\begin{lemma}
If $f\mcol [a,b] \to \RR^n$ is absolutely continuous, then $f$ is rectifiable.
\end{lemma}
\begin{proof}  Here we flesh out a proof of Royden \ycite{Royden}*{p.\ 105}.
Since $f$ is absolutely continuous, there exists some $\delta > 0$ such that
for every collection of pairwise-disjoint intervals $(x_i,y_i)$ with $\sum |x_i
- y_i| < \delta$, we have $\sum \|f(x_i) - f(y_i)\| < 1$.  Let $m =
\left\lfloor \frac{b-a}{\delta} \right\rfloor$.  Given any partition $a < t_0 <
t_1 < \dotsb < t_N = b$, we can split the intervals $(t_i,t_{i+1})$ (if
necessary) and group them in such a way that we have $m+1$ groups: $m$ whose
total length falls in the range $(\delta - \frac{1}{m} \delta, \delta)$ and 1
with total length less than $m \frac{1}{m} \delta = \delta$.  

By the triangle inequality, splitting intervals in a partition will never
decrease the value of 
$$\sum_{j=1}^N \|f(t_j) - f(t_{j-1})\|$$
so we know that for any partition, that sum is at most $(m + 1) (1) = m+1$,
so the length is also bounded above by $m+1 < \infty$.  
\end{proof}
Wait a moment.  How can it be true that rectifiable implies Lipschitz implies
absolutely continuous implies rectifiable and yet absolutely continuous is 
not equivalent to Lipschitz?

The answer is that rectifiability is an attribute of the trace of a curve in
space, while Lipschitz and absolute continuity are attributes of the function
mapping out that trace.  An absolutely continuous function maps out a curve
which admits a Lipschitz parameterization, but it may not itself be that
parameterization.  In our example, the trace of the function $x \mapsto
\sqrt{x}$ is the interval $[0,\infty)$.  And this certainly admits a Lipschitz
parameterization, namely the function $x \mapsto  x$ (for $x \in [0,\infty)$),
which is Lipschitz with Lipschitz constant 1.

What else do we know?

\begin{theorem}
\label{theorem:dae}
If $f$ is absolutely continuous, then $f$ has a derivative almost everywhere.
\end{theorem}
\begin{proof}
See, for example, \ocite{Royden}*{p.\ 105}.
\end{proof}

So we get derivatives almost everywhere.  That allows us to look at the
(Lebesgue) integral
$$L(\V) = \int_a^b \|\V'(s)\| \ds.$$
Then we have 
\begin{theorem}
The function $L(\V)$ gives the arclength of $\V$ (as defined using inscribed
polygons) from $a$ to~$b$.
\end{theorem}
\begin{proof}
See \ocite{MR1835418}*{p.\ 57}.
\end{proof}

We are finally ready to look for those vector fields which preserve length.  We
need $\V$ to remain rectifiable throughout, so we will restrict ourselves to
vector fields which are absolutely continuous.  

Let's summarize: $\V$ is a simple, rectifiable curve with a finite number of
connected components and total length $\ell < \infty$ and $V\mcol [0,\ell] \to
\RR^n$ is an absolutely continuous vector field that moves $\V$ while locally
preserving $\V$'s arclength to first order.  What can we say about $V$?

\begin{prop}
\label{prop:orthoprime}
If $\V(s)$ is a simple, rectifiable, arclength parameterized curve
with length $\ell$ and $V(s)$ is an absolutely continuous variation of $\V(s)$
which locally preserves the arclength of $\V(s)$ to first order, then $\langle
V'(s), \V'(s) \rangle = 0$ almost everywhere.
\end{prop}

We'll jump into the proof in a moment, but before we start we'll need a 
result about Lebesgue integrals.
\begin{theorem}[Dominated Convergence Theorem]
\IndexDef{Dominated Convergence Theorem}[3!12 ]
Let $\{ f_n \}$ be a sequence of integrable functions such that
\begin{enumerate}[label=(\alph*)]
\item $f_n \to f$ almost everywhere and
\item there exists a nonnegative, integrable function $g$ such that $|f_n| \le
g$ almost everywhere for all~$n$.
\end{enumerate}
Then $f$ is integrable and $\int f = \lim_{n \to \infty} \int f_n$.
\end{theorem}
\begin{proof}
See, for example, \ocite{MR1681462}*{p.\ 54}
\end{proof}

Now we're ready for that proof.
\begin{proof}[Proof of Proposition \ref{prop:orthoprime}]
For any $a,b \in [0,\ell]$ and $\varepsilon > 0$, we know that the length of
the section of $\V(s) + \varepsilon V(s)$ between $a$ and $b$ is given by 
$$\int_a^b \left\| \frac{d}{ds} \left( \V(s) + \varepsilon V(s) \right) \right\|
\ds = \int_a^b \left\| \V'(s) + \varepsilon V'(s) \right\| \ds.$$
If $V$ is to be locally arclength preserving to first order, we must have that
$$\left[ \frac{d}{d\varepsilon} \int_a^b \left\| \V'(s) + \varepsilon V'(s)
\right\| \,ds \right]_{\varepsilon = 0} = 0$$
for every choice of $a,b \in [0,\ell]$.

Applying the definition of the derivative, we get
\begin{align*}
0 & = \lim_{\varepsilon \to 0} \frac{\int_a^b \|\V'(s) + \varepsilon V'(s)\| \ds
- \int_a^b \|\V'(s) - (0) V'(s)\| \ds}{\varepsilon - 0} \\
& = \lim_{\varepsilon \to 0} \int_a^b \frac{\|\V'(s) + \varepsilon V'(s)\| 
- \|\V'(s)\|}{\varepsilon} \ds
\end{align*}

Now if we let $\{ \varepsilon_n \}$ be a sequence of positive values
approaching $0$, we see that the functions
$$f_n(s) = \frac{\|\V'(s) + \varepsilon_n V'(s)\| - \|\V'(s)\|}{\varepsilon_n}$$
approach
$$f(s) = \left[\frac{\partial}{\partial \varepsilon} \|\V'(s) + \varepsilon
V'(s)\|\right]_{\varepsilon = 0}$$
wherever $\V'(s)$ and $V'(s)$ are both defined (which is almost everywhere).
Now either $\|\V'(s) + \varepsilon_n V'(s)\| \ge \|\V'(s)\|$ or not.
If it is, then
\begin{align*}
0 & \le \|\V'(s) + \varepsilon_n V'(s)\| - \|\V'(s)\| \\
  & \le \|\V'(s) + \varepsilon_n V'(s) - \V'(s)\| 
    = \|\varepsilon_n V'(s)\| = \varepsilon_n \|V'(s)\|.
\end{align*}
Otherwise, we have
\begin{align*}
0 & \le \|\V'(s)\| - \|\V'(s) + \varepsilon_n V'(s)\| \\
  & \le \|\V'(s) - \V'(s) - \varepsilon_n V'(s)\| 
    = \|- \varepsilon_n V'(s)\|
    = \varepsilon_n \|V'(s)\|.
\end{align*}
So, by the triangle inequality, 
$$\left|\frac{\|\V'(s) + \varepsilon_n V'(s)\| -
\|\V'(s)\|}{\varepsilon_n}\right|
\le \left|\frac{\varepsilon_n\|V'(s)\|}{\varepsilon_n}\right| = \|V'(s)\|$$
for all~$n$.  Since $\|V'(s)\|$ is nonnegative and integrable, the Dominated
Convergence Theorem\Index{Dominated Convergence Theorem}[3!12 ] tells us that
$$\left[ \frac{d}{d\varepsilon} \int_a^b \|\V'(s) + \varepsilon V'(s)\| \ds
\right]_{\varepsilon = 0} = \int_a^b \left[\frac{\partial}{\partial
\varepsilon} \|\V'(s) + \varepsilon V'(s)\| \right]_{\varepsilon=0} \ds.$$

Since this is true for \emph{any} choice of $a$ and $b$, it must be the case
that 
$$ \left[ \frac{\partial}{\partial\varepsilon} \left\| \V'(s) + \varepsilon
V'(s) \right\| \right]_{\varepsilon = 0} = 0 $$
almost everywhere.

Now
\begin{align*}
\left[ \frac{\partial}{\partial\varepsilon} \left\|\V'(s) + \varepsilon
V'(s)\right\| \right]_{\varepsilon = 0} 
& = \left[ \frac{\partial}{\partial\varepsilon} \left\langle \V'(s) +
    \varepsilon V'(s), \V'(s) + \varepsilon V'(s)
    \right\rangle^{\nicefrac{1}{2}} \right]_{\varepsilon=0} \\
& = \left[ \frac{2 \left\langle V'(s), \V'(s) + \varepsilon V'(s)
    \right\rangle}{2\left\| \V'(s) + \varepsilon V'(s) \right\|} 
    \right]_{\varepsilon=0} \\
& = \frac{\left\langle V'(s), \V'(s)
    \right\rangle}{\left\| \V'(s) \right\|} \\
& = \langle V'(s), \V'(s) \rangle,
\end{align*}
so we need $V'(s)$ and $\V'(s)$ to be orthogonal almost everywhere.
\end{proof}

Our $\X$, then, is the space of absolutely continuous vector fields $V$ on
$\gamma$ with $V'(s) \perp \V'(s)$ for almost all $s \in [0,\ell]$.

We will need to check that $\X$ is a subspace of $\Vf(\V)$.
\begin{prop}
\label{prop:Xsubspace}
$\X$ is a subspace of $\Vf(\V)$.
\end{prop}
\begin{proof}
Let $V,W \in \X$ and $\alpha,\beta \in \RR$.  Then $\alpha V + \beta W$ is
absolutely continuous, and here's why.  Let $\varepsilon > 0$.  Then (by the
absolute continuity of $V$ and $W$) there exists $\delta_V$ and $\delta_W$
such that whenever $(x_i,y_i)$ is a finite collection of pairwise-disjoint
intervals we have
$$\sum |x_i - y_i| < \delta_V \Rightarrow \sum \|V(x_i) - V(y_i)\| <
\frac{\varepsilon}{2|\alpha|}$$
and
$$\sum |x_i - y_i| < \delta_W \Rightarrow \sum \|W(x_i) - W(y_i)\| <
\frac{\varepsilon}{2|\beta|}.$$
But then, whenever $\sum |x_i - y_i| < \min\{\delta_V,\delta_W\}$, we
have 
\begin{align*}
\sum \|(\alpha V + \beta W)(x_i) & - (\alpha V + \beta W)(y_i) \| \\
& = \sum \|\alpha (V(x_i) - V(y_i)) + \beta(W(x_i) - W(y_i)) \| \\
& \le \sum |\alpha| \|V(x_i) - V(y_i)\| + |\beta| \|W(x_i) - W(y_i)\| \\
& < |\alpha| \frac{\varepsilon}{2|\alpha|} + |\beta|
\frac{\varepsilon}{2|\beta|} = \varepsilon.
\end{align*}
That gives us absolute continuity.  We need to check the orthogonality
condition.  But 
$$(\alpha V + \beta W)'(s) \cdot \V'(s) = \alpha V'(s) \cdot \V'(s) + \beta
W'(s) \cdot V'(s) = 0$$
almost everywhere, so $\alpha V + \beta W \in \X$, and $\X$ is a subspace.
\end{proof}

\subsection{The Circle of Struts}
\Index[|(]{circle of struts}[ ] Now that we have a grasp on $\X$, let's return
to our examples.  For both examples, $\V'$ is made up of unit vectors tangent
to the circle and pointing counter-clockwise around it.  Any vector field $V
\in \X$, then, must have $V'$ pointing radially outward or inward on the
circle.  So everything in $\X$ will be of the nature 
$$V(s) = \int_0^s (f(t) \cos t, f(t) \sin t) \dt$$ 
\label{Vpoints}
(where $f\mcol \qrp \to \RR$ is not necessarily continuous and only need be
defined almost everywhere).  But why should we believe that every such $V(s)$
is absolutely continuous?
\begin{theorem}
\label{theorem:iiac}
A function $f$ is an indefinite integral if and only if it is absolutely
continuous.
\end{theorem}
\begin{proof}
See, for example, \cite{Royden}*{p.\ 106}.
\end{proof}
Since we also need absolute continuity at $0$, we'll have the additional 
requirement that $V(2\pi) = V(0)$.  That is,
$$\int_0^{2\pi} (f(t) \cos t, f(t) \sin t) \dt = \zero.$$ 
 
For the circle-of-struts example, let's try building a stress.  The simplest
option would be to take the measure $d\theta$, which is uniform on
$\qrp$ (that is, on the struts), and integrates to $\pi$.  Let's suppose that
$V$ is some variation in $\X$ and see what happens.  

Now $YV(\theta) = (V(\theta) - V(\theta + \pi)) \cdot (\cos \theta, \sin
\theta)$.  So we get
\begin{align*}
\int_0^{\pi} YV(\theta) \, d\theta 
& = \int_0^\pi (V(\theta) - V(\theta + \pi)) \cdot (\cos \theta, \sin \theta) \, d\theta \\ 
& = \int_0^\pi V(\theta) \cdot (\cos \theta, \sin \theta) \, d\theta +
    \int_\pi^{2\pi} V(\theta) \cdot (\cos \theta, \sin \theta) \, d\theta \\
& = \int_0^{2\pi} V(\theta) \cdot (\cos \theta, \sin \theta) \, d\theta, \\
\intertext{which we can integrate by parts to get}
& = V(\theta) \cdot (\sin \theta, -\cos \theta) \big|_0^{2\pi} -
    \int_0^{2\pi} V'(\theta) \cdot (\sin \theta, -\cos \theta) \, d\theta.
\intertext{Since we require that $V(0) = V(2\pi)$, that first term is zero, 
so we get}
& = \int_0^{2\pi} V'(\theta) \cdot (-\sin \theta, \cos \theta) \, d\theta.
\end{align*}
On the other hand, $(-\sin \theta, \cos \theta) = \V'(\theta)$ and since 
$V \in \X$, we know that $V'(\theta) \cdot \V'(\theta) = 0$ almost
everywhere, so we have $\int V'(\theta) \cdot \V'(\theta) = 0$.  Hence 
$d\theta \in Y(\X)^\perp$.

So $d\theta$ is a stress.  But taking any nonempty, open set of edges and
integrating by $d\theta$ will give a positive value, so $d\theta$ is a strictly
positive stress.  The circle of struts is $\X$-bar equivalent.  If we were to
apply Theorem \ref{theorem:nonnegativeMeasure}, we would have:
\begin{prop}
\label{prop:COSone}
\statementCOSone
\end{prop}
But Theorem \ref{theorem:positiveImplies} gives us the stronger:
\begin{prop}
\label{prop:COStwo}
\statementCOStwo
\end{prop}
\Index[|)]{circle of struts}[ ]

\subsection{Almost half a circle of struts}
Now let's consider our other example.  In this case, we'd like to show that 
there is no semipositive stress.  

Let's choose a ``good'' element of $\X$, one in $T(p) \setminus I(p)$.
Remember that our tensegrity has less than half the struts that the
circle-of-struts example had, so we should be able to expand along those and
contract elsewhere, keeping the curve of vertices unstretched.  Let's try
$$V_g = \left(\frac{\sin \theta}{2} + \frac{\sin 3\theta}{6},
              \frac{\cos \theta}{2} - \frac{\cos 3\theta}{6}\right),$$
which is shown in \vref{fig:Vg}.
\begin{figure}[ht]
\hspace*{\fill}
\begin{tikzpicture}[scale=2]
\foreach \x in {5,10,...,86} {
  \draw[strut,black!20!white] (\x:1cm) -- (\x:-1cm);
}
\draw[vcurve] (0,0) circle (1cm);
\draw[->] ( 1.000, 0.000) -- ( 1.000, 0.333);
\draw[->] ( 0.996, 0.087) -- ( 1.083, 0.424);
\draw[->] ( 0.985, 0.174) -- ( 1.155, 0.522);
\draw[->] ( 0.966, 0.259) -- ( 1.213, 0.624);
\draw[->] ( 0.940, 0.342) -- ( 1.255, 0.729);
\draw[->] ( 0.906, 0.423) -- ( 1.279, 0.833);
\draw[->] ( 0.866, 0.500) -- ( 1.283, 0.933);
\draw[->] ( 0.819, 0.574) -- ( 1.267, 1.026);
\draw[->] ( 0.766, 0.643) -- ( 1.232, 1.109);
\draw[->] ( 0.707, 0.707) -- ( 1.179, 1.179);
\draw[->] ( 0.643, 0.766) -- ( 1.109, 1.232);
\draw[->] ( 0.574, 0.819) -- ( 1.026, 1.267);
\draw[->] ( 0.500, 0.866) -- ( 0.933, 1.283);
\draw[->] ( 0.423, 0.906) -- ( 0.833, 1.279);
\draw[->] ( 0.342, 0.940) -- ( 0.729, 1.255);
\draw[->] ( 0.259, 0.966) -- ( 0.624, 1.213);
\draw[->] ( 0.174, 0.985) -- ( 0.522, 1.155);
\draw[->] ( 0.087, 0.996) -- ( 0.424, 1.083);
\draw[->] ( 0.000, 1.000) -- ( 0.333, 1.000);
\draw[->] (-0.087, 0.996) -- ( 0.250, 0.909);
\draw[->] (-0.174, 0.985) -- ( 0.174, 0.815);
\draw[->] (-0.259, 0.966) -- ( 0.106, 0.719);
\draw[->] (-0.342, 0.940) -- ( 0.044, 0.624);
\draw[->] (-0.423, 0.906) -- (-0.013, 0.534);
\draw[->] (-0.500, 0.866) -- (-0.067, 0.449);
\draw[->] (-0.574, 0.819) -- (-0.121, 0.371);
\draw[->] (-0.643, 0.766) -- (-0.176, 0.300);
\draw[->] (-0.707, 0.707) -- (-0.236, 0.236);
\draw[->] (-0.766, 0.643) -- (-0.300, 0.176);
\draw[->] (-0.819, 0.574) -- (-0.371, 0.121);
\draw[->] (-0.866, 0.500) -- (-0.449, 0.067);
\draw[->] (-0.906, 0.423) -- (-0.534, 0.013);
\draw[->] (-0.940, 0.342) -- (-0.624,-0.044);
\draw[->] (-0.966, 0.259) -- (-0.719,-0.106);
\draw[->] (-0.985, 0.174) -- (-0.815,-0.174);
\draw[->] (-0.996, 0.087) -- (-0.909,-0.250);
\draw[->] (-1.000, 0.000) -- (-1.000,-0.333);
\draw[->] (-0.996,-0.087) -- (-1.083,-0.424);
\draw[->] (-0.985,-0.174) -- (-1.155,-0.522);
\draw[->] (-0.966,-0.259) -- (-1.213,-0.624);
\draw[->] (-0.940,-0.342) -- (-1.255,-0.729);
\draw[->] (-0.906,-0.423) -- (-1.279,-0.833);
\draw[->] (-0.866,-0.500) -- (-1.283,-0.933);
\draw[->] (-0.819,-0.574) -- (-1.267,-1.026);
\draw[->] (-0.766,-0.643) -- (-1.232,-1.109);
\draw[->] (-0.707,-0.707) -- (-1.179,-1.179);
\draw[->] (-0.643,-0.766) -- (-1.109,-1.232);
\draw[->] (-0.574,-0.819) -- (-1.026,-1.267);
\draw[->] (-0.500,-0.866) -- (-0.933,-1.283);
\draw[->] (-0.423,-0.906) -- (-0.833,-1.279);
\draw[->] (-0.342,-0.940) -- (-0.729,-1.255);
\draw[->] (-0.259,-0.966) -- (-0.624,-1.213);
\draw[->] (-0.174,-0.985) -- (-0.522,-1.155);
\draw[->] (-0.087,-0.996) -- (-0.424,-1.083);
\draw[->] (-0.000,-1.000) -- (-0.333,-1.000);
\draw[->] ( 0.087,-0.996) -- (-0.250,-0.909);
\draw[->] ( 0.174,-0.985) -- (-0.174,-0.815);
\draw[->] ( 0.259,-0.966) -- (-0.106,-0.719);
\draw[->] ( 0.342,-0.940) -- (-0.044,-0.624);
\draw[->] ( 0.423,-0.906) -- ( 0.013,-0.534);
\draw[->] ( 0.500,-0.866) -- ( 0.067,-0.449);
\draw[->] ( 0.574,-0.819) -- ( 0.121,-0.371);
\draw[->] ( 0.643,-0.766) -- ( 0.176,-0.300);
\draw[->] ( 0.707,-0.707) -- ( 0.236,-0.236);
\draw[->] ( 0.766,-0.643) -- ( 0.300,-0.176);
\draw[->] ( 0.819,-0.574) -- ( 0.371,-0.121);
\draw[->] ( 0.866,-0.500) -- ( 0.449,-0.067);
\draw[->] ( 0.906,-0.423) -- ( 0.534,-0.013);
\draw[->] ( 0.940,-0.342) -- ( 0.624, 0.044);
\draw[->] ( 0.966,-0.259) -- ( 0.719, 0.106);
\draw[->] ( 0.985,-0.174) -- ( 0.815, 0.174);
\draw[->] ( 0.996,-0.087) -- ( 0.909, 0.250);
\end{tikzpicture}
\hspace*{\fill}
\capt[The almost-half-a-circle-of-struts with the ``good'' motion~$V_g$.]{The
almost-half-a-circle-of-struts example with the ``good'' motion~$V_g$.}
\label{fig:Vg}
\end{figure}

Now, $YV_g$ is nonnegative on all of~$\E$.  In specific, 
\begin{align*}
p(\theta) - p(\theta + \pi) & =
(\cos \theta, \sin \theta) - (\cos(\theta + \pi), \sin(\theta + \pi)) \\
& = (\cos \theta, \sin \theta) - (-\cos \theta, -\sin \theta) \\
& = (2\cos \theta, 2 \sin \theta)
\end{align*}
and similarly,
\begin{align*}
V_g(\theta) - V_g(\theta + \pi) 
& = \left(\sin \theta + \frac{\sin 3\theta}{3},
          \cos \theta - \frac{\cos 3\theta}{3}\right),
\end{align*}
so
\begin{align*}
YV_g(\theta) & = (V_g(\theta) - V_g(\theta + \pi)) \cdot (p(\theta) - p(\theta + \pi)) \\
& = 4 \sin \theta \cos \theta + \frac{2 \sin 3 \theta \cos \theta
                   - 2 \cos 3 \theta \sin \theta}{3} \\
& = 2 \sin 2 \theta + \frac{2 \sin 2 \theta}{3} = \frac{8}{3} \sin 2 \theta,
\end{align*}
which is strictly positive on the edgeset, that is on $\theta \in
[\varepsilon,\nicefrac{\pi}{2}-\varepsilon]$.
But that strict positivity means that for any semipositive $\mu$, we would have
$\mu Y V_g$ greater than zero.  So there is no semipositive $\mu \in
Y(\X)^\perp$.

\subsection{Any open set will do}
\label{subsection:anyOpenSet}
Our success in showing that the almost-half-a-circle-of-struts example is not
$\X$-bar equivalent leads to the next question: If we remove \emph{any}
arbitrary open set of edges from the circle of struts, is the result no longer
$\X$-bar equivalent?  The answer is yes.\Index{circle of struts}[ ]
\begin{prop}
\label{prop:minBarEq}
The circle of struts with any nontrivial open set of edges removed is not
$\X$-bar equivalent.
\end{prop}
\begin{proof}
As before, we consider the set $\E$ to be parameterized by $\theta \in \qrp$.
Since open sets on $\E$ and $\V$ are (in this case) unions of open intervals,
removing a nontrivial open set (neither the empty set nor all of $\E$) from
$\E$ will mean that some open interval is free of edges and some other open
interval is not.  We'll create a variation in which the vertex curve moves
inward along that open interval and outward elsewhere, lengthing all of the
struts: a motion.

First, we want to show that for any real number $k > 1$, there is a variation
that is tangent to the circle at $0$ and $\nicefrac{\pi}{k}$, pointing into the
interval $(0,\nicefrac{\pi}{k})$; that points away from the center of the
circle at all points between $0$ and $\nicefrac{\pi}{k}$ and whose derivative
is radial (see \vref{fig:BuildingBlock}).  
\begin{figure}
\hfil
\begin{tikzpicture}[scale=4]
\draw[vcurve] (1,0) arc (0:180:1cm);
\draw[->] (  1.000,  0.000) -- (  1.000,  0.125);
\draw[->] (  0.999,  0.044) -- (  1.043,  0.170);
\draw[->] (  0.996,  0.087) -- (  1.082,  0.216);
\draw[->] (  0.991,  0.131) -- (  1.119,  0.264);
\draw[->] (  0.985,  0.174) -- (  1.151,  0.313);
\draw[->] (  0.976,  0.216) -- (  1.178,  0.363);
\draw[->] (  0.966,  0.259) -- (  1.199,  0.413);
\draw[->] (  0.954,  0.301) -- (  1.215,  0.463);
\draw[->] (  0.940,  0.342) -- (  1.223,  0.512);
\draw[->] (  0.924,  0.383) -- (  1.226,  0.559);
\draw[->] (  0.906,  0.423) -- (  1.221,  0.605);
\draw[->] (  0.887,  0.462) -- (  1.209,  0.648);
\draw[->] (  0.866,  0.500) -- (  1.191,  0.688);
\draw[->] (  0.843,  0.537) -- (  1.166,  0.723);
\draw[->] (  0.819,  0.574) -- (  1.134,  0.755);
\draw[->] (  0.793,  0.609) -- (  1.097,  0.782);
\draw[->] (  0.766,  0.643) -- (  1.055,  0.804);
\draw[->] (  0.737,  0.676) -- (  1.008,  0.820);
\draw[->] (  0.707,  0.707) -- (  0.957,  0.832);
\draw[->] (  0.676,  0.737) -- (  0.903,  0.839);
\draw[->] (  0.643,  0.766) -- (  0.846,  0.840);
\draw[->] (  0.609,  0.793) -- (  0.788,  0.837);
\draw[->] (  0.574,  0.819) -- (  0.728,  0.829);
\draw[->] (  0.537,  0.843) -- (  0.668,  0.818);
\draw[->] (  0.500,  0.866) -- (  0.608,  0.804);
\begin{scope}[rotate=90,black!50!white]
\draw[->] (  1.000,  0.000) -- (  1.000, -0.125);
\draw[->] (  0.999,  0.044) -- (  0.956, -0.082);
\draw[->] (  0.996,  0.087) -- (  0.910, -0.042);
\draw[->] (  0.991,  0.131) -- (  0.864, -0.003);
\draw[->] (  0.985,  0.174) -- (  0.819,  0.034);
\draw[->] (  0.976,  0.216) -- (  0.775,  0.070);
\draw[->] (  0.966,  0.259) -- (  0.733,  0.105);
\draw[->] (  0.954,  0.301) -- (  0.693,  0.139);
\draw[->] (  0.940,  0.342) -- (  0.656,  0.172);
\draw[->] (  0.924,  0.383) -- (  0.622,  0.206);
\draw[->] (  0.906,  0.423) -- (  0.592,  0.240);
\draw[->] (  0.887,  0.462) -- (  0.565,  0.276);
\draw[->] (  0.866,  0.500) -- (  0.541,  0.312);
\draw[->] (  0.843,  0.537) -- (  0.521,  0.351);
\draw[->] (  0.819,  0.574) -- (  0.504,  0.392);
\draw[->] (  0.793,  0.609) -- (  0.489,  0.436);
\draw[->] (  0.766,  0.643) -- (  0.477,  0.482);
\draw[->] (  0.737,  0.676) -- (  0.467,  0.531);
\draw[->] (  0.707,  0.707) -- (  0.457,  0.582);
\draw[->] (  0.676,  0.737) -- (  0.448,  0.636);
\draw[->] (  0.643,  0.766) -- (  0.439,  0.692);
\draw[->] (  0.609,  0.793) -- (  0.430,  0.750);
\draw[->] (  0.574,  0.819) -- (  0.419,  0.809);
\draw[->] (  0.537,  0.843) -- (  0.406,  0.869);
\draw[->] (  0.500,  0.866) -- (  0.392,  0.929);
\end{scope}
\end{tikzpicture}
\hfil
\capt[The building-block variation.]{The building-block variation for $k = 3$
(shown on the interval $\theta \in (0,\nicefrac{\pi}{3})$).  It is tangential
to the circle at the ends of the interval $(0,\nicefrac{\pi}{k})$, pointing
toward the interval.  Throughout the interval it points out of the circle.  Its
derivative is everywhere radial.  For interest's sake, a variation that also
has $k=3$, but for which $m = -1$ instead of $m=1$ is shown in grey on the
interval $\theta \in (\nicefrac{\pi}{2},\nicefrac{5\pi}{6})$.}
\label{fig:BuildingBlock}
\end{figure}

Consider the variation (with $m > 0$) defined by
\begin{align*}
V(\theta) & = \frac{m}{2(k^2-1)} 
  \Big( (k+1)\sin((k-1)\theta) + (k-1)\sin((k+1)\theta), \\
& \hspace*{1.5in}
  (k+1)\cos((k-1)\theta) - (k-1)\cos((k+1)\theta) \Big).
\end{align*}

Certainly, $V(0) =  \frac{m}{k^2-1} \left( 0, 1 \right)$, which is tangential
to the circle and points toward the interval.  Also, $V(\nicefrac{\pi}{k}) =
\frac{m}{k^2-1} \left(\sin(\nicefrac{\pi}{k}), -\cos(\nicefrac{\pi}{k})
\right)$, which is also tangential and points toward the interval.  Both of
these vectors have length $\frac{m}{k^2-1}$.  If we take the dot product
$V(\theta) \cdot (\cos \theta, \sin \theta)$, we get $\frac{k}{k^2 - 1} \sin(k
\theta)$ which is positive in the interval $(0,\nicefrac{\pi}{k})$.  Finally,
direct calculation gives us 
$$V'(\theta) = m \cos(k\theta) (\cos \theta, \sin \theta)$$
which points radially outward on the interval $(0,\nicefrac{\pi}{2k})$ and
radially inward on the interval $(\nicefrac{\pi}{2k},\nicefrac{\pi}{k})$.

Of course, if $m < 0$, the derivative remains radial and the variation
tangential at the endpoints.  The only change is that the variation now points
away from the interval at the ends and has negative dot product with the 
outward pointing normals in the interval.

We construct a variation in the following fashion.  Given some open set $U$ of
edges which are to be removed from $\E$, we find an open interval
$(\theta_1,\theta_2)$ which contains all of $\E \setminus U$.  As $U$ is open
and thus $\E \setminus U$ is closed, $\theta_1$ and $\theta_2$ are not in 
$\E \setminus U$.  So we can use the above construction to create four partial
variations, one for each of the intervals $[\theta_1,\theta_2]$, $[\theta_2,
\pi + \theta_1]$, $[\pi + \theta_1, \pi + \theta_2]$ and $[\pi + \theta_2, 2\pi
+ \theta_1]$.  We will create them so that at $\theta_1$, $\theta_2$, $\pi +
\theta_1$ and $\pi + \theta_2$, the vectors point into the intervals
$(\theta_1,\theta_2)$ and $(\pi + \theta_1, \pi + \theta_2)$ and expand
all of the struts.  By adjusting the respective $m$'s, we can match
the lengths of the vectors at the endpoints.  
\end{proof}

\begin{figure}
\hfil
\begin{tikzpicture}[scale=2]
\draw[vcurve] (0,0) circle (1cm);
\draw[strut] (55:1cm) -- (55:-1cm);
\draw[strut] (65:1cm) -- (65:-1cm);
\draw[strut] (70:1cm) -- (70:-1cm);
\draw[strut] (75:1cm) -- (75:-1cm);
\draw[strut] (80:1cm) -- (80:-1cm);
\draw[strut] (85:1cm) -- (85:-1cm);
\draw[strut] (90:1cm) -- (90:-1cm);
\draw[strut] (95:1cm) -- (95:-1cm);
\draw[strut] (105:1cm) -- (105:-1cm);
\draw[strut] (110:1cm) -- (110:-1cm);
\draw[strut] (115:1cm) -- (115:-1cm);
\draw[strut] (120:1cm) -- (120:-1cm);
\draw[strut] (125:1cm) -- (125:-1cm);
\draw[strut] (130:1cm) -- (130:-1cm);
\draw[strut] (145:1cm) -- (145:-1cm);
\draw[strut] (150:1cm) -- (150:-1cm);
\draw[ultra thin] (30:1.1cm) -- (30:1.5cm) 
  node[above right] {$\nicefrac{\pi}{6}$};
\draw[ultra thin] (165:1.1cm) -- (165:1.5cm) 
  node[above left] {$\nicefrac{11\pi}{12}$};
\draw[ultra thin] (34:1.3cm) arc (34:161:1.3cm);
\end{tikzpicture}
\hfil
\begin{tikzpicture}[scale=2]
\draw[vcurve] (0,0) circle(1cm);
\begin{scope}
\draw[strut] (55:1cm) -- (55:-1cm);
\draw[strut] (65:1cm) -- (65:-1cm);
\draw[strut] (70:1cm) -- (70:-1cm);
\draw[strut] (75:1cm) -- (75:-1cm);
\draw[strut] (80:1cm) -- (80:-1cm);
\draw[strut] (85:1cm) -- (85:-1cm);
\draw[strut] (90:1cm) -- (90:-1cm);
\draw[strut] (95:1cm) -- (95:-1cm);
\draw[strut] (105:1cm) -- (105:-1cm);
\draw[strut] (110:1cm) -- (110:-1cm);
\draw[strut] (115:1cm) -- (115:-1cm);
\draw[strut] (120:1cm) -- (120:-1cm);
\draw[strut] (125:1cm) -- (125:-1cm);
\draw[strut] (130:1cm) -- (130:-1cm);
\draw[strut] (145:1cm) -- (145:-1cm);
\draw[strut] (150:1cm) -- (150:-1cm);
\end{scope}
\begin{scope}[rotate=-15]
\draw[->] (  1.000,  0.000) -- (  1.000, -0.133);
\draw[->] (  0.999,  0.044) -- (  0.912, -0.092);
\draw[->] (  0.996,  0.087) -- (  0.825, -0.054);
\draw[->] (  0.991,  0.131) -- (  0.742, -0.019);
\draw[->] (  0.985,  0.174) -- (  0.665,  0.014);
\draw[->] (  0.976,  0.216) -- (  0.596,  0.044);
\draw[->] (  0.966,  0.259) -- (  0.537,  0.075);
\draw[->] (  0.954,  0.301) -- (  0.489,  0.107);
\draw[->] (  0.940,  0.342) -- (  0.454,  0.141);
\draw[->] (  0.924,  0.383) -- (  0.431,  0.179);
\draw[->] (  0.906,  0.423) -- (  0.421,  0.222);
\draw[->] (  0.887,  0.462) -- (  0.421,  0.271);
\draw[->] (  0.866,  0.500) -- (  0.433,  0.327);
\draw[->] (  0.843,  0.537) -- (  0.453,  0.390);
\draw[->] (  0.819,  0.574) -- (  0.480,  0.461);
\draw[->] (  0.793,  0.609) -- (  0.511,  0.538);
\draw[->] (  0.766,  0.643) -- (  0.546,  0.622);
\draw[->] (  0.737,  0.676) -- (  0.580,  0.710);
\draw[->] (  0.707,  0.707) -- (  0.613,  0.801);
\end{scope}
\begin{scope}[rotate=30]
\draw[->] (  1.000,  0.000) -- (  1.000,  0.133);
\draw[->] (  0.999,  0.044) -- (  1.004,  0.177);
\draw[->] (  0.996,  0.087) -- (  1.005,  0.221);
\draw[->] (  0.991,  0.131) -- (  1.005,  0.265);
\draw[->] (  0.985,  0.174) -- (  1.003,  0.309);
\draw[->] (  0.976,  0.216) -- (  0.998,  0.352);
\draw[->] (  0.966,  0.259) -- (  0.992,  0.396);
\draw[->] (  0.954,  0.301) -- (  0.984,  0.439);
\draw[->] (  0.940,  0.342) -- (  0.974,  0.481);
\draw[->] (  0.924,  0.383) -- (  0.962,  0.523);
\draw[->] (  0.906,  0.423) -- (  0.948,  0.565);
\draw[->] (  0.887,  0.462) -- (  0.932,  0.606);
\draw[->] (  0.866,  0.500) -- (  0.914,  0.646);
\draw[->] (  0.843,  0.537) -- (  0.894,  0.685);
\draw[->] (  0.819,  0.574) -- (  0.873,  0.723);
\draw[->] (  0.793,  0.609) -- (  0.849,  0.760);
\draw[->] (  0.766,  0.643) -- (  0.824,  0.795);
\draw[->] (  0.737,  0.676) -- (  0.797,  0.830);
\draw[->] (  0.707,  0.707) -- (  0.769,  0.863);
\draw[->] (  0.676,  0.737) -- (  0.739,  0.895);
\draw[->] (  0.643,  0.766) -- (  0.707,  0.925);
\draw[->] (  0.609,  0.793) -- (  0.674,  0.954);
\draw[->] (  0.574,  0.819) -- (  0.640,  0.981);
\draw[->] (  0.537,  0.843) -- (  0.604,  1.006);
\draw[->] (  0.500,  0.866) -- (  0.567,  1.029);
\draw[->] (  0.462,  0.887) -- (  0.530,  1.051);
\draw[->] (  0.423,  0.906) -- (  0.491,  1.070);
\draw[->] (  0.383,  0.924) -- (  0.451,  1.088);
\draw[->] (  0.342,  0.940) -- (  0.410,  1.104);
\draw[->] (  0.301,  0.954) -- (  0.369,  1.117);
\draw[->] (  0.259,  0.966) -- (  0.326,  1.129);
\draw[->] (  0.216,  0.976) -- (  0.284,  1.139);
\draw[->] (  0.174,  0.985) -- (  0.241,  1.146);
\draw[->] (  0.131,  0.991) -- (  0.198,  1.151);
\draw[->] (  0.087,  0.996) -- (  0.154,  1.154);
\draw[->] (  0.044,  0.999) -- (  0.110,  1.155);
\draw[->] (  0.000,  1.000) -- (  0.067,  1.154);
\draw[->] ( -0.044,  0.999) -- (  0.023,  1.151);
\draw[->] ( -0.087,  0.996) -- ( -0.020,  1.145);
\draw[->] ( -0.131,  0.991) -- ( -0.063,  1.138);
\draw[->] ( -0.174,  0.985) -- ( -0.106,  1.128);
\draw[->] ( -0.216,  0.976) -- ( -0.148,  1.116);
\draw[->] ( -0.259,  0.966) -- ( -0.190,  1.103);
\draw[->] ( -0.301,  0.954) -- ( -0.231,  1.087);
\draw[->] ( -0.342,  0.940) -- ( -0.271,  1.070);
\draw[->] ( -0.383,  0.924) -- ( -0.310,  1.050);
\draw[->] ( -0.423,  0.906) -- ( -0.348,  1.029);
\draw[->] ( -0.462,  0.887) -- ( -0.386,  1.006);
\draw[->] ( -0.500,  0.866) -- ( -0.422,  0.981);
\draw[->] ( -0.537,  0.843) -- ( -0.457,  0.955);
\draw[->] ( -0.574,  0.819) -- ( -0.491,  0.927);
\draw[->] ( -0.609,  0.793) -- ( -0.523,  0.898);
\draw[->] ( -0.643,  0.766) -- ( -0.555,  0.867);
\draw[->] ( -0.676,  0.737) -- ( -0.584,  0.835);
\draw[->] ( -0.707,  0.707) -- ( -0.613,  0.801);
\end{scope}
\begin{scope}[rotate=165]
\draw[->] (  1.000,  0.000) -- (  1.000, -0.133);
\draw[->] (  0.999,  0.044) -- (  0.912, -0.092);
\draw[->] (  0.996,  0.087) -- (  0.825, -0.054);
\draw[->] (  0.991,  0.131) -- (  0.742, -0.019);
\draw[->] (  0.985,  0.174) -- (  0.665,  0.014);
\draw[->] (  0.976,  0.216) -- (  0.596,  0.044);
\draw[->] (  0.966,  0.259) -- (  0.537,  0.075);
\draw[->] (  0.954,  0.301) -- (  0.489,  0.107);
\draw[->] (  0.940,  0.342) -- (  0.454,  0.141);
\draw[->] (  0.924,  0.383) -- (  0.431,  0.179);
\draw[->] (  0.906,  0.423) -- (  0.421,  0.222);
\draw[->] (  0.887,  0.462) -- (  0.421,  0.271);
\draw[->] (  0.866,  0.500) -- (  0.433,  0.327);
\draw[->] (  0.843,  0.537) -- (  0.453,  0.390);
\draw[->] (  0.819,  0.574) -- (  0.480,  0.461);
\draw[->] (  0.793,  0.609) -- (  0.511,  0.538);
\draw[->] (  0.766,  0.643) -- (  0.546,  0.622);
\draw[->] (  0.737,  0.676) -- (  0.580,  0.710);
\draw[->] (  0.707,  0.707) -- (  0.613,  0.801);
\end{scope}
\begin{scope}[rotate=210]
\draw[->] (  1.000,  0.000) -- (  1.000,  0.133);
\draw[->] (  0.999,  0.044) -- (  1.004,  0.177);
\draw[->] (  0.996,  0.087) -- (  1.005,  0.221);
\draw[->] (  0.991,  0.131) -- (  1.005,  0.265);
\draw[->] (  0.985,  0.174) -- (  1.003,  0.309);
\draw[->] (  0.976,  0.216) -- (  0.998,  0.352);
\draw[->] (  0.966,  0.259) -- (  0.992,  0.396);
\draw[->] (  0.954,  0.301) -- (  0.984,  0.439);
\draw[->] (  0.940,  0.342) -- (  0.974,  0.481);
\draw[->] (  0.924,  0.383) -- (  0.962,  0.523);
\draw[->] (  0.906,  0.423) -- (  0.948,  0.565);
\draw[->] (  0.887,  0.462) -- (  0.932,  0.606);
\draw[->] (  0.866,  0.500) -- (  0.914,  0.646);
\draw[->] (  0.843,  0.537) -- (  0.894,  0.685);
\draw[->] (  0.819,  0.574) -- (  0.873,  0.723);
\draw[->] (  0.793,  0.609) -- (  0.849,  0.760);
\draw[->] (  0.766,  0.643) -- (  0.824,  0.795);
\draw[->] (  0.737,  0.676) -- (  0.797,  0.830);
\draw[->] (  0.707,  0.707) -- (  0.769,  0.863);
\draw[->] (  0.676,  0.737) -- (  0.739,  0.895);
\draw[->] (  0.643,  0.766) -- (  0.707,  0.925);
\draw[->] (  0.609,  0.793) -- (  0.674,  0.954);
\draw[->] (  0.574,  0.819) -- (  0.640,  0.981);
\draw[->] (  0.537,  0.843) -- (  0.604,  1.006);
\draw[->] (  0.500,  0.866) -- (  0.567,  1.029);
\draw[->] (  0.462,  0.887) -- (  0.530,  1.051);
\draw[->] (  0.423,  0.906) -- (  0.491,  1.070);
\draw[->] (  0.383,  0.924) -- (  0.451,  1.088);
\draw[->] (  0.342,  0.940) -- (  0.410,  1.104);
\draw[->] (  0.301,  0.954) -- (  0.369,  1.117);
\draw[->] (  0.259,  0.966) -- (  0.326,  1.129);
\draw[->] (  0.216,  0.976) -- (  0.284,  1.139);
\draw[->] (  0.174,  0.985) -- (  0.241,  1.146);
\draw[->] (  0.131,  0.991) -- (  0.198,  1.151);
\draw[->] (  0.087,  0.996) -- (  0.154,  1.154);
\draw[->] (  0.044,  0.999) -- (  0.110,  1.155);
\draw[->] (  0.000,  1.000) -- (  0.067,  1.154);
\draw[->] ( -0.044,  0.999) -- (  0.023,  1.151);
\draw[->] ( -0.087,  0.996) -- ( -0.020,  1.145);
\draw[->] ( -0.131,  0.991) -- ( -0.063,  1.138);
\draw[->] ( -0.174,  0.985) -- ( -0.106,  1.128);
\draw[->] ( -0.216,  0.976) -- ( -0.148,  1.116);
\draw[->] ( -0.259,  0.966) -- ( -0.190,  1.103);
\draw[->] ( -0.301,  0.954) -- ( -0.231,  1.087);
\draw[->] ( -0.342,  0.940) -- ( -0.271,  1.070);
\draw[->] ( -0.383,  0.924) -- ( -0.310,  1.050);
\draw[->] ( -0.423,  0.906) -- ( -0.348,  1.029);
\draw[->] ( -0.462,  0.887) -- ( -0.386,  1.006);
\draw[->] ( -0.500,  0.866) -- ( -0.422,  0.981);
\draw[->] ( -0.537,  0.843) -- ( -0.457,  0.955);
\draw[->] ( -0.574,  0.819) -- ( -0.491,  0.927);
\draw[->] ( -0.609,  0.793) -- ( -0.523,  0.898);
\draw[->] ( -0.643,  0.766) -- ( -0.555,  0.867);
\draw[->] ( -0.676,  0.737) -- ( -0.584,  0.835);
\draw[->] ( -0.707,  0.707) -- ( -0.613,  0.801);
\end{scope}
\end{tikzpicture}
\capt[All struts lie in $(\nicefrac{\pi}{6},\nicefrac{11\pi}{12})$.]{All struts
lie in $(\nicefrac{\pi}{6},\nicefrac{11\pi}{12})$, so we can expand all the
struts by pulling in where there are no struts.}
\label{fig:StrutsInInterval}
\end{figure}
For example, if we can establish that the interval
$(\nicefrac{\pi}{6},\nicefrac{11\pi}{12})$ contains all of $\E \setminus U$
(see \vref{fig:StrutsInInterval}), then we can patch together a variation which
pulls the circle inward in the intervals
$(-\nicefrac{\pi}{12},\nicefrac{\pi}{6})$ and
$(\nicefrac{11\pi}{12},\nicefrac{7\pi}{6})$ and pushes it outward on
$(\nicefrac{\pi}{6},\nicefrac{11\pi}{12})$ and
$(\nicefrac{7\pi}{6},\nicefrac{23\pi}{12})$.  To do the patching, we could use
$m_1 = -2$ for the ``non-strut'' portion of the variation and $m_2 =
\nicefrac{14}{135}$ for remainder so the vectors at the transition points are
length: 
$$
\frac{2}{15} = \frac{m_1}{4^2 - 1} = \frac{m_2}{\left(\nicefrac{4}{3}\right)^2
- 1} = \frac{9}{7} \cdot \frac{14}{135} = \frac{2}{15}.
$$

\section{\texorpdfstring{Minimal $\X$-bar Equivalence}{Minimal X-bar
  Equivalence}}
\label{section:MBE}
\subsection{Definition and example}
All we lack to complete our results now is a theorem which says that any
$\X$-bar equivalent tensegrity has a strictly positive measure.  We don't have
such a theorem for the general case, but there is a class of tensegrities for
which we do get the desired result.

\begin{definition}
We term a tensegrity $G(p)$ \emph{minimally $\X$-bar
equivalent}\Index{minimally $\X$-bar@X-bar equivalent}[23!1 ] if it is $\X$-bar
equivalent but no nonempty subtensegrity formed by removing an open set of
edges from it is $\X$-bar equivalent.
\end{definition}

\begin{lemma}
\label{lemma:notJustSemi}
A minimally $\X$-bar equivalent tensegrity $G(p)$ admits no semipositive stress
that is not also strictly positive.
\end{lemma}
\begin{proof}
By Corollary \vref{corollary:partiallyBE}, if the tensegrity admits a 
semipositive (but not strictly positive) stress $\mu$, then the subtensegrity
whose edgeset is $\supp \mu$ is $\X$-bar equivalent, and $\E - \supp \mu$ is
nonempty and open, so $G(p)$ is not minimally $\X$-bar equivalent.
\end{proof}
\begin{theorem}[Third Main Theorem]
\label{theorem:minimallyBarEquivalent}
\statementMinimallyBarEquivalent
\end{theorem}
\begin{proof}
By Theorem \ref{theorem:nonnegativeMeasure}, $G(p)$ has a semipositive stress,
and by Lemma \ref{lemma:notJustSemi}, that stress must be strictly positive.
\end{proof}

Having seen that a minimally $\X$-bar equivalent tensegrity cannot admit a
(solely) semipositive stress, it seems tempting to conjecture:
\begin{conjecture}
\label{conj:minXbarConj}
\statementminXbarConj
\end{conjecture}

\renewcommand{\V}{\gamma}
What sort of tensegrity is minimally $\X$-bar equivalent?
\vref{subsection:anyOpenSet} showed that the circle of struts is.\Index{circle
of struts}[ ]

Interestingly enough, if we make the circle into a square filled with either
antipodal struts or with struts parallel to the edges (the left and right
pictures in \vref{fig:squareOfStruts}), the tensegrity ceases to be $\X$-bar
equivalent ($\X$ is still the motions that are local isometries on $\V$).  

The center picture of \ref{fig:squareOfStruts} shows a motion that induces a
semipositive load on~$\E$.  In the first case, that load is zero only on the
diagonals and in the second case it is zero only on the outermost struts.
\begin{figure}[ht]
\hspace*{\fill}
\begin{tikzpicture}
{
  \clip (0,0) -- (2,0) -- (2,2) -- (0,2) -- cycle;
  \foreach \x in {0,15,...,181} {
    \draw[strut,xshift=1cm,yshift=1cm] (\x:2cm) -- (\x:-2cm);
  }
}
\draw[vcurve] (0,0) -- (2,0) -- (2,2) -- (0,2) -- cycle;
\end{tikzpicture}
\hspace*{\fill}
\begin{tikzpicture}
\foreach \x/\y in {0.2/0.1, 0.4/0.2, 0.6/0.3, 0.8/0.4, 1/0.5, 1.2/0.4, 1.4/0.3,
                1.6/0.2, 1.8/0.1} {
  \draw[->] (0,\x) -- (-\y,\x);
  \draw[->] (\x,0) -- (\x,-\y);
  \draw[->,yshift=2cm] (\x,0) -- (\x,\y);
  \draw[->,xshift=2cm] (0,\x) -- (\y,\x);
}
\draw[vcurve] (0,0) -- (2,0) -- (2,2) -- (0,2) -- cycle;
\end{tikzpicture}
\hspace*{\fill}
\begin{tikzpicture}
\draw[vcurve] (0,0) -- (2,0) -- (2,2) -- (0,2) -- cycle;
\foreach \x in {0.1,0.4,...,1.95} {
  \draw[strut] (0,\x) -- (2,\x) (\x,0) -- (\x,2);
}
\end{tikzpicture}
\hspace*{\fill}
\capt[Two squares of struts and a semipositive motion.]{Two squares of struts
and a semipositive motion that works for either.}
\label{fig:squareOfStruts}
\end{figure}
So any strictly positive stress $\mu$ will give $\mu YV > 0$ for this motion.
On the other hand, the semipositive stress that just assigns an atom to each
diagonal (resp.\ each outermost strut) annihilates $Y(\X)$, as we will show
below.

Now $\X$ contains those vector fields for which $V'(v) \cdot \V'(v) = 0$
except on a set of measure zero.  So along the top and bottom lines of
vertices, the horizontal component of $V(v)$ cannot vary, while along the
sides, its vertical component cannot vary.  That means that 
\begin{align*}
(V((1,0)) - V((0,0))) \cdot (1,0) & = 0 &
(V((0,1)) - V((0,0))) \cdot (0,1) & = 0 \\
(V((1,1)) - V((0,1))) \cdot (1,0) & = 0 &
(V((1,1)) - V((1,0))) \cdot (0,1) & = 0 
\end{align*}

So taking $\mu$ to be an atom of size $1$ on each outermost strut
of the right hand tensegrity, we get
\begin{align*}
\mu YV = & (V((1,0)) - V((0,0))) \cdot (1,0) + 
           (V((1,1)) - V((1,0))) \cdot (0,1) + \\
         & (V((0,1)) - V((0,0))) \cdot (0,1) +
           (V((1,1)) - V((0,1))) \cdot (1,0) = 0.
\end{align*}
On the other hand, taking $\mu$ to give an atom of size $1$ on 
each diagonal strut of the left hand tensegrity, we get
\begin{align*}
\mu YV & = (V((1,1)) - V((0,0))) \cdot (1,1) +
           (V((1,0)) - V((0,1))) \cdot (1,-1) \\
       & = (V((1,1)) - V((1,0)) + V((1,0)) - V((0,0))) \cdot (1,1) \\
       & \hspace*{1cm} +
           (V((1,0)) - V((1,1)) + V((1,1)) - V((0,1))) \cdot (1,-1) \\
       & = (V((1,1)) - V((1,0))) \cdot (1,0)
         + (V((1,0)) - V((0,0))) \cdot (0,1) \\
       & \hspace*{1cm}
         + (V((1,0)) - V((1,1))) \cdot (1,0)
         + (V((1,1)) - V((0,1))) \cdot (0,-1) \\
       & = (V((1,0)) - V((0,0))) \cdot (0,1) 
         + (V((1,1)) - V((0,1))) \cdot (0,-1) \\
       & = (V((1,0)) - V((0,0))
         +  V((0,1)) - V((1,1))) \cdot (0,1) \\
       & = (V((1,0)) - V((1,1))) \cdot (0,1) 
         + (V((0,1)) - V((0,0))) \cdot (0,1) = 0 
\end{align*}
So in each case we have a semipositive stress and the two are partially
$\X$-bar equivalent.  Alternatively, that semipositive stress is a strictly
positive stress on $\supp \mu$, so the subtensegrities (shown in
\vref{fig:minimalSquares}) consisting of $\V$ and the mentioned edges are
$\X$-bar equivalent.
\begin{figure}[ht]
\hfil
\begin{tikzpicture}
\draw[vcurve] (0,0) rectangle (2,2);
\draw[strut] (0,0) -- (2,2) (0,2) -- (2,0);
\end{tikzpicture}
\hfil
\begin{tikzpicture}
\draw[vcurve,scale=1.04] (-1,-1) rectangle (1,1);
\draw[strut] (-1,-1) rectangle (1,1);
\end{tikzpicture}
\hfil
\capt{The $\X$-bar equivalent subtensegrities of the squares of struts.}
\label{fig:minimalSquares}
\end{figure}
\renewcommand{\V}{\mathscr{V}}

\subsection{Finite is not minimal}
\label{subsection:minimal}
It is worth a moment to remember, prior to jumping into the next section, that
not all finite tensegrities are minimal.  The crossed square\Index[|(]{crossed
square}[ ] with one of the struts replaced with a bar (\vref*{fig:notMinimal})
is just such an example.  We saw, in \ref{subsection:TheOtherDirection}, that
the stresses of that tensegrity look like $\begin{bmatrix} a & a & a & a & b &
a & (b-a) \end{bmatrix}^\top$ where $a \in [0,b]$.

If we set $a = 0$, we get a stress which is strictly positive on the bar but 
zero everywhere else.  If we set $a = b$, we get a stress which is strictly
positive on the usual crossed square and zero on that new ``cable''
(see \vref{fig:minSub}).  But these two subtensegrities are the only minimally
bar equivalent subtensegrities in this example.\footnote{We can see this by
looking at the vertices.  If any edge is removed from vertex $1$ or vertex $3$,
all three edges must be removed.  That leaves the bar to deal with.  If we 
remove the strut portion of the bar, we must remove all of the other edges
at $4$ and $2$ (and hence also at $1$ and $3$).}
\begin{figure}[ht]
\hspace*{\fill}
\begin{tikzpicture}
\begin{scope}[xshift=-5cm]
\draw[strut] (0,0) node [vertex] {1} (2,2) node [vertex] {3};
\draw[bar] (0,2) node [vertex] {4} -- (2,0) node [vertex] {2};
\end{scope}
\draw[->] (-0.5,1) .. controls (-1,1.5) and (-2,1.5) .. 
  node[above] {$a \to 0$} (-2.5,1);
\draw[cable] (0,0) -- (2,0) -- (2,2) -- (0,2) -- (0,0) -- (2,0);
\draw[strut] (0,0) node [vertex] {1} -- (2,2) node [vertex] {3};
\draw[bar] (0,2) node [vertex] {4} -- (2,0) node [vertex] {2};
\draw[->] (2.5,1) .. controls (3,1.5) and (4,1.5) .. 
  node[above] {$a \to 3$} (4.5,1);
\begin{scope}[xshift=5cm]
\draw[cable] (0,0) -- (2,0) -- (2,2) -- (0,2) -- (0,0) -- (2,0);
\draw[strut] (0,0) node [vertex] {1} -- (2,2) node [vertex] {3}
             (0,2) node [vertex] {4} -- (2,0) node [vertex] {2};
\end{scope}
\end{tikzpicture}
\hspace*{\fill}
\capt{The non-minimal crossed square with its minimal subtensegrities.}
\label{fig:minSub}
\end{figure}
\Index[|)]{crossed square}[ ]

We note that while the tensegrity itself is not minimally bar equivalent, 
every edge in it belongs to a minimally bar equivalent subtensegrity.  Perhaps
that is key to having a strictly positive stress.

\section{Covered Tensegrities}
\label{section:Covered}
Let's consider tensegrities which are covered by minimal subtensegrities.
Better yet, let's consider tensegrities which are almost covered by them.
\begin{definition}
We say that a tensegrity is \emph{countably covered}\IndexDef{countably
covered}[ ] by some class of subtensegrities, if there exists a countable,
dense subset $\{e_n\} \subset \E$ such that every $e_n$ belongs to a 
subtensegrity of that class.
\end{definition}
Having a countable, dense subset is not an unreasonable thing to ask, as $\E$
is compact and metrizable and hence second countable.\footnote{This is an
exercise (with hint) in \ocite{Munkres}*{p.\ 194}.  If $X$ is metrizable and
compact, then for each positive integer $n$, we can cover $X$ with balls of
size~$\nicefrac{1}{n}$.  Since $X$ is compact, each such cover has a finite
subcover,~$\mathscr{A}_n$.  The union of all of the $\mathscr{A}_n$ forms a
countable basis for the (metric) topology on $X$ and the centers of the balls
give a countable, dense subset of~$X$.}

\begin{theorem}[Fourth Main Theorem]
\label{theorem:CountablyCovered}
\statementCountablyCovered
\end{theorem}
\begin{proof}
Let $\{e_n\}$ be the countable, dense subset of edges from the definition of
``countably covered''.  For each $e_n$ select a subtensegrity $G_n$ with
edgeset $\E_n$ of which $e_n$ is a member and which has a strictly positive
stress (in the case where $e_n$ belongs to infinitely many subtensegrities,
this may involve the Axiom of Choice\Index{Axiom of Choice}[ ]).  

Scale that stress so that its norm is $\frac{1}{2^n}$ and call it~$\muh_n$.
Finally, extend $\muh_n$ to be a stress on the entire $G(p)$ by setting, for
each open $U \subset \E$, $\mu_n(U) = \muh_n(U \cap \E_n)$.

Now we can construct
$$\mu = \sum_{n=1}^\infty \mu_n.$$
What do we know about $\mu$?  Well, if $V \in \X$, then every $\mu_n(YV) =
0$ since $\mu_n$ is $0$ outside $\E_n$ and takes $YV$ to zero on $\E_n$, so 
$\mu(YV) = 0$.  That makes $\mu$ a stress.  Furthermore, if $U \subset \E$ is
open and nonempty, then it must contain at least one of the~$e_n$.  But then,
by construction, $\mu_n(U) > 0$ and so $\mu(U) > 0$.
\end{proof}
So if we have a tensegrity which is countably covered either by minimally 
$\X$-bar equivalent subtensegrities or by finite, $\X$-bar equivalent 
subtensegrities (just use the stress from \ocite{MR610958} to build an atomic
measure) or by some combination of the two, then our tensegrity has a 
strictly positive stress and thus is $\X$-bar equivalent.

What might such a creature look like?  One simple example would be to move the
circle-of-struts example ``up a dimension''.  A ``cylinder of struts'' is
represented in \vref{fig:SCyl}.
\begin{figure}[ht]
\hspace*{\fill}
\begin{tikzpicture}
\foreach \x in {-1,-0.9,...,1.1} {
  \begin{scope}[cm={1,0,0,1,(\x,0)},cm={1,0.5,0,0.866,(\x,-\x)}]
    \path[draw=black,fill=white] (0,0) circle(1.5cm);
  \end{scope}
}
\begin{scope}
  \clip[cm={1,0.5,0,0.855,(2,-1)}] (0,0) circle(1.5cm);
  \foreach \x in {-0.85,-0.5,...,1.1} {
    \begin{scope}[cm={1,0,0,1,(\x,0)},cm={1,0.5,0,0.866,(\x,-\x)}]
      \draw[draw=white,double=black,line width=1.8pt,double distance=1.6pt]
        (0,0) -- (  0:1.5cm)
        (0,0) -- ( 30:1.5cm)
        (0,0) -- ( 60:1.5cm)
        (0,0) -- ( 90:1.5cm)
        (0,0) -- (120:1.5cm)
        (0,0) -- (150:1.5cm)
        (0,0) -- (180:1.5cm)
        (0,0) -- (210:1.5cm)
        (0,0) -- (240:1.5cm)
        (0,0) -- (270:1.5cm)
        (0,0) -- (300:1.5cm)
        (0,0) -- (330:1.5cm);
    \end{scope}
  }
\end{scope}
\end{tikzpicture}
\hspace*{\fill}
\capt{A cylinder filled with antipodal struts.}
\label{fig:SCyl}
\end{figure}
As with the circle of struts, our design variations are those which are local
isometries on the surface of the cylinder.  

Of course, any given strut in the cylinder of struts is part of a circle of
struts\Index{circle of struts}[ ].  So every strut is contained in a minimally
$\X$-bar equivalent subtensegrity and our tensegrity is (more than) countably
covered in minimally $\X$-bar equivalent subtensegrities.  By Theorem
\ref{theorem:CountablyCovered}, the tensegrity is $\X$-bar equivalent.  

\ref{chapter:Examples} contains a number of tensegrities which are countably
covered by finite subtensegrities.

Now, the ``countably covered tensegrities'' are a subset of all the $\X$-bar
equivalent tensegrities, but how large a subset is it?  That's not a settled
question, but it may be all of them:
\newlength{\spacelength}
\settowidth{\spacelength}{~}
\begin{conjecture}
\label{conj:CountCov}
\statementCountCov\hspace*{-\spacelength}\footnote{The matching conjecture,
that every $\X$-bar equivalent continuous tensegrity is countably covered by
finite, $\X$-bar equivalent subtensegrities, is false.  One counterexample is
the ``on a circle'' example of \vref{section:onACircle} when the cable ends are
an irrational fraction of the circle circumference apart.}
\end{conjecture}

Of course, if that turns out not to be true, we can retreat to the weaker
claim:
\begin{conjecture}
\label{conj:StrictlyPositive}
\statementStrictlyPositive
\end{conjecture}
\end{chapter}
% 
% \hrule
% 
% \begin{lemma}
% \label{lemma:CiCC}
% If there exists a dense subset $\E_m \subset \E$ such that every edge in 
% $\E_m$ belongs to a given class of subtensegrities, then there exists a
% countable dense subset for which that is also true.
% \end{lemma}
% \begin{proof}
% Cover $\E$ by balls of size $\nicefrac{1}{n}$ centered at points of $\E_m$.
% Since $\E_m$ is dense, this must be possible.  That cover admits a finite
% subcover.  Call it $\mathscr{A}_n$.  The collection of $\mathscr{A}_n$ for all
% positive integers $n$ is a countable union of finite sets and hence countable.
% Let $e \in \E$ and let $U$ be an open set containing $e$.  Then there is some
% open ball centered at $e$ of radius $\nicefrac{1}{m}$ for some positive
% integer $m$.  But then $\mathscr{A}_m$ must contain an open set which
% contains $e$.  The center of that open set must lie inside $U$.  So the
% centers of the collection of $\mathscr{A}_n$ form a countable dense subset of
% $\E$.
% \end{proof}

\chapter{Examples}
\label{chapter:Examples}
Here we have a collection of examples that shed light on various aspects of
the subject.  Unless stated otherwise, $\X = \Vf(\V)$ for these examples.

\section{Rectangle}
\label{section:rectangle}
We start with an example that is simple enough to analyze pretty thoroughly.
As we are primarily interested in bar equivalence in this paper, we usually
give only passing attention to infinitesimal rigidity.  But this example is
infinitesimally rigid and we take the time to prove it.

The rectangle example lies in $\RR^2$ and its has vertices arranged along two
intervals, say $(0,0)$ to $(2,0)$ and $(0,1)$ to~$(2,1)$.  Each vertex is
connected to one directly opposite it by a strut.  Every vertex is also
connected to the two opposite endpoints by cables (that means that the corner
vertices are connected vertically by strut-cable pairs, which is to say, bars).
Finally, the two upper corner vertices and the two lower corner vertices are
connected by horizontal struts (see \vref{fig:rectangle}).
\begin{figure}[ht]
\hspace*{\fill}
\begin{tikzpicture}[scale=4]
\foreach \x in {0.1, 0.3, ..., 1.99} {
  \draw[strut] (\x,0) -- (\x,1);
  \draw[cable] (0,0) -- (\x,1) -- (2,0) 
                    (0,1) -- (\x,0) -- (2,1);
}
\draw[strut] (0,0) -- (2,0) 
  (0,1) node [anchor=south east] {$(0,1)$} -- 
  (2,1) node [anchor=south west] {$(2,1)$};
\draw (1.3,0) node [anchor=north] {$(x,0)$};
\draw[bar,cap=round] (0,0) -- (0,1) (2,0) -- (2,1);
\end{tikzpicture}
\hspace*{\fill}
\capt{The rectangle example.}
\label{fig:rectangle}
\end{figure}
Remember that the horizontal struts only connect the corner vertices.  They
pass through, but do not affect, the other vertices.

We can think of the possible vector fields as four continuous functions from
$[0,2]$ to $\RR$, one for the vertical direction and one for the horizontal on
each interval.  Let's call them $f_{vt}$ (vertical top), $f_{vb}$, $f_{ht}$
and~$f_{hb}$.  Then $\im Y$ is as shown in \ref{table:imY} (where we've
omitted scaling factors since we will only be concerned with signs).
\begin{table}[ht]
\capt{The values of $YV$ for the different edges of the rectangle tensegrity.}
\label{table:imY}
$$\begin{array}{c|c}
\textbf{Edge} & \textbf{Value} \\
(x,0) \to (x,1) \text{ (strut)} & f_{vt}(x) - f_{vb}(x) \\
(x,0) \to (0,1) \text{ (cable)} & 
  -x(f_{hb}(x) - f_{ht}(0)) + f_{vb}(x) - f_{vt}(0) \\
(x,0) \to (2,1) \text{ (cable)} & 
  (2-x)(f_{hb}(x) - f_{ht}(2)) + f_{vb}(x) - f_{vt}(2) \\
(x,1) \to (0,0) \text{ (cable)} & 
  -x(f_{ht}(x) - f_{hb}(0)) - f_{vt}(x) + f_{vb}(0) \\
(x,1) \to (2,0) \text{ (cable)} & 
  (2-x)(f_{ht}(x) - f_{hb}(2)) - f_{vt}(x) + f_{vb}(2) \\
(0,0) \to (2,0) \text{ (strut)} & f_{hb}(0) - f_{hb}(2) \\
(0,1) \to (2,1) \text{ (strut)} & f_{ht}(0) - f_{ht}(2) \\
\end{array}$$
\end{table}

So what are the implications?  We'll first examine the associated bar
framework.  As we work through the families of constraints above, we'll call
them struts and cables, but in every case we'll be treating them as equalities,
as if they were bars.  
  
Because of the vertical struts, we must have $f_{vt}(x) = f_{vb}(x)$.  
We can eliminate the Euclidean motions\Index{Euclidean motion} by
requiring that $f_{vb}(0) = f_{hb}(0) = f_{ht}(0) = 0$.  That results 
in having $f_{vt}(0) = 0$, since $f_{vt}(0) = f_{vb}(0)$.  Our first family of
cables then gives us that 
\begin{equation}
f_{vb}(x) = x f_{hb}(x)
\label{eq:CablesOne}
\end{equation}
and the third set means that 
\begin{equation}
f_{vt}(x) = - x f_{ht}(x).
\label{eq:CablesThree}
\end{equation}
Since we have changed all edges to bars, the final two struts give us 
$f_{hb}(2) = f_{ht}(2) = 0$.  Now using the second and fourth cable families
with $x = 0$, we get 
$$2 (f_{hb}(0) - f_{ht}(2)) + f_{vb}(0) - f_{vt}(2) \Rightarrow f_{vt}(2) = 0$$
and also $f_{vb}(2) = 0$.  For other values of $x$ these give
\begin{equation}
f_{vb}(x) = -(2-x)f_{hb}(x) \text{\quad and \quad} 
f_{vt}(x) = (2-x)f_{ht}(x).
\label{eq:CablesTwoFour}
\end{equation}
So if we combine \reftextlabelrange{eq:CablesOne}{eq:CablesTwoFour}, for all
values of $x$ we get
$$-x f_{ht}(x) = (2-x) f_{ht}(x) \Rightarrow 2 f_{ht}(x) = 0 \Rightarrow
f_{ht}(x) = 0$$
and similarly $f_{hb}(x) = 0$ and thus immediately that $f_{vt}(x) = f_{vb}(x)
= 0$.

So the bar framework is infinitesimally rigid\Index{infinitesimally rigid}[2!1
].  Its only motions are the Euclidean motions\Index{Euclidean motion}
of~$\RR^2$.  If we can exhibit a strictly positive stress, we'll have that the
tensegrity is bar equivalent and hence also infinitesimally rigid.  Define
$\mu$ as follows.

\begin{figure}[ht]
\hspace*{\fill}
\begin{tikzpicture}[scale=2]
\draw[cable] 
  (0,1) -- node [pos=0.35,anchor=north east] {$e_1$} 
  (0.5,0) -- node [anchor=north west] {$e_3$} (2,1);
\draw[strut] (0.5,0) -- node [pos=0.7,anchor=west] {$e_2$} (0.5,1);
\draw[yshift=-0.05 cm,snake=brace,segment amplitude=0.1 cm] (0.5,0) -- 
  node [below] {$x$} (0,0);
\draw[yshift=-0.05 cm,snake=brace,segment amplitude=0.1 cm] (2,0) -- 
  node [below] {$2-x$} (0.5,0);
\end{tikzpicture}
\hspace*{\fill}
\capt{The edges connecting to an internal vertex in the rectangle.}
\label{fig:ThreeEdges}
\end{figure}
At any internal point, there are only three edges touching, so we must have
something that looks like \vref{fig:ThreeEdges}.  Using $e_1$, $e_2$ and $e_3$
as shown in the figure, we get
\begin{align}
\mu(e_2) & = \frac{\mu(e_1)}{\sqrt{1 + x^2}} + \frac{\mu(e_3)}{\sqrt{1 +
(2-x)^2}} && \text{from the vertical stresses, and} 
\label{eq:fromVertical} \\
\frac{x}{\sqrt{1 + x^2}}\mu(e_1) & = \frac{(2-x)}{\sqrt{1 + (2-x)^2}}\mu(e_3)
&& \text{from the horizontal ones.}
\label{eq:fromHorizontal}
\end{align}

If we let $m(e)$ be Lebesgue measure on the edges, inherited from Lebesgue
measure on the vertices, and decide that $\mu(e_2) = m(e_2)$, we get 
\begin{align*}
\mu(e_1) & = \sqrt{1 + x^2} \left(m(e_2) - \frac{\mu(e_3)}{\sqrt{1 +
(2-x)^2}}\right) & \text{from \ref{eq:fromVertical}} \\
\intertext{and}
\mu(e_1) & = \frac{(2-x) \sqrt{1 + x^2}}{x \sqrt{1 + (2-x)^2}} \mu(e_3)
& \text{from \ref{eq:fromHorizontal}.}
\end{align*}
Setting these expressions equal results in
\begin{align*}
\frac{(2-x) \mu(e_3)}{x \sqrt{1 + (2-x)^2}} = m(e_2) - \frac{\mu(e_3)}{\sqrt{1
+ (2-x)^2}} 
& \Rightarrow \mu(e_3) = \frac{x\sqrt{1 + (2-x)^2}}{2} m(e_2)
\end{align*}
and directly that
$$\mu(e_1) = \frac{(2-x)\sqrt{1+x^2}}{2} m(e_2).$$
That takes care of all but the endpoints.  There are infinitely many cables 
ending at each corner, but we can integrate to find out how much force is being
exerted on those corners.  

The forces on the upper left-hand corner all come from those $e_1$'s above.  The
$e_1$ forces resolve into horizontal 
$$\frac{x}{\sqrt{1-x^2}} \mu(e_1) = \frac{(2-x)x}{2} m(e_2)$$
and vertical
$$\frac{1}{\sqrt{1-x^2}} \mu(e_1) = \frac{2-x}{2} m(e_2).$$
By symmetry, these are the horizontal and vertical forces at every corner.

Integrating these gives us
$$\int_0^2 \frac{x (2-x)}{2} \, dm$$ 
for a total horizontal force of $\nicefrac{2}{3}$ at each corner and
$$\int_0^2 \frac{2-x}{2} \, dm$$ 
for a total vertical force of $1$ at each corner.  So to balance things
out, there is an atom of size $\nicefrac{4}{3}$ on the two horizontal struts
and an atom of size $2$ on the two bars.

We have a strictly positive stress and the rectangle is infinitesimally rigid.
Alternatively, we could state it thus:
\begin{prop}
Any (non-Euclidean) motion of two parallel line segments does at least one of
the following: \begin{enumerate}
\item Moves two corresponding points on the segments closer together,
\item Shortens one or both of the segments, or
\item Moves a point on one of the segments farther from one of the endpoints of
the opposite segment.
\end{enumerate}
\end{prop}

\section{On a Circle}
\label{section:onACircle}
\begin{figure}[ht]
\hspace*{\fill}
\begin{tikzpicture}
\draw[vcurve] (0,0) circle(2cm);
\draw[cable] (270:2cm) -- (30:2cm) -- (150:2cm) -- (270:2cm)
             (90:2cm) -- (210:2cm) -- (330:2cm) -- (90:2cm);
\draw[strut] (270:2cm) -- (90:2cm) (150:2cm) -- (330:2cm) (30:2cm) -- (210:2cm);
\end{tikzpicture}
\hspace*{1in}
\begin{tikzpicture}
\fill [strut] (0,0) circle (2cm); % Struts
\foreach \x in {0,3,...,361} {
  \draw[cable,black!52!white] (\x:2cm) -- (\x+107:2cm);
}
\end{tikzpicture}
\hspace*{\fill}
\capt{Two examples of the on-a-circle example.}
\label{fig:onACircle}
\end{figure}

We next consider a family of examples.  Each of these has a unit circle for its 
vertex set with antipodal points connected by struts.  Cables are
connected from every point to the two points $h$ units in arc length around the
circle from it.  If the cable connection distance is a rational fraction of the
circumference, the tensegrity consists of infinitely many unconnected finite
subtensegrities.  If the distance is irrational, it is a single connected
tensegrity (see \vref{fig:onACircle}).

We let $m$ be Lebesgue measure on the edges.  Then we create a measure $\mu$ on
the edges by letting $\mu = m$ for all struts and $\mu = \alpha m$ for all
cables.  Let's look at the angles involved (see \vref{fig:AnglesInCircle}).
\begin{figure}[ht]
\hspace*{\fill}
\begin{tikzpicture}
\draw[vcurve] (0,0) circle (2cm);
\draw[strut] (0,-2) -- (0,2);
\draw[cable] (-2,-0.268) -- (0,-2) -- (2,-0.268);
\draw[dotted,black!70!white] (-2,-0.268) -- (0,0) -- (2,-0.268);
\draw (0,0) node [anchor=north west] {$h$};
\draw (0,-2) node [anchor=south west,yshift=0.42cm] {$\frac{\pi - h}{2}$};
\end{tikzpicture}
\hspace*{\fill}
\capt{The angles involved in the on-a-circle example.}
\label{fig:AnglesInCircle}
\end{figure}
If the cables stretch $h$ units around the circle (where $0 < h < \pi$), 
then the cables make an angle of $\frac{\pi - h}{2}$ with the strut meeting
at the same vertex.  

We want equilibrium at every vertex.  The situation is symmetrical across the
strut, so we need only check the sum along the strut.  If we set $\alpha =
\frac{1}{\cos \frac{\pi - h}{2}}$, we get
$$2 \, dm - 2 \alpha \left(\cos \frac{\pi-h}{2}\right) \, dm = 0 \, dm$$
at every vertex.

So we have
$$\mu(YV) = \int_\E YV \, d\mu = \int_{v \in \V} V(v) \cdot (0,0) \, dm = 0.$$
Thus $\mu$ is a strictly positive stress for the tensegrity.  Every one of
these tensegrities is bar equivalent.

It is worthy of note that if $h \to 0$, $\alpha \to \infty$, so we cannot turn
this into the example from \vref{section:TwoExamples} by a simple limiting
process.  \label{example:StarsOfDavid}

\section{On a Sphere}
Let's move that example up into three dimensions.  We start by examining an 
octahedral tensegrity that has three struts and twelve cables (see
\vref{fig:Octasegrity}).
\begin{figure}[ht]
\hspace*{\fill}
\begin{tikzpicture}
\draw[strut] (-1,0,0) -- (1,0,0) (0,-1,0) -- (0,1,0) (0,0,-1) -- (0,0,1);
\draw[cable] 
  (-1,0,0) -- (0,-1,0) -- (0,0,-1) -- (-1,0,0) -- (0,1,0) -- (0,0,1) -- cycle 
  (1,0,0) -- (0,-1,0) -- (0,0,1) -- (1,0,0) -- (0,1,0) -- (0,0,-1) -- cycle;
\end{tikzpicture}
\label{fig:Octasegrity}
\hspace*{\fill}
\capt{The octahedron tensegrity.}
\end{figure}

If we place a compression\Index{compression} of unit size on each strut and a
tension\Index{tension} of $\nicefrac{\sqrt{2}}{2}$ on each cable, then the net
force at vertex is zero.  So this octahedron is bar equivalent.  Now, imagine
this octahedron embedded within a sphere as represented in
\vref{fig:OctInSphere}.
\begin{figure}[ht]
\hspace*{\fill}
\begin{tikzpicture}
\pgfdeclareradialshading{sphere}{\pgfpoint{0.5cm}{0.5cm}}%
 {rgb(0cm)=(0.9,0.9,0.9);
  rgb(0.7cm)=(0.7,0.7,0.7);
  rgb(1cm)=(0.5,0.5,0.5);
  rgb(1.05cm)=(1,1,1)}
\pgfuseshading{sphere}
\draw[strut] (-1,0,0) -- (1,0,0) (0,-1,0) -- (0,1,0) (0,0,-1) -- (0,0,1);
\draw[cable] 
  (-1,0,0) -- (0,-1,0) -- (0,0,-1) -- (-1,0,0) -- (0,1,0) -- (0,0,1) -- cycle 
  (1,0,0) -- (0,-1,0) -- (0,0,1) -- (1,0,0) -- (0,1,0) -- (0,0,-1) -- cycle;
\end{tikzpicture}
\hspace*{\fill}
\capt{Octahedra cover the spherical tensegrity.}
\label{fig:OctInSphere}
\end{figure}
If we connect every point on a sphere with its antipode by a strut, and connect
each point to every point on its associated equator by a cable, we'll have a
three-dimensional tensegrity that is covered by octahedra as above and which
is thus bar equivalent.

\section{On a General Curve}
Both of the previous two examples have had a lot of symmetry on which we could 
depend.  What if we were to consider something less symmetric?  Let's
take as our vertex set any simple closed curve in $\RR^2$ and perform the same
construction we did in \ref{example:StarsOfDavid}.  It will be easier, though
to think of the cable skip ($h$ in that section) not as a length, but rather as
a fraction of the arclength of the vertex curve (see 
\vref{fig:OACGeneral} for two examples with a skip of $\nicefrac{1}{5}$ of the
total arclength).
\begin{figure}[ht]
\hspace*{\fill}
\begin{tikzpicture}
\draw[vcurve] (0,2) arc (90:270:2cm) arc (270:360:1cm) -- (1,1) arc (0:90:1cm);
\draw[strut] (0,2) -- (-0.563,-1.919)
             (-1.081,1.682) -- (0.541,-1.841)
             (-1.819,0.831) -- (1,-0.857)
             (-1.980,-0.285) -- (1,0.286)
             (-1.511,-1.310) -- (0.910,1.415);
\draw[cable] 
  (0,2) -- (-1.819,0.831) -- (-1.511,-1.310) -- (0.541,-1.841) -- 
  (1,0.286) -- cycle 
  (-1.081,1.682) -- (-1.980,-0.285) -- (-0.563,-1.919) -- (1,-0.857) -- 
  (0.910,1.415) -- cycle;
\end{tikzpicture}
\hspace*{1in}
\begin{tikzpicture}[scale=0.667]
\draw[vcurve] (0,-3) arc (-90:0:3cm) -- (3,2) arc (0:90:1cm) -- 
              (0,3) arc (90:270:1cm) arc (90:-90:1cm) arc(90:270:1cm);
\draw[strut] (1.168,-2.763) -- (-0.697,2.717) (2.613,-1.475) -- (-0.389,1.079) 
 (3,0.429) -- (1,0) (2.921,2.389) -- (-0.389,-1.079) (1.2,3) -- (-0.697,-2.717);
\draw[cable] 
  (1.168,-2.763) -- (3,0.429) -- (1.2,3) -- (-0.389,1.079) -- 
  (-0.389,-1.079) -- cycle 
  (2.613,-1.475) -- (2.921,2.389) -- (-0.697,2.717) -- (1,0) -- 
  (-0.697,-2.717) -- cycle;
\end{tikzpicture}
\hspace*{\fill}
\capt{The on-a-circle example with more general curves.}
\label{fig:OACGeneral}
\end{figure}

It seems complicated.  Perhaps we should set our sights a little lower at first
and see what we can learn.

\section{Understanding a Quadrilateral}
Suppose we have some quadrilateral that has struts for the diagonals and
cables for the edges.  In the calculations that follow, we'll use the symbol
$\ev{e_{ij}}$ to mean the unit vector pointing from vertex $i$ toward vertex
$j$ (so $\ev{e_{ij}} = - \ev{e_{ji}}$) and $w_{ij}$ to mean the value of the
stress on the edge connecting vertices $i$ and~$j$.  What do we know?
\begin{figure}[ht]
\hspace*{\fill}
\begin{tikzpicture}
\draw[strut] (0,0) -- node [above left] {$w_{13}$} (5,1) 
             (5,0) -- node [pos=0.7,below left] {$w_{24}$} (-1,3);
\draw[cable] (0,0) node [vertex] {$1$} -- node [below] {$w_{12}$}
             (5,0) node [vertex] {$2$} -- node [right] {$w_{23}$}
             (5,1) node [vertex] {$3$} -- node [above] {$w_{34}$}
             (-1,3) node [vertex] {$4$} -- node [left] {$w_{14}$}
             (0,0);
\end{tikzpicture}
\hspace*{\fill}
\capt{A general quadrilateral tensegrity.}
\end{figure}

For one thing, we know that at any vertex, the weighted sum of the cable
vectors must equal the weighted strut vector.  So the negative of the strut
vector must lie in the cone generated by the cable vectors. That is, the angle
between the cables must be less than $\pi$ on the side where the strut lies and
the strut must lie between them.  So the quadrilateral must be convex.

Putting the above equalities into symbols, we get
\begin{align*}
w_{12} \ev{e_{12}} + w_{14} \ev{e_{14}} + w_{13} \ev{e_{31}} & = 0 &
w_{12} \ev{e_{21}} + w_{23} \ev{e_{23}} + w_{24} \ev{e_{42}} & = 0 \\
w_{23} \ev{e_{32}} + w_{34} \ev{e_{34}} + w_{13} \ev{e_{13}} & = 0 &
w_{14} \ev{e_{41}} + w_{34} \ev{e_{43}} + w_{24} \ev{e_{24}} & = 0
\end{align*}
Let's change the signs so that the subscripts are in ascending order:
\begin{align}
  w_{12} \ev{e_{12}} + w_{14} \ev{e_{14}} - w_{13} \ev{e_{13}} & = 0 
\label{eq:quad1} \\
- w_{12} \ev{e_{12}} + w_{23} \ev{e_{23}} - w_{24} \ev{e_{24}} & = 0
\label{eq:quad2} \\
- w_{23} \ev{e_{23}} + w_{34} \ev{e_{34}} + w_{13} \ev{e_{13}} & = 0
\label{eq:quad3} \\
- w_{14} \ev{e_{14}} - w_{34} \ev{e_{34}} + w_{24} \ev{e_{24}} & = 0
\label{eq:quad4}
\end{align}

Since $\ev{e_{13}}$ must lie between $\ev{e_{12}}$ and $\ev{e_{14}}$, it
splits the plane into two half-planes, each containing one of those vectors.
There are two candidates for ``the vector normal to $\ev{e_{13}}$''. We'll
choose the one that lies in the same half-plane as~$\ev{e_{14}}$.

\ref{eq:quad1} can be thought of as two equations in three variables:
\begin{align}
w_{12} \ev{e_{12}} \cdot \ev{e_{13}} + w_{14} \ev{e_{14}} \cdot \ev{e_{13}} -
w_{13} & = 0 \label{eq:quad5} \\
w_{12} \ev{e_{12}} \cdot \ev{e_{13}}^\perp + w_{14} \ev{e_{14}} \cdot
\ev{e_{13}}^\perp & = 0. \label{eq:quad6}  
\end{align}
We note that the two dot products in \ref{eq:quad5} are both positive, while
our choice of $\ev{e_{13}}^\perp$ makes $\ev{e_{14}} \cdot \ev{e_{31}}^\perp$
positive and $\ev{e_{12}} \cdot \ev{e_{31}}^\perp$ negative.

We can solve this system of equations:
\begin{align*}
w_{12} & = \frac{\ev{e_{14}} \cdot \ev{e_{13}}^\perp}
              {-\ev{e_{12}} \cdot \ev{e_{13}}^\perp} w_{14} \Rightarrow
w_{13} = \left[\left(\ev{e_{12}} \cdot \ev{e_{13}}\right) 
         \frac{\ev{e_{14}} \cdot \ev{e_{13}}^\perp}
              {-\ev{e_{12}} \cdot \ev{e_{13}}^\perp} + 
              \ev{e_{14}} \cdot \ev{e_{13}} \right] w_{14}.
\end{align*}
Letting $\theta_{jik}$ be the angle between $\ev{e_{ij}}$ and $\ev{e_{ik}}$
we can rewrite these in terms of sines and cosines:
\begin{align*}
w_{12} & = \frac{\sin \theta_{314}}{\sin \theta_{213}} w_{14} \\
w_{13} & = \left[\cos \theta_{213} 
           \frac{\sin \theta_{314}}{\sin \theta_{213}} +
              \cos \theta_{314}\right] w_{14}
         = \frac{\sin \theta_{214}}{\sin \theta_{213}} w_{14}
\end{align*}
But since $0 < \theta_{213}, \theta_{314}, \theta_{214} < \pi$ we know that
$w_{12}$, $w_{13}$ and $w_{14}$ will all have the same sign (and can 
certainly be positive).

Similar calculations on \ref{eq:quad3} give us that
\begin{align*}
w_{23} & = \frac{\sin \theta_{134}}{\sin \theta_{132}} w_{34} 
\text{\quad and \quad}
w_{13} = \frac{\sin \theta_{234}}{\sin \theta_{132}} w_{34},
\end{align*}
again with the same sign.  So we have a one-parameter family of solutions:
\begin{align*}
w_{12} & = \frac{\sin \theta_{314}}{\sin \theta_{214}} w_{13} &
w_{14} & = \frac{\sin \theta_{213}}{\sin \theta_{214}} w_{13} &
w_{23} & = \frac{\sin \theta_{134}}{\sin \theta_{234}} w_{13} &
w_{34} & = \frac{\sin \theta_{132}}{\sin \theta_{234}} w_{13}
\end{align*}
There are some points worth noting about this.  First, the ``double arrow''
consisting of a strut and the four cables to which it connects (like the ones
shown in \vref{fig:doubleArrows}) seems a good basic unit to use in the
analysis.  Once $w_{13}$ is set, all of the cable weights follow.  Second, the
weight on a given cable is a function not only of the angle between it and its
strut but also the angle between its neighbor cable and their shared strut.
That will turn out to be significant later.

Now both \ref{eq:quad2} and \vref{eq:quad4} claim to be able to tell us what
$w_{24}$ is.  Suppose we call the \ref{eq:quad2} answer $w_{24}$ and the 
\ref{eq:quad4} answer~$\hat{w}_{24}$.  How do they relate?
\begin{align*}
\hat{w}_{24} \ev{e_{24}} & = w_{14} \ev{e_{14}} + w_{34} \ev{e_{34}} 
& \text{(\ref{eq:quad4})} \\
& = w_{13} \ev{e_{13}} - w_{12} \ev{e_{12}} + w_{23} \ev{e_{23}} - w_{13}
\ev{e_{13}} & \text{(\ref{eq:quad1} \& \ref{eq:quad3})} \\
& = w_{23} \ev{e_{23}} - w_{12} \ev{e_{12}} = w_{24} \ev{e_{24}}
& \text{(\ref{eq:quad2})} 
\end{align*}

That's rather exciting news.\footnote{News that can be established also using
the law of sines on the quadrilateral.}  That compatibility shows that we can
establish a strictly positive stress on any such convex quadrilateral.

Suppose that we have any convex curve of vertices.  If we put struts between
every pair of antipodal (by arc length) points and cables between every pair of
points whose distance is $\nicefrac{1}{4}$ of the length of the curve, the
resulting tensegrity is covered by convex quadrilaterals and so it has a
strictly positive stress and is bar equivalent.

But there are two hypotheses that seem like they could be removed.  First, 
surely \nicefrac{1}{4} isn't the only fraction of the curve length that works.
Secondly, is it really necessary to have a convex curve?

\section{Other Fractions}
Let's start by trying other fractions.  We started with \nicefrac{1}{4} because
the quadrilateral is pretty straightforward.  Of course, if we use the fraction
$\nicefrac{1}{2n}$ for any $n > 1$, we get a closed $2n$-gon of cables with
struts connecting opposite points.  As before, we certainly need each strut to
fall in the cone generated by its cables.  Also as before, we can generate
families of positive solutions for each ``double arrow''.  We need only
establish that those solutions can be compatible with each other.
\vref{fig:Hexasegrity} shows a curve with one of the hexagonal tensegrities
coming from a skip of $\nicefrac{1}{6}$ of the curve length.
\begin{figure}[ht]
\hspace*{\fill}
\begin{tikzpicture}[scale=0.7]
\draw[vcurve] (0,-1) arc (270:135:1cm) -- (2.293,3.707) arc (135:0:1cm) --
              (4,0) arc (0:-90:1cm) -- cycle;
\draw[strut] (3.891,-0.454) -- (0.571,1.985) (1.345,-1) -- (2.55,3.893) 
 (-0.987,-0.161) -- (4,2.283);
\draw[cable] 
  (3.891,-0.454) node [vertex] {1} -- (1.345,-1) node [vertex] {2} --
  (-0.987,-0.161) node [vertex] {3} -- (0.571,1.985) node [vertex] {4} --
  (2.55,3.893) node [vertex] {5} -- (4,2.283) node [vertex] {6} -- cycle;
\end{tikzpicture}
\hspace*{\fill}
\capt{A more general hexagon.}
\label{fig:Hexasegrity}
\end{figure}

Now suppose that we have two ``adjacent'' struts (shown in
\vref{fig:doubleArrows}) along with their associated cables.  
\begin{figure}[ht]
\hspace*{\fill}
\begin{tikzpicture}[scale=0.7]
\draw[black!15!white,line width=2mm]
  (1.345,-1) -- (-0.987,-0.161) (2.55,3.893) -- (4,2.283);
\draw[strut] (-0.987,-0.161) -- (4,2.283);
\draw[cable] 
  (1.345,-1) node [vertex] {2} -- (-0.987,-0.161) node [vertex] {3} --
  (0.571,1.985) node [vertex] {4}
  (2.55,3.893) node [vertex] {5} -- (4,2.283) node [vertex] {6} --
  (3.891,-0.454) node [vertex] {1};
\end{tikzpicture}
\hspace*{1in}
\begin{tikzpicture}[scale=0.7]
\draw[black!15!white,line width=2mm]
  (1.345,-1) -- (-0.987,-0.161) (2.55,3.893) -- (4,2.283);
\draw[strut] (1.345,-1) -- (2.55,3.893);
\draw[cable] 
  (3.891,-0.454) node [vertex] {1} -- (1.345,-1) node [vertex] {2} --
  (-0.987,-0.161) node [vertex] {3}
  (0.571,1.985) node [vertex] {4} -- (2.55,3.893) node [vertex] {5} --
  (4,2.283) node [vertex] {6};
\end{tikzpicture}
\hspace*{\fill}
\capt[Two struts with their associated cables.]{Two of the struts from the
hexagonal tensegrity above, along with their associated cables.  The two cables
they share are shaded.}
\label{fig:doubleArrows}
\end{figure}
These two struts will share one cable on each end.  Certainly we can scale
the answers for the two ``double arrows'' so that the weights match on one of
the shared cables, but it's not at all clear that we can make both ends match.
And it turns out that in general it doesn't work.  
\vref{fig:HexWithVect} shows a set of vectors that increase all strut lengths
and decrease all cable lengths for this particular example, and since none of
the angles are straight, there will be an open set of neighbor hexagons that
also have a motion.  So this tensegrity is not bar equivalent.
\begin{figure}[ht]
\hspace*{\fill}
\begin{tikzpicture}
\draw[vcurve] (0,-1) arc (270:135:1cm) -- (2.293,3.707) arc (135:0:1cm) --
              (4,0) arc (0:-90:1cm) -- cycle;
\draw[strut] (3.891,-0.454) -- (0.571,1.985) (1.345,-1) -- (2.55,3.893) 
 (-0.987,-0.161) -- (4,2.283);
\draw[cable] 
  (3.891,-0.454) node [vertex] {1} -- (1.345,-1) node [vertex] {2} --
  (-0.987,-0.161) node [vertex] {3} -- (0.571,1.985) node [vertex] {4} --
  (2.55,3.893) node [vertex] {5} -- (4,2.283) node [vertex] {6} -- cycle;
% Vector: -1.2454,0.4376  -0.4984,-1.2140  0.3672,0  -1.3138,0.7544  0,-1.1324
%          0.5677,0
\draw[->] (3.891,-0.454) -- (2.6456,-0.0164);
\draw[->] (1.345,-1) -- (0.8466,-2.2140);
\draw[->] (-0.987,-0.161) -- (-0.6198,-0.161);
\draw[->] (0.571,1.985) -- (-0.7428,2.7394);
\draw[->] (2.55,3.893) -- (2.55,2.7606);
\draw[->] (4,2.283) -- (4.5677,2.283);
\end{tikzpicture}
\hspace*{\fill}
\capt{The hexagon with a motion.}
\label{fig:HexWithVect}
\end{figure}

What does work?  Well, certainly any \emph{regular} figure (where the
construction at each vertex is identical to that at all of the others).
\vref{fig:TwoRegulars} shows two regular figures that could arise.
\begin{figure}[ht]
\hspace*{\fill}
\begin{tikzpicture}[scale=0.75]
\draw[vcurve] (0,0) circle(2cm);
\draw[strut] (2,0) -- (-2,0) 
 (1.618,1.176) -- (-1.618,-1.176) 
 (0.618,1.902) -- (-0.618,-1.902) 
 (-0.618,1.902) -- (0.618,-1.902) 
 (-1.618,1.176) -- (1.618,-1.176);
% 0,5,1,6,2,7,3,8,4,9
\draw[cable] (2,0) -- (-1.618,1.176) -- (0.618,-1.902) -- (0.618,1.902) -- (-1.618,-1.176) -- cycle;
% 0,2,4,6,8
\draw[cable] (1.618,1.176) -- (-2,0) -- (1.618,-1.176) -- (-0.618,1.902) -- (-0.618,-1.902) -- cycle;
% 1,3,5,7,9
\end{tikzpicture}
\hspace*{1cm}
\begin{tikzpicture}[scale=1.5]
\draw[vcurve] (0,0) -- (2,0) -- (2,2) -- (0,2) -- cycle;
\draw[strut] (0,1.5) -- (2,0.5) (1.5,2) -- (0.5,0);
\draw[cable] (0,1.5) -- (1.5,2) -- (2,0.5) -- (0.5,0) --
  cycle;
\end{tikzpicture}
\hspace*{\fill}
\capt{Two regular figures that arise with this construction.}
\label{fig:TwoRegulars}
\end{figure}

But other figures can provide bar equivalence.  Suppose, for example, that we
have the non-regular hexagon shown in \vref{fig:NonRegHex}.
\begin{figure}[ht]
\hspace*{\fill}
\begin{tikzpicture}
\draw[strut] (2,0) -- (-2,0) (1,1) -- (-1,-1) (-1,1) -- (1,-1);
\draw[cable] (2,0) -- (1,1) -- (-1,1) -- (-2,0) -- (-1,-1) -- (1,-1) -- cycle;
\draw (2,0) node[anchor=west,scale=0.5] {$(2,0)$}
      (1,1) node[anchor=south west,scale=0.5] {$(1,1)$}
     (-1,1) node[anchor=south east,scale=0.5] {$(-1,1)$}
     (-2,0) node[anchor=east,scale=0.5] {$(-2,0)$}
    (-1,-1) node[anchor=north east,scale=0.5] {$(-1,-1)$}
     (1,-1) node[anchor=north west,scale=0.5] {$(1,-1)$};
\end{tikzpicture}
\hspace*{\fill}
\capt{A bar equivalent tensegrity that is a nonregular hexagon.}
\label{fig:NonRegHex}
\end{figure}
It turns out that a stress of 2 units on each of the cables and 1 unit each on
the struts will put this in equilibrium, so it \emph{is} bar equivalent, and 
this suggests a family of bar equivalent hexagons.
\vref{fig:AnotherNonRegHex} gives another, with the weights shown on the edges:
\begin{figure}[ht]
\hspace*{\fill}
\begin{tikzpicture}[scale=2]
\draw[strut] 
  (0.6,0.6) -- node [pos=0.25,weight] {$1$} (-0.6,-0.6) 
  (1,0)     -- node [pos=0.25,weight] {$3$} (-1,0) 
  (0,1)     -- node [pos=0.25,weight] {$3$} (0,-1);
\draw[cable] (0.6,0.6) -- 
  node [weight] {$6$} (1,0) -- 
  node [weight] {$3.6$} (0,-1) -- 
  node [weight] {$6$} (-0.6,-0.6) -- 
  node [weight] {$6$} (-1,0) -- 
  node [weight] {$3.6$} (0,1) -- 
  node [weight] {$6$} (0.6,0.6);
\draw (1,0) node[anchor=west,scale=0.5] {$(1,0)$}
      (0.6,0.6) node[anchor=south west,scale=0.5] {$(0.6,0.6)$}
      (-1,0) node[anchor=east,scale=0.5] {$(-1,0)$}
      (0,1) node[anchor=south east,scale=0.5] {$(0,1)$}
      (-0.6,-0.6) node[anchor=north east,scale=0.5] {$(-0.6,-0.6)$}
      (0,-1) node[anchor=north west,scale=0.5] {$(0,-1)$};
\end{tikzpicture}
\hspace*{\fill}
\capt{Another in the family of nonregular, bar equivalent hexagons.}
\label{fig:AnotherNonRegHex}
\end{figure}
Now these are simply example hexagons.  We have not addressed how much stranger
they can get, whether they can actually arise while building tensegrities from
curves in that fashion, nor whether, if they should arise, there must also be
non-bar-equivalent subtensegrities in the same construction.  This seems a
field worth more exploration in the future.  For now, let's see if we can
settle something about non-convexity.
\label{section:exampleToDo}

\section{Non-convexity}
\label{section:nonConv}
\vref{fig:TwoNonConvex} shows two non-convex vertex curves with
antipodal struts and fixed-skip cables.  The skip in the first case is
\nicefrac{1}{4} of the curve length, and in the second it is \nicefrac{1}{6}.
\begin{figure}[ht]
\hspace*{\fill}
\begin{tikzpicture}[scale=0.8]
\foreach \x in {0,90,180,270} {
  \draw[vcurve,rotate=\x] (-4,2) arc (180:90:2cm) arc (90:0:1cm) arc
    (-180:0:1cm) arc (180:90:1cm);
}
\draw[strut] (-3,-1) -- (3,1) 
% (-2.076,-0.383) -- (2.076,0.383) 
% (-2.293,0.707) -- (2.293,-0.707) 
% (-3.383,1.076) -- (3.383,-1.076) 
% (-4,2) -- (4,-2) 
% (-3.663,3.111) -- (3.663,-3.111) 
% (-2.765,3.848) -- (2.765,-3.848) 
% (-1.617,3.924) -- (1.617,-3.924) 
 (-1,3) -- (1,-3);
% (-0.383,2.076) -- (0.383,-2.076) 
% (0.707,2.293) -- (-0.707,-2.293) 
% (1.076,3.383) -- (-1.076,-3.383) 
% (2,4) -- (-2,-4) 
% (3.111,3.663) -- (-3.111,-3.663) 
% (3.848,2.765) -- (-3.848,-2.765) 
% (3.924,1.617) -- (-3.924,-1.617);
\draw[cable] 
(-3,-1) -- (-1,3) -- (3,1) -- (1,-3) -- cycle;
% (-2.076,-0.383) -- (-0.383,2.076) -- (2.076,0.383) -- (0.383,-2.076) -- cycle
% (-2.293,0.707) -- (0.707,2.293) -- (2.293,-0.707) -- (-0.707,-2.293) -- cycle
% (-3.383,1.076) -- (1.076,3.383) -- (3.383,-1.076) -- (-1.076,-3.383) -- cycle
% (-4,2) -- (2,4) -- (4,-2) -- (-2,-4) -- cycle
% (-3.663,3.111) -- (3.111,3.663) -- (3.663,-3.111) -- (-3.111,-3.663) -- cycle
% (-2.765,3.848) -- (3.848,2.765) -- (2.765,-3.848) -- (-3.848,-2.765) -- cycle
% (-1.617,3.924) -- (3.924,1.617) -- (1.617,-3.924) -- (-3.924,-1.617) -- cycle;
\end{tikzpicture}
\hspace*{1cm}
\begin{tikzpicture}
\foreach \x in {0,60,120,180,240,300} {
  \draw[vcurve,rotate=\x] (-1,1.732) -- (0,3.464) -- (1,1.732);
\draw[strut]
 (-0.75,2.165) -- (0.75,-2.165) 
% (0.25,3.031) -- (-0.25,-3.031) 
 (1.5,1.732) -- (-1.5,-1.732) 
% (2.75,1.299) -- (-2.75,-1.299) 
 (2.25,-0.433) -- (-2.25,0.433);
% (2.5,-1.732) -- (-2.5,1.732);
% 0,6,1,7,2,8,3,9,4,10,5,11
\draw[cable] (-0.75,2.165) -- (1.5,1.732) -- (2.25,-0.433) -- (0.75,-2.165) -- (-1.5,-1.732) -- (-2.25,0.433) -- cycle;
% 0,2,4,6,8,10
% \draw[cable] (0.25,3.031) -- (2.75,1.299) -- (2.5,-1.732) -- (-0.25,-3.031) -- (-2.75,-1.299) -- (-2.5,1.732) -- cycle;
% 1,3,5,7,9,11
% \draw[strut] (0.5,2.598) -- (-0.5,-2.598) 
%  (-0.5,2.598) -- (0.5,-2.598) 
%  (2.5,0.866) -- (-2.5,-0.866) 
%  (2,1.732) -- (-2,-1.732) 
%  (2,-1.732) -- (-2,1.732) 
%  (2.5,-0.866) -- (-2.5,0.866);
% 0,6,1,7,2,8,3,9,4,10,5,11
% \draw[cable] (0.5,2.598) -- (2.5,0.866) -- (2,-1.732) -- (-0.5,-2.598) -- (-2.5,-0.866) -- (-2,1.732) -- cycle;
% 0,2,4,6,8,10
% \draw[cable] (-0.5,2.598) -- (2,1.732) -- (2.5,-0.866) -- (0.5,-2.598) -- (-2,-1.732) -- (-2.5,0.866) -- cycle;
% 1,3,5,7,9,11
\draw[dashed] (-0.75,2.165) -- (2.25,-0.433) -- (-1.5,-1.732) -- cycle
                       (1.5,1.732) -- (0.75,-2.165) -- (-2.25,0.433) -- cycle;
}
\end{tikzpicture}
\hspace*{\fill}
\capt[Two nonconvex examples.]{Here are two examples of bar-equivalent
tensegrities that arise from nonconvex curves.  One subtensegrity is shown for
each.  The dashed lines are explained in the text.}
\label{fig:TwoNonConvex}
\end{figure}
The resulting tensegrity is, in each case, covered with regular
subtensegrities, all of which are bar equivalent.  So each tensegrity is itself
bar equivalent (though neither is infinitesimally rigid, since a vector field
with the same rotational symmetry will move the subtensegrities with respect to
each other).

What's going on here?  As we have seen above, what we need is for the two
cables meeting at each vertex of each subtensegrity to form a cone that
contains the strut.  That will happen if, for a given vertex, the strut
intersects the line that connects the two ``distant'' cable ends (the dashed
lines in the \ref{fig:TwoNonConvex} and \ref{fig:Relationship}).   Or to put it
differently, if the curve has length $L$ and the cables connect points $\alpha
L$ apart, then the tensegrity will be bar equivalent if all segments connecting
points $2 \alpha L$ apart lie entirely within the curve (see
\vref{fig:Relationship}).
\begin{figure}[ht]
\hspace*{\fill}
\begin{tikzpicture}
\draw[vcurve] 
  (-1.732,6) arc (90:30:1cm) arc (210:330:1cm) arc (150:90:1cm);
\draw[strut] (0,5) -- (0,3);
\draw[cable] (-3,4.5) -- (0,5) -- (3,4.5);
\draw[dashed] (-3,4.5) -- (3,4.5);
\draw (1.5,3.5) node {\parbox{1in}{\footnotesize Possibly \\ bar equivalent}};
\end{tikzpicture}
\quad
\begin{tikzpicture}
\draw[vcurve] 
  (-1.732,5) arc (90:30:1cm) arc (210:330:1cm) arc (150:90:1cm);
\draw[strut] (0,4) -- (0,2);
\draw[cable] (-3,4.5) -- (0,4) -- (3,4.5);
\draw[dashed] (-3,4.5) -- (3,4.5);
\draw (2.0,2.5) node {\footnotesize Not bar equivalent};
\end{tikzpicture}
\hspace*{\fill}
\capt[``Double skip'' and bar equivalence.]{How the relationship between the
vertex curve and the ``double skip'' tells on the bar equivalence of the
subtensegrities.}
\label{fig:Relationship}
\end{figure}

For a given non-convex curve, then, there is a minimum distance below which the
cables must not go.  This distance can be found in the following fashion.  For
each vertex $v$, have two points start at the antipode and move, one in each
direction at equal speeds around the curve toward~$v$.

If, at some point, the segment between them no longer lies entirely within the
curve, stop.  Call that distance~$d_v$.  Otherwise, $d_v=0$.  Now, to get a 
bar equivalent tensegrity, the cable skip distance must be strictly greater
than $\max_{v \in \V} d_v$.  

\begin{figure}[ht]
\hspace*{\fill}
\begin{tikzpicture}[scale=0.5]
\draw[vcurve] (1,-1) arc (0:180:1cm) arc (0:-180:1cm) --
  (-3,1) arc (0:180:3cm) arc (180:360:1cm) arc (180:0:1cm) --
  (-5,-1) arc (180:360:3cm);
\draw[strut] (0,0) -- (-8,0);
\end{tikzpicture}
\hspace*{\fill}
\capt[Vertex curve and strut that cannot be balanced.]{A vertex curve and a
strut that cannot balanced by fixed-skip cables.}
\label{fig:cant}
\end{figure}
This is not always possible.  \vref{fig:cant} shows a vertex curve with
one of its struts.  If cables are attached with a skip \emph{anything} less
than \nicefrac{1}{2} of the overall length, the two associated cables will both
end on the same side of the strut, so the strut will not lie in the cone they
create.

\section{Affine Transformations}
\label{section:CrossedMore}
\subsection{Deeper into the Crossed Square}
In \vref{section:WorkingRWex}, we found that the crossed
square\Index[|(]{crossed square}[ ] (also shown in
\ref{fig:CSNum}\subref{subfig:csen}) is bar equivalent.
\begin{figure}[ht]
\hspace*{\fill}
\subfloat[The crossed square with edges numbered.]{
\begin{tikzpicture}[scale=1.25]
\useasboundingbox (-1,-1) rectangle (3,3);
\draw[cable] 
  (0,0) node [scale=0.75,below left] {$(0,0)$} -- node [weight] {1}
  (2,0) node [scale=0.75,below right] {$(1,0)$} -- node [weight] {4}
  (2,2) node [scale=0.75,above right] {$(1,1)$} -- node [weight] {6}
  (0,2) node [scale=0.75,above left] {$(0,1)$} -- node [weight] {3}
  (0,0) -- (2,0);
\draw[strut] 
  (0,0) node [vertex] {1} -- node [weight,pos=0.25] {2}
  (2,2) node [vertex] {3}
  (0,2) node [vertex] {4} -- node [weight,pos=0.75] {5}
  (2,0) node [vertex] {2};
\end{tikzpicture}
\label{subfig:csen}
}
\subfloat[The crossed square after transformation.]{
\begin{tikzpicture}[scale=1.25]
\begin{scope}[cm={1,-0.25,-0.25,1.25,(5cm,0cm)}]
\useasboundingbox (-1,-0.1) rectangle (3,2.1);
\draw[cable] 
  (0,0) node[scale=0.75,below left] {$(0,0)$} -- node [weight] {1}
  (2,0) node[scale=0.75,below right] {$(a,b)$} -- node [weight] {4}
  (2,2) node[scale=0.75,above right] {$(a+b,b+c)$} -- node [weight] {6}
  (0,2) node[scale=0.75,above left] {$(b,c)$} -- node [weight] {3}
  (0,0) -- (2,0);
\draw[strut] (0,0) node [vertex] {1} -- node [weight,pos=0.25] {2}
  (2,2) node [vertex] {3}
  (0,2) node [vertex] {4} -- node [weight,pos=0.75] {5}
  (2,0) node [vertex] {2};
\end{scope}
\end{tikzpicture}
\label{subfig:csat}
}
\hspace*{\fill}
\capt[The crossed-square, original and transformed.]{The crossed-square
example, both the original and after being transformed by the matrix
$L^\top L = \left[\begin{smallmatrix}a & b \\ b & c\end{smallmatrix}\right]$.}
\label{fig:CSNum}
\end{figure}
Let's find how its set of loads, $Y(\X)$, looks.  

Suppose $V = (x_1,y_1,x_2,y_2,x_3,y_3,x_4,y_4) \in \Vf(\V)$.  Then, letting the
edges be $e_1$ through $e_6$ (as numbered in the figure), we have
\begin{align*}
YV(e_1) & = -(x_2 - x_1) & YV(e_4) & = -(y_3 - y_2) \\
YV(e_2) & = x_3 - x_1 + y_3 - y_1 & YV(e_5) & = x_2 - x_4 - (y_2 - y_4) \\
YV(e_3) & = -(y_4 - y_1) & YV(e_6) & = -(x_3 - x_4)
\end{align*}
But then 
\begin{align*}
\sum_{i=1}^6 YV(e_i) & = (-x_2 + x_1) + (x_3 - x_1 + y_3 - y_1) + (-y_4 + y_1)\\
& \hspace*{2em} + (-y_3 + y_2) + (x_2 - x_4 - y_2 + y_4) + (-x_3 + x_4) = 0\\
\end{align*}

So $\im Y$ is contained in the hyperplane $\sum_{i=1}^6 YV(e_i) = 0$.  On the
other hand, given any six real numbers $r_1,\dotsc,r_6$ with $r_1 + \dotsb +
r_6 = 0$, the vector field 
$$(r_1,r_3,0,r_4,r_1+r_2+r_3,0,-r_4-r_5,0)$$ 
is one of (infinitely) many for which $YV(e_i) = r_i, i \in \{1, \dotsc, 6\}$.  

So $\im Y$ is the hyperplane $\sum_{i=1}^6 YV(e_i) = 0$ and from that
equation, it is clear to see that $\im Y \cap C(\E)^+ = \{\zero\}$ (whenever
one of the $YV(e_i)$ is positive, another must be negative to offset it).  
We note that the hyperplane is 5-dimensional, a reasonable thing since
$\Vf(\V)$ is 8-dimensional and there are three degrees of freedom given to
rigid motions of the plane (one rotational and two translational) that, of
course, lie in the kernel of~$Y$.
% \cite{MR1300410}*{p.\ 409}. for the dimension of SE(2)

Now suppose that we apply the linear transformation whose matrix is $L$ to our
example.  How will $\im Y$ change?  Well, $Y$ is defined by taking the dot
product of vectors.  We note that 
\begin{equation}
\begin{split}
(LV(v_1) - LV(v_2)) & \cdot (Lp(v_1) - Lp(v_2)) \\
 & = L(V(v_1) - V(v_2)) \cdot L (p(v_1) - p(v_2)) \\
 & = (V(v_1) - V(v_2))^\top L^\top L (p(v_1) - p(v_2)) \\
 & = (V(v_1) - V(v_2)) \cdot L^\top L (p(v_1) - p(v_2)).
\end{split}
\label{eq:dotProd}
\end{equation}
So we can accomplish the same purpose by transforming $p$ by $L^\top L$ and 
leaving $V$ alone.  A simple calculation shows that $L^\top L$ is always
symmetric (with nonnegative entries on the diagonal), so we'll write it as 
$$L^\top L = \begin{bmatrix} a & b \\ b & c \end{bmatrix}$$
and see what effect it has (one possibility is shown in
\ref{fig:CSNum}\subref{subfig:csat}).

Now we get (using $\hat{Y}$ for the rigidity operator after transformation)
\begin{align*}
\hat{Y}V(e_1) & = - a (x_2 - x_1) - b (y_2 - y_1) \\
\hat{Y}V(e_2) & = (a + b)(x_3 - x_1) + (b + c)(y_3 - y_1) \\
\hat{Y}V(e_3) & = - b (x_4 - x_1) - c (y_4 - y_1) \\
\hat{Y}V(e_4) & = - b (x_3 - x_2) - c (y_3 - y_2) \\
\hat{Y}V(e_5) & = (a - b)(x_2 - x_4) + (b - c)(y_2 - y_4) \\
\hat{Y}V(e_6) & = - a (x_3 - x_4) - b (y_3 - y_4).
\end{align*}
The arithmetic is only touch more complicated, but the result is the same,
$\im Y$ is the hyperplane $\sum_{i=1}^6 YV(e_i) = 0$.  So the crossed square
remains bar equivalent under linear transformation.  

This is far from obvious.  Linear transformation does \emph{not}, in general,
preserve dot product.  So a vector field that is in the kernel of $Y$ before
transformation is not necessarily in the kernel of $\hat{Y}$ afterward (see
\vref{fig:kernYtrans} for an example).  This bears more exploration.
\begin{figure}[ht]
\hspace*{\fill}
\begin{tikzpicture}[scale=1.25]
\begin{scope}
  \draw[cable] (0,0) -- (2,0) -- (2,2) -- (0,2) -- cycle;
  \draw[strut] (0,0) node [vertex] {1} -- (2,2) node [vertex] {3}
               (0,2) node [vertex] {4} -- (2,0) node [vertex] {2};
  \draw[->] (0,0) -- +(-0.3, 0.0);
  \draw[->] (2,0) -- +(-0.3,-0.6);
  \draw[->] (2,2) -- +( 0.3,-0.6);
  \draw[->] (0,2) -- +( 0.3, 0.0);
  \fill[blue] (0,1) circle (1pt);
\end{scope}
\begin{scope}[cm={0.75,0.125,0.25,0.875,(4cm,0cm)}]
  \draw[cable] (0,0) -- (2,0) -- (2,2) -- (0,2) -- cycle;
  \draw[strut] (0,0) node [vertex] {1} -- (2,2) node [vertex] {3}
               (0,2) node [vertex] {4} -- (2,0) node [vertex] {2};
  \draw[->] (0,0) -- +(-0.3, 0.0);
  \draw[->] (2,0) -- +(-0.3,-0.6);
  \draw[->] (2,2) -- +( 0.3,-0.6);
  \draw[->] (0,2) -- +( 0.3, 0.0);
\end{scope}
\end{tikzpicture}
\hspace*{\fill}
\capt[A motion that is in $\ker Y$ only before transformation.]{The
crossed-square, original and transformed, with a motion (rotation around the
marked point) that is in $\ker Y$ before transformation by $L =
\left[\begin{smallmatrix} \nicefrac{3}{4} & \nicefrac{1}{4} \\ \nicefrac{1}{8}
& \nicefrac{7}{8} \end{smallmatrix}\right]$ and not in $\ker \hat{Y}$
afterward.}
\label{fig:kernYtrans}
\end{figure}
\Index[|)]{crossed square}[ ]

\subsection{Triangles again}
Let's look at one more example and see if we can figure out what is going on.
In \vref{fig:triangle} we find a tensegrity that looks like the one in
\vref{fig:hasAmotion}, but we are going to treat it differently.
\begin{figure}[ht]
\hspace*{\fill}
\begin{tikzpicture}
\foreach \x/\y in {0/0cm,1/6cm} {
  \begin{scope}[xshift=\y,cm={1,0,\x,1,(0,0)}]
    \draw[bar] (2,0) -- (0,2);
    \draw[cable] (2,0) node [vertex] {2} -- (0,0) node [vertex] {1} -- 
      (0,2) node [vertex] {3};
    \draw[->] (0,2) -- (-0.5,2);
  \end{scope}
}
\draw[->] (1.5,1.5) .. controls +(45:1cm) and +(135:1cm) ..
  node [below] {$\left[\begin{smallmatrix} 1 & 1 \\ 0 & 1
  \end{smallmatrix}\right]$} (6.5,1.5);
\end{tikzpicture}
\hspace*{\fill}
\capt[$\X$-bar equivalent before transformation but not after.]{A tensegrity
that is $\X$-bar equivalent before transformation by $\left[\begin{smallmatrix}
1 & 1 \\ 0 & 1 \end{smallmatrix}\right]$ but not afterward.  The motion shown,
moving only vertex 3, results in a load that has negative components for the
lefthand tensegrity, but is semipositive for the right.}
\label{fig:triangle}
\end{figure}
This time we are going to define $\X$ to be those variations that do not change
the length of the cables (to first order).  That means
$$\X = \{(x_1,y_1,x_2,y_2,x_3,y_3) \in \Vf(\V) : x_1 = x_2, y_1 = y_3\}$$

The bar functions as a strut and a cable, so $C(\E)$ is a 4-dimensional 
space.  Taking the edges in the order ``1-2 cable, 1-3 cable, 2-3 cable, 2-3
strut'', the elements of $Y(\X)$ turn out to all be of the form
$$(0,0,-(y_3-y_2)-(x_2-x_3),(y_3-y_2)+(x_2-x_3))$$
so $Y(\X) \cap C(\E)^+ = \{ \zero \}$ and the tensegrity is $\X$-bar equivalent.

On the other hand, if we transform the tensegrity with $L =
\left[\begin{smallmatrix} 1 & 1 \\ 0 & 1 \end{smallmatrix}\right]$, then we get
$$L\X = \{(x_1,y_1,x_2,y_2,x_3,y_3) \in \Vf(\V) : x_1 + y_1 = x_2 + y_2, y_1
= y_3\}$$
and so the elements of $\hat{Y}(L\X)$ are 
$$(y_2 - y_1,x_1 - x_3, y_1 - y_2, y_2 - y_1).$$
Now $\hat{Y}(L\X) \cap C(\E)^+$ contains functions that are zero except on the
``newly diagonal'' cable, the edge running from vertex 1 to vertex 3.  On that
edge, they can be strictly positive.  

Since there is a semipositive stress, namely $(1,0,2,1)$, our transformed 
tensegrity is partially $L\X$-bar equivalent, but it is not fully so.

\subsection{Conclusion}
So what is the difference between the two preceeding examples?  $\X$.  

To be a little clearer, let's think back to \vref{eq:dotProd} in which we 
showed that for two vectors $v$ and $p$, $Lv \cdot Lp = v \cdot L^\top L p$.
We could have gone the other direction and we would have gotten
$$(LV(v_1) - LV(v_2)) \cdot (Lp(v_1) - Lp(v_2)) = L^\top L (V(v_1) - V(v_2))
\cdot (p(v_1) - p(v_2)).$$
That equation gives the answer to the question.  In our first example, $L^\top
L$ was an automorphism of $\X$ (which, after all, was all of $\Vf(\V)$).  In the
second it was not.  So that suggests a proposition.
\begin{prop}
\label{prop:affineOne}
Let $G(p)$ be an $\X$-bar equivalent tensegrity in $\RR^n$, $L$ a linear
transformation on $\RR^n$ and $t \in \RR^n$.  Define $A$ by $Ax = Lx + t$ for
all $x \in \RR^n$.  Then $G(Ap)$ is $A\X$-bar equivalent if $L^\top L \X
\subset \X$.
\end{prop}
\begin{proof}
Let $AV \in A\X$.  Then, for any $\{v_1,v_2\} \in \E$ we have 
\begin{align*}
\hat{Y}(AV)(\{v_1,v_2\}) & = \pm (LV(v_1)+t - LV(v_2)-t) \cdot (Lp(v_1)+t -
Lp(v_2)-t) \\
& = \pm (LV(v_1) - LV(v_2)) \cdot L(p(v_1) - p(v_2)) \\
& = \pm L^\top (LV(v_1) - LV(v_2)) \cdot (p(v_1) - p(v_2)) \\
& = \pm (L^\top LV(v_1) - L^\top L V(v_2)) \cdot (p(v_1) - p(v_2)) \\
& = Y(L^\top L V)(\{v_1,v_2\})
\end{align*}
(where the choice of sign depends on the type of edge).  But since
$L^\top L \X \subset \X$, we have $\hat{Y}(AV) \in Y(\X)$.  So $\hat{Y}(A\X)
\subset Y(\X)$, which means that 
$$\hat{Y}(A\X) \cap C(\E)^+ \subset Y(\X) \cap C(\E)^+ = \{ \zero \}$$
and $G(Ap)$ is $A\X$-bar equivalent.

Note that the proposition and proof still stand if ``$\X$-bar equivalent'' is
replaced by ``partially $\X$-bar equivalent'', ``$C(\E)^+$'' is replaced with
``$\INT C(\E)^+$'' and ``$\{\zero\}$'' by ``$\varnothing$''.
\end{proof}

Remembering that ``bar equivalent'' means ``$\Vf(\V)$-bar equivalent'', we 
immediately get a corollary.
\begin{corollary}
\label{corollary:affine}
If $G(p)$ is a bar equivalent tensegrity in $\RR^n$ and $A$ is an affine
transformation on $\RR^n$, then $G(Ap)$ is bar equivalent.
\end{corollary}
\begin{proof}
Note that for any linear transformation $L$, $L^\top L\Vf(\V) \subset \Vf(\V)$,
and apply Proposition \ref{prop:affineOne}.
\end{proof}

In the case of an invertible transformation (for which we only need that $L$ is
invertible, since $A^{-1}x = L^{-1}(x-t)$), we can say something more.
\begin{prop}
Let $G(p)$ be an $\X$-bar equivalent tensegrity in $\RR^n$, $L$ an invertible
linear transformation and $t \in \RR^n$.  Define $A$ by $Ax = Lx + t$.  Then
$G(Ap)$ is $(L^\top)^{-1}\X$-bar equivalent.
\end{prop}
\begin{proof}
Let $(L^\top)^{-1}V \in (L^\top)^{-1}\X$.  Then, for any $\{v_1,v_2\} \in \E$,
we have
\begin{align*}
\hat{Y}((L^\top)^{-1}V)(\{v_1,v_2\}) 
& = ((L^\top)^{-1}V(v_1) - (L^\top)^{-1}V(v_2)) \cdot (Ap(v_1) - Ap(v_2)) \\
& = (L^\top)^{-1}(V(v_1) - V(v_2)) \cdot (Lp(v_1) + t - Lp(v_2) - t) \\
& = (L^\top)^{-1}(V(v_1) - V(v_2)) \cdot L(p(v_1) - p(v_2)) \\
& = L^\top(L^\top)^{-1}(V(v_1) - V(v_2)) \cdot (p(v_1) - p(v_2)) \\
& = (V(v_1) - V(v_2)) \cdot (p(v_1) - p(v_2)) \\
& = YV(\{v_1,v_2\})
\end{align*}
So $\hat{Y}((L^\top)^{-1}V) = YV \in Y(\X)$ and we already know that $Y(\X)
\cap C(\E)^+ = \{\zero\}$.

Once again, the proposition and proof are true, \textit{mutatis mutandis} for
partially $(L^\top)^{-1}\X$-bar equivalent.
\end{proof}

\subsection{More Circles}
So what does this say about the on-a-circle example of
\ref{section:onACircle}?  Since $\X = \Vf(\V)$, Corollary
\ref{corollary:affine} says that any affine transformation of the on-a-circle
example is still bar equivalent.  Note that this is different from what
one would get by first transforming the circle of vertices and then applying
the on-a-circle construction (see figure \vref{fig:OACaffine}).
\begin{figure}[ht]
\hspace*{\fill}
\begin{tikzpicture}[scale=0.5]
\draw[dashed] (-3,-4) -- (-3,4);
\begin{scope}[cm={1,0,0,2,(0,0)}]
\draw[vcurve] (0,0) circle(2cm);
\draw[cable] 
% (90:2cm) -- (210:2cm) -- (330:2cm) -- cycle 
% (270:2cm) -- (30:2cm) -- (150:2cm) -- cycle;
  (30:2cm) -- (90:2cm) -- (150:2cm) -- (210:2cm) -- (270:2cm) -- (330:2cm) --
  cycle;
\draw[strut] (90:2cm) -- (270:2cm) (210:2cm) -- (30:2cm) 
                   (330:2cm) -- (150:2cm);
\end{scope}
\begin{scope}[cm={1,0,0,2,(-12,0)}]
\draw[vcurve] (0,0) circle(2cm);
\draw[cable] 
  (10:2cm) -- (70:2cm) -- (130:2cm) -- (190:2cm) -- (250:2cm) -- (310:2cm) --
  cycle;
\draw[strut] (10:2cm) -- (190:2cm) (130:2cm) -- (310:2cm) 
                   (250:2cm) -- (70:2cm);
\end{scope}
\begin{scope}[xshift=6cm]
\draw[vcurve] (0,0) ellipse (2cm and 4cm);
\draw[cable] 
% (0,-4) -- (1.827,1.627) -- (-1.827,1.627) -- cycle
% (0,4) -- (-1.827,-1.627) -- (1.827,-1.627) -- cycle;
  (0,-4) -- (1.827,-1.627) -- (1.827,1.627) -- (0,4) -- (-1.827,1.627) --
  (-1.827,-1.627) -- cycle;
\draw[strut] (0,4) -- (0,-4) (1.827,1.627) -- (-1.827,-1.627) 
                   (-1.827,1.627) -- (1.827,-1.627);
\end{scope}
\begin{scope}[xshift=-6cm]
\draw[vcurve] (0,0) ellipse (2cm and 4cm);
\draw[cable] 
% (1.982,0.538) -- (-1.508,2.628) -- (-0.937,-3.534) -- cycle
% (-1.982,-0.538) -- (1.508,-2.628) -- (0.937,3.534) -- cycle;
  (1.982,0.538) -- (0.937,3.534) -- (-1.508,2.628) -- (-1.982,-0.538) --
  (-0.937,-3.534) -- (1.508,-2.628) -- cycle;
\draw[strut] 
  (1.982,0.538) -- (-1.982,-0.538) 
 (-1.508,2.628) -- (1.508,-2.628) 
 (-0.937,-3.534) -- (0.937,3.534);
\end{scope}
\end{tikzpicture}
\hspace*{\fill}
\capt[Two ellipse constructions.]{Two ways to build elliptical tensegrites.
Two different subtensegrities have been shown.  In both cases the picture on
the left is the result of building a circular tensegrity and then transforming
it and the picture on the right is the result of building the tensegrity on 
an ellipse.}
\label{fig:OACaffine}
\end{figure}

\section{A Corner Case}
\label{section:corner}
What might it look like to have an example where not every subtensegrity is 
bar equivalent?  Let's try to build such a creature.  Consider the family of
tensegrities represented in \vref{fig:triangles}.
\begin{figure}[ht]
\hspace*{\fill}
\begin{tikzpicture}
\foreach \x/\s in {1/0,0.5/3.5,0/7,-0.2/10.5} {
  \draw[xshift=\s cm,strut] (0,\x) -- (2.828,\x) -- (1.414,-1.732) -- cycle;
  \draw[xshift=\s cm,cable] (0,\x) node [vertex] {} -- 
                            (1.414,0) node [vertex] {} -- 
                            (1.414,-1.732) node [vertex] {}
                            (2.828,\x) node [vertex] {} -- 
                            (1.414,0) node [vertex] {};
}
\end{tikzpicture}
\hspace*{\fill}
\capt{Representatives of our family of triangular tensegrities.}
\label{fig:triangles}
\end{figure}
For the first triangle (where each side is length 2), a stress which has weight
3 on all of the cables and weight 1 on each of the struts will suffice to show
that it is bar equivalent.  As the top comes down to meet the midpoint, the
weights on the upper two cables must rise (or all of the others fall) until,
when the top strut either touches or runs below the midpoint, the tensegrity is
no longer bar equivalent.

Let's group together a lot of the bar equivalent ones and slip one that isn't
bar equivalent in amongst them.  We can choose a family of such tensegrities
with the top strut approaching the center point, as shown in \vref{fig:bigY}
(we've removed the struts from the picture for visibility's sake).
\begin{figure}[ht]
\hspace*{\fill}
\begin{tikzpicture}[xscale=3,yslant=-0.5]
\draw[cable]
   (-1.000, 1.000, 2.000) -- (-1.000, 0.000, 1.000)
   (-1.000, 1.000, 0.000) -- (-1.000, 0.000, 1.000) -- (-1.000,-1.732, 1.000);
\draw[cable]
   (-0.967, 0.934, 2.000) -- (-0.967, 0.000, 1.000)
   (-0.967, 0.934, 0.000) -- (-0.967, 0.000, 1.000) -- (-0.967,-1.732, 1.000);
\draw[cable]
   (-0.933, 0.871, 2.000) -- (-0.933, 0.000, 1.000)
   (-0.933, 0.871, 0.000) -- (-0.933, 0.000, 1.000) -- (-0.933,-1.732, 1.000);
\draw[cable]
   (-0.900, 0.810, 2.000) -- (-0.900, 0.000, 1.000)
   (-0.900, 0.810, 0.000) -- (-0.900, 0.000, 1.000) -- (-0.900,-1.732, 1.000);
\draw[cable]
   (-0.867, 0.751, 2.000) -- (-0.867, 0.000, 1.000)
   (-0.867, 0.751, 0.000) -- (-0.867, 0.000, 1.000) -- (-0.867,-1.732, 1.000);
\draw[cable]
   (-0.833, 0.694, 2.000) -- (-0.833, 0.000, 1.000)
   (-0.833, 0.694, 0.000) -- (-0.833, 0.000, 1.000) -- (-0.833,-1.732, 1.000);
\draw[cable]
   (-0.800, 0.640, 2.000) -- (-0.800, 0.000, 1.000)
   (-0.800, 0.640, 0.000) -- (-0.800, 0.000, 1.000) -- (-0.800,-1.732, 1.000);
\draw[cable]
   (-0.767, 0.588, 2.000) -- (-0.767, 0.000, 1.000)
   (-0.767, 0.588, 0.000) -- (-0.767, 0.000, 1.000) -- (-0.767,-1.732, 1.000);
\draw[cable]
   (-0.733, 0.538, 2.000) -- (-0.733, 0.000, 1.000)
   (-0.733, 0.538, 0.000) -- (-0.733, 0.000, 1.000) -- (-0.733,-1.732, 1.000);
\draw[cable]
   (-0.700, 0.490, 2.000) -- (-0.700, 0.000, 1.000)
   (-0.700, 0.490, 0.000) -- (-0.700, 0.000, 1.000) -- (-0.700,-1.732, 1.000);
\draw[cable]
   (-0.667, 0.444, 2.000) -- (-0.667, 0.000, 1.000)
   (-0.667, 0.444, 0.000) -- (-0.667, 0.000, 1.000) -- (-0.667,-1.732, 1.000);
\draw[cable]
   (-0.633, 0.401, 2.000) -- (-0.633, 0.000, 1.000)
   (-0.633, 0.401, 0.000) -- (-0.633, 0.000, 1.000) -- (-0.633,-1.732, 1.000);
\draw[cable]
   (-0.600, 0.360, 2.000) -- (-0.600, 0.000, 1.000)
   (-0.600, 0.360, 0.000) -- (-0.600, 0.000, 1.000) -- (-0.600,-1.732, 1.000);
\draw[cable]
   (-0.567, 0.321, 2.000) -- (-0.567, 0.000, 1.000)
   (-0.567, 0.321, 0.000) -- (-0.567, 0.000, 1.000) -- (-0.567,-1.732, 1.000);
\draw[cable]
   (-0.533, 0.284, 2.000) -- (-0.533, 0.000, 1.000)
   (-0.533, 0.284, 0.000) -- (-0.533, 0.000, 1.000) -- (-0.533,-1.732, 1.000);
\draw[cable]
   (-0.500, 0.250, 2.000) -- (-0.500, 0.000, 1.000)
   (-0.500, 0.250, 0.000) -- (-0.500, 0.000, 1.000) -- (-0.500,-1.732, 1.000);
\draw[cable]
   (-0.467, 0.218, 2.000) -- (-0.467, 0.000, 1.000)
   (-0.467, 0.218, 0.000) -- (-0.467, 0.000, 1.000) -- (-0.467,-1.732, 1.000);
\draw[cable]
   (-0.433, 0.188, 2.000) -- (-0.433, 0.000, 1.000)
   (-0.433, 0.188, 0.000) -- (-0.433, 0.000, 1.000) -- (-0.433,-1.732, 1.000);
\draw[cable]
   (-0.400, 0.160, 2.000) -- (-0.400, 0.000, 1.000)
   (-0.400, 0.160, 0.000) -- (-0.400, 0.000, 1.000) -- (-0.400,-1.732, 1.000);
\draw[cable]
   (-0.367, 0.134, 2.000) -- (-0.367, 0.000, 1.000)
   (-0.367, 0.134, 0.000) -- (-0.367, 0.000, 1.000) -- (-0.367,-1.732, 1.000);
\draw[cable]
   (-0.333, 0.111, 2.000) -- (-0.333, 0.000, 1.000)
   (-0.333, 0.111, 0.000) -- (-0.333, 0.000, 1.000) -- (-0.333,-1.732, 1.000);
\draw[cable]
   (-0.300, 0.090, 2.000) -- (-0.300, 0.000, 1.000)
   (-0.300, 0.090, 0.000) -- (-0.300, 0.000, 1.000) -- (-0.300,-1.732, 1.000);
\draw[cable]
   (-0.267, 0.071, 2.000) -- (-0.267, 0.000, 1.000)
   (-0.267, 0.071, 0.000) -- (-0.267, 0.000, 1.000) -- (-0.267,-1.732, 1.000);
\draw[cable]
   (-0.233, 0.054, 2.000) -- (-0.233, 0.000, 1.000)
   (-0.233, 0.054, 0.000) -- (-0.233, 0.000, 1.000) -- (-0.233,-1.732, 1.000);
\draw[cable]
   (-0.200, 0.040, 2.000) -- (-0.200, 0.000, 1.000)
   (-0.200, 0.040, 0.000) -- (-0.200, 0.000, 1.000) -- (-0.200,-1.732, 1.000);
\draw[cable]
   (-0.167, 0.028, 2.000) -- (-0.167, 0.000, 1.000)
   (-0.167, 0.028, 0.000) -- (-0.167, 0.000, 1.000) -- (-0.167,-1.732, 1.000);
\draw[cable]
   (-0.133, 0.018, 2.000) -- (-0.133, 0.000, 1.000)
   (-0.133, 0.018, 0.000) -- (-0.133, 0.000, 1.000) -- (-0.133,-1.732, 1.000);
\draw[cable]
   (-0.100, 0.010, 2.000) -- (-0.100, 0.000, 1.000)
   (-0.100, 0.010, 0.000) -- (-0.100, 0.000, 1.000) -- (-0.100,-1.732, 1.000);
\draw[cable]
   (-0.067, 0.004, 2.000) -- (-0.067, 0.000, 1.000)
   (-0.067, 0.004, 0.000) -- (-0.067, 0.000, 1.000) -- (-0.067,-1.732, 1.000);
\draw[cable]
   (-0.033, 0.001, 2.000) -- (-0.033, 0.000, 1.000)
   (-0.033, 0.001, 0.000) -- (-0.033, 0.000, 1.000) -- (-0.033,-1.732, 1.000);
\draw[cable]
   (-0.000, 0.000, 2.000) -- (-0.000, 0.000, 1.000)
   (-0.000, 0.000, 0.000) -- (-0.000, 0.000, 1.000) -- (-0.000,-1.732, 1.000);
\draw[cable]
   ( 0.033, 0.001, 2.000) -- ( 0.033, 0.000, 1.000)
   ( 0.033, 0.001, 0.000) -- ( 0.033, 0.000, 1.000) -- ( 0.033,-1.732, 1.000);
\draw[cable]
   ( 0.067, 0.004, 2.000) -- ( 0.067, 0.000, 1.000)
   ( 0.067, 0.004, 0.000) -- ( 0.067, 0.000, 1.000) -- ( 0.067,-1.732, 1.000);
\draw[cable]
   ( 0.100, 0.010, 2.000) -- ( 0.100, 0.000, 1.000)
   ( 0.100, 0.010, 0.000) -- ( 0.100, 0.000, 1.000) -- ( 0.100,-1.732, 1.000);
\draw[cable]
   ( 0.133, 0.018, 2.000) -- ( 0.133, 0.000, 1.000)
   ( 0.133, 0.018, 0.000) -- ( 0.133, 0.000, 1.000) -- ( 0.133,-1.732, 1.000);
\draw[cable]
   ( 0.167, 0.028, 2.000) -- ( 0.167, 0.000, 1.000)
   ( 0.167, 0.028, 0.000) -- ( 0.167, 0.000, 1.000) -- ( 0.167,-1.732, 1.000);
\draw[cable]
   ( 0.200, 0.040, 2.000) -- ( 0.200, 0.000, 1.000)
   ( 0.200, 0.040, 0.000) -- ( 0.200, 0.000, 1.000) -- ( 0.200,-1.732, 1.000);
\draw[cable]
   ( 0.233, 0.054, 2.000) -- ( 0.233, 0.000, 1.000)
   ( 0.233, 0.054, 0.000) -- ( 0.233, 0.000, 1.000) -- ( 0.233,-1.732, 1.000);
\draw[cable]
   ( 0.267, 0.071, 2.000) -- ( 0.267, 0.000, 1.000)
   ( 0.267, 0.071, 0.000) -- ( 0.267, 0.000, 1.000) -- ( 0.267,-1.732, 1.000);
\draw[cable]
   ( 0.300, 0.090, 2.000) -- ( 0.300, 0.000, 1.000)
   ( 0.300, 0.090, 0.000) -- ( 0.300, 0.000, 1.000) -- ( 0.300,-1.732, 1.000);
\draw[cable]
   ( 0.333, 0.111, 2.000) -- ( 0.333, 0.000, 1.000)
   ( 0.333, 0.111, 0.000) -- ( 0.333, 0.000, 1.000) -- ( 0.333,-1.732, 1.000);
\draw[cable]
   ( 0.367, 0.134, 2.000) -- ( 0.367, 0.000, 1.000)
   ( 0.367, 0.134, 0.000) -- ( 0.367, 0.000, 1.000) -- ( 0.367,-1.732, 1.000);
\draw[cable]
   ( 0.400, 0.160, 2.000) -- ( 0.400, 0.000, 1.000)
   ( 0.400, 0.160, 0.000) -- ( 0.400, 0.000, 1.000) -- ( 0.400,-1.732, 1.000);
\draw[cable]
   ( 0.433, 0.188, 2.000) -- ( 0.433, 0.000, 1.000)
   ( 0.433, 0.188, 0.000) -- ( 0.433, 0.000, 1.000) -- ( 0.433,-1.732, 1.000);
\draw[cable]
   ( 0.467, 0.218, 2.000) -- ( 0.467, 0.000, 1.000)
   ( 0.467, 0.218, 0.000) -- ( 0.467, 0.000, 1.000) -- ( 0.467,-1.732, 1.000);
\draw[cable]
   ( 0.500, 0.250, 2.000) -- ( 0.500, 0.000, 1.000)
   ( 0.500, 0.250, 0.000) -- ( 0.500, 0.000, 1.000) -- ( 0.500,-1.732, 1.000);
\draw[cable]
   ( 0.533, 0.284, 2.000) -- ( 0.533, 0.000, 1.000)
   ( 0.533, 0.284, 0.000) -- ( 0.533, 0.000, 1.000) -- ( 0.533,-1.732, 1.000);
\draw[cable]
   ( 0.567, 0.321, 2.000) -- ( 0.567, 0.000, 1.000)
   ( 0.567, 0.321, 0.000) -- ( 0.567, 0.000, 1.000) -- ( 0.567,-1.732, 1.000);
\draw[cable]
   ( 0.600, 0.360, 2.000) -- ( 0.600, 0.000, 1.000)
   ( 0.600, 0.360, 0.000) -- ( 0.600, 0.000, 1.000) -- ( 0.600,-1.732, 1.000);
\draw[cable]
   ( 0.633, 0.401, 2.000) -- ( 0.633, 0.000, 1.000)
   ( 0.633, 0.401, 0.000) -- ( 0.633, 0.000, 1.000) -- ( 0.633,-1.732, 1.000);
\draw[cable]
   ( 0.667, 0.444, 2.000) -- ( 0.667, 0.000, 1.000)
   ( 0.667, 0.444, 0.000) -- ( 0.667, 0.000, 1.000) -- ( 0.667,-1.732, 1.000);
\draw[cable]
   ( 0.700, 0.490, 2.000) -- ( 0.700, 0.000, 1.000)
   ( 0.700, 0.490, 0.000) -- ( 0.700, 0.000, 1.000) -- ( 0.700,-1.732, 1.000);
\draw[cable]
   ( 0.733, 0.538, 2.000) -- ( 0.733, 0.000, 1.000)
   ( 0.733, 0.538, 0.000) -- ( 0.733, 0.000, 1.000) -- ( 0.733,-1.732, 1.000);
\draw[cable]
   ( 0.767, 0.588, 2.000) -- ( 0.767, 0.000, 1.000)
   ( 0.767, 0.588, 0.000) -- ( 0.767, 0.000, 1.000) -- ( 0.767,-1.732, 1.000);
\draw[cable]
   ( 0.800, 0.640, 2.000) -- ( 0.800, 0.000, 1.000)
   ( 0.800, 0.640, 0.000) -- ( 0.800, 0.000, 1.000) -- ( 0.800,-1.732, 1.000);
\draw[cable]
   ( 0.833, 0.694, 2.000) -- ( 0.833, 0.000, 1.000)
   ( 0.833, 0.694, 0.000) -- ( 0.833, 0.000, 1.000) -- ( 0.833,-1.732, 1.000);
\draw[cable]
   ( 0.867, 0.751, 2.000) -- ( 0.867, 0.000, 1.000)
   ( 0.867, 0.751, 0.000) -- ( 0.867, 0.000, 1.000) -- ( 0.867,-1.732, 1.000);
\draw[cable]
   ( 0.900, 0.810, 2.000) -- ( 0.900, 0.000, 1.000)
   ( 0.900, 0.810, 0.000) -- ( 0.900, 0.000, 1.000) -- ( 0.900,-1.732, 1.000);
\draw[cable]
   ( 0.933, 0.871, 2.000) -- ( 0.933, 0.000, 1.000)
   ( 0.933, 0.871, 0.000) -- ( 0.933, 0.000, 1.000) -- ( 0.933,-1.732, 1.000);
\draw[cable]
   ( 0.967, 0.934, 2.000) -- ( 0.967, 0.000, 1.000)
   ( 0.967, 0.934, 0.000) -- ( 0.967, 0.000, 1.000) -- ( 0.967,-1.732, 1.000);
\draw[cable]
   ( 1.000, 1.000, 2.000) -- ( 1.000, 0.000, 1.000)
   ( 1.000, 1.000, 0.000) -- ( 1.000, 0.000, 1.000) -- ( 1.000,-1.732, 1.000);
\draw[vcurve] (-1,0,1) -- (1,0,1) (-1,-1.732,1) -- (1,-1.732,1);
\draw[vcurve] plot [smooth] coordinates {
  (-1,1,0) (-0.967,0.934,0) (-0.933,0.871,0) (-0.900,0.810,0) (-0.867,0.751,0)
  (-0.833,0.694,0) (-0.800,0.640,0) (-0.767,0.588,0) (-0.733,0.538,0)
  (-0.700,0.490,0) (-0.667,0.444,0) (-0.633,0.401,0) (-0.600,0.360,0)
  (-0.567,0.321,0) (-0.533,0.284,0) (-0.500,0.250,0) (-0.467,0.218,0)
  (-0.433,0.188,0) (-0.400,0.160,0) (-0.367,0.134,0) (-0.333,0.111,0)
  (-0.300,0.090,0) (-0.267,0.071,0) (-0.233,0.054,0) (-0.200,0.040,0)
  (-0.167,0.028,0) (-0.133,0.018,0) (-0.100,0.010,0) (-0.067,0.004,0)
  (-0.033,0.001,0) (-0,0,0) (0.033,0.001,0) (0.067,0.004,0) (0.100,0.010,0)
  (0.133,0.018,0) (0.167,0.028,0) (0.200,0.040,0) (0.233,0.054,0)
  (0.267,0.071,0) (0.300,0.090,0) (0.333,0.111,0) (0.367,0.134,0)
  (0.400,0.160,0) (0.433,0.188,0) (0.467,0.218,0) (0.500,0.250,0)
  (0.533,0.284,0) (0.567,0.321,0) (0.600,0.360,0) (0.633,0.401,0)
  (0.667,0.444,0) (0.700,0.490,0) (0.733,0.538,0) (0.767,0.588,0)
  (0.800,0.640,0) (0.833,0.694,0) (0.867,0.751,0) (0.900,0.810,0)
  (0.933,0.871,0) (0.967,0.934,0) (1,1,0)
};
\draw[vcurve] plot [smooth] coordinates {
  (-1,1,2) (-0.967,0.934,2) (-0.933,0.871,2) (-0.900,0.810,2) (-0.867,0.751,2)
  (-0.833,0.694,2) (-0.800,0.640,2) (-0.767,0.588,2) (-0.733,0.538,2)
  (-0.700,0.490,2) (-0.667,0.444,2) (-0.633,0.401,2) (-0.600,0.360,2)
  (-0.567,0.321,2) (-0.533,0.284,2) (-0.500,0.250,2) (-0.467,0.218,2)
  (-0.433,0.188,2) (-0.400,0.160,2) (-0.367,0.134,2) (-0.333,0.111,2)
  (-0.300,0.090,2) (-0.267,0.071,2) (-0.233,0.054,2) (-0.200,0.040,2)
  (-0.167,0.028,2) (-0.133,0.018,2) (-0.100,0.010,2) (-0.067,0.004,2)
  (-0.033,0.001,2) (-0,0,2) (0.033,0.001,2) (0.067,0.004,2) (0.100,0.010,2)
  (0.133,0.018,2) (0.167,0.028,2) (0.200,0.040,2) (0.233,0.054,2)
  (0.267,0.071,2) (0.300,0.090,2) (0.333,0.111,2) (0.367,0.134,2)
  (0.400,0.160,2) (0.433,0.188,2) (0.467,0.218,2) (0.500,0.250,2)
  (0.533,0.284,2) (0.567,0.321,2) (0.600,0.360,2) (0.633,0.401,2)
  (0.667,0.444,2) (0.700,0.490,2) (0.733,0.538,2) (0.767,0.588,2)
  (0.800,0.640,2) (0.833,0.694,2) (0.867,0.751,2) (0.900,0.810,2)
  (0.933,0.871,2) (0.967,0.934,2) (1,1,2)
};
\end{tikzpicture}
\hspace*{\fill}
\capt[Continuous tensegrity formed from triangle subtensegrities.]{Triangle
subtensegrities with varying top angle forming a continuous tensegrity.  The
outer strut triangles are not shown.}
\label{fig:bigY}
\end{figure}

Now in three dimensions, even those subtensegrities that previously were
infinitesimally rigid are no longer so, as they have no way to resist
variations normal to their planes of definition.  However, as Proposition
\vref{prop:Ambient} tells us, they remain bar equivalent.
\begin{prop}
\label{prop:Ambient}
A finite tensegrity that is bar equivalent in $\RR^n$ is bar equivalent in
Euclidean space of all dimensions higher than~$n$.
\index{ambient space!increasing dimension}
\end{prop}
\begin{proof}
Let $G(p)$ be a finite, bar equivalent tensegrity in~$\RR^n$.  By Lemma 
\vref{lemma:NewFiveOne}, that tensegrity has a strictly positive stress, that
is, a positive linear dependence among the row vectors of the rigidity matrix.
Moving to a higher dimension will change those vectors, but only by adding 
columns of zeros for the new dimension, so it doesn't change that linear
dependence.  Since $G(p)$ continues to have a positive stress, it continues
to be bar equivalent.
\end{proof}

So our continuous family is countably covered by bar equivalent subtensegrities
and hence must be bar equivalent itself.  That means, of course, that 
$Y(\X) \cap C(\E)^+ = \{ \zero \}$, which is to say, no matter which $V$
we pick, $YV$ must have some negative values.

That's true, but there's something a little strange going on here.  It turns
out that the image of $Y$, while intersecting the nonnegative orthant only
at the origin, gets arbitrarily close to it elsewhere.  Let's see if we can
watch that happen.

Take the vertical cable in each subtensegrity.  We can shrink this cable by 
putting a length-1 vector pointing directly downward at the central node and 
balancing it with a length $\nicefrac{1}{3}$ upward-pointing vector at each 
of the other three nodes (see \vref{fig:balancingForce}).  
\begin{figure}[ht]
\hspace*{\fill}
\begin{tikzpicture}
\draw[strut] (0,0) -- (2,0) -- (1,-1.732) -- cycle;
\draw[cable] (0,0) -- (1,-0.577) -- (2,0) (1,-0.577) --
  (1,-1.732);
\draw[->] (0,-0.333) -- (0,0);
\draw[->] (2,-0.333) -- (2,0) ;
\draw[->] (1,-2.065) -- (1,-1.732);
\draw[->] (1,0.423) -- (1,-0.577);
\end{tikzpicture}
\hspace*{\fill}
\capt{Balancing forces placed on a subtensegrity.}
\label{fig:balancingForce}
\end{figure}

This set of vectors will have no effect whatsoever on the struts.  It will 
produce a value of $\nicefrac{4}{3}$ on the vertical cable, and it will 
produce a force on the other two cables that depends on their angle (see
\vref{fig:armForces}).  
\begin{figure}[ht]
\hspace*{\fill}
\begin{tikzpicture}[scale=3,>=stealth]
\draw[cable] (1,-0.577) -- (2,0);
\draw[xshift=0.05 cm,snake=brace,segment amplitude=0.1 cm] (2,0) -- 
  node [right] {$y$} (2,-0.577);
\draw[yshift=0.05 cm,snake=brace,segment amplitude=0.1 cm] (1,0) -- 
  node [above] {$1$} (2,0);
\draw[->] (2,-0.33) -- (2,0);
\draw[->] (1,0.423) -- (1,-0.577);
\end{tikzpicture}
\hspace*{\fill}
\capt{The forces present on one ``arm'' of the subtensegrity.}
\label{fig:armForces}
\end{figure}

In every case, the horizontal length of each of those cables is 1.  If the
vertical length is $y$, then the rigidity operator will take our set of
vectors to the value $-\nicefrac{2}{3} \frac{y}{\sqrt{1 + y^2}}$.  
Now when $y = 0$, this value is~$0$.  That makes sense, since the forces are
now acting orthogonally to the edge, but we can use that fact to our advantage.

Suppose, now, that we identify our family of subtensegrities with the interval
$[-1,1]$ in such a way that if $x \in [-1,1]$ gives a certain subtensegrity,
then $y = x^2$ for that tensegrity (this is function used in \ref{fig:bigY}).
Now we take a family of vector fields on our vertices.  For every
value $i \in \ZZ$, we let $V_i(x)$ be the zero vector for $x \notin
[-\frac{1}{i},\frac{1}{i}]$ and then have vectors in the ratios we've discussed
increasing linearly from $0$ at $x = -\frac{1}{i}$ to $1$ at $x=0$ and back
down to $0$ at $x = \frac{1}{i}$ (the case $i = 2$ is shown in
\vref{fig:bigYwithMotion}).  
\begin{figure}[ht]
\hfil
\begin{tikzpicture}[xscale=3,yslant=-0.5,>=stealth]
\draw[cable]
   (-1.000, 1.000, 2.000) -- (-1.000, 0.000, 1.000)
   (-1.000, 1.000, 0.000) -- (-1.000, 0.000, 1.000) -- (-1.000,-1.732, 1.000);
\draw[cable]
   (-0.967, 0.934, 2.000) -- (-0.967, 0.000, 1.000)
   (-0.967, 0.934, 0.000) -- (-0.967, 0.000, 1.000) -- (-0.967,-1.732, 1.000);
\draw[cable]
   (-0.933, 0.871, 2.000) -- (-0.933, 0.000, 1.000)
   (-0.933, 0.871, 0.000) -- (-0.933, 0.000, 1.000) -- (-0.933,-1.732, 1.000);
\draw[cable]
   (-0.900, 0.810, 2.000) -- (-0.900, 0.000, 1.000)
   (-0.900, 0.810, 0.000) -- (-0.900, 0.000, 1.000) -- (-0.900,-1.732, 1.000);
\draw[cable]
   (-0.867, 0.751, 2.000) -- (-0.867, 0.000, 1.000)
   (-0.867, 0.751, 0.000) -- (-0.867, 0.000, 1.000) -- (-0.867,-1.732, 1.000);
\draw[cable]
   (-0.833, 0.694, 2.000) -- (-0.833, 0.000, 1.000)
   (-0.833, 0.694, 0.000) -- (-0.833, 0.000, 1.000) -- (-0.833,-1.732, 1.000);
\draw[cable]
   (-0.800, 0.640, 2.000) -- (-0.800, 0.000, 1.000)
   (-0.800, 0.640, 0.000) -- (-0.800, 0.000, 1.000) -- (-0.800,-1.732, 1.000);
\draw[cable]
   (-0.767, 0.588, 2.000) -- (-0.767, 0.000, 1.000)
   (-0.767, 0.588, 0.000) -- (-0.767, 0.000, 1.000) -- (-0.767,-1.732, 1.000);
\draw[cable]
   (-0.733, 0.538, 2.000) -- (-0.733, 0.000, 1.000)
   (-0.733, 0.538, 0.000) -- (-0.733, 0.000, 1.000) -- (-0.733,-1.732, 1.000);
\draw[cable]
   (-0.700, 0.490, 2.000) -- (-0.700, 0.000, 1.000)
   (-0.700, 0.490, 0.000) -- (-0.700, 0.000, 1.000) -- (-0.700,-1.732, 1.000);
\draw[cable]
   (-0.667, 0.444, 2.000) -- (-0.667, 0.000, 1.000)
   (-0.667, 0.444, 0.000) -- (-0.667, 0.000, 1.000) -- (-0.667,-1.732, 1.000);
\draw[cable]
   (-0.633, 0.401, 2.000) -- (-0.633, 0.000, 1.000)
   (-0.633, 0.401, 0.000) -- (-0.633, 0.000, 1.000) -- (-0.633,-1.732, 1.000);
\draw[cable]
   (-0.600, 0.360, 2.000) -- (-0.600, 0.000, 1.000)
   (-0.600, 0.360, 0.000) -- (-0.600, 0.000, 1.000) -- (-0.600,-1.732, 1.000);
\draw[cable]
   (-0.567, 0.321, 2.000) -- (-0.567, 0.000, 1.000)
   (-0.567, 0.321, 0.000) -- (-0.567, 0.000, 1.000) -- (-0.567,-1.732, 1.000);
\draw[cable]
   (-0.533, 0.284, 2.000) -- (-0.533, 0.000, 1.000)
   (-0.533, 0.284, 0.000) -- (-0.533, 0.000, 1.000) -- (-0.533,-1.732, 1.000);
\draw[cable]
   (-0.500, 0.250, 2.000) -- (-0.500, 0.000, 1.000)
   (-0.500, 0.250, 0.000) -- (-0.500, 0.000, 1.000) -- (-0.500,-1.732, 1.000);
\draw[cable]
   (-0.467, 0.218, 2.000) -- (-0.467, 0.000, 1.000)
   (-0.467, 0.218, 0.000) -- (-0.467, 0.000, 1.000) -- (-0.467,-1.732, 1.000);
\draw[cable]
   (-0.433, 0.188, 2.000) -- (-0.433, 0.000, 1.000)
   (-0.433, 0.188, 0.000) -- (-0.433, 0.000, 1.000) -- (-0.433,-1.732, 1.000);
\draw[cable]
   (-0.400, 0.160, 2.000) -- (-0.400, 0.000, 1.000)
   (-0.400, 0.160, 0.000) -- (-0.400, 0.000, 1.000) -- (-0.400,-1.732, 1.000);
\draw[cable]
   (-0.367, 0.134, 2.000) -- (-0.367, 0.000, 1.000)
   (-0.367, 0.134, 0.000) -- (-0.367, 0.000, 1.000) -- (-0.367,-1.732, 1.000);
\draw[cable]
   (-0.333, 0.111, 2.000) -- (-0.333, 0.000, 1.000)
   (-0.333, 0.111, 0.000) -- (-0.333, 0.000, 1.000) -- (-0.333,-1.732, 1.000);
\draw[cable]
   (-0.300, 0.090, 2.000) -- (-0.300, 0.000, 1.000)
   (-0.300, 0.090, 0.000) -- (-0.300, 0.000, 1.000) -- (-0.300,-1.732, 1.000);
\draw[cable]
   (-0.267, 0.071, 2.000) -- (-0.267, 0.000, 1.000)
   (-0.267, 0.071, 0.000) -- (-0.267, 0.000, 1.000) -- (-0.267,-1.732, 1.000);
\draw[cable]
   (-0.233, 0.054, 2.000) -- (-0.233, 0.000, 1.000)
   (-0.233, 0.054, 0.000) -- (-0.233, 0.000, 1.000) -- (-0.233,-1.732, 1.000);
\draw[cable]
   (-0.200, 0.040, 2.000) -- (-0.200, 0.000, 1.000)
   (-0.200, 0.040, 0.000) -- (-0.200, 0.000, 1.000) -- (-0.200,-1.732, 1.000);
\draw[cable]
   (-0.167, 0.028, 2.000) -- (-0.167, 0.000, 1.000)
   (-0.167, 0.028, 0.000) -- (-0.167, 0.000, 1.000) -- (-0.167,-1.732, 1.000);
\draw[cable]
   (-0.133, 0.018, 2.000) -- (-0.133, 0.000, 1.000)
   (-0.133, 0.018, 0.000) -- (-0.133, 0.000, 1.000) -- (-0.133,-1.732, 1.000);
\draw[cable]
   (-0.100, 0.010, 2.000) -- (-0.100, 0.000, 1.000)
   (-0.100, 0.010, 0.000) -- (-0.100, 0.000, 1.000) -- (-0.100,-1.732, 1.000);
\draw[cable]
   (-0.067, 0.004, 2.000) -- (-0.067, 0.000, 1.000)
   (-0.067, 0.004, 0.000) -- (-0.067, 0.000, 1.000) -- (-0.067,-1.732, 1.000);
\draw[cable]
   (-0.033, 0.001, 2.000) -- (-0.033, 0.000, 1.000)
   (-0.033, 0.001, 0.000) -- (-0.033, 0.000, 1.000) -- (-0.033,-1.732, 1.000);
\draw[cable]
   (-0.000, 0.000, 2.000) -- (-0.000, 0.000, 1.000)
   (-0.000, 0.000, 0.000) -- (-0.000, 0.000, 1.000) -- (-0.000,-1.732, 1.000);
\draw[cable]
   ( 0.033, 0.001, 2.000) -- ( 0.033, 0.000, 1.000)
   ( 0.033, 0.001, 0.000) -- ( 0.033, 0.000, 1.000) -- ( 0.033,-1.732, 1.000);
\draw[cable]
   ( 0.067, 0.004, 2.000) -- ( 0.067, 0.000, 1.000)
   ( 0.067, 0.004, 0.000) -- ( 0.067, 0.000, 1.000) -- ( 0.067,-1.732, 1.000);
\draw[cable]
   ( 0.100, 0.010, 2.000) -- ( 0.100, 0.000, 1.000)
   ( 0.100, 0.010, 0.000) -- ( 0.100, 0.000, 1.000) -- ( 0.100,-1.732, 1.000);
\draw[cable]
   ( 0.133, 0.018, 2.000) -- ( 0.133, 0.000, 1.000)
   ( 0.133, 0.018, 0.000) -- ( 0.133, 0.000, 1.000) -- ( 0.133,-1.732, 1.000);
\draw[cable]
   ( 0.167, 0.028, 2.000) -- ( 0.167, 0.000, 1.000)
   ( 0.167, 0.028, 0.000) -- ( 0.167, 0.000, 1.000) -- ( 0.167,-1.732, 1.000);
\draw[cable]
   ( 0.200, 0.040, 2.000) -- ( 0.200, 0.000, 1.000)
   ( 0.200, 0.040, 0.000) -- ( 0.200, 0.000, 1.000) -- ( 0.200,-1.732, 1.000);
\draw[cable]
   ( 0.233, 0.054, 2.000) -- ( 0.233, 0.000, 1.000)
   ( 0.233, 0.054, 0.000) -- ( 0.233, 0.000, 1.000) -- ( 0.233,-1.732, 1.000);
\draw[cable]
   ( 0.267, 0.071, 2.000) -- ( 0.267, 0.000, 1.000)
   ( 0.267, 0.071, 0.000) -- ( 0.267, 0.000, 1.000) -- ( 0.267,-1.732, 1.000);
\draw[cable]
   ( 0.300, 0.090, 2.000) -- ( 0.300, 0.000, 1.000)
   ( 0.300, 0.090, 0.000) -- ( 0.300, 0.000, 1.000) -- ( 0.300,-1.732, 1.000);
\draw[cable]
   ( 0.333, 0.111, 2.000) -- ( 0.333, 0.000, 1.000)
   ( 0.333, 0.111, 0.000) -- ( 0.333, 0.000, 1.000) -- ( 0.333,-1.732, 1.000);
\draw[cable]
   ( 0.367, 0.134, 2.000) -- ( 0.367, 0.000, 1.000)
   ( 0.367, 0.134, 0.000) -- ( 0.367, 0.000, 1.000) -- ( 0.367,-1.732, 1.000);
\draw[cable]
   ( 0.400, 0.160, 2.000) -- ( 0.400, 0.000, 1.000)
   ( 0.400, 0.160, 0.000) -- ( 0.400, 0.000, 1.000) -- ( 0.400,-1.732, 1.000);
\draw[cable]
   ( 0.433, 0.188, 2.000) -- ( 0.433, 0.000, 1.000)
   ( 0.433, 0.188, 0.000) -- ( 0.433, 0.000, 1.000) -- ( 0.433,-1.732, 1.000);
\draw[cable]
   ( 0.467, 0.218, 2.000) -- ( 0.467, 0.000, 1.000)
   ( 0.467, 0.218, 0.000) -- ( 0.467, 0.000, 1.000) -- ( 0.467,-1.732, 1.000);
\draw[cable]
   ( 0.500, 0.250, 2.000) -- ( 0.500, 0.000, 1.000)
   ( 0.500, 0.250, 0.000) -- ( 0.500, 0.000, 1.000) -- ( 0.500,-1.732, 1.000);
\draw[cable]
   ( 0.533, 0.284, 2.000) -- ( 0.533, 0.000, 1.000)
   ( 0.533, 0.284, 0.000) -- ( 0.533, 0.000, 1.000) -- ( 0.533,-1.732, 1.000);
\draw[cable]
   ( 0.567, 0.321, 2.000) -- ( 0.567, 0.000, 1.000)
   ( 0.567, 0.321, 0.000) -- ( 0.567, 0.000, 1.000) -- ( 0.567,-1.732, 1.000);
\draw[cable]
   ( 0.600, 0.360, 2.000) -- ( 0.600, 0.000, 1.000)
   ( 0.600, 0.360, 0.000) -- ( 0.600, 0.000, 1.000) -- ( 0.600,-1.732, 1.000);
\draw[cable]
   ( 0.633, 0.401, 2.000) -- ( 0.633, 0.000, 1.000)
   ( 0.633, 0.401, 0.000) -- ( 0.633, 0.000, 1.000) -- ( 0.633,-1.732, 1.000);
\draw[cable]
   ( 0.667, 0.444, 2.000) -- ( 0.667, 0.000, 1.000)
   ( 0.667, 0.444, 0.000) -- ( 0.667, 0.000, 1.000) -- ( 0.667,-1.732, 1.000);
\draw[cable]
   ( 0.700, 0.490, 2.000) -- ( 0.700, 0.000, 1.000)
   ( 0.700, 0.490, 0.000) -- ( 0.700, 0.000, 1.000) -- ( 0.700,-1.732, 1.000);
\draw[cable]
   ( 0.733, 0.538, 2.000) -- ( 0.733, 0.000, 1.000)
   ( 0.733, 0.538, 0.000) -- ( 0.733, 0.000, 1.000) -- ( 0.733,-1.732, 1.000);
\draw[cable]
   ( 0.767, 0.588, 2.000) -- ( 0.767, 0.000, 1.000)
   ( 0.767, 0.588, 0.000) -- ( 0.767, 0.000, 1.000) -- ( 0.767,-1.732, 1.000);
\draw[cable]
   ( 0.800, 0.640, 2.000) -- ( 0.800, 0.000, 1.000)
   ( 0.800, 0.640, 0.000) -- ( 0.800, 0.000, 1.000) -- ( 0.800,-1.732, 1.000);
\draw[cable]
   ( 0.833, 0.694, 2.000) -- ( 0.833, 0.000, 1.000)
   ( 0.833, 0.694, 0.000) -- ( 0.833, 0.000, 1.000) -- ( 0.833,-1.732, 1.000);
\draw[cable]
   ( 0.867, 0.751, 2.000) -- ( 0.867, 0.000, 1.000)
   ( 0.867, 0.751, 0.000) -- ( 0.867, 0.000, 1.000) -- ( 0.867,-1.732, 1.000);
\draw[cable]
   ( 0.900, 0.810, 2.000) -- ( 0.900, 0.000, 1.000)
   ( 0.900, 0.810, 0.000) -- ( 0.900, 0.000, 1.000) -- ( 0.900,-1.732, 1.000);
\draw[cable]
   ( 0.933, 0.871, 2.000) -- ( 0.933, 0.000, 1.000)
   ( 0.933, 0.871, 0.000) -- ( 0.933, 0.000, 1.000) -- ( 0.933,-1.732, 1.000);
\draw[cable]
   ( 0.967, 0.934, 2.000) -- ( 0.967, 0.000, 1.000)
   ( 0.967, 0.934, 0.000) -- ( 0.967, 0.000, 1.000) -- ( 0.967,-1.732, 1.000);
\draw[cable]
   ( 1.000, 1.000, 2.000) -- ( 1.000, 0.000, 1.000)
   ( 1.000, 1.000, 0.000) -- ( 1.000, 0.000, 1.000) -- ( 1.000,-1.732, 1.000);
\draw[vcurve] (-1,0,1) -- (1,0,1) (-1,-1.732,1) -- (1,-1.732,1);
\draw[vcurve] plot [smooth] coordinates {
  (-1,1,0) (-0.967,0.934,0) (-0.933,0.871,0) (-0.900,0.810,0) (-0.867,0.751,0)
  (-0.833,0.694,0) (-0.800,0.640,0) (-0.767,0.588,0) (-0.733,0.538,0)
  (-0.700,0.490,0) (-0.667,0.444,0) (-0.633,0.401,0) (-0.600,0.360,0)
  (-0.567,0.321,0) (-0.533,0.284,0) (-0.500,0.250,0) (-0.467,0.218,0)
  (-0.433,0.188,0) (-0.400,0.160,0) (-0.367,0.134,0) (-0.333,0.111,0)
  (-0.300,0.090,0) (-0.267,0.071,0) (-0.233,0.054,0) (-0.200,0.040,0)
  (-0.167,0.028,0) (-0.133,0.018,0) (-0.100,0.010,0) (-0.067,0.004,0)
  (-0.033,0.001,0) (-0,0,0) (0.033,0.001,0) (0.067,0.004,0) (0.100,0.010,0)
  (0.133,0.018,0) (0.167,0.028,0) (0.200,0.040,0) (0.233,0.054,0)
  (0.267,0.071,0) (0.300,0.090,0) (0.333,0.111,0) (0.367,0.134,0)
  (0.400,0.160,0) (0.433,0.188,0) (0.467,0.218,0) (0.500,0.250,0)
  (0.533,0.284,0) (0.567,0.321,0) (0.600,0.360,0) (0.633,0.401,0)
  (0.667,0.444,0) (0.700,0.490,0) (0.733,0.538,0) (0.767,0.588,0)
  (0.800,0.640,0) (0.833,0.694,0) (0.867,0.751,0) (0.900,0.810,0)
  (0.933,0.871,0) (0.967,0.934,0) (1,1,0)
};
\draw[vcurve] plot [smooth] coordinates {
  (-1,1,2) (-0.967,0.934,2) (-0.933,0.871,2) (-0.900,0.810,2) (-0.867,0.751,2)
  (-0.833,0.694,2) (-0.800,0.640,2) (-0.767,0.588,2) (-0.733,0.538,2)
  (-0.700,0.490,2) (-0.667,0.444,2) (-0.633,0.401,2) (-0.600,0.360,2)
  (-0.567,0.321,2) (-0.533,0.284,2) (-0.500,0.250,2) (-0.467,0.218,2)
  (-0.433,0.188,2) (-0.400,0.160,2) (-0.367,0.134,2) (-0.333,0.111,2)
  (-0.300,0.090,2) (-0.267,0.071,2) (-0.233,0.054,2) (-0.200,0.040,2)
  (-0.167,0.028,2) (-0.133,0.018,2) (-0.100,0.010,2) (-0.067,0.004,2)
  (-0.033,0.001,2) (-0,0,2) (0.033,0.001,2) (0.067,0.004,2) (0.100,0.010,2)
  (0.133,0.018,2) (0.167,0.028,2) (0.200,0.040,2) (0.233,0.054,2)
  (0.267,0.071,2) (0.300,0.090,2) (0.333,0.111,2) (0.367,0.134,2)
  (0.400,0.160,2) (0.433,0.188,2) (0.467,0.218,2) (0.500,0.250,2)
  (0.533,0.284,2) (0.567,0.321,2) (0.600,0.360,2) (0.633,0.401,2)
  (0.667,0.444,2) (0.700,0.490,2) (0.733,0.538,2) (0.767,0.588,2)
  (0.800,0.640,2) (0.833,0.694,2) (0.867,0.751,2) (0.900,0.810,2)
  (0.933,0.871,2) (0.967,0.934,2) (1,1,2)
};
\draw[->,ultra thin] (-0.367,0.046,0) -- (-0.367,0.134,0);
\draw[->,ultra thin] (-0.333,0,0) -- (-0.333,0.111,0);
\draw[->,ultra thin] (-0.300,-0.043,0) -- (-0.300,0.090,0);
\draw[->,ultra thin] (-0.267,-0.084,0) -- (-0.267,0.071,0);
\draw[->,ultra thin] (-0.233,-0.123,0) -- (-0.233,0.054,0);
\draw[->,ultra thin] (-0.200,-0.160,0) -- (-0.200,0.040,0);
\draw[->,ultra thin] (-0.167,-0.194,0) -- (-0.167,0.028,0);
\draw[->,ultra thin] (-0.133,-0.227,0) -- (-0.133,0.018,0);
\draw[->,ultra thin] (-0.100,-0.257,0) -- (-0.100,0.010,0);
\draw[->,ultra thin] (-0.067,-0.284,0) -- (-0.067,0.004,0);
\draw[->,ultra thin] (-0.033,-0.310,0) -- (-0.033,0.001,0);
\draw[->,ultra thin] (-0,-0.333,0) -- (-0,0,0);
\draw[->,ultra thin] (0.033,-0.310,0) -- (0.033,0.001,0);
\draw[->,ultra thin] (0.067,-0.284,0) -- (0.067,0.004,0);
\draw[->,ultra thin] (0.100,-0.257,0) -- (0.100,0.010,0);
\draw[->,ultra thin] (0.133,-0.227,0) -- (0.133,0.018,0);
\draw[->,ultra thin] (0.167,-0.194,0) -- (0.167,0.028,0);
\draw[->,ultra thin] (0.200,-0.160,0) -- (0.200,0.040,0);
\draw[->,ultra thin] (0.233,-0.123,0) -- (0.233,0.054,0);
\draw[->,ultra thin] (0.267,-0.084,0) -- (0.267,0.071,0);
\draw[->,ultra thin] (0.300,-0.043,0) -- (0.300,0.090,0);
\draw[->,ultra thin] (0.333,-0,0) -- (0.333,0.111,0);
\draw[->,ultra thin] (0.367,0.046,0) -- (0.367,0.134,0);
\draw[->,ultra thin] (0.400,0.093,0) -- (0.400,0.160,0);
\draw[->,ultra thin] (-0.367,0.046,2) -- (-0.367,0.134,2);
\draw[->,ultra thin] (-0.333,0,2) -- (-0.333,0.111,2);
\draw[->,ultra thin] (-0.300,-0.043,2) -- (-0.300,0.090,2);
\draw[->,ultra thin] (-0.267,-0.084,2) -- (-0.267,0.071,2);
\draw[->,ultra thin] (-0.233,-0.123,2) -- (-0.233,0.054,2);
\draw[->,ultra thin] (-0.200,-0.160,2) -- (-0.200,0.040,2);
\draw[->,ultra thin] (-0.167,-0.194,2) -- (-0.167,0.028,2);
\draw[->,ultra thin] (-0.133,-0.227,2) -- (-0.133,0.018,2);
\draw[->,ultra thin] (-0.100,-0.257,2) -- (-0.100,0.010,2);
\draw[->,ultra thin] (-0.067,-0.284,2) -- (-0.067,0.004,2);
\draw[->,ultra thin] (-0.033,-0.310,2) -- (-0.033,0.001,2);
\draw[->,ultra thin] (-0,-0.333,2) -- (-0,0,2);
\draw[->,ultra thin] (0.033,-0.310,2) -- (0.033,0.001,2);
\draw[->,ultra thin] (0.067,-0.284,2) -- (0.067,0.004,2);
\draw[->,ultra thin] (0.100,-0.257,2) -- (0.100,0.010,2);
\draw[->,ultra thin] (0.133,-0.227,2) -- (0.133,0.018,2);
\draw[->,ultra thin] (0.167,-0.194,2) -- (0.167,0.028,2);
\draw[->,ultra thin] (0.200,-0.160,2) -- (0.200,0.040,2);
\draw[->,ultra thin] (0.233,-0.123,2) -- (0.233,0.054,2);
\draw[->,ultra thin] (0.267,-0.084,2) -- (0.267,0.071,2);
\draw[->,ultra thin] (0.300,-0.043,2) -- (0.300,0.090,2);
\draw[->,ultra thin] (0.333,-0,2) -- (0.333,0.111,2);
\draw[->,ultra thin] (0.367,0.046,2) -- (0.367,0.134,2);
\draw[->,ultra thin] (0.400,0.093,2) -- (0.400,0.160,2);
\draw[->,ultra thin] (-0.367,-1.821,1) -- (-0.367,-1.732,1);
\draw[->,ultra thin] (-0.333,-1.843,1) -- (-0.333,-1.732,1);
\draw[->,ultra thin] (-0.300,-1.865,1) -- (-0.300,-1.732,1);
\draw[->,ultra thin] (-0.267,-1.888,1) -- (-0.267,-1.732,1);
\draw[->,ultra thin] (-0.233,-1.910,1) -- (-0.233,-1.732,1);
\draw[->,ultra thin] (-0.200,-1.932,1) -- (-0.200,-1.732,1);
\draw[->,ultra thin] (-0.167,-1.954,1) -- (-0.167,-1.732,1);
\draw[->,ultra thin] (-0.133,-1.976,1) -- (-0.133,-1.732,1);
\draw[->,ultra thin] (-0.100,-1.999,1) -- (-0.100,-1.732,1);
\draw[->,ultra thin] (-0.067,-2.021,1) -- (-0.067,-1.732,1);
\draw[->,ultra thin] (-0.033,-2.043,1) -- (-0.033,-1.732,1);
\draw[->,ultra thin] (-0,-2.065,1) -- (-0,-1.732,1);
\draw[->,ultra thin] (0.033,-2.043,1) -- (0.033,-1.732,1);
\draw[->,ultra thin] (0.067,-2.021,1) -- (0.067,-1.732,1);
\draw[->,ultra thin] (0.100,-1.999,1) -- (0.100,-1.732,1);
\draw[->,ultra thin] (0.133,-1.976,1) -- (0.133,-1.732,1);
\draw[->,ultra thin] (0.167,-1.954,1) -- (0.167,-1.732,1);
\draw[->,ultra thin] (0.200,-1.932,1) -- (0.200,-1.732,1);
\draw[->,ultra thin] (0.233,-1.910,1) -- (0.233,-1.732,1);
\draw[->,ultra thin] (0.267,-1.888,1) -- (0.267,-1.732,1);
\draw[->,ultra thin] (0.300,-1.865,1) -- (0.300,-1.732,1);
\draw[->,ultra thin] (0.333,-1.843,1) -- (0.333,-1.732,1);
\draw[->,ultra thin] (0.367,-1.821,1) -- (0.367,-1.732,1);
\draw[->,ultra thin] (0.400,-1.799,1) -- (0.400,-1.732,1);
\draw[->,ultra thin] (-0.467,0.067,1) -- (-0.467,0,1);
\draw[->,ultra thin] (-0.433,0.133,1) -- (-0.433,0,1);
\draw[->,ultra thin] (-0.400,0.200,1) -- (-0.400,0,1);
\draw[->,ultra thin] (-0.367,0.267,1) -- (-0.367,0,1);
\draw[->,ultra thin] (-0.333,0.333,1) -- (-0.333,0,1);
\draw[->,ultra thin] (-0.300,0.400,1) -- (-0.300,0,1);
\draw[->,ultra thin] (-0.267,0.467,1) -- (-0.267,0,1);
\draw[->,ultra thin] (-0.233,0.533,1) -- (-0.233,0,1);
\draw[->,ultra thin] (-0.200,0.600,1) -- (-0.200,0,1);
\draw[->,ultra thin] (-0.167,0.667,1) -- (-0.167,0,1);
\draw[->,ultra thin] (-0.133,0.733,1) -- (-0.133,0,1);
\draw[->,ultra thin] (-0.100,0.800,1) -- (-0.100,0,1);
\draw[->,ultra thin] (-0.067,0.867,1) -- (-0.067,0,1);
\draw[->,ultra thin] (-0.033,0.933,1) -- (-0.033,0,1);
\draw[->,ultra thin] (-0,1,1) -- (-0,0,1);
\draw[->,ultra thin] (0.033,0.933,1) -- (0.033,0,1);
\draw[->,ultra thin] (0.067,0.867,1) -- (0.067,0,1);
\draw[->,ultra thin] (0.100,0.800,1) -- (0.100,0,1);
\draw[->,ultra thin] (0.133,0.733,1) -- (0.133,0,1);
\draw[->,ultra thin] (0.167,0.667,1) -- (0.167,0,1);
\draw[->,ultra thin] (0.200,0.600,1) -- (0.200,0,1);
\draw[->,ultra thin] (0.233,0.533,1) -- (0.233,0,1);
\draw[->,ultra thin] (0.267,0.467,1) -- (0.267,0,1);
\draw[->,ultra thin] (0.300,0.400,1) -- (0.300,0,1);
\draw[->,ultra thin] (0.333,0.333,1) -- (0.333,0,1);
\draw[->,ultra thin] (0.367,0.267,1) -- (0.367,0,1);
\draw[->,ultra thin] (0.400,0.200,1) -- (0.400,0,1);
\draw[->,ultra thin] (0.433,0.133,1) -- (0.433,0,1);
\draw[->,ultra thin] (0.467,0.067,1) -- (0.467,0,1);
\end{tikzpicture}
\hfil
\capt{The tensegrity with the vector field described.}
\label{fig:bigYwithMotion}
\end{figure}

Now $YV_i$ is a function that is strictly positive on the center
subtensegrity. After all, it shrinks that center cable, and since it acts 
orthogonally to them, it doesn't change the lengths of the other cables, to
first order.  Of course, $YV_i$ is negative on the upper cables of those 
tensegrities surrounding the center one--we knew it had to be negative 
somewhere since the tensegrity is bar equivalent.

For each $V_i$, we can define $f_i(e) = \max \{ YV_i(e), 0 \}$.  Then $f_i \in
C(\E)^+$, and as $i$ gets large, the distance $\|YV_i - f_i\|$ gets arbitrarily
small.  On the other hand $\|YV_i\| = 1 = \|f_i\|$ for all $i$, so we have a
subspace, $\im Y$ that, at a distance $1$ from the origin, gets arbitrarily
close to the nonnegative orthant and yet which intersects that orthant only at
the origin.

\chapter{Conclusion and Future Directions}
\label{chapter:Future}
So where are we and where do we go from here?  

We set out to take the theorem ``A tensegrity is bar equivalent if and only if
it has a strictly positive stress'' into the new world of continuous
tensegrities.  We've come close to accomplishing that.  We have:
\begin{enumerate}[label=\roman*)]
\item A tensegrity is partially bar equivalent if and only if it has a semipositive stress.
\item A tensegrity is bar equivalent if it has a strictly positive stress.
\item A tensegrity is minimally bar equivalent only if it has a strictly
positive stress.
\item If a tensegrity is countably covered with subtensegrities that have
strictly positive stresses, then it has a strictly positive stress.
\end{enumerate}
If we want to show that a tensegrity has a motion, showing that it has no
semipositive stress gives a very strong result.  On the other hand, if we wish
to show that a tensegrity has no motion, showing that it has a strictly
positive stress gives a strong result and showing that it has a semipositive
stress still gives a valuable result.  We can handle all of these situations
--- but only when $\E$ is compact.  

That is not always a given.  For example, Connelly et al.  \ycite{CDR} find a
flow of any nonconvex polygon in $\RR^2$ that ``convexifies'' the polygon and
in which all non-adjacent vertices are constantly moving away from each other.
\begin{figure}[ht] \hfil \begin{tikzpicture} \draw[cable] (1,0) node[vertex]
(a) {} -- (3,0) node[vertex] (b) {} -- (4.5,0.5) node[vertex] (c) {} -- (5,2)
node[vertex] (d) {} -- (3,3) node[vertex] (e) {} -- (2,4) node[vertex] (f) {}
-- (0,2) node[vertex] (g) {} -- cycle; \draw[strut] (a) -- (c) -- (e) -- (g) --
(b) -- (d) -- (f) -- (a); \draw[strut] (a) -- (d) -- (g) -- (c) -- (f) -- (b)
-- (e) -- (e); \end{tikzpicture} \hfil \begin{tikzpicture} \draw[vcurve]
plot[smooth cycle] coordinates {(1,0) (3,0) (4.5,0.5) (5,2) (3,3) (2,4) (0,2)};
\clip plot[smooth cycle] coordinates {(1,0) (3,0) (4.5,0.5) (5,2) (3,3) (2,4)
(0,2)}; \foreach \x in {-5,5,...,180} { \draw[strut] (1,0) -- +(\x:6cm); }
\end{tikzpicture} \hfil \capt[Moving the ``convexifying'' problem in
to the continuous world.]{Moving the ``convexifying polygons'' problem in to
the continuous world.  In the continuous picture, the struts from only one
vertex are shown.} \label{fig:toContRealm} \end{figure} Moving into the
continuous realm (see \vref{fig:toContRealm}), we find that for each vertex, we
get a ``fan'' of struts.  

Within that fan, we can find a sequence of struts which limits to the
``zero-length strut'' connecting the vertex to itself.  If we exclude those
``zero-length struts'', the edgeset is not compact.  If we include them, the
load on those struts will have to be zero, so any motion is semipositive at
best.  On the other hand, perhaps it could be shown that the motion, though
semipositive, is strictly positive on all of the struts connecting distinct
points.

In the work here, we have not addressed the question of infinitesimal rigidity.
A bar equivalent continuous tensegrity is infinitesimally rigid if the
associated bar framework is, but determining whether a continuous bar framework
is infinitesimally rigid is not always simple.   Roth \ycite{MR619413} gives a
``rigidity predictor'' for bar frameworks.  Perhaps it could be extended to
continuous bar frameworks and thence to continuous tensegrities.  

Any work on infinitesimal rigidity of continuous bar frameworks is likely to be
informed by the literature on continuous families of linear equalities.  That
is a realization which may well be exceedingly fruitful.  What we have done 
here is valuable in its own right in solving problems with continuous 
tensegrities, but the statement that $YV \in C(\E)+$ is really just an infinite
family of linear inequalities.  By modifying $Y$, the theory could be applied
to other situations.  For example, Ashton et al. \ycite{RidgeRunner} use both
struts and kinks (places where the vertex curve has reached a curvature 
constraint) to build the rigidity operator.  A similar idea could perhaps be
employed to provide an alternate proof of Schur's Theorem (page
\pageref{SchursTheorem}).

As another example, consider moving the differential version of the
``convexifying'' task up one dimension.  In this realm, we could still seek to
have all vertex-vertex distances increase, or we could instead build our
rigidity operator using the volumes of the tetrahedra determined by sets of 4
points on the surface, and we could require that those volumes be strictly
increasing.  Or, with careful choice of rigidity operator, perhaps our results
can provide new proofs for various of the Theorems of the Alternative (see
\ref{section:ToA} and \ref{appendix:Motzkin}).\Index{Theorem of the
Alternative}[1!234] 

Speaking of Theorems of the Alternative, at the end of \ref{appendix:Motzkin},
we relate stresses of a tensegrity to variations lying in $\X^\perp$.  It would
be interesting to know how those variations relate to the ``curvature force''
of \ocite{CFKSW}.  For that matter, there is plenty of work yet to be done to
understand exactly which variations thus arise.  In fact, at this point, the
circle-of-struts example is the only bar equivalent tensegrity we have whose
design variations are local isometries of the vertex curve.

In the world where all variations are design variations, there is exploration
that could be done into which tensegrities that arise from the antipodal
strut/fixed-skip cable technique are bar equivalent (see
\vref{section:exampleToDo}).  Also, the continuous analog of the Gr\"unbaum
polygons from page \pageref{GrunbaumPolygons} waits to be found.

Finally, there are those conjectures from pages \pageref{conj:minXbarConj} and
\pageref{conj:CountCov}:
\begin{ConjThree}
\statementminXbarConj
\end{ConjThree}
\begin{ConjOne}
\statementCountCov
\end{ConjOne}
\begin{ConjTwo}
\statementStrictlyPositive
\end{ConjTwo}
With regard to that last one, having the strictly positive loads in the
interior of $C(\E)^+$ and the semipositive ones on the boundary seems
appropriate but that $C^*(\E)^+$ has no interior is surprising.  Perhaps a
different topology on $C(\E)$ would give an interior to $C^*(\E)^+$ and place
the strictly positive stresses in that interior and the semipositive ones along
the boundary.  Or perhaps one of the other Theorems of the Alternative, such as
those of Tucker \ycite{MR0089112} or Slater \ycite{MR0061636} could solve the
problem with the topology as it stands.

\nocite{MR1300410}
\nocite{MR1129886}
\nocite{MR1129886}
\nocite{MR1129886}
\nocite{MR1129886}
\nocite{MR1453120}
\nocite{MR1202431}
\nocite{MR1307623}
\nocite{MR0447465}
\nocite{MR1453120}
\nocite{MR0252038}
\nocite{Motzkin}
\nocite{MR515723}
\nocite{MR515723}
\nocite{Pardon}
\nocite{MR1681462}
\nocite{MR1998947}
\nocite{Royden}
\nocite{MR1681462}
\nocite{CFKSW}

\pseudochapter{Bibliography}

\appendix
\labelformat{chapter}{Appendix~#1}
\chapter{Euclidean Motions: Infinitesimal Rigid Motions of Space}
\label{appendix:rigidMotions}
\section{The Special Orthogonal Group: \texorpdfstring{$SO_n$}{SO(n)} and
         \texorpdfstring{$\mathfrak{so}_n$}{so(n)}}
As we mentioned on page \pageref{rigidMotionMention}, the vector fields in
$\Vf(\V)$ vary in nature.  \vref{fig:differentVF} shows three different vector
fields on the same tensegrity.  We are here interested in identifying which
vector fields, like the first two of those in the figure, can arise from rigid
motions of space.
\IndDefBeg{Euclidean motion}
\begin{figure}[ht]
\hspace*{\fill}
\tikzstyle{every node}=[vertex]
\begin{tikzpicture}[scale=2]
\draw[cable] 
  (0,0) node {} -- (1,0) node {} -- (1,1) node {} -- (0,1) node {} -- cycle;
\draw[->] (0,0) -- +(0.3,0.1);
\draw[->] (1,0) -- +(0.3,0.1);
\draw[->] (1,1) -- +(0.3,0.1);
\draw[->] (0,1) -- +(0.3,0.1);
\begin{scope}[xshift=2cm]
\draw[cable] 
  (0,0) node {} -- (1,0) node {} -- (1,1) node {} -- (0,1) node {} -- cycle;
\draw[->] (0,0) -- +(0.3,-0.1);
\draw[->] (1,0) -- +(0.3,0.3);
\draw[->] (1,1) -- +(-0.1,0.3);
\draw[->] (0,1) -- +(-0.1,-0.1);
\fill[blue] (0.25,0.75) circle (1pt);
\end{scope}
\begin{scope}[xshift=4cm]
  \draw[cable] 
    (0,0) node {} -- (1,0) node {} -- (1,1) node {} -- (0,1) node {} -- cycle;
  \draw[->] (0,0) -- +(0.3,0.3);
  \draw[->] (1,0) -- +(-0.5,0.3);
  \draw[->] (1,1) -- +(-0.1,-0.1);
  \draw[->] (0,1) -- +(0.3,-0.5);
\end{scope}
\end{tikzpicture}
\hspace*{\fill}
\capt[Three different vector fields on a square made of cables.]{Three different
vector fields on a square made of cables.  The left hand one would result in a
simple translation of the tensegrity.  The center one indicates a rotation
around the marked point.  The right field is the only one which would change
the shape of the tensegrity.  The first two are elements of $T(p)$, the third
is in $I(p)$ but not in $\Ib(p)$ or~$T(p)$.}
\label{fig:differentVF}
\end{figure}

In the work we do here we will look for inspiration to Murray et al.
\ycite{MR1300410}, who much of the same work, though primarily in~$\RR^3$.  

Intuitively, a \emph{rigid motion of space}\Index{rigid motion of
space}[2!1,34 1!234] is a continuous movement of $\RR^n$ which preserves all
distances and angles.  More formally, $H_t$ is a rigid motion of space if 
$$H_t\mcol [0,1] \to \aut \RR^n$$
(where $\aut \RR^n$ is the set of automorphisms of $\RR^n$) such that $H_0$ is
the identity map on $\RR^n$ and for all $v,w \in \RR^n$ and $t \in [0,1]$,
$$\langle H_t(v), H_t(w) \rangle = \langle v, w \rangle.$$ We will call a map
$g\mcol \RR^n \to \RR^n$ a \emph{rigid transformation}\IndexDef{rigid
transformation}[2,1 1!2] if there is some rigid motion of space, $H_t$ such
that $g = H_1$.

Clearly, any such transformation must take an orthonormal basis for $\RR^n$ to
an orthonormal basis for $\RR^n$, so we can describe it by telling where the
origin goes and how the new set of coordinates is oriented with respect to the
old.  That is, it will consist of a translation and a ``reorientation''.  That
reorientation may be a rotation (in fact, must be in $\RR^2$ and $\RR^3$)
or a composition of rotations \cite{MR1129886}*{p.\ 125}.

We'll look first at the reorientation and then see what it takes to include the 
translation.  Taking as our original basis for $\RR^n$ the traditional 
orthonormal basis consisting of unit vectors in the coordinate directions,
$e_1, \dotsc, e_n$, we can describe a reorientation by taking the vectors
$\hat{e}_1, \dotsc, \hat{e}_n$ of the new coordinate axes and building a 
reorientation matrix
$$R = [\hat{e}_1  \cdots  \hat{e}_n].$$
Now since the $\hat{e}_1, \dotsc, \hat{e}_n$ form an orthonormal basis
and hence are mutually orthogonal and of unit length, we have
$$R^\top R = R R^\top = \begin{bmatrix}
1 & 0 & \cdots & 0 \\
0 & 1 & \cdots & 0 \\
\vdots & \vdots & \ddots & \vdots \\
0 & 0 & \cdots & 1
\end{bmatrix} = I\text{ (the identity matrix).}$$

Then, since 
\begin{equation}
1 = \det I = \det R^\top R = \det R^\top \det R = (\det R)^2,
\label{eq:detSquared}
\end{equation}
we have that
$\det R = \pm 1$.  Now if we think of the set of all $n \cross n$ matrices
as the space $\RR^{n^2}$, then the determinant is a continuous function from
that space into~$\RR$.  But the set $\{-1,1\}$ has two connected components,
so the set of matrices described by \ref{eq:detSquared} must have at least two
connected components.\footnote{It turns out to be exactly two (see
\ocite{MR1129886}*{p.\ 124}), but we don't need that information here.}
Furthermore, since the rigid motions provide a path between the identity
transformation and the reorientation in question, all of the reorientations
which play into our rigid transformations must lie in the connected component
in which the identity has its home.  

That means that we have $\det R = 1$ for every reorientation $R$.  So our
reorientations are contained in the set of ``special orthogonal matrices'',
which we can write formally as
$$SO_n = \left\{R \in \RR^{n \cross n} : R^\top R = I\text{ and } \det R =
1\right\}$$\IndexDef{$SO_n$@SOn (Special Orthogonal Group)}[ ]
(see, for example, \ocite{MR1129886}*{p.\ 124}).

On the other hand, if $R \in SO_n$, then for any $v,w \in \RR^n$,
\begin{align*}
\langle Rv, Rw \rangle = v^\top R^\top R w = v^\top w = \langle v, w \rangle,
\end{align*}
and we can produce a rigid motion of space by taking the rotations which make
up $R$ and varying their angles from $0$ through the angle desired.  Hence our
set of reorientations is $SO_n$.

However, we're not really interested in $SO_n$, \textit{per se}.  We want the
``infinitesimal reorientations'', the directions in the space of $n \cross n$
matrices in which a rigid motion of space can leave the origin -- the tangent
space to $SO_n$ at the identity (more on that in a moment).

Now the $n \cross n$ matrices form a normed linear space using the Frobenius
norm (see, for example, \ocite{MR1129886}*{p.\ 153}), which is essentially the
Euclidean norm with the matrices thought of as vectors in~$\RR^{n^2}$.

$SO_n$ is a subset of that space of $n \cross n$ matrices (not a subspace,
since, for example, scaling a matrix changes its determinant), but it is more
than that.  It is a Lie group\Index{Lie group}[2!1 ],\Index{special orthogonal
group}[3!12 ] which is to say, a group that is also a smooth manifold.  So its
tangent space at the identity is well-defined.  This tangent space to $SO_n$ at
the identity is usually called $\mathfrak{so}_n$\IndexDef{$\mathfrak{so}_n$@son
(tangent space to $SO_n$)}[ ] \cite{MR1453120}*{pp.\ 12,64}.

Now the elements of a tangent space to a surface at a given point provide us
the best linear approximation to that surface at that point.  So if our surface
is defined as the level curve of some function $f(x) = c$, the elements of the
tangent space are precisely those vectors $v$ for which $f(x + \varepsilon v) -
c$ is proportional to $\varepsilon^2$ rather than to~$\varepsilon$.

In our case, $SO_n$ is defined by the equation\footnote{There is also the determinant requirement, but that simply selects which connected component we are concerned with.  As we are asking for the tangent space \emph{at the identity}, we
will get the same results if we think of just the component in hand or of
both components together.} $R^\top R = I$, so we can see it as the level curve
$f(R) = 0$ for $f(R) = \|R^\top R - I\|$.

Then $dR \in \mathfrak{so}_n$ if and only if $\frac{d}{d\varepsilon}
f(I + \varepsilon dR)\big|_{\varepsilon=0} = 0$.  That is,
\begin{align*}
0 & = \left. \frac{d}{d\varepsilon} f(I + \varepsilon dR)\right|_{\varepsilon =
      0} \\
& = \left. \frac{d}{d\varepsilon} \|(I + \varepsilon dR)^\top (I +
    \varepsilon dR) - I\| \right|_{\varepsilon = 0} \\
& = \left. \frac{d}{d\varepsilon} \|I^\top I + \varepsilon dR^\top I +
    I^\top \varepsilon dR + \varepsilon^2 dR^\top dR - I\| \right|_{\varepsilon =
    0} \\
& = \left. \frac{d}{d\varepsilon} \|\varepsilon dR^\top + \varepsilon dR +
    \varepsilon^2 dR^\top dR\| \right|_{\varepsilon = 0}. \\
& = \left. \frac{d}{d\varepsilon} \varepsilon \|dR^\top + dR +
    \varepsilon dR^\top dR\| \right|_{\varepsilon = 0}. \\
& = \left. \|dR^\top + dR + \varepsilon dR^\top dR\| +
    \varepsilon \frac{d}{d\varepsilon} \|dR^\top + dR +
    \varepsilon dR^\top dR\| \right|_{\varepsilon = 0}. \\
\intertext{If we knew that $\frac{d}{d\varepsilon} \|dR^\top + dR +
\varepsilon dR^\top dR\|$ was finite, we'd know that that entire term is zero
at $\varepsilon = 0$.  Well, denoting the entry of $dR$ in the $i^\text{th}$
row and $j^\text{th}$ column by $d_{ij}$ and the equivalent entry in $dR^\top
dR$ by $c_{ij}$ gives us}
& = \left. \|dR^\top + dR + \varepsilon dR^\top dR\| +
    \varepsilon \frac{\sum_{i=1}^n
   \sum_{j=1}^n (d_{ji} + d_{ij} + \varepsilon c_{ij}) c_{ij}}{\|dR^\top + dR +
   \varepsilon dR^\top dR\|} \right|_{\varepsilon = 0}. \\
& = \|dR^\top + dR\| \Longrightarrow dR^\top + dR = \zero \Longrightarrow 
dR = - dR^\top.
\end{align*}

So $dR \in \mathfrak{so}_n$ if and only if $dR$ is a skew-symmetric matrix.
That means that the entries above the diagonal may be arbitrary but then the
rest of the matrix is determined, so $\mathfrak{so}_n$ is
$\frac{n(n-1)}{2}$-dimensional.  Since the dimension of a tangent space is the
same as the dimension of the manifold (see, for example, \ocite{MR1202431}*{p.
7}), that tells us also the dimension of~$SO_n$.  However, there's an intuitive
argument for the same value.  

Thinking back, we remember that the columns of any element of $SO_n$ form an
orthonormal basis for~$\RR^n$.  So we start with $n^2$ degrees of freedom
and apply $n$ constraints to make the column vectors unit length.  Then, for
any pair of vectors, we add another constraint to make them orthogonal.  That's
$\binom{n}{2} = \frac{n(n-1)}{2}$ more constraints for a total of 
$\frac{n(n+1)}{2}$ degrees of freedom gone and $\frac{n(n-1)}{2}$ remaining.

\section{The Special Euclidean Group: \texorpdfstring{$SE_n$}{SE(n)} and
         \texorpdfstring{$\mathfrak{se}_n$}{se(n)}}
Now that we've seen what the reorientations are like, we'd like to add in the
translations.  Of course, since any linear transformation on $\RR^n$ takes the
origin to the origin, translations are outside that realm.  But by moving up
one dimension, we can return to the world of linear transformations.  

To do so we will take a point in $\RR^n$, $p = [\begin{smallmatrix} p_1 & p_2 &
\cdots & p_n\end{smallmatrix}]^\top$ and move it into $\RR^{n+1}$ by adding a
coordinate and setting its value to $1$, thus: $[\begin{smallmatrix} p_1 &
\cdots & p_n & 1 \end{smallmatrix}]^\top$.  So then, if $R \in SO_n$ tells how
we wish to reorient the axes and $t \in \RR^n$ tells where the origin should
go, we can find the new location of the point $p$ by the linear transformation:
\begin{equation}
\begin{bmatrix}\hat{p} \\ 1\end{bmatrix} = \begin{bmatrix} R & t \\ \zero & 1
\end{bmatrix} \begin{bmatrix} p \\ 1 \end{bmatrix} = \begin{bmatrix} Rp + t \\
1 \end{bmatrix}
\label{eq:sen}
\end{equation}

The matrices $\left[\begin{smallmatrix} R & t \\ \zero & 1
\end{smallmatrix}\right]$ form the \emph{Special Euclidean
Group}\IndexDef{special Euclidean group}[3!12 2!3 ]
$SE_n$\IndexDef{$SE_n$@SEn (Special Euclidean Group)}[ ] (sometimes called just
the \emph{Euclidean Group}, $Euc_n$).  If we label the $(n+1) \cross (n+1)$
matrices in the fashion $\left[\begin{smallmatrix} R & t \\ c & d
\end{smallmatrix}\right]$ (where $R$ is an $n \cross n$ matrix, $t$ is an
$n$-coordinate column vector, $c$ an $n$-coordinate row vector and $d$ a
scalar), then $SE_n$ is the subset of that space defined by $R^\top R = I$, $c
= \zero$ and $d = 1$. 

Since we can pair any reorientation with any translation, $SE_n$ appears to be
isomorphic to $SO_n \cross \RR^n$ and so would be of dimension
$\frac{n(n+1)}{2}$, and that's close to true.  The dimension is right, but the
group law is wrong.

Reorientations and translations don't commute.  From \vref{eq:sen} we can see
that an element of $SE_n$ first applies the reorientation and then the
translation, so suppose we have two elements of $SE_n$, say 
\begin{enumerate}[label=\ensuremath{\alph* = (R_\alph*,t_\alph*)}:,leftmargin=*]
\item ``rotate around the origin by $\nicefrac{\pi}{3}$ and then move right $1$
unit'', and
\item ``rotate around the origin by $\nicefrac{\pi}{6}$ and then move right $1$
unit''.
\end{enumerate}
\begin{figure}[ht]
\hfil
\subfloat[The effect of applying $(R_b,t_b)$ and then $(R_a,t_a)$.]{
\begin{tikzpicture}[scale=0.65]
% Smile
\draw[ultra thin,blue!40,<->] (-2,0) -- (8,0);
\draw[ultra thin,blue!40,<->] (0,-2) -- (0,5);
\draw (-1,1) circle (1pt) (1,1) circle (1pt) (-1,0) arc (-180:0:1cm);
\begin{scope}[xshift=4cm]
  \begin{scope}[rotate=30]
    \draw (-1,1) circle (1pt) (1,1) circle (1pt) (-1,0) arc (-180:0:1cm);
  \end{scope}
\end{scope}
\begin{scope}[xshift=4cm]
  \begin{scope}[rotate=60]
    \begin{scope}[xshift=4cm]
      \begin{scope}[rotate=30]
        \draw (-1,1) circle (1pt) (1,1) circle (1pt) (-1,0) arc (-180:0:1cm);
      \end{scope}
    \end{scope}
  \end{scope}
\end{scope}
\draw[dashed,->] (0,2) .. controls +(45:1cm) and +(135:1cm) .. 
  node [sloped,above,scale=0.7] {$(R_b,t_b)$} (3,1.732);
\draw[dashed,->] (3,1.732) .. controls +(90:1cm) and +(210:1cm) ..
  node [sloped,above,scale=0.7] {$(R_a,t_a)$} (4,3.464);
\end{tikzpicture}
}
\hfil
\subfloat[The effect of applying $(R_a R_b,t_a + t_b)$.]{
\begin{tikzpicture}[scale=0.65]
\draw[ultra thin,blue!40,<->] (-2,0) -- (10,0);
\draw[ultra thin,blue!40,<->] (0,-2) -- (0,5);
% Smile
\draw (-1,1) circle (1pt) (1,1) circle (1pt) (-1,0) arc (-180:0:1cm);
\begin{scope}[xshift=8cm]
  \begin{scope}[rotate=90]
    \draw (-1,1) circle (1pt) (1,1) circle (1pt) (-1,0) arc (-180:0:1cm);
  \end{scope}
\end{scope}
\draw[dashed,->] (0,2) .. controls +(45:2cm) and +(180:2cm) .. 
  node [sloped,above,scale=0.7] {$(R_a R_b,t_a + t_b)$} (6,0);
\end{tikzpicture}
}
\hfil
\capt[$SE_n \not\cong SO_n \cross \RR^n$]{$SE_n$ is not isomorphic to $SO_n
\cross \RR^n$, since the operations of rotating around the origin and
translating do not commute.}
\label{fig:smileys}
\end{figure}
Then, as \vref{fig:smileys} shows, the operations $(R_a,t_a)(R_b,t_b)$ and 
$(R_a R_b, t_a + t_b)$ are not, in general, the same.  However, the reorientations are automorphisms of $\RR^n$ and the translations are $\RR^n$.  So we can
think of the reorientations acting on the translations and define a group 
operation thus:
$$(R_1, t_1) (R_2, t_2) = (R_1 R_2, t_1 + (R_1 t_2)).$$
The group thus formed is known as the \emph{semidirect
product}\Index_{semidirect product} of $SO_n$ and $\RR^n$, denoted $SO_n
\ltimes \RR^n$ (see, for example, \ocite{MR1307623}*{p.\ 167} or
\ocite{MR0447465}*{p.\ 236}, for more information about semidirect products).

\begin{prop}
$SO_n \ltimes \RR^n$ is a group and is isomorphic to~$SE_n$.
\end{prop}
\begin{proof}
Let $R_1,R_2,R_3 \in SO_n$ and $t_1,t_2,t_3 \in \RR^n$.  Then
\begin{align*}
(R_1,t_1)((R_2,t_2)(R_3,t_3)) 
& = (R_1,t_1)(R_2 R_3,t_2 + R_2 t_3) \\
& = (R_1 R_2 R_3, t_1 + R_1 t_2 + R_1 R_2 t_3) \\
& = (R_1 R_2, t_1 + R_1 t_2)(R_3, t_3) \\
& = ((R_1,t_1)(R_2,t_2))(R_3, t_3) 
\end{align*}
so we have associativity.  The element $(I,\zero)$ is clearly the identity, so
we only need inverses.  Since $SO_n$ is a group, we have inverses for all of
its elements.   Using these we can build the inverses we need:
\begin{align*}
(R_1,t_1)(R_1^{-1},-R_1^{-1} t_1) 
& = (R_1 R_1^{-1},t_1 + R_1 (-R_1^{-1} t_1)) \\
& = (I,t_1 - t_1) = (I,\zero).
\end{align*}
So $SO_n \ltimes \RR^n$ is a group.  Consider the map $f\mcol SO_n \ltimes
\RR^n \to SE_n$ by $(R,t) \mapsto \left[\begin{smallmatrix} R & t \\ \zero & 1 
\end{smallmatrix} \right]$.  Then 
\begin{align*}
f((R_1,t_1)) f((R_2,t_2)) 
& = \begin{bmatrix} R_1 & t_1 \\ \zero & 1 \end{bmatrix}
    \begin{bmatrix} R_2 & t_2 \\ \zero & 1 \end{bmatrix} \\
& = \begin{bmatrix} R_1 R_2 & R_1 t_2 + t_1 \\ \zero & 1 \end{bmatrix} 
  = f((R_1,t_1)(R_2,t_2)).
\end{align*}
So the two groups are homomorphic.  Furthermore, $f$ is onto and has trivial 
kernel, by inspection, so $SO_n \ltimes \RR^n$ and $SE_n$ are isomorphic.
\end{proof}

It turns out that $SE_n$ is a Lie group\Index{Lie group}[2!1 ], just as $SO_n$
was \cite{MR1453120}*{p.\ 64}, which is good, as we are still interested in
the tangent space, $\mathfrak{se}_n$.\IndexDef{$\mathfrak{se}_n$@sen (tangent
space to $SE_n$)}[ ]  Since the elements of $SE_n$ look like
$\left[\begin{smallmatrix} R & t \\ c & d \end{smallmatrix}\right]$,
the elements of $\mathfrak{se}_n$ look like
$\left[\begin{smallmatrix} dR & dt \\ dc & dd \end{smallmatrix}\right]$.
We already know that $dR$ must be in $\mathfrak{so}_n$ and we can simply take
differentials of the other two equations to give us $dc = \zero$
and $dd = 0$.  So the elements of $\mathfrak{se}_n$ look like
$\left[\begin{smallmatrix} dR & dt \\ \zero & 0 \end{smallmatrix}\right]$
where $dR \in \mathfrak{so}_n$ and $dt \in \RR^n$.

All that remains is, then, is to figure out which vector fields correspond to 
the elements of $\mathfrak{se}_n$.  To find those, we take the points in
$p(\V)$ and apply the elements of $\mathfrak{se}_n$ to them.  In the next
subsection, we give an example.

\subsection{Examples}
Consider the crossed-square\Index{crossed square}[ ] example from
\vref{fig:RWexample}, shown also in \vref{fig:repeatVF}.  We've shown that the
elements of $\mathfrak{se}_2$ look like 
$\left[\begin{smallmatrix} 0 & a & x \\ -a & 0 & y \\ 0 & 0 & 0
\end{smallmatrix}\right].$
For the crossed square, $p(2) = (1,0)$, so we multiply
$$\begin{bmatrix} 0 & a & x \\ -a & 0 & y \\ 0 & 0 & 0 \end{bmatrix} 
\begin{bmatrix} 1 \\ 0 \\ 1 \end{bmatrix} = \begin{bmatrix} x \\ -a + y \\ 0 
\end{bmatrix}.$$
Thus we have $V(2) = (x,y-a)$.  Following through the rest of the arithmetic
gives us that the vector fields in $T(p)$ look like
$$V(v) = \begin{cases}
(x,y), & v=1 \\
(x,y-a), & v=2 \\
(x+a,y-a), & v=3 \\
(x+a,y), & v=4.
\end{cases}$$
\begin{figure}[ht]
\hfil
\begin{tikzpicture}[scale=2]
\draw[strut] (0,0) -- (1,1) (1,0) -- (0,1);
\draw[cable] 
  (0,0) node [vertex] {1} -- (1,0) node [vertex] {2} -- 
  (1,1) node [vertex] {3} -- (0,1) node [vertex] {4} -- cycle;
\draw[->] (0,0) -- +(0.3,0.1);
\draw[->] (1,0) -- +(0.3,0.1);
\draw[->] (1,1) -- +(0.3,0.1);
\draw[->] (0,1) -- +(0.3,0.1);
\draw (0.5,-0.25) node {$\{a = 0, x = 0.3, y = 0.1\}$};
\begin{scope}[xshift=3cm]
\draw[strut] (0,0) -- (1,1) (1,0) -- (0,1);
\draw[cable] 
  (0,0) node [vertex] {1} -- (1,0) node [vertex] {2} -- 
  (1,1) node [vertex] {3} -- (0,1) node [vertex] {4} -- cycle;
\draw[->] (0,0) -- +(0.3,-0.1);
\draw[->] (1,0) -- +(0.3,0.3);
\draw[->] (1,1) -- +(-0.1,0.3);
\draw[->] (0,1) -- +(-0.1,-0.1);
\draw (0.5,-0.25) node {$\{a = -0.4, x = 0.3, y = -0.1\}$};
\end{scope}
\end{tikzpicture}
\hfil
\capt{Two vector fields in $T(p)$ on the crossed square.}
\label{fig:repeatVF}
\end{figure}
For example, the two vector fields in \vref{fig:repeatVF} (which are the same
as the first two in \vref{fig:differentVF}) are given by $\{a = 0, x = 0.3, y =
0.1\}$ and $\{a = -0.4, x = 0.3, y = -0.1\}$ respectively.
\IndDefEnd{Euclidean motion}

\chapter{Motzkin's Theorem}
\label{appendix:Motzkin}
\section{Preparing to use Motzkin's Theorem}
\label{section:MotSetup}
\subsection{Background}
In \vref{section:ToA}, we talked about Theorems of the
Alternative\Index{Theorem of the Alternative}[1!234] and how our main theorem
relates to them.  Here we will use one of them to give an alternate proof for
the first main theorem.  The new version actually has stronger hypotheses, so
it's not as strong a theorem, but we include it here because the technique is
interesting.

The theorems of Stiemke\Index{Stiemke Theorem}[2!1] (\ref{theorem:Stiemke}) and
Gordan\Index{Gordan Theorem}[2!1] (\ref{theorem:Gordan}) in \ref{section:ToA}
are for finite-dimensional Euclidean space.  We need them extended to more
general spaces.  Unfortunately, a literature search failed to locate 
generalized versions.

However, Mangasarian\Index_,{Olvi Mangasarian} \ycite{MR0252038} proves
Stiemke's theorem using a theorem from Theodore Motzkin's\Index_,{Theodore S.
Motzkin} dissertation \ycite{Motzkin} commonly known as Motzkin's Transposition
Theorem.  And B.\ D.\ Craven\Index_,{Bruce D. Craven} \ycite{MR515723} gives a
version of Motzkin's Theorem for normed linear spaces.

The theorem reads as follows (here $L(X,W)$ denotes the space of continuous
linear maps from $X$ to $W$ and $A^T$ is Craven's notation for the adjoint of
$A$):
\begin{theorem}[Motzkin Transposition Theorem]
\label{theorem:MAT}
\IndexDef{Motzkin Transposition Theorem}[3!12 ]
Let $X$, $W$, $Z$ be normed spaces; let $S \subset W$ be a convex cone, with
$\INT S \ne \varnothing$; let $T \subset Z$ be a closed convex cone; let $A \in
L(X,Z)$ and $B \in L(X,W)$.  If the convex cone $A^T(T^*)$ is weak-* closed,
then exactly one of the two following systems has a solution:
\begin{enumerate}[label=(\Roman*)]
\item $-Ax \in T$, $-Bx \in \INT S$ $(x \in X)$;
\item $p \circ B + q \circ A = 0, q \in T^*, 0 \ne p \in S^*$.
\label{item:caseTwo}
\end{enumerate}
\end{theorem}
\begin{proof}
See \ocite{MR515723}*{p.\ 32}.
\end{proof}
\begin{figure}[ht]
$$
\xymatrix{
& X^*  & X \ar[dl]_A \ar[dr]^B \\
T^* \subset Z^* \ar[ur]^{A^T} & T \subset Z & & W \supset S
}
$$
\capt{The diagram for the Motzkin Transposition Theorem.}
\label{fig:MotzDiag}
\end{figure}
\ref{fig:MotzDiag} shows the diagram for Motzkin's theorem.  We know that
$C(\E)^+$ is a convex cone with nonempty interior, so it seems reasonable to
set $S = C(\E)^+$ and, of course, $W = C(\E)$.  In that case, if we made $B$ be
the negative of our rigidity operator and set $X = \Vf(\V)$, then $-Bx \in \INT
S$ would mean $YV \in \INT C(\E)^+$, that is, $V$ is a strictly positive
motion.  

That seems promising, but what are $T$ and $Z$?  Well, that $-Ax \in T$ could
allow us to implement the design variations.  If $T = \X$ (which, of course,
means that $Z \in \Vf(\V)$) and $A = -\mathrm{Id}$, the negative of the 
identity map, then $-Ax \in T$ becomes $V \in \X$.  

How does that work out with the other alternative?  Using the definitions
we have put forth, we get Case \ref{item:caseTwo} as $-pY - q = 0$, where
$0 \ne p \in C^*(\E)^+$ and $q \in \X^*$.  So we could take $p$ as our stress $\mu$ and (so long as $\X$ is a subspace so that $\X^* = \X^\perp$), we get
the statement $Y^*\mu \in \X^\perp$.  That is, there exists a semipositive
stress.  So this seems a good match for our needs.  Our only concerns are
about having $\X^\perp$ be closed and the weak-* closure of $T^*$.  To make
those easier to answer, we'll make $\Vf(\V)$ into a Hilbert space.

\renewcommand{\V}{\gamma}
\subsection{$\V$ and $\Vf(\V)$}
\label{subsection:newV}
In what follows, we will work exclusively with tensegrities whose vertices lie
along a simple, rectifiable (in fact, arclength parametrized) curve (again
called $\V$) with finitely many connected components and total length $\ell$.

We want $\Vf(\V)$ to be a Hilbert space, so we need to endow it with an inner
product.  In doing this, we'll follow the lead of John Pardon \ycite{Pardon}.
We'll redefine $\Vf(\V)$ to be a subset of the absolutely continuous functions
and establish the inner product $\langle V_1, V_2 \rangle =
\int V_1' \cdot V_2' \, ds$ on it.

For that inner product to be finite, we'll need $\int |V'|^2 \ds < \infty$
for all $V \in \Vf(\V)$, and that to get completeness, we'll want the
derivative map on $\Vf(\V)$ to be 1-1.  

Pardon accomplishes that last task by requiring that $V(0) = \zero$, but we
take a different route.  We would like $\int_0^\ell V(s) \ds = \zero$, that is,
we want the ``center of mass'' of the variation to be at~$\zero$.  This will
have the effect that $V(s)$ as a variation will give $\V$ no net translation. 
In summary:
\begin{definition}
\label{definition:newVfV}
The space of variations on $\V$, $\Vf(\V)$ will consist of all absolutely
continuous vector fields $V$ on $\V$ which satisfy $\int |V'|^2 < \infty$
and $\int V = \zero$.
\end{definition}

It may help to have a more explicit definition for elements of~$\Vf(\V)$.  
\begin{prop}
$\Vf(\V)$, as defined in Definition \vref{definition:newVfV} consists of those
vector fields on $\V$ which satisfy $\int |V'|^2 < \infty$ and 
\begin{equation}
V(s) = \int_0^s V'(t) \dt - \frac{1}{\ell} \int_0^\ell \int_0^t V'(r) \, dr
\dt.
\label{eq:Vs}
\end{equation}
\end{prop}
\begin{proof}
By Theorem \vref{theorem:iiac}, vector fields which satisfy \vref{eq:Vs} are
absolutely continuous, so we need only show that this is the equation which
gives us $\int V = \zero$.  We know that $V(s) = \int_0^s V'(t) \dt + V(0)$. 
So 
\begin{align*}
\int_0^\ell V(s) \ds 
& = \int_0^\ell \int_0^s V'(t) \dt \ds + \int_0^\ell V(0) \ds
& = \int_0^\ell \int_0^s V'(t) \dt \ds + \ell V(0)
\end{align*}
To get $\int_0^\ell V(s) \ds = 0$, then, we need 
$$V(0) = -\frac{1}{\ell} \int_0^\ell \int_0^s V'(t) \dt \ds.$$
\end{proof}

Now we are ready to establish the nature of our new~$\Vf(\V)$.
\begin{prop}
$\Vf(\V)$ is a Hilbert space under the inner product $\langle V_1, V_2 \rangle
= \int V_1' \cdot V_2' \ds$.
\end{prop}
\begin{proof}
First we need to establish that $\langle V_1, V_2 \rangle$ as we've defined it
really \emph{is} an inner product.  We have three things to check:
\begin{enumerate}[label=(\roman*)]
\item Symmetry:
$$\langle V_1, V_2 \rangle = \int V_1' \cdot V_2' \ds = \int V_2' \cdot V_1'
\ds = \langle V_2, V_1 \rangle$$
\item Positive definiteness:
$$\langle V, V \rangle = \int V' \cdot V' \ds \ge 0$$
We note that $\langle V, V \rangle = 0$ only when $V' = \zero$ almost
everywhere.  So $V$ must be constant.   But because $\int V = \zero$, this only
happens when $V = \zero$.  Also, by our definition of $\Vf(\V)$, $\int \|V'\|^2 < \infty$, so $\langle V, V \rangle$ is finite.
\item Linearity in the first position:
\begin{align*}
\langle aV_1 + bV_2, V_3 \rangle & = \int (aV_1' + bV_2') \cdot V_3 \ds \\
& = a\int V_1' \cdot V_3' + b\int V_2' \cdot V_3' \ds 
  = a\langle V_1, V_3 \rangle + b \langle V_2, V_3 \rangle
\end{align*}
(for all $a,b \in \RR$ and $V_1, V_2, V_3 \in \Vf(\V)$).
\end{enumerate}
So $\langle V_1, V_2 \rangle$ is an inner product.  

Next we need to show that our space is complete with the norm $\|V\| = \langle
V, V \rangle^\frac{1}{2}$.

Suppose that there is a Cauchy sequence $\left\{ V_n \right\}$ of elements
of~$\Vf(\V)$.  That is, for any $\varepsilon > 0$, there exists some
positive integer $N$ such that whenever $m,n > N$, we have 
$$\varepsilon > \|V_m - V_n\| = \left(\int \|V'_m -
V'_n\|^2\right)^{\nicefrac{1}{2}}.$$
But we recognize that final term as $\|V'_m - V'_n\|_2$, the $L^2$ norm of
$V'_m - V'_n$.  $L^2$ is complete under its norm (see, for example,
\ocite{MR1681462}*{p.\ 183}), so there is some $V' \in L^2$ such that $V'_n \to
V'$.  With that in hand we can calculate $V(s)$ from \vref{eq:Vs}.  The result
is an absolutely continuous (by Theorem \vref{theorem:iiac}) vector field for
which $\int |V'_n|^2 < \infty$ and $\int V = 0$. Thus $V$ is an element
of~$\Vf(\V)$.  And so $\Vf(\V)$ with our inner product is a Hilbert space.%
% \footnote{The space we have constructed bears some resemblence to a Sobolev
% Space.  In fact, if we had demanded that $\|V\|_2 = 0$ instead of $\int V =
% \zero$, we would have a subspace of Sobolev Space $W^{1,2}$.  The interested
% reader is encouraged to seek further information in \ocite{SobolevSpaces}.}
% 
% \begin{color}[rgb]{0.2,0.5,0}
% \begin{prop}
% If the center of mass of $V$ is $\zero$, then norms $\|V'\|_2$ and
% $\|V\|_2+\|V'\|_2$ are equivalent.
% \end{prop}
% \begin{proof}
% We need to show that there exists some numbers $c_1,c_2 > 0$ such that 
% \begin{equation}
% \label{eq:eqNorms}
% c_1 \|V'\|_2 \le \|V\|_2 + \|V'\|_2 \le c_2 \|V'\|_2
% \end{equation}
% for all $V \in \Vf(\V)$ (see, for example, \ocite{MR1681462}*{p.\ 152}).  The
% left side of \ref{eq:eqNorms} is trivially true, so we go in search of
% $c_2$.  In fact, it will suffice to show that there exists some $k > 0$ such
% that $\|V\|_2 \le k \|V'\|_2$.  But
% \begin{align*}
% \|V\|_2 
% & = \left(\int_0^\ell |V(s)|^2 \ds\right)^\frac{1}{2} \\
% & = \left(\int_0^\ell \left|\int_0^s V'(t) \dt + V(0)\right|^2
% \ds\right)^\frac{1}{2} \\
% & \le \left(\int_0^\ell \int_0^s |V'(t)|^2 \dt \int_0^s \one
% \dt\ds\right)^\frac{1}{2} \text{by Cauchy-Schwarz} \\
% & = \left(\int_0^\ell s \int_0^s |V'(t)|^2 \dt \ds\right)^\frac{1}{2} \\
% & \le \left(\int_0^\ell \ell \int_0^\ell |V'(t)|^2 \dt \ds\right)^\frac{1}{2} \\
% & = \left(\ell^2 \|V'(t)\|^2\right)^\frac{1}{2} \\
% & = \ell \|V'\|_2
% \end{align*}
% \end{proof}
% \end{color}
\end{proof}

One delightful aspect of Hilbert spaces is that every linear functional on a
Hilbert space is given by inner product with some element of the space.  This
is the result of another theorem also called the ``Riesz Representation
Theorem'', which we give below.  This allows us to think of elements of
a Hilbert space and its dual as lying in the same space.
\begin{theorem}[Riesz Representation Theorem for Hilbert Spaces]
\label{theorem:RRHilbert}
Let $X$ be a Hilbert space.  Then, for every $f \in X^*$, there exists a unique
element $p \in X$ such that $f(x) = \langle x, p \rangle$ for all $x \in X$ and
$\|f\| = \|p\|$\IndexDef{Riesz Representation Theorem {for Hilbert
Spaces}}[3!12!4 123!4].
\end{theorem}
\begin{proof}
See \ocite{MR1998947}*{p.\ 430}.
\end{proof}

\subsection{$\X$ and $\X^\perp$}
\label{subsection:X}
Having established a new $\Vf(\V)$, we choose our ``design variations'' $\X$ as
follows:
\begin{equation}
\X = \left\{ V \in \Vf(\V): V'(s) \cdot \V'(s) = 0 \text{ almost
everywhere } \right\}.
\label{eq:defM}
\end{equation}
We will need to check that $\X$ is a subspace of our new $\Vf(\V)$.
Furthermore, we'll need $\X$ to be closed.  In proving that it is, we'll make
use of the Cauchy-Schwarz inequality:
\begin{theorem}[Cauchy-Schwarz Inequality]
\Index{Cauchy-Schwarz Inequality}[ ]
For all $x,y$ in a Hilbert space, $\langle x,y \rangle \le \|x\| \|y\|$.
\end{theorem}
\begin{proof}
See \ocite{Royden}*{p.\ 210}.
\end{proof}

\begin{lemma}
\label{lemma:MVclosed}
$\X$, as defined in \ref{eq:defM}, is a closed subspace of~$\Vf(\V)$.
\end{lemma}
\begin{proof}
First we'll dispose of the subspace portion.  Let $V,W \in \X$ and
$\alpha,\beta \in \RR$.  Then 
$$(\alpha V + \beta W)'(s) \cdot \V'(s) = \alpha V'(s) \cdot \V'(s) + \beta
W'(s) \cdot V'(s) = 0$$
almost everywhere, putting $\alpha V + \beta W \in \X$.  Now for closure.

Let $\{ V_n \} \in \X$ be a sequence of vector fields that limits to some
vector field $V \in \Vf(\V)$.  Then
\begin{align*} 
\int \left(V'(s) \cdot \V'(s)\right)^2 \ds 
  & = \int \left( V'(s) \cdot \V'(s) - V'_n(s) \cdot \V'(s) \right)^2 \ds 
      \text{ since } V_n \in \X \\ 
  & = \int \left( (V'(s) - V'_n(s)) \cdot \V'(s) \right)^2 \ds \\ 
  & \le \int \|V'(s) - V'_n(s)\|^2 \|\V'(s)\|^2 \ds \text{ by Cauchy-Schwarz}\\ 
  & = \int \|V'(s) - V'_n(s)\|^2 \ds \text{ since $\V$ is unit speed } \\
  & = \|V - V_n\|^2 
\end{align*} 
But $V_n \to V$, so for any $\varepsilon > 0$, there is a positive integer $N$
such that $\|V - V_n\| < \sqrt{\varepsilon} \Rightarrow \|V - V_n\|^2 <
\varepsilon$.  That means that $\int ( V'(s) \cdot \V'(s))^2$ is nonnegative,
but less than every positive number.  Hence,
$$\int (V'(s) \cdot \V'(s))^2 \, ds = 0$$ 
and thus $V'(s) \cdot \V'(s) = 0$ almost everywhere.  So $V \in \X$, and
$\X$ is closed.
\end{proof}

\begin{lemma}
\label{lemma:Xpsubspace}
The dual cone\Index{dual cone}[2,1 ] of a subspace $X$ is its annihilator.
That is, if $X$ is a subspace, $X^* = X^\perp$.
\end{lemma}
\begin{proof}
Let $\mu \in X^*$.  Then $\mu \cdot x \ge 0$ by definition.  However, since
$X$ is a subspace, $x \in X \Rightarrow -x \in X$, so $\mu \cdot (-x) = -\mu
\cdot x$ must also be nonnegative.  But that means that $\mu \cdot x = 0$.
That is, $\mu \in X^\perp$.

Conversely, let $\eta \in X^\perp$.  Then, for all $x \in X$, we have
$\eta \cdot x = 0 \ge 0$, so $\eta \in X^*$.
\end{proof}

\begin{lemma}
\label{lemma:Subspacen}
Since $\X$ is a subspace, $\X^\perp$ and $Y(\X)$ are subspaces.
\end{lemma}
\begin{proof}
Let $\alpha, \beta \in \RR$ and $V_1, V_2 \in \X^\perp$.  Then, for all $V \in
\X$, we have
$$(\alpha V_1 + \beta V_2) \cdot V = \alpha V_1 \cdot V + \alpha V_2 \cdot V =
0,$$
so $\X^\perp$ is a subspace.

Similarly, if $YV,YW \in Y(\X)$, then $\alpha V + \beta W \in \X$ (since $\X$
is a subspace) and $\alpha YV + \alpha YW = Y(\alpha V + \alpha W)$.  So
$Y(\X)$ is a subspace.
\end{proof}

\section{Motzkin's Theorem Applied}
\label{section:MT}
We are now ready to apply the Motzkin Transposition Theorem to our setup.
For clarity, \vref{fig:TwoDiagrams} shows the Motzkin Diagram and our
situation given as a parallel diagrams.
\begin{figure}[ht]
$$
\xymatrix{
  & X \ar[dl]_A \ar[dr]^B \\
T \subset Z & & Y \supset S
}
\quad
\xymatrix{
  & \Vf(\V) \ar[dl]_{-\mathrm{Id}} \ar[dr]^{-Y} \\
\X \subset \Vf(\V) & & C(\E) \supset C(\E)^+
}
$$
\capt{The Motzkin diagram and our matching situation.}
\label{fig:TwoDiagrams}
\end{figure}

We need to understand what ``$A^T(T^*)$ is weak-* closed'' means in our
situation.  For us, $A$ and $A^T$ are both negative identity functions, one on
$\Vf(\V)$ and the other on $\Vf^*(\V)$, which, thanks to the Riesz
Representation Theorem for Hilbert Spaces\Index{Riesz Representation Theorem
{for Hilbert Spaces}}[3!12!4 123!4] (Theorem \vref{theorem:RRHilbert}), we can
identify with~$\Vf(\V)$.  

We showed, in Lemma \ref{lemma:Xpsubspace}, that $\X^* = \X^\perp$, so
$A^T(T^*)$ is simply $-\X^\perp$, which (by Lemma \ref{lemma:Subspacen}) is
just $\X^\perp$.  Now for the theorem.

\begin{theorem}
\Index{semipositive stress $\Leftrightarrow$ partially $\X$-bar@X-bar
equivalent equiv.}[1!2!3457 2!1!3457 456!312 56!4!312]
A tensegrity $G(p)$ is partially $\X$-bar equivalent if and only if $G(p)$ has
a semipositive stress.
\end{theorem}
\begin{proof}
$\Vf(\V)$ and $C(\E)$ are normed spaces, $C(\E)$ with the $\sup$ norm and
$\Vf(\V)$ with the norm arising from its inner product.  So we have the right
kind of spaces for Motzkin's Theorem (\vref*{theorem:MAT}).  For its
hypotheses, we need the following to be true:\begin{enumerate}[label=(\alph*)]
\item $C(\E)^+$ is a convex cone with nonempty interior, 
\item $\X$ is a closed convex cone and 
\label{item:Xccc}
\item $\X^\perp$ is closed.
\label{item:Xpc}
\end{enumerate} 
If those are true, then exactly one of the follow systems has a solution:
\begin{enumerate}[label=(\Roman*)]
\item $V \in \X$, $YV \in \INT C(\E)^+$;
\label{item:SPM}
\item $Y^* \mu = V$, $-V \in \X^\perp$, $\zero \ne \mu \in C^*(\E)^+$.
\label{item:SPS}
\end{enumerate}
That is, $G(p)$ has either a strictly positive motion \ref{item:SPM} or a
semipositive stress \ref{item:SPS}.

$\X$ is not only a closed convex cone, but moreover a closed subspace, by Lemma
\ref{lemma:MVclosed}, so Item \ref{item:Xccc} is satisfied.  Since $\Vf(\V)$ is
a Hilbert space, the subspace $\X^\perp$ (since it is defined by orthogonality
to a set) is closed (see, for example, \ocite{MR1681462}*{p.\ 173}).  That
takes care of Item \ref{item:Xpc}.

It remains only to show that $C(\E)^+$ is a convex cone with nonempty interior.
Certainly, if $f$ is a nonnegative continuous function from $\E$ to $\RR$, then
$\alpha f$ is as well, for any $\alpha > 0$.  Likewise, if $f,g \in C(\E)^+$,
then any convex combination of $f$ and $g$ is also in $C(\E)^+$, so $C(\E)^+$
is a convex cone.  Furthermore, by Lemma \vref{lemma:intCEp}, $\INT C(\E)^+$
contains all of the strictly positive functions (and thus is nonempty).

By applying Motzkin's theorem, we have our result.
\end{proof}

\section{Motzkin on the Other Hand}
\label{section:MotOthHand}
Having had such success with that direction, it seems only natural to see if we
can get anything more out of Motzkin's Theorem.  What we've proven is 
$$\text{there is no strictly positive motion} \Leftrightarrow 
  \text{there is a semipositive stress}.$$
Can we use the Motzkin Transposition Theorem to show that 
$$\text{there is no semipositive motion} \Leftrightarrow
  \text{there is a strictly positive stress}?$$

The natural thing would be to try to require that our measures fall in the 
interior of~$C^*(\E)^+$.  But we saw in Theorem \vref{theorem:noInterior}
that when $\E$ is infinite, $C^*(\E)^+$ has no interior.  So any success along
these lines must be limited to the finite case.  Let's give it a try.
\renewcommand{\V}{\mathscr{V}}

We need $\INT S$ to be $\INT C^*(\E)^+$, but in our basic setup,
$Y^*$ maps from $C^*(\E)^+$ and nothing maps to it.  So that's where
we'll send our identity function (see \ref{fig:OneMoreDiagram}).
\begin{figure}[ht]
$$
\xymatrix{
  & C^*(\E) \ar[dr]^{-\mathrm{Id}} \ar[dl]_{Y^*} \\
\X^\perp \subset \Vf^*(\V) & & C^*(\E) \supset C^*(\E)^+
}
$$
\capt{Our new use of the Motzkin theorem.}
\label{fig:OneMoreDiagram}
\end{figure}

Now, we need $C^*(\E)^+$ to have non-empty interior, which it does.  We need
$\X^\perp = \{ \zero \}$ to be a closed convex cone, which it is.  Finally, we
need 
$$Y^{**}((\X^\perp)^*) = Y(\Vf(\V))$$ 
to be weak-* closed.  But $Y$ is a finite-dimensional linear map and so its
image is a subspace of the finite-dimensional $C(\E)$ and thus it is closed.

Now we can apply Motzkin's theorem and get that exactly one of these two
systems has a solution:
\begin{enumerate}[label=(\Roman*)]
\item $-Y^*\mu = \zero$, $\mu \in \INT C^*(\E)^+$
\item $YV = f$, $V \in \Vf(\V)$, $0 \ne f \in C(\E)^+$.
\end{enumerate}
That is, ``either there is a strictly positive stress or else there is a
semipositive motion''. Thus we have a new proof of Roth and Whiteley's theorem.

\renewcommand{\V}{\gamma}
\section{Stresses as Variations}
One somewhat surprising effect of making $\Vf(\V)$ into a Hilbert space is that
any given stress $\mu$ has a variation, $V_\mu \in \X^\perp$ associated with
it.  That is because $Y^*\mu \in \X^\perp \subset \Vf^*(\V)$ and $\Vf(\V)$ and
$\Vf^*(\V)$ can be identified with each other, thanks to the Riesz
Representation Theorem for Hilbert Spaces (Theorem \vref{theorem:RRHilbert}).
\Index{Riesz Representation Theorem {for Hilbert Spaces}}[3!12!4 123!4]

For example, in the case of the circle of struts, we found a stress $d\theta$,
which was uniform on the struts.  $Y^* \, d\theta$, then corresponds to some
variation $V_{d\theta}$ in such way that
\begin{equation}
\begin{split}
\langle V_{d\theta}, V \rangle & = (Y^* \, d\theta) V = d\theta(YV) =
\int_0^\pi YV
\, d\theta \\
& = \int_0^\pi (V(\theta)-V(\theta+\pi)) \cdot ((\cos \theta,\sin \theta)
- (\cos(\theta+\pi),\sin(\theta+\pi)))\,d\theta \\
& = \int_0^\pi (V(\theta)-V(\theta+\pi)) \cdot 2(\cos \theta,\sin \theta)
\,d\theta \\
& = \int_0^\pi V(\theta) \cdot 2(\cos \theta,\sin \theta) \,d\theta 
  - \int_0^\pi V(\theta+\pi) \cdot 2(\cos \theta,\sin \theta) \,d\theta \\
& = \int_0^\pi V(\theta) \cdot 2(\cos \theta,\sin \theta) \,d\theta 
  - \int_\pi^{2\pi} V(\theta) \cdot 2(\cos(\theta-\pi),\sin(\theta-\pi))
    \,d\theta \\
& = \int_0^\pi V(\theta) \cdot 2(\cos \theta,\sin \theta) \,d\theta 
  + \int_\pi^{2\pi} V(\theta) \cdot 2(\cos \theta,\sin(\theta) \,d\theta \\
& = \int_0^{2\pi} V(\theta) \cdot 2(\cos \theta,\sin \theta) \,d\theta. \\
\end{split}
\label{eq:IPOne}
\end{equation}
On the other hand, our inner product on $\Vf(\V)$ says that
$$\langle V_{d\theta}, V \rangle = \int_0^{2\pi} V'_{d\theta}(s) \cdot V'(s)
\,ds.$$
We can integrate this by parts to get
\begin{equation}
\langle V_{d\theta}, V \rangle = V'_{d\theta}(s) \cdot V(s) \big|_0^{2\pi}
- \int_0^{2\pi} V''_{d\theta}(s) \cdot V(s) \,ds.
\label{eq:IPTwo}
\end{equation}
Now, $0$ and $2\pi$ are the same point on our domain, so that first term is
zero.  We can combine \ref{eq:IPOne} and \ref{eq:IPTwo} to get
$$ - \int_0^{2\pi} V''_{d\theta}(s) \cdot V(s) \,ds.
= \int_0^{2\pi} V(\theta) \cdot 2(\cos \theta,\sin \theta) \,d\theta,$$
which must be true for all $V \in \Vf(\V)$.  Since $s = \theta$ on the unit
circle, that gives us that
$-V''_{d\theta}(\theta) = 2(\cos \theta,\sin \theta)$
or 
$$V_{d\theta}(\theta) = 2(\cos \theta, \sin \theta),$$
which is shown in \vref{fig:withVm}.
\begin{figure}[ht]
\hspace*{\fill}
\begin{tikzpicture}[scale=0.5]
\foreach \x in {0,20,...,181} {
  \draw[strut] (\x:2cm) -- (\x+180:2cm);
  \draw[->] (\x:2cm) -- (\x:4cm);
  \draw[->] (\x+180:2cm) -- (\x+180:4cm);
}
\draw[vcurve] (0,0) circle (2cm);
\end{tikzpicture}
\hspace*{\fill}
\capt{The circle of struts with the vectors of $V_{d\theta}$.}
\label{fig:withVm}
\end{figure}
\Index{circle of struts}[ ]

Thinking back to our original definition of stress, this might seem the obvious
thing.  After all, it just looks to be the weighted edge vectors at each
vertex.  However, on second glance it seems less obvious.  After all, this 
example has a high degree of symmetry.  We were able to transform an integral
on the edges into an integral on the vertices in the form $\int V \cdot
\text{(something)}$ which we could then relate to our inner product.  

This is an area which bears further investigation, both in identifying the 
variation for a given stress and in understanding the interplay between that
and the ``curvature force''\Index{curvature force}[2,1 ] of \ocite{CFKSW}.

\clearpage
\phantomsection
\addcontentsline{toc}{pseudochapter}{Index}
\renewcommand{\V}{\mathscr{V}}
\printindex
% \listoffixmes
% \chapter{Stuff that probably should go in somewhere}
% \begin{prop}
% Any bar-equivalent tensegrity that cannot be embedded in $\RR^1$ has at least
% 6 edges.
% \end{prop}
% \begin{proof}
% Since the tensegrity cannot be embedded in $\RR^1$, it must have at least one
% vertex at that edges meet at a non-straight angle.  At that vertex at least
% three edges must meet, since the weighted edgevector sum at that vertex must
% be~$\zero$.  So now we have 4 vertices.  The original vertex and 3 remote
% ones.  But no line containing two of the remote vertices can contain the
% original vertex as well, and since increasing the vertex count will only
% increase the edge count, the best we can do is to add the three remaining
% edges to make it a complete graph.  That will give three edges at each vertex
% and a total of 6 edges.
% \end{proof}
\end{document}